\documentclass[11pt]{article}
\usepackage[titletoc]{appendix}
\usepackage[labelsep=period]{caption}
\usepackage{bm}
\usepackage{fullpage}
\usepackage{epsfig}
\usepackage{graphics}
\usepackage{latexsym}
\usepackage{amsmath}
\usepackage{amsfonts}
\usepackage{amssymb}
\usepackage{mathrsfs}
\usepackage{pifont}
\usepackage{yhmath}
\usepackage{undertilde}
\usepackage{bbm}
\usepackage{dsfont}
\usepackage{eufrak}
\usepackage{wasysym}
\usepackage{underscore}
\usepackage{epstopdf}
\usepackage{fontenc}
\usepackage{amsthm}
\usepackage{xcolor}
\newtheorem{theorem}{Theorem}[section]

\newtheorem{corollary}[theorem]{Corollary}
\newtheorem{definition}[theorem]{Definition}
\newtheorem{lemma}[theorem]{Lemma}

\newtheorem{conjecture}{Conjecture}[section]
\newtheorem{claim}{Claim}[section]

\newcommand{\CC}{\mathcal{C}}

\newcommand{\TT}{\mathcal {T}}

\newcommand{\sigmabar}{\overline{\sigma}}
\newcommand{\mubar}{\overline{\mu}}
\newcommand{\rhobar}{\overline{\varrho}}
\newcommand{\phibar}{\overline{\varphi}}
\newcommand{\phiv}{\varphi}

\def\bfm#1{\mbox{\boldmath$#1$}}
\def\qed{\hfill \rule{4pt}{7pt}}
\usepackage{fullpage,graphicx}

\DeclareMathAlphabet{\mathpzc}{OT1}{pzc}{m}{it}
\title{\bf Proof of the Goldberg-Seymour Conjecture on Edge-Colorings of Multigraphs}
\author{\vspace{2mm} Guantao Chen$^{a}$\thanks{Supported in part by NSF grant DMS-1855716. E-mail: gchen@gsu.edu.} 
\quad Guangming Jing$^{b}$\thanks{Corresponding author. Supported in part by NSF grant DMS-2001130. E-mail: gjing@augusta.edu.} \quad
Wenan Zang$^{c}$\thanks{Supported in part by the Research Grants Council of Hong Kong. E-mail: wzang@maths.hku.hk.}\\
$\stackrel{a}{}$ Department of Mathematics and Statistics, Georgia State University\\
Atlanta, GA 30303, USA \smallskip\\ 
$\stackrel{b}{}$ Department of Mathematics, Augusta University\\Augusta, GA 30912, USA \smallskip\\ 
$\stackrel{c}{}$ Department of Mathematics, The University of Hong Kong\\Hong Kong, China }
\begin{document}
\date{}
\maketitle

\begin{abstract}
Given a multigraph $G=(V,E)$, the edge-coloring problem (ECP) is to color the edges of $G$ with the minimum number 
of colors so that no two adjacent edges have the same color. This problem can be naturally formulated as an integer
program, and its linear programming relaxation is referred to as the fractional edge-coloring problem (FECP). 
The optimal value of ECP (resp. FECP) is called the chromatic index (resp. fractional chromatic 
index) of $G$, denoted by $\chi'(G)$ (resp. $\chi^*(G)$). Let $\Delta(G)$ be the maximum degree of $G$ and let 
$\Gamma(G)$ be the density of $G$, defined by  
\[\Gamma(G)=\max \left\{\frac{2|E(U)|}{|U|-1}:\,\, U \subseteq V, \,\, |U|\ge 3 \hskip 2mm
{\rm and \hskip 2mm odd} \right\},\]
where $E(U)$ is the set of all edges of $G$ with both ends in $U$.  Clearly, $\max\{\Delta(G), \, 
\lceil \Gamma(G) \rceil \}$ is a lower bound for $\chi'(G)$. As shown by Seymour, $\chi^*(G)=\max\{\Delta(G), \, 
\Gamma(G)\}$. In the early 1970s Goldberg and Seymour independently conjectured that $\chi'(G) \le \max\{\Delta(G)+1, 
\, \lceil \Gamma(G) \rceil\}$. Over the past five decades this conjecture, a cornerstone in modern edge-coloring, 
has been a subject of extensive research, and has stimulated an important body of work. In this paper we present 
a proof of this conjecture. Our result implies that, first, there are only two possible values for 
$\chi'(G)$, so an analogue to Vizing's theorem on edge-colorings of simple graphs holds for multigraphs; second, 
although it is $NP$-hard in general to determine $\chi'(G)$, we can approximate it within one of its true value, and 
find it exactly in polynomial time when $\Gamma(G)>\Delta(G)$; third, every multigraph $G$ satisfies $\chi'(G)-\chi^*(G) 
\le 1$, and thus FECP has a fascinating integer rounding property.  
\end{abstract}

\newpage

\section{Introduction}

Given a multigraph $G=(V,E)$, the {\em edge-coloring problem} (ECP)\label{ECP} is to color the edges of $G$ with the minimum 
number of colors so that no two adjacent edges have the same color. Its optimal value is called the {\em chromatic 
index}\label{cindex} of $G$,  denoted by $\chi'(G)$. In addition to its great theoretical interest, ECP arises in a variety 
of applications, so it has attracted tremendous research efforts in several fields, such as discrete mathematics, combinatorial 
optimization, and theoretical computer science. Holyer \cite{H} proved that it is $NP$-hard in general 
to determine $\chi'(G)$, even when restricted to a simple cubic graph, so there is no efficient algorithm for 
solving ECP exactly unless $NP=P$, and hence the focus of extensive research has been on near-optimal solutions to 
ECP or good estimates of $\chi'(G)$.  

Let $\Delta(G)$ be the maximum degree of $G$. Clearly, $\chi'(G)\ge \Delta(G)$. There are two classical upper bounds
on the chromatic index: the first of these,  $\chi'(G) \le \lfloor \frac{3 \Delta(G)}{2} \rfloor$, was established 
by Shannon \cite{Sh} in 1949, and the second,  $\chi'(G) \le \Delta(G) +\mu(G)$, where $\mu(G)$ is the maximum multiplicity 
of edges in $G$, was proved independently by Vizing \cite{V} and Gupta \cite{G} in the 1960s. This second result is widely known as 
Vizing's theorem, which is particularly appealing when applied to a simple graph $G$, because it reveals that $\chi'(G)$
can take only two possible values: $\Delta(G)$ and $\Delta(G)+1$. Nevertheless, in the presence of multiple edges, 
the gap between $\chi'(G)$ and these three bounds can be arbitrarily large. Therefore it is desirable 
to find some other graph theoretic parameters connected to the chromatic index.

Observe that each color class in an edge-coloring of $G$ is a matching, so it contains at most $(|U|-1)/2$ edges 
in $E(U)$ for any $U \subseteq V$ with $|U|$ odd, where $E(U)$ is the set of all edges of $G$ with both ends in $U$.   
Hence the {\em density} of $G$\label{density}, defined by 
\[\Gamma(G)=\max \left\{\frac{2|E(U)|}{|U|-1}:\,\, U \subseteq V, \,\, |U|\ge 3 \hskip 2mm
{\rm and \hskip 2mm odd} \right\},\]
provides another lower bound for $\chi'(G)$. It follows that $\chi'(G) \ge \max\{\Delta(G), \, \Gamma(G)\}$.

To facilitate better understanding of the parameter $\max\{\Delta(G), \, \Gamma(G)\}$, let $A$ be the edge-matching 
incidence matrix of $G$ (that is, each column of $A$ is the incidence vector of a matching). Then ECP can be naturally 
formulated as an integer program, whose linear programming (LP) relaxation is exactly as given below:
\begin{center}
\begin{tabular}{ll}
\hbox{Minimize} \ \ \ & $ {\bm 1}^T {\bm x}$  \\
\hbox{\hskip 0.2mm subject to} & $ A{\bm x} = {\bm 1}$  \\
& \hskip 3mm ${\bm x} \ge {\bm 0}.$
\end{tabular}
\end{center}
This linear program is called the {\em fractional edge-coloring problem} (FECP)\label{FECP}, and its optimal 
value is called the {\em fractional chromatic index}\label{findex} of $G$, denoted by  $\chi^*(G)$. As shown by Seymour \cite{Se}
using Edmonds' matching polytope theorem \cite{e65}, it is always true that $\chi^*(G)=\max\{\Delta(G), \,\Gamma(G) \}$. 
Thus the preceding inequality actually states that $\chi'(G) \ge \chi^*(G)$.

As $\chi'(G)$ is integer-valued, we further obtain $\chi'(G) \ge \max\{\Delta(G), \, \lceil \Gamma(G) \rceil \}$. 
How good is this bound? In the early 1970s Goldberg \cite{G73} and Seymour \cite{Se} independently made the following 
conjecture.

\begin{conjecture} \label{GS}
Every multigraph $G$ satisfies $\chi'(G)\le \max\{\Delta(G)+1, \, \lceil \Gamma(G) \rceil\}$.
\end{conjecture} 

Over the past five decades Conjecture \ref{GS} has been a subject of extensive research, and has stimulated an
important body of work; see McDonald \cite{Mc} for a survey on this conjecture and Stiebitz et al. \cite{SSTF} for 
a comprehensive account of edge-colorings.   

Several approximate results state that $\chi'(G)\le \max\{\Delta(G)+\tau(G), \, \lceil \Gamma(G) \rceil\}$, where 
$\tau(G)$ is a positive number depending on $G$. Asymptotically, Kahn \cite{K} showed that $\tau(G)=o(\Delta(G))$.  
Scheide \cite{S} and Chen, Yu, and Zang \cite{CYZ} independently proved that $\tau(G)\le \sqrt{\Delta(G)/2}$. 
Chen et al. \cite{GT} improved this to $\tau(G)\le \sqrt[3]{\Delta(G)/2}$. Recently, Chen and Jing \cite{GG} 
further strengthened this as  $\tau(G)\le \sqrt[3]{\Delta(G)/4}$. 

There is another family of results, asserting that $\chi'(G)\le \max\{\frac{m \Delta(G)+(m-3)}{m-1}, \, \lceil \Gamma(G) 
\rceil\}$, for increasing values of $m$. Such results have been obtained by Andersen \cite{An77} and  Goldberg \cite{G73} 
for $m=5$, Andersen \cite{An77} for $m=7$, Goldberg \cite{G84} and Hochbaum, Nishizeki, and Shmoys \cite{HNS86} for $m=9$, 
Nishizeki and Kashiwagi \cite{NK} and Tashkinov \cite{T} for $m=11$, Favrholdt, Stiebitz, and Toft \cite{FST} for $m=13$, 
Scheide \cite{S} for $m=15$,  Chen et al. \cite{GT} for $m=23$, and Chen and Jing \cite{GG} for $m=39$. It is
worthwhile pointing out that, when $\Delta(G) \le 39$, the validity of Conjecture \ref{GS} follows instantly from  
Chen and Jing's result \cite{GG}, because $\frac{39 \Delta(G)+36}{38}<\Delta(G)+2$.  

Haxell and McDonald \cite{HM} obtained a different sort of approximation to Conjecture \ref{GS}: $\chi'(G)\le \max\{\Delta(G)+
2 \sqrt{\mu(G) \log \Delta(G)}, \, \lceil \Gamma(G) \rceil\}$. Another way to obtain approximations for Conjecture \ref{GS}
is to incorporate the order $n$ of $G$ (that is, number of vertices) into a bound. In this direction, Plantholt \cite{P99}
proved that $\chi'(G)\le \max\{\Delta(G),  \, \lceil \Gamma(G) \rceil+1+ \sqrt{n \log (n/6)}\}$ for any multigraph
$G$ with ever order $n\ge 572$. In \cite{P13}, he established an improved result that is applicable to all multigraphs.   

Marcotte \cite{M} showed that $\chi'(G)=\max\{\Delta(G), \, \lceil \Gamma(G) \rceil\}$ for any multigraph $G$ with no 
$K_5^-$-minor, thereby confirming Conjecture \ref{GS} for this multigraph class. Recently, Haxell, Krivelevich, and 
Kronenberg \cite{HKK} established Conjecture \ref{GS} for random multigraphs. 

The purpose of this paper is to present a proof of the Goldberg-Seymour conjecture. 

\begin{theorem} \label{ThmGS}
	Every multigraph $G$ satisfies $\chi'(G)\le \max\{\Delta(G)+1, \, \lceil \Gamma(G) \rceil\}$.
\end{theorem} 

Let $r$ be a positive integer. A multigraph $G=(V,E)$ is called an $r$-{\em graph} if $G$
is regular of degree $r$ and for every $X \subseteq V$ with $|X|$ odd, the number of edges
between $X$ and $V-X$ is at least $r$. Note that if $G$ is an $r$-graph, then $|V(G)|$ is even
and $\Gamma(G)=r$. Seymour \cite{Se} also proposed the following weaker version of Conjecture 
\ref{GS}, which amounts to saying that $\chi'(G)\le \max\{\Delta(G), \, \lceil \Gamma(G) \rceil\}+1$ for 
any multigraph $G$.

\begin{conjecture} \label{Sey} 
Every $r$-graph $G$ satisfies $\chi'(G)\le r+1$. 
\end{conjecture}

The following conjecture was posed by Gupta \cite{G} in 1967 and can be deduced from Conjecture \ref{GS}, as 
shown by Scheide \cite{S07}. 

\begin{conjecture} \label{Gupta67} 
Let $G$ be a multigraph such that $\Delta(G)$ cannot be expressed in the form $2p\mu(G)-q$,
for any two integers $p$ and $q$ satisfying $p > \lfloor (q+1)/2 \rfloor$ and $q \ge 0$. Then 
$\chi'(G)\le \Delta(G)+\mu(G)-1$. 
\end{conjecture}

A multigraph $G$ is called {\em critical}\label{critical} if $\chi'(H)<\chi'(G)$ for any proper subgraph $H$ of $G$. 
In edge-coloring theory, critical multigraphs are of special interest, because they have much more 
structural properties than arbitrary multigraphs. The following two conjectures are due to Jakobsen \cite{J0,J} 
and were proved by Andersen \cite{An77} to be weaker than Conjecture \ref{GS}. 

\begin{conjecture} \label{And0} 
Let $G$ be a critical multigraph with $\chi'(G) \ge \Delta(G)+2$. Then $G$ contains an odd number of
vertices.
\end{conjecture}

\begin{conjecture} \label{And} 
Let $G$ be a critical multigraph with $\chi'(G)>\frac{m \Delta(G)+(m-3)}{m-1}$ for an odd 
integer $m \ge 3$. Then $G$ has at most $m-2$ vertices.
\end{conjecture}

Motivated by Conjecture \ref{GS}, Hochbaum, Nishizeki, and Shmoys \cite{HNS86} formulated the following 
conjecture concerning the approximability of ECP.

\begin{conjecture} \label{HNS} 
There is a polynomial-time algorithm that colors the edges of any multigraph $G$ using at most $\chi'(G)+1$
colors. 
\end{conjecture}

Since Conjectures \ref{Sey}-\ref{And} are all weaker than the Goldberg-Seymour conjecture, the truth of them follows 
from Theorem \ref{ThmGS} as corollaries.

\begin{theorem} \label{ThmSey} 
Every $r$-graph $G$ satisfies $\chi'(G)\le r+1$. 
\end{theorem}

\begin{theorem} \label{ThmGupta67} 
Let $G$ be a multigraph such that $\Delta(G)$ cannot be expressed in the form $2p\mu(G)-q$,
for any two integers $p$ and $q$ satisfying $p > \lfloor (q+1)/2 \rfloor$ and $q \ge 0$. Then 
$\chi'(G)\le \Delta(G)+\mu(G)-1$. 
\end{theorem}

\begin{theorem} \label{ThmAnd0} 
Let $G$ be a critical multigraph with $\chi'(G) \ge \Delta(G)+2$. Then $G$ contains an odd number of
vertices.
\end{theorem}

\begin{theorem} \label{ThmAnd} 
Let $G$ be a critical multigraph with $\chi'(G)>\frac{m \Delta(G)+(m-3)}{m-1}$ for an odd 
integer $m \ge 3$. Then $G$ has at most $m-2$ vertices.
\end{theorem}

We have seen that FECP is intimately tied to ECP. For any multigraph $G$, the fractional chromatic index 
$\chi^*(G)=\max\{\Delta(G), \,\Gamma(G) \}$ can be determined in polynomial time 
by combining the Padberg-Rao separation algorithm for $b$-matching polyhedra \cite{PR} (see also \cite{LRT, PW}) with 
binary search. In \cite{CZZ}, Chen, Zang, and Zhao designed a combinatorial polynomial-time algorithm for finding 
the density $\Gamma(G)$ of any multigraph $G$, thereby resolving a problem posed in both Stiebitz et al. \cite{SSTF} 
and Jensen and Toft \cite{JT}. Nemhauser and Park \cite{NP} observed that FECP can be solved in polynomial 
time by an ellipsoid algorithm, because the separation problem of its LP dual is exactly the maximum-weight 
matching problem (see also Schrijver \cite{Sc}, Theorem 28.6 on page 477). In \cite{CZZ}, 
Chen, Zang, and Zhao also came up with a combinatorial polynomial-time algorithm for FECP.
 
Although our proof of Theorem \ref{ThmGS} is not algorithmic in nature, we believe that it can be adapted to 
yield a polynomial-time algorithm for finding an edge-coloring of any multigraph $G$ with at most 
$\max\{\Delta(G)+1, \, \lceil \Gamma(G) \rceil\}$ colors. A successful implementation would lead to an 
affirmative answer to Conjecture \ref{HNS} as well. 

Some remarks may help to put Theorem \ref{ThmGS} in proper perspective. 

First, by Theorem \ref{ThmGS}, there are only two possible values for the chromatic index of a multigraph $G$: 
$\max\{\Delta(G), \, \lceil \Gamma(G) \rceil\}$ and  $\max\{\Delta(G)+1, \, \lceil \Gamma(G) \rceil\}$. Thus an 
analogue to Vizing's theorem on edge-colorings of simple graphs, a fundamental result in discrete mathematics, holds for 
multigraphs.

Second, Theorem \ref{ThmGS} exhibits a dichotomy on edge-coloring: While Holyer's theorem \cite{H} tells us that it is 
$NP$-hard to determine $\chi'(G)$, we can approximate it within one of its true value, because $\max\{\Delta(G)+1, \, 
\lceil \Gamma(G) \rceil\} - \chi'(G) \le 1$. Furthermore, if $\Gamma(G)>\Delta(G)$, then $\chi'(G)=\lceil \Gamma(G) \rceil$,
so it can be found in polynomial time \cite{CZZ,PR}.  

Third, by Theorem \ref{ThmGS} and aforementioned Seymour's theorem on fractional chromatic index, every multigraph $G=(V,E)$ satisfies 
$\chi'(G)-\chi^*(G) \le 1$, which can be naturally extended to the weighted case. Let $w(e)$ be a 
nonnegative integral weight on each edge $e\in E$ and let ${\bm w}=(w(e): e\in E)$. The {\em chromatic index} 
of $(G,{\bfm w})$, denoted by $\chi'_w(G)$, is the minimum number of matchings in $G$ such that
each edge $e$ is covered exactly $w(e)$ times by these matchings, and the {\em fractional chromatic index} of 
$(G, {\bfm w})$, denoted by $\chi^*_w(G)$, is the optimal value of the following linear program:    
\begin{center}
\begin{tabular}{ll}
\hbox{Minimize} \ \ \ & $ {\bm 1}^T {\bm x}$  \\
\hbox{\hskip 0.2mm subject to} & $ A{\bm x} = {\bm w}$  \\
& \hskip 3mm ${\bm x} \ge {\bm 0},$
\end{tabular}
\end{center}
where $A$ is again the edge-matching incidence matrix of $G$. Clearly, $\chi'_w(G)$ is the optimal value of
the corresponding integer program. Let $G_w$ be obtained from $G$ by replacing each edge $e$ with $w(e)$ 
parallel edges between the same ends. It is then routine to check that $\chi'_w(G)=\chi'(G_w)$ and 
$\chi^*_w(G)=\chi^*(G_w)$. So the inequality $\chi'_w(G)-\chi^*_w(G) \le 1$ holds for all nonnegative
integral weight functions ${\bm w}$, and hence FECP has a fascinating integer rounding property (see Schrijver 
\cite{Sc86,Sc}). (The LP relaxation (LP) of a combinatorial optimization problem (IP) is said to
have an integer rounding property if there exists an absolute positive constant $c$, such that the optimal 
value of LP differs from that of IP by at most $c$, for all weight functions. This property is of great interest
in integer programming and combinatorial optimization.) 

So far the most powerful and sophisticated technique for multigraph edge-coloring is the method of Tashkinov trees
\cite{T}, which generalizes the earlier methods of Vizing fans \cite{V} and Kierstead paths \cite{Ki}. (These
methods are named after the authors who invented them, respectively.) Most recent results described above 
Theorem \ref{ThmGS} were obtained by using the method of Tashkinov trees. As remarked by McDonald \cite{Mc}, the 
Goldberg-Seymour conjecture and ideas culminating in this method are two cornerstones in modern edge-coloring. Nevertheless, 
this method suffers some theoretical limitation when applied to prove the conjecture; see 
Asplund and McDonald \cite{AM} for detailed information. Despite various attempts to extend the Tashkinov trees (see, for 
instance, \cite{GT,GG,CYZ,S,SSTF}), the difficulty encountered by the method remains unresolved and, undesirably, another
problem emerges: it becomes very difficult to preserve the structure of an extended Tashkinov tree under Kempe changes 
(the most useful tool in edge-coloring theory). In this paper we introduce a new type of extended Tashkinov trees and 
develop an effective control mechanism over Kempe changes, which can overcome all the aforementioned difficulties.   
The reader is referred to Chen and Jing \cite{GG} for a prototype of this control mechanism and its role
in the derivation of the best approximate result on the Goldberg-Seymour conjecture presently available.

\vskip 2mm 

The remainder of this paper is organized as follows. In Section 2, we introduce some basic concepts and 
techniques of edge-coloring theory, and exhibit some important properties of stable colorings. 
In Section 3, we define the extended Tashkinov trees to be employed in subsequent proof, and give an outline 
of our proof strategy. In Section 4, we establish some auxiliary results on the extended Tashkinov trees
and stable colorings, which ensure that this type of trees is preserved under some restricted Kempe changes.  In Section 5, 
we develop an effective control mechanism over Kempe changes, the so-called good hierarchy of an extended 
Tashkinov tree, based on its prototype introduced by Chen and Jing in \cite{GG} (see Condition R2 therein).
In Section 6, we derive some properties satisfied by the good hierarchies introduced in the preceding section. 
In Section 7, we present the last step of our proof based on these good hierarchies. 

\section{Preliminaries}

This section presents some basic definitions, terminology, and notation used in our paper, along with some important 
properties and results.

\subsection{Terminology and Notation}

Let $G=(V,E)$ be a multigraph.  For each $X \subseteq V$, let $G[X]$ denote the subgraph of $G$ induced by $X$, and let 
$G-X$ denote $G[V-X]$; we write $G-x$ for $G-X$ if $X=\{x\}$. Moreover, we use $\partial(X)$ to denote the set of all 
edges with precisely one end in $X$, and write $\partial(x)$ for $\partial(X)$ if $X=\{x\}$. For each pair $x, y 
\in V$, let $E(x,y)$ denote the set of all edges between $x$ and $y$. As it is no longer appropriate to represent 
an edge $f$ between $x$ and $y$ by $xy$ in a multigraph, we write $f \in E(x,y)$ instead. For each subgraph $H$ of $G$, 
let $V(H)$ and $E(H)$ denote the vertex set and edge set of $H$, respectively, let $|H|=|V(H)|$, and let $G[H]=G[V(H)]$ and $\partial(H) =\partial(V(H))$.     

Let $e$ be an edge of $G$. A {\em tree-sequence}\label{tsequence} with respect to $G$ and $e$ is a sequence 
$T=(y_0,e_1,y_1, \ldots, e_p, y_p)$ with $p\ge 1$, consisting of distinct edges $e_1,e_2, \ldots, e_p$ and 
distinct vertices $y_0,y_1, \ldots, y_p$, such that $e_1=e$ and each edge $e_j$ with $1\le j \le p$ is between 
$y_j$ and some $y_i$ with $0\le i <j$.  Given a tree-sequence $T=(y_0,e_1,y_1, \ldots, e_p, y_p)$, we can naturally 
associate a linear order $\prec$ with its vertices, such that $y_i \prec y_j$ if $i<j$. We
write $y_i \preceq y_j$ if $i\le j$.  This linear order will be used repeatedly in subsequent sections.
For each vertex $y_j$ of $T$ with $j \ge 1$, let $T(y_j)$ denote $(y_0,e_1,y_1, \ldots, e_j, y_j)$. Clearly, 
$T(y_j)$ is also a tree-sequence with respect to $G$ and $e$. We call $T(y_j)$ the {\em segment}\label{segment} of $T$ induced
by $y_j$. Let $T_1$ and $T_2$ be two tree-sequences with respect to $G$ and $e$. We write $T_2-T_1$ for 
$E[T_2]-E[T_1]$, write $T_1 \subseteq T_2$ if $T_1$ is a segment of $T_2$, and write $T_1 \subset T_2$ if $T_1$ 
is a proper segment of $T_2$; that is, $T_1 \subseteq T_2$ and $T_1 \ne T_2$. 

A {\em $k$-edge-coloring}\label{kedgecoloring} of $G$ is an assignment of $k$ colors, $1,2, \ldots, k$, to the edges of $G$ 
so that no two adjacent edges have the same color. By definition, the chromatic index $\chi'(G)$ of $G$ is the minimum $k$ 
for which $G$ has a $k$-edge-coloring.  We use $[k]$ to denote the color set $\{1,2, \ldots, k\}$,
and use ${\cal C}^k(G)$ to denote the set of all $k$-edge-colorings of $G$. Note that every $k$-edge-coloring
of $G$ is a mapping from $E$ to $[k]$.

Let $\varphi$ be a $k$-edge-coloring of $G$. For each $\alpha \in [k]$, the edge set $E_{\varphi, \alpha}=\{e\in E:\, 
\varphi(e)=\alpha\}$ is called a {\em color class}\label{colorclass}, which is a matching in $G$.  For any two distinct colors 
$\alpha$ and $\beta$ in $[k]$, let $H$ be the spanning subgraph of $G$ with $E(H)=E_{\varphi, \alpha} 
\cup E_{\varphi, \beta}$. Then each component of $H$ is either a path or an even cycle; we refer to such a component 
as an $(\alpha, \beta)$-{\em chain} with respect to $\varphi$, and also call it an $(\alpha, \beta)$-{\em path} \label{alphabetapath}
(resp. $(\alpha, \beta)$-{\em cycle}) if it is a path (resp. cycle). Possibly a component of $H$ is an isolated
vertex. We use $P_v(\alpha, \beta, \varphi)$ to denote the unique $(\alpha, \beta)$-chain containing the vertex $v$. 
Clearly, for any two distinct vertices $u$ and $v$, $P_u(\alpha, \beta, \varphi)$ and $P_v(\alpha, \beta, \varphi)$ 
are either identical or vertex-disjoint. Let $C$ be an $(\alpha, \beta)$-chain with respect to $\varphi$, and let 
$\varphi'$ be the $k$-edge-coloring arising from $\varphi$ by interchanging $\alpha$ and $\beta$ on $C$. We say 
that $\varphi'$ is obtained from $\varphi$ by {\em recoloring} $C$, and write $\varphi'=\varphi /C$. This 
operation is called a {\em Kempe change}\label{kempe}. We also say that this Kempe change is {\em rooted} at $v$ if 
it has degree one in $C$. 

Let $F$ be an edge subset of $G$. As usual, $G-F$ stands for the multigraph obtained from $G$ by deleting all
edges in $F$; we write $G-f$ for $G-F$ if $F=\{f\}$. Let $\pi \in {\cal C}^k(G-F)$. For each $K \subseteq E$, 
define $\pi(K)=\cup_{e\in K-F} \, \pi(e)$. For each $v \in V$, define 
\[\pi(v)=\pi(\partial(v)) \hskip 2mm {\rm and} \hskip 2mm \overline{\pi}(v)=[k]-\pi(v).\] 
We call $\pi(v)$ the set of colors {\em present} at $v$ and call $\overline{\pi}(v)$ the set of colors {\em missing}\label{missingcolor}
at $v$.  For each $X\subseteq V$, define 
\[\overline{\pi}(X)= \cup_{v\in X} \, \overline{\pi}(v).\] 
We call $X$ {\em elementary}\label{elementaryset} with respect to $\pi$ if $\overline{\pi}(u) \cap \overline{\pi}(v)
=\emptyset$ for any two distinct vertices $u, v\in X$. We call $X$ {\em closed}\label{closedsets} with respect to $\pi$ if
$\pi(\partial(X))\cap \overline{\pi}(X)=\emptyset$; that is, no missing color of $X$ appears on the 
edges in $\partial(X)$.  Furthermore, we call $X$ {\em strongly closed}\label{sclosed} with respect to $\pi$ if $X$ is closed
with respect to $\pi$ and $\pi(e) \ne \pi(f)$ for any two distinct colored edges $e, f \in 
\partial(X)$.  For each subgraph $H$ of $G$, write $\overline{\pi}(H)$ for $\overline{\pi}(V(H))$,
and write ${\pi}\langle H \rangle$ for ${\pi}(E(H))$. Moreover, define 
\[\partial_{\pi, \alpha}(H)=\{e\in \partial(H): \pi(e)=\alpha\},\] 
and define 
\[I[\partial_{\pi, \alpha}(H)]=\{v\in V(H): v \hskip 2mm \mbox{is incident with an edge in} \hskip 2mm
\partial_{\pi, \alpha}(H)\}.\] 
For an edge $e\in \partial(H)$, we call its end in (resp. outside) $H$ the {\em in-end} (resp. {\em out-end})
relative to $H$.  Thus $I[\partial_{\pi, \alpha}(H)]$ consists of all in-ends (relative to $H$) of edges in 
$\partial_{\pi, \alpha}(H)$. If $|\partial_{\pi, \alpha}(H)| \ge 2$, we call $\alpha$ a {\em defective color}\label{defectivecolor} 
of $H$ with respect to $\pi$, call each edge in $\partial_{\pi, \alpha}(H)$ a {\em defective edge}\label{defectiveedge} of $H$ with 
respect to $\pi$, and call each vertex in $I[\partial_{\pi, \alpha}(H)]$ a {\em defective vertex}\label{defectivevertex} of $H$ with 
respect to $\pi$. A color $\alpha \in \overline{\pi}(H)$ is called {\em closed} in $H$ under $\pi$ 
if $\partial_{\pi, \alpha}(H)=\emptyset$.  For convenience, we say that $H$ is {\em closed} (resp. {\em strongly closed}) 
with respect to $\pi$ if $V(H)$ is closed (resp. strongly closed) with respect to $\pi$. Let $\alpha$ and $\beta$ be 
two colors that are not assigned to $\partial(H)$ under $\pi$. We use $\pi/(G-H, \alpha, \beta)$ to denote the 
coloring $\pi'$ obtained from $\pi$ by interchanging $\alpha$ and $\beta$ in $G-V(H)$. Since 
$\pi$ belongs to ${\cal C}^k(G-F)$, so does $\pi'$.

\subsection{Elementary Multigraphs}

Let $G=(V,E)$ be a multigraph. We call $G$ an {\em elementary multigraph}\label{elementarymulti} if $\chi'(G)=\lceil \Gamma(G) \rceil$.
With this notion, Conjecture \ref{GS} can be rephrased as follows.

\begin{conjecture} \label{GS2}
Every multigraph $G$ with $\chi'(G) \ge \Delta(G)+2$ is elementary.
\end{conjecture} 

Recall that $G$ is critical if $\chi'(H)<\chi'(G)$ for any proper subgraph $H$ of $G$. As pointed out by Stiebitz 
et al. \cite{SSTF} (see page 7), for a proof of Conjecture \ref{GS2}, it suffices to consider critical multigraphs. To 
see this, let $G$ be an arbitrary multigraph with $\chi'(G) \ge \Delta(G)+2$. Then $G$ contains a critical multigraph 
$H$ with $\chi'(H)=\chi'(G)$, which implies that $\chi'(H) \ge \Delta(H)+2$. Note that if $H$ is elementary, 
then so is $G$, because $\lceil \Gamma(G) \rceil \le \chi'(G) = \chi'(H) = \lceil \Gamma(H) \rceil \le 
\lceil \Gamma(G) \rceil$. Thus both inequalities hold with equalities, and hence $\chi'(G)=\lceil \Gamma(G) \rceil$.

To prove Conjecture \ref{GS}, we shall actually establish the following statement.

\begin{theorem} \label{ThmGS2} 
Every critical multigraph $G$ with $\chi'(G) \ge \Delta(G)+2$ is elementary.
\end{theorem}

In our proof we shall appeal to the following theorem, which reveals some intimate connection between elementary 
multigraphs and elementary sets. This result is implicitly contained in Andersen \cite{An77} and Goldberg \cite{G84}, 
and explicitly stated in Stiebitz et al. \cite{SSTF}  (see Theorem 1.4 on page 8).

\begin{theorem} \label{egraph2}
Let $G=(V,E)$ be a multigraph with $\chi'(G)=k+1$ for an integer $k \ge \Delta(G)+1$. If $G$ is critical, then the 
following conditions are equivalent:
\begin{itemize}
\vspace{-2mm}
\item[(i)] $G$ is an elementary multigraph.
\vspace{-2mm}
\item[(ii)] For each edge $e\in E$ and each coloring $\varphi \in {\cal C}^k(G-e)$, the vertex set $V$ is elementary
with respect to $\varphi$.
\vspace{-2mm}
\item[(iii)] There exists an edge $e\in E$ and a coloring $\varphi \in {\cal C}^k(G-e)$, such that the vertex set $V$ 
is elementary with respect to $\varphi$.
\vspace{-2mm}
\item[(iv)] There exists an edge $e\in E$, a coloring $\varphi \in {\cal C}^k(G-e)$, and a subset $X$ of $V$, 
such that $X$ contains both ends of $e$, and $X$ is elementary as well as strongly closed with respect to $\varphi$.
\end{itemize}
\end{theorem} 

\subsection{Stable Colorings}

\vskip 2mm
In this subsection, we assume that $T$ is a tree-sequence with respect to a multigraph $G=(V,E)$ and an edge $e$, 
$C$ is a subset of $[k]$, and $\varphi$ is a coloring in ${\cal C}^k(G-e)$, where $k \ge \Delta(G)+1$. We say
that an edge $f$ of $G$ is {\em incident} to $T$ if at least one end of $f$ is contained in $T$; this definition
applies to edges of $T$ as well.  Since our proof consists of a sophisticated sequence of Kempe changes, 
the concept of stable coloring introduced below will be employed to preserve some important coloring properties 
of $T$, such as, among others, the color on each edge and the set of colors missing at each vertex. Usually, 
$C$ is the set of colors assigned to $E(T)$ but not missing at any vertex of $T$. 

To be specific, a coloring $\pi \in {\cal C}^k(G-e)$ is called a $(T, C, \varphi)$-{\em stable coloring}\label{tcstable} 
if the following two conditions are satisfied:
\begin{itemize}
\vspace{-2mm}
\item[$(i)$] $\pi(f) = \varphi(f)$ for any $f\in E$ incident to $T$ with $\varphi(f)\in \overline{\varphi}(T)\cup C$; and
\vspace{-2mm}
\item[$(ii)$] $\overline{\pi} (v) = \overline{\varphi}(v)$ for any $v\in V(T)$.
\end{itemize}

\vspace{-2mm}
By convention, $\pi(e)=\varphi(e)=\emptyset$. The following lemma gives an equivalent definition of stable colorings.

\begin{lemma} \label{sc1}
A coloring $\pi \in {\cal C}^k(G-e)$ is $(T, C, \varphi)$-stable iff $\pi(f)=\varphi(f)$ for any $f\in E$ incident to $T$ with 
$\varphi(f)\in \overline{\varphi}(T)\cup C$ or $\pi(f)\in \overline{\varphi}(T)\cup C$.
\end{lemma}

{\bf Proof.}  Let $(i')$ stand for the condition specified in the ``if" part. We propose to show that $(i')$ holds iff
both $(i)$ and $(ii)$ hold. 

Trivially, $(i')$ implies $(i)$. If there exists $v\in V(T)$ such that $\overline{\pi} (v) \ne \overline{\varphi}(v)$, 
then some edge $f$ incident to $v$ satisfies $\pi(f) \in \overline{\varphi}(v)$ because $|\overline{\pi} (v)|=|\overline{\varphi}(v)|$.
From $(i')$ we deduce that $\pi(f)=\varphi(f)$ and hence $\varphi(f) \in \overline{\varphi}(v)$, a contradiction. So $(i')$ implies 
$(ii)$ as well.  

Conversely, let $f \in E$ be an arbitrary edge incident to $T$ with $\pi(f)\in \overline{\varphi}(T)\cup C$. 
We claim that $\varphi(f) = \pi(f)$. Assume the contrary: $\varphi(f) \ne \pi(f)$. Let $v\in V(T)$ be an end of $f$. By $(ii)$, 
we have $\overline{\pi} (v) = \overline{\varphi}(v)$. So ${\pi} (v) = {\varphi}(v)$ and hence there exists an edge 
$g \in \partial(v)-\{f\}$ with $\varphi(g) = \pi(f)$. It follows that $\varphi(g) \in \overline{\varphi}(T)\cup C$.
By $(i)$, we obtain $\pi(g)=\varphi(g)$, which implies $\pi(f)=\pi(g)$, contradicting the hypothesis that 
$\pi \in {\cal C}^k(G-e)$. Our claim asserts that $\varphi(f)=\pi(f)$ for any $f\in E$ incident to $T$ with $\pi(f)\in 
\overline{\varphi}(T)\cup C$. Combining this with $(i)$, we conclude that $(i')$ holds. \qed

\vskip 3mm

From the definition and Lemma \ref{sc1} we see that the following statements hold for a $(T, C, \varphi)$-stable 
coloring $\pi$: 

$\bullet$ if $T' \subseteq T$ and $\overline{\varphi}(T')\cup C' \subseteq \overline{\varphi}(T)\cup C$, then $\pi$ is also $(T', C', \varphi)$-stable;

$\bullet$ if a color $\alpha \in \overline{\varphi}(T)$ is closed in $T$ under $\varphi$, then it is also closed in $T$ under $\pi$; and

$\bullet$ if ${\varphi}\langle T \rangle \subseteq \overline{\varphi}(T)\cup C$, then $\pi(f) = \varphi(f)$ for all edges $f$ on $T$.    

\vskip 2mm
Let us derive some further properties satisfied by stable colorings.

\begin{lemma} \label{sc2}
Being $(T, C, \cdot )$-stable is an equivalence relation on ${\cal C}^k(G-e)$. 
\end{lemma}

{\bf Proof.} From Lemma \ref{sc1} and the above condition $(ii)$,  it is clear that being $(T, C, \cdot )$-stable is reflexive,
symmetric, and transitive. So it defines an equivalence relation on ${\cal C}^k(G-e)$. \qed 

\vskip 3mm

\begin{lemma}\label{zang1}
Suppose $T$ is closed but not strongly closed with respect to $\varphi$, with $|V(T)|$ odd. If $\pi$ is 
a $(T, C, \phiv)$-stable coloring, then $T$ is also closed but not strongly closed with respect to $\pi$.
\end{lemma}

{\bf Proof.} Let $X=V(T)$ and let $t$ be the size of the set $[k]-\overline{\varphi}(X)$. By hypotheses, $|V(T)|$ 
is odd and $T$ is not strongly closed with respect to $\varphi$. Thus under the coloring $\varphi$ each color in $[k]-
\overline {\varphi}(X)$ is assigned to at least one edge in $\partial(T)$, and some color in $[k]-\overline{\varphi}(X)$ is assigned 
to at least two edges in $\partial(T)$. It follows that $|\partial(T)|\ge t+1$. Since $\pi$ is a $(T, C, \phiv)$-stable coloring, 
from Lemma \ref{sc1} and the above condition $(ii)$ we deduce that $T$ is closed with respect to $\pi$ and that $\overline{\pi}(X)=\overline{\varphi}(X)$ 
(so $[k]-\overline{\pi}(X)$ is also of size $t$). As only colors in $[k]-\overline{\pi}(X)$ can be assigned to 
edges in $\partial(T)$ under $\pi$, some of these colors is used at least twice by the Pigeonhole Principle. Hence $T$ is not strongly 
closed with respect to $\pi$. \qed 

\vskip 3mm
Let $P$ be a path in $G$ whose edges are colored alternately by $\alpha$ and $\beta$ in $\varphi$,
with $|P|\ge 2$, and let $u$ and $v$ be the ends of $P$ with $v \in V(T)$. We say that $P$ is a $T$-{\em exit path}\label{exitpath} 
with respect to $\varphi$ if $V(T) \cap V(P)=\{v\}$ and $\overline{\varphi}(u) \cap \{\alpha,\beta\} \ne \emptyset$; 
in this case, $v$ is called a $(T,\varphi,\{\alpha,\beta\})$-{\em exit}\label{exit} and $P$ is also called a 
$(T,\varphi,\{\alpha,\beta\})$-{\em exit path}\label{tphiexitpath}. Note that possibly $\overline{\varphi}(v) \cap \{\alpha,\beta\}=\emptyset$;
now $P$ is a proper subpath of an $(\alpha,\beta)$-path.    

\begin{lemma}\label{LEM:extable}
Suppose $T$ is closed with respect to $\varphi$, and $f\in E(u,v)$ is an edge in $\partial (T)$ with $v \in V(T)$. If there 
exists a $(T, C\cup \{\phiv(f)\}, \phiv)$-stable coloring $\pi$, such that  $\overline{\pi}(u) \cap \overline{\pi}(T)
\ne \emptyset$, then for any $\alpha \in \phibar(v)$ there exists a $(T, C\cup \{\phiv(f)\}, \phiv)$-stable coloring $\sigma$, 
such that $v$ is a $(T, \sigma, \{\alpha, \phiv(f)\})$-exit. 
\end{lemma}
 
{\bf Proof.} Let $\beta \in \overline{\pi}(u) \cap \overline{\pi}(T)$. By the definition of stable coloring, $\beta \in \phibar(T)$.
Since both $\alpha$ and $\beta$ are closed in $T$ under $\varphi$, they are also closed in $T$ under $\pi$ by Lemma \ref{sc1}. Define 
$\sigma=\pi/ (G-T, \alpha,\beta)$. Clearly, $\sigma$ is a $(T, C\cup \{\phiv(f)\}, \pi)$-stable coloring. By Lemma \ref{sc2}, $\sigma$ 
is also a $(T, C\cup \{\phiv(f)\}, \phiv)$-stable coloring. Since $P_v(\alpha, \phiv(f), \sigma)$ consists of a single edge $f$, it 
is a $T$-exit path with respect to $\sigma$. Hence $v$ is a $(T, \sigma, \{\alpha, \phiv(f)\})$-exit. \qed

\subsection{Tashkinov Trees}

A multigraph $G$ is called $k$-{\em critical}\label{kcritical} if it is critical and $\chi'(G)=k+1$. Throughout this paper, by 
a $k$-{\em triple}\label{ktriple} we mean a $k$-critical multigraph $G=(V,E)$, where $k \ge \Delta(G)+1$,  together with an 
uncolored edge $e\in E$ and a coloring $\varphi \in {\cal C}^k(G-e)$; we denote it by $(G,e, \varphi)$.  

Let $(G,e, \varphi)$ be a $k$-triple. A {\em Tashkinov tree}\label{tashkinovtree} with respect to $e$ and $\varphi$ is a 
tree-sequence $T=(y_0,e_1,y_1, \ldots, e_p, y_p)$ with respect to $G$ and $e$, such that 
for each edge $e_j$ with $2\le j \le p$, there is a vertex $y_i$ with $0 \le i <j$ satisfying
$\varphi(e_j) \in  \overline{\varphi}(y_i)$.

\vskip 2mm
The following theorem is due to Tashkinov \cite{T}; its proof can also be found in Stiebitz et al. 
\cite{SSTF}  (see Theorem 5.1 on page 116).

\begin{theorem} \label{TashTree}
Let $(G,e, \varphi)$ be a $k$-triple and let $T$ be a Tashkinov tree with respect to $e$ and $\varphi$.
Then $V(T)$ is elementary with respect to $\varphi$. 
\end{theorem} 

Let $G=(V,E)$ be a $k$-critical multigraph $G$ with $k \ge \Delta(G)+1$. For each edge $e\in E$ and each 
coloring $\varphi \in {\cal C}^k(G-e)$, there is at least one Tashkinov tree $T$ with respect to $e$ and $\varphi$. 
The {\em Tashkinov order}\label{torder} of $G$, denoted by $t(G)$, is the largest number of vertices contained in 
such a Tashkinov tree over all $e$ and $\varphi \in {\cal C}^k(G-e)$.  Scheide \cite{S} (see Proposition 4.5) 
has established the following result, which will be employed in our proof.

\begin{theorem} \label{ThmScheide}
Let $G$ be a critical multigraph $G$ with $\chi'(G) \ge \Delta(G)+2$. If $t(G)<11$, then $G$ is an elementary
multigraph. 
\end{theorem}

Tashkinov trees have been used successfully to establish various approximate results on Conjecture \ref{GS}. 
The crux of this approach is to capture the density $\Gamma(G)$ by exploring a sufficiently large Tashkinov tree 
(see Theorem \ref{TashTree}). However, this target may become unreachable when the upper bound on $\chi'(G)$
(one wishes to derive) gets close to $\chi'(G)$, even if we allow for an unlimited number of Kempe changes; such 
an example has been found by Asplund and McDonald \cite{AM}. To carry out a proof of Conjecture \ref{GS}, we 
introduce a type of extended Tashkinov trees in this paper by  the procedure described below.

\begin{definition} \label{TAA}
{\rm Given a $k$-triple $(G,e, \varphi)$ and a tree-sequence $T$ with respect to $G$ and $e$, we say that a tree-sequence $(T,f,y)$ is obtained from $T$
by a {\em Tashkinov augmentation}\label{TA} (TA) under $\varphi$ if $\varphi(f)\in\phibar(T)$, one end $x$ of $f$ is contained in $T$, and 
the other end $y$ of $f$ is outside $T$. A {\em Tashkinov augmentation algorithm} (TAA)\label{TAAalgorithm} consists of a sequence of TAs under the same
edge coloring. We call a tree-sequence $T'$ a {\em closure}\label{closure} of $T$ under $\varphi$ if $T'$ arises from $T$ by 
TAA and cannot grow further by TA under $\varphi$ (equivalently, $T'$ has become closed).} 
\end{definition}

So a Tashkinov tree with respect to $e$ and $\varphi$ is a tree-sequence obtained from $(y_0,e,y_1)$ by TAA, where
$y_0$ and $y_1$ are two ends of $e$. We point out that, although there might be several ways to construct a closure of $T$ under $\varphi$, the 
vertex set of these closures is unique.   
 
In the next section we shall give a detailed description of an algorithm for constructing the aforementioned extended Tashkinov trees, and present the main
result of this paper, which implies Theorem~\ref{ThmGS2}. In view of Theorem \ref{egraph2}, to prove Conjecture \ref{GS}, we may turn to finding a strongly 
closed elementary tree-sequence. Thus in our algorithm we introduce three types of extensions from a {\em closed} elementary tree-sequence $T$ (say, 
a closed Tashkinov tree), which has defective edges under a coloring $\varphi$ (so $T$ is not strongly closed). Specifically, we consider the maximum defective vertex $v$ in the order $\prec$ 
over all $(T,C,\varphi)$-stable colorings (here $C$ contains all colors used to construct $T$ but are not missing at vertices in $T$, so $C=\emptyset$ when 
$T$ is a Tashkinov tree). Let $f$ be a defective edge incident to $v$ under $\varphi$, with $f\in E(u,v)$ and $\varphi(f)=\delta$. 

If $T\cup \{u\}$ is elementary under all $(T,C\cup\{\delta\},\varphi)$-stable colorings, we add this edge $f$ to $T$ and extend the resulting tree-sequence 
using TAA. This first type of extensions is called series extension (SE) in our algorithm.

Otherwise, we pick a color $\gamma$ in $\phibar(v)$; Lemma \ref{LEM:extable} guarantees the existence of a stable coloring such that $v$ is the only common 
vertex of a $(\gamma,\delta)$-path $P$ and $T$. We then perform the second type of extensions, called parallel extension (PE). Each PE is followed by a sequence of
the third type of extensions, called revisiting extension (RE) and performed whenever possible. During PE, we switch colors along $P$ and apply TAA to $T$, 
which is no longer closed now. During REs, we repeatedly add an edge on some $(\gamma,\delta)$-cycle intersecting $T$ (the original $T$ before the previous 
PE) and apply TAA, until the vertices of all $(\gamma,\delta)$-cycles intersecting $T$ are contained in the resulting tree-sequence.  
We may view each PE and its succeeding REs as a whole, in which PE is the primary extension while REs are auxiliary extensions.   

Applying the above three extensions to closed elementary tree-sequences recursively, we shall end up with a strongly closed elementary tree-sequence, 
thereby proving the Goldberg-Seymour conjecture.  The elementary property of such tree-sequences will be established in Sections 5, 6 and 7 by
further developing techniques introduced in Chen and Jing \cite{GG}, which essentially allow us to prove that if the elementary property of a closed tree-sequence $T$
is preserved under adding a vertex using Algorithm 3.1, then this property can be extended even further to a closure of the resulting tree-sequence. To the 
best of our knowledge, no previous investigation has ever employed these three extensions; the weaker extension used in~\cite{GG} (see Condition R1 therein) 
can only lead to some partial results despite several new techniques originating from this study. 

It is worthwhile pointing out that,  when encountering a non-elementary tree-sequence, almost all previous methods proceed by reducing the size of the non-elementary tree-sequence 
and eventually reach a contradiction by coloring the uncolored critical edge $e$.  However, our algorithm goes the other way around: when $T\cup\{u\}$ is 
not elementary, it modifies the coloring and employs PE to construct a larger elementary tree-sequence while avoiding the edge $f$. Essentially our proof only 
requires the edge $e$ to be critical as potential non-critical edges like $f$ are bypassed.

\section{Extended Tashkinov Trees}

The purpose of this section is to present extended Tashkinov trees to be used in our proof and to give an outline of
our proof strategy.

Given a $k$-triple $(G, e, \phiv)$, a {\em Tashkinov series}\label{taskinovseries} constructed from it is a series of tuples $(T_n, \phiv_{n-1}, \\ S_{n-1}, F_{n-1}, 
\Theta_{n-1})$ for $n=1,2, \ldots $ output by the following algorithm, where $T_n+f_n$ stands for the tree-sequence 
augmented from $T_n$ by adding an edge $f_n$, and the definition of $(T,\varphi,\{\alpha,\beta\})$-{\em exit}  can be found in the paragraph right above Lemma~\ref{LEM:extable}.

To help gain a clearer picture of our algorithm, we intentionally use a descriptive language and include Iteration 1, although it is contained in the general Iteration $n$.

\newpage
\vskip 4mm
\noindent {\bf Algorithm 3.1}

\vskip 2mm

\noindent {\bf Iteration 0.} Let $(T_1, \phiv_0, S_0, F_0, \Theta_0)$ be the initial tuple, such that $\varphi_0=\varphi$, $T_1$ is a closure of $e$ under 
$\varphi_0$ (which is a closed Tashkinov tree with respect to $e$ and $\varphi_0$), and $S_0 = F_0 = \Theta_0=\emptyset$.  
\vskip 2mm
\noindent {\bf Iteration 1.}
If $T_1$ is strongly closed with respect to $\phiv_{0}$, stop. Else, we construct the tuple $(T_2, \phiv_1, S_1, F_1, \Theta_1)$ as follows. 
Set $D_{0}=\emptyset$.
\begin{itemize}
\vspace{-2mm}	
\item Let $v_1$ be the maximum defective vertex\footnote{For each  $(T_1, D_{0}, \phiv_{0})$-stable coloring $\pi$, let $v_{\pi}$ be the 
largest defective vertex of $T_1$ with respect to $\pi$ in the order $\prec$. Then $v_1$ is the largest vertex among all these vertices 
$v_{\pi}$ in the order $\prec$. By definition, $v_1=v_{\pi_0}$. } in the order $\prec$ over all $(T_1, D_{0}, \phiv_{0})$-stable colorings, let 
$\pi_{0}$ be a corresponding coloring, let $f_1$ be a defective edge (of $T_1$ with respect to $\pi_0$) incident to $v_1$, let $u_1$ be 
the other end of $f_1$, and let $\delta_1=\pi_{0}(f_1)$. 
	\begin{itemize}
	\vspace{-2mm}	
		\item If for every $(T_1, D_{0}\cup\{\delta_1\}, \pi_{0})$-stable coloring $\pi$,  we have $\overline{\pi}(u_1) \cap \overline{\pi} (T_1)=\emptyset$, apply {\bf SE} with $n=1$.
		\item Else, let $\gamma_1$ be an arbitrary color in $\overline{\pi}_0(v_1)$ and let $\pi_{0}'$ be a $(T_1, D_{0}\cup\{\delta_1\}, 
		\pi_{0})$-stable coloring such that $v_1$ is a $(T_1, \pi_{0}', \{\gamma_1, \delta_1\})$-exit (such $\pi_{0}'$ exists by   
		Lemma \ref{LEM:extable}), apply {\bf PE} with $n=1$.
	\end{itemize}
\end{itemize}

\vspace{-1mm}

\noindent {\bf Iteration} ${\bm n}$. If $T_n$ is strongly closed with respect to $\phiv_{n-1}$, stop. Else, we construct the tuple 
$(T_{n+1}, \phiv_n, S_n, F_n, \Theta_n)$ as follows. Set $D_{n-1}=\cup_{i \le n-1} \, S_i - \phibar_{n-1}(T_{n-1})$ (so $D_0=\emptyset$).
\begin{itemize}
	
	\vspace{-2mm}	
	\item If there is a subscript $h \le n-1$ with $\Theta_h =PE$ and $S_h =\{\delta_h, \gamma_h\}$, such that 
	some $(\gamma_h, \delta_h)$-cycle $O$ with respect to $\phiv_{n-1}$ intersects both $V(T_h)$ and $V(G)-V(T_n)$, 
	apply {\bf RE}. (Note that the existence of cycle $O$ implies $\Theta_i=RE$ for all $i$ with $h+1 \le i \le n-1$, if any, and that
	this case cannot occur when $n=1$.)
	\vspace{-2mm}	
	\item Else, let $v_n$ be the maximum defective vertex\label{maximumdefective} in the order $\prec$ over all
	$(T_n, D_{n-1}, \phiv_{n-1})$-stable colorings, let $\pi_{n-1}$ be a corresponding coloring, let $f_n$ be a defective edge (of $T_n$ with respect to $\pi_{n-1}$)
	incident to $v_n$, 	let $u_n$ be the other end of $f_n$, and let $\delta_n=\pi_{n-1}(f_n)$. 
	\begin{itemize}
		\vspace{-2mm}	
		\item If for every $(T_n, D_{n-1}\cup\{\delta_n\}, \pi_{n-1})$-stable coloring $\pi$,  we have $\overline{\pi}(u_n) \cap \overline{\pi} 
		(T_n)=\emptyset$, apply {\bf SE}.
		\item Else, pick a color $\gamma_n$ in $\overline{\pi}_{n-1}(v_n)$ as follows. If $v_n=v_i$ for some $1\leq i<n$ with $\Theta_i={PE}$, let $n'$ be the 
		largest such $i$ and let $\gamma_n=\delta_{n'}$. Otherwise, let $\gamma_n$ be an arbitrary color in $\overline{\pi}_{n-1}(v_n)$.  Let $\pi_{n-1}'$ be an arbitrary
		$(T_n, D_{n-1}\cup\{\delta_n\}, \pi_{n-1})$-stable coloring so that $v_n$ is a $(T_n, \pi_{n-1}', \{\gamma_n, \delta_n\})$-exit (such $\pi_{n-1}'$ exists 
		by Lemma \ref{LEM:extable}), apply {\bf PE}. 
		
	\end{itemize}
\end{itemize}

\noindent {\bf RE.}\label{REextension} Let $f_n$ be an  edge in $O \cap \partial(T_n)$ such that $O$ contains a path $L$
connecting $f_n$ and $V(T_h)$ with $V(L) \subseteq V(T_n)$. Let $\phiv_n = \phiv_{n-1}$ and $T_{n+1}$ be a closure of 
$T_n+f_n$ under $\phiv_{n}$.  Set $\delta_n=\delta_h$, $\gamma_n=\gamma_h$, $S_n=\{\delta_n,\gamma_n\}$, $F_n = \{ f_n\}$, and $\Theta_n = RE$. 
We call this extension a {\bf revisiting extension} (RE), call $f_n$ an {\em RE connecting edge}\label{connectingedge}, and call $\delta_n$ and $\gamma_n$ 
{\em connecting colors}\label{connectingcolor}. Let $v_n$ be the end of $f_n$ in $T_n$. Note that $v_n$ here is neither called an extension vertex nor called a supporting vertex. 

\vskip 2mm
\noindent {\bf SE.}\label{SEextension} Let $\phiv_n = \pi_{n-1}$ and let $T_{n+1}$ be a closure of $T_n+f_n$ under $\phiv_{n}$. Set $S_n= \{ \delta_n \}$,  $F_n=\{f_n\}$, and 
$\Theta_n = SE$. We call this extension a {\bf series extension} (SE), call $f_n$ an {\em SE connecting edge}, call $\delta_n$ a {\em connecting 
	color}, and call $v_n$ an {\em extension vertex}\label{extensionvertex}.

\vskip 2mm
\noindent {\bf PE.}\label{PEextension}  Let $\phiv_n  = \pi_{n-1}'/ P_{v_n}(\gamma_n, \delta_n, \pi_{n-1}')$. Note that $P_{v_n}(\gamma_n, \delta_n, \pi_{n-1}')\cap V(T_n)=\{v_n\}$ 
and $\delta_n\in\phibar_n(v_n)$ is a defective color of $T_n$. So $T_n$ is not closed under $\varphi_n$. Let $T_{n+1}$ be a closure of $T_n$ under $\phiv_n$. Set $S_n = \{\delta_n, \gamma_n\}$, $F_n=\{f_n\}$, and $\Theta_n =PE$. We call this extension a {\bf parallel extension} (PE), call $f_n$ a {\em PE 
connecting edge}, call $\delta_n$ and $\gamma_n$ {\em connecting colors}, and call $v_n$ a {\em supporting vertex}\label{supportingvertex}. As 
$f_n$ is the first edge along $P_{v_n}(\gamma_n, \delta_n, \pi_{n-1}')$ and is colored by $\gamma_n$ under $\varphi_n$, 
it is not necessarily contained in $T_{n+1}$.

\begin{figure}[htpb]
	\vspace{-1mm}
	\centerline{\includegraphics[width=8.58cm]{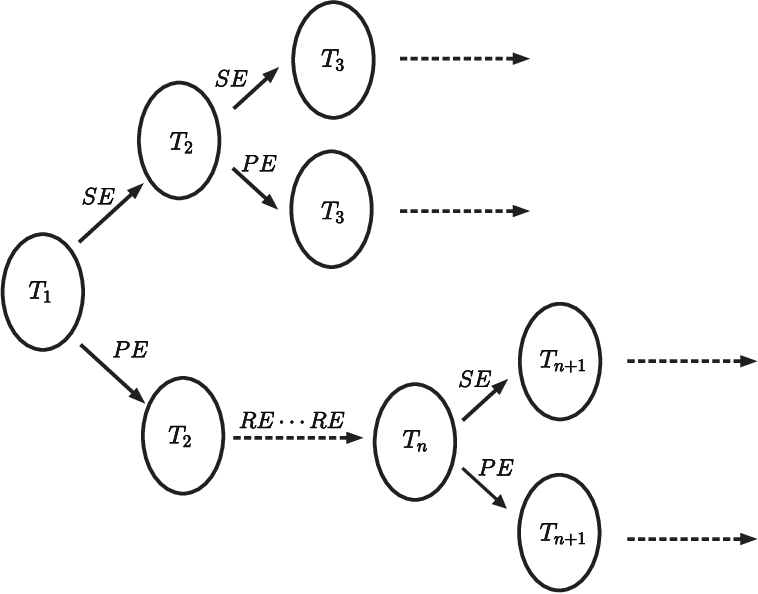}}
	\vspace{-1mm}\caption{Tashkinov series constructed by Algorithm 3.1}
		\label{abc}
\end{figure}

\vskip 3mm
Figure~\ref{abc} shows the possible choices of extensions used in the construction of a Tashkinov series.
\vspace{-3mm}

\begin{definition} \label{wz1}
{\rm Let $(G,e, \varphi)$ be a $k$-triple and let $\TT=\{(T_i, \phiv_{i-1}, S_{i-1}, F_{i-1}, \Theta_{i-1}): 	1\le i \le n+1\}$  be a Tashkinov 
series constructed from $(G,e, \varphi)$.  A tree-sequence $T$ is called an {\em extended Tashkinov 
tree} (ETT)\label{ETT} constructed from $\TT$ under $\varphi_n$ if $T_n\subset T\subseteq T_{n+1}$. 
We say that $\varphi_n$ is the {\em generating coloring}\label{generatingcoloring} of $T$. If the Tashkinov series $\TT$ is clear 
from the context, we may simply say that $T$ is an ETT under $\varphi_n$.}
\end{definition}

We shall mainly work on an ETT $T$ as defined above in the remainder of this paper. Such $T$ is not necessarily closed under 
$\varphi_n$, while $T_{n+1}$ is closed.

Throughout we reserve all symbols used for the same usage as in the algorithm. In particular, $D_i=\cup_{h\leq i} \, S_h-\phibar_i(T_i)$
for $i \ge 0$. To help understand the algorithm and ETTs better, let us make a few remarks and offer some simple observations.

{\bf (3.1)}  In our proof we shall always restrict our attention to the case when $|T_n|$ is odd. Suppose $T_n$ is closed but not 
strongly closed with respect to $\phiv_{n-1}$. Then, by Lemma \ref{zang1}, the same property holds for $T_n$ with respect to 
any $(T_n, D_{n-1}, \phiv_{n-1})$-stable coloring $\pi$. Let $v_{\pi}$ denote the largest defective vertex of $T_n$ with respect to $\pi$ 
in the order $\prec$. Note that $v_n$ involved in SE and PE is the largest vertex among all these vertices $v_{\pi}$ in the order $\prec$ 
and so it is uniquely determined by the triple $(T_n, D_{n-1}, \phiv_{n-1})$, while $f_n$ involved in each extension might be selected in several ways.  

{\bf (3.2)}  In the algorithm $\delta_n$ is a defective color of $T_n$ with respect to $\varphi_n$ when $\Theta_n =SE$ 
or $PE$ (as $|\partial_{ \pi_{n-1}, \delta_n}(T_n)|\ge 3$ when $|T_n|$ is odd), while $\gamma_n$ is a defective color 
of $T_n$ with respect to $\varphi_n$ when $\Theta_n =RE$.  Unlike PE or SE, $v_n$ involved in RE may not be a maximum defective 
vertex. Moreover, the set $D_{n-1}=\cup_{i \le n-1} \, S_i - \phibar_{n-1}(T_{n-1})$ is used to store colors employed in the construction of 
$T_n$ but not missing at any vertex of $T_{n-1}$ under $\varphi_{n-1}$.

{\bf (3.3)}  As described in the algorithm, after performing each PE, we grow the Tashkinov series by using RE, whenever possible. 
So revisiting extension (RE) has priority over both series and parallel extensions (SE and PE). If $\Theta_n=RE$, then all edges
in $O \cap \partial(T_h)$ are colored with $\delta_h$ and all edges in $O \cap \cup \partial(T_n)$ are colored with $\gamma_h$, 
because $\delta_h$ is a color missing at $v_h$ under $\phiv_{h}=\varphi_{n-1}$ (thereby $v_h$ is outside $O$), the only edge
$f_h$ in $\partial_{\varphi_h, \gamma_h}(T_h)$ (see Lemma \ref{hku}(v) to be proved) is adjacent to $v_h$, and $T_n$ is closed with 
respect to $\varphi_{n-1}$. Hence $O$ contains at least one edge colored with $\gamma_h$ in $G[T_h]$, at least two boundary edges of 
$T_h$ colored with $\delta_h$, and at least two boundary edges of $T_n$ colored with $\gamma_h$. 

RE is illustrated by the following figure, in which $O \cap \partial(T_n)$ contains four edges colored with $\gamma_h$: both the top and
bottom edges are good candidates for $f_n$, but neither of the middle edges can serve this purpose, because $O$ contains no path $L$ 
connecting them and $T_h$, with $V(L) \subseteq V(T_n)$.

\begin{figure}[htpb]
	\vspace{-1mm}
	\centerline{\includegraphics[width=7cm]{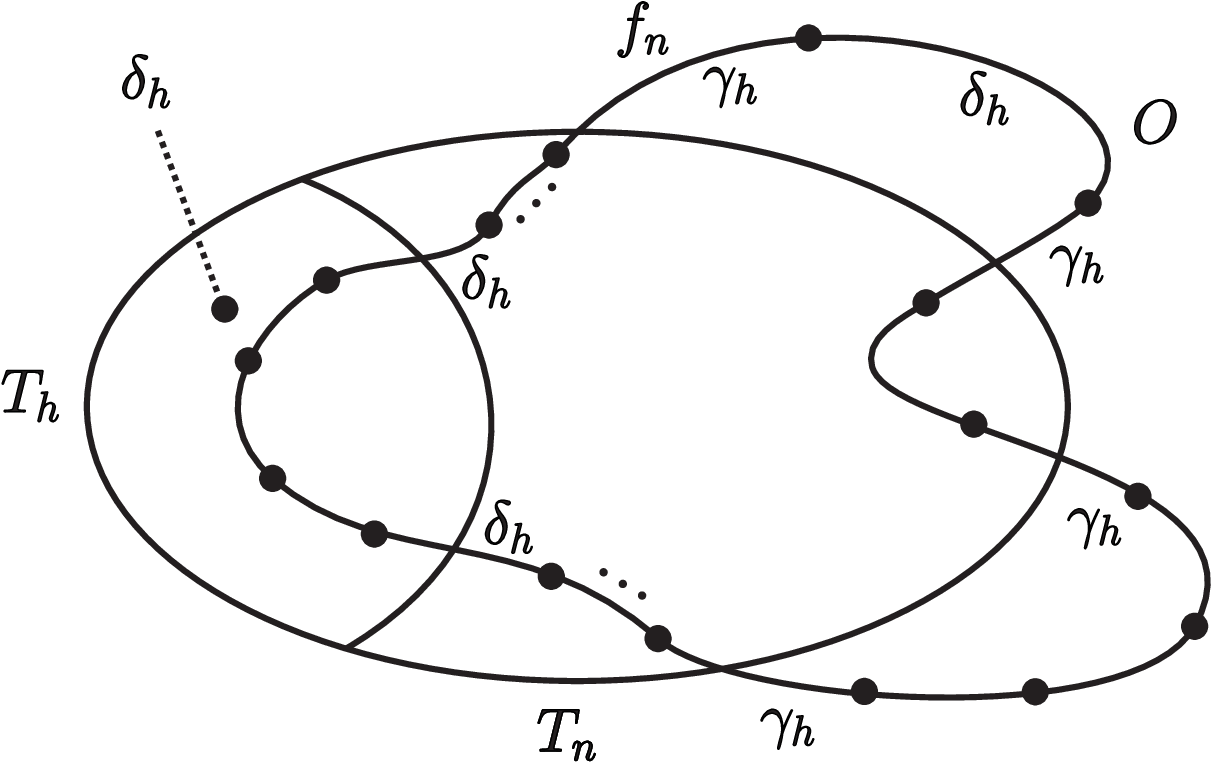}}
	\vspace{-1mm}\caption{Revisiting extension (RE)}
\end{figure}

In all previous approaches to the Goldberg-Seymour conjecture using the method of Tashkinov trees,
the trees involved in the proofs were constructed under a fix coloring. In sharp contrast, Algorithm 3.1 constructs tree-sequences and 
edge-colorings simultaneously as it progresses. Therefore, the structural property of an ETT embodied in the extension type and the corresponding 
coloring might be very fragile (see the next paragraph for details), even under stable colorings. This could cause a serious problem when we try to prove that an 
ETT $T$ with $T_n\subset T\subseteq T_{n+1}$ is elementary.

For example, if $T$ is not elementary, then we would apply Kempe changes to reduce the size of this counterexample to reach a contradiction while keeping 
the coloring $(T_n,D_n,\varphi_n)$-stable. However, if $\Theta_n=PE$ and the $(\gamma_n,\delta_n)$-path starting at $v_n$ has evolved to 
contain two or more vertices from $T_n$ during the process, then the resulting tree-sequence may no longer be an ETT under the new coloring, because PE requires that 
the $(\gamma_n,\delta_n)$-path starting at $v_n$ share exactly one vertex with $T_n$ in order to get $T_{n+1}$. Moreover, if $\Theta_n=RE$ and the edge $f_n$ 
no longer belongs to any $(\gamma_n,\delta_n)$-cycle during the process, then the resulting tree-sequence may not be an ETT under the new coloring anymore, because RE 
requires $f_n$ to be an edge of a $(\gamma_n,\delta_n)$-cycle.  

To circumvent this problem, we introduce the concept of mod coloring (see Definition 3.7) and impose a maximum property (see Definition 3.8) on ETT, 
and then we can ensure that the ETT structure is preserved under stable colorings (see Theorem 3.10(vi) and Lemma~\ref{td}). REs play an important role in the 
proof of Lemma~\ref{extension base} (Theorem 3.10(iv)), which in turn leads to Lemma~\ref{td} (Theorem 3.10(vi)). These technical results are essential to deriving the elementary 
property of the ETT we consider.

\vskip 2mm

Let us look back at Algorithm 3.1. Clearly, PE is the only extension that involves a non-stable coloring (in which one missing color at the 
supporting vertex has been changed). Based on this observation, we can exhibit some basic coloring properties (Lemmas 3.2-3.6) satisfied by 
ETTs. Recall that $D_{n-1}=\cup_{i \le n-1} \, S_i - \phibar_{n-1}(T_{n-1})$ is used to store colors employed in the construction of 
$T_n$ but not missing at any vertex of $T_{n-1}$ under $\varphi_{n-1}$.

\begin{lemma}\label{hku}
For $n \ge 1$, the following statements hold:
\begin{itemize}
\vspace{-1.5mm}
\item[(i)] $\phibar_{n-1}(T_n)\cup D_{n-1}\subseteq \phibar_{n}(T_n)\cup D_n \subseteq \phibar_{n}(T_{n+1})\cup D_n$. 
\vspace{-1.5mm}
\item[(ii)] For any $i\leq n$ and $v\in V(T_i)$, we have $\phibar_{i-1}(v)=\phibar_n(v)$ if $v$ is not used as a supporting vertex at any 
iteration $j$ with $i \leq j\leq n$ and $\Theta_j=PE$.
\vspace{-1.5mm}
\item[(iii)] For any edge $f$ incident to $T_n$, if $\varphi_{n-1}(f) \in \phibar_{n-1}(T_n)\cup D_{n-1}$,
then $\varphi_{n}(f)=\varphi_{n-1}(f)$, unless $\Theta_n=PE$ and $f=f_n$. So $\varphi_{n}(f) \in \phibar_{n}(T_n)\cup D_n$
provided that $\varphi_{n-1}(f) \in \phibar_{n-1}(T_n)\cup D_{n-1}$.
\vspace{-1.5mm}
\item[(iv)] $\varphi_{n-1} \langle T_n \rangle \subseteq \phibar_{n-1}(T_n)\cup D_{n-1}$
and $\varphi_n \langle T_n \rangle \subseteq \phibar_n(T_n)\cup D_n$. So $\sigma_n(f)=\varphi_n(f)$ for 
any $(T_n, D_n, \varphi_n)$-stable coloring $\sigma_n$ and any edge $f$ on $T_n$.  
\vspace{-1.5mm}
\item[(v)] If $\Theta_n=PE$,  then $\partial_{\varphi_n, \gamma_n}(T_n)=\{f_n\}$, and edges in 
$\partial_{\varphi_n, \delta_n}(T_n)$ are all incident to $V(T_n(v_n)-v_n)$. Furthermore, 
each color in  $\overline{\varphi}_n(T_n)-\{\delta_n\}$ is closed in $T_n$ under $\varphi_n$.
\end{itemize} 
\end{lemma}	

{\bf Proof.} By definition, $D_{n-1}=\cup_{i\leq n-1}S_i-\phibar_{n-1}(T_{n-1})$. So $\phibar_{n-1}(T_n)\cup D_{n-1}=\phibar_{n-1}(T_n) \cup [\cup_{i\leq n-1}S_i-\phibar_{n-1}(T_{n-1})]$. Since $\phibar_{n-1}(T_{n-1}) \subseteq \phibar_{n-1}(T_n)$, we obtain 

(1) $\phibar_{n-1}(T_n)\cup D_{n-1}=\phibar_{n-1}(T_n) \cup (\cup_{i\leq n-1}S_i)$.

\noindent Similarly, we can prove that 

(2) $\phibar_{n}(T_n)\cup D_n= \phibar_{n}(T_n)\cup (\cup_{i\leq n} S_i)$. 

(i) For any $\alpha \in \phibar_{n-1}(T_n)$, from Algorithm 3.1 and definition of stable colorings we see that $\alpha \in \phibar_{n}(T_n)$, unless $\Theta_n=PE$ and $\alpha=\gamma_n$; in this exceptional case, $\alpha  \in S_n$. 
So $\phibar_{n-1}(T_n) \subseteq \phibar_{n}(T_n) \cup S_n$ and hence $\phibar_{n-1}(T_n) \cup (\cup_{i\leq n-1}S_i) \subseteq \phibar_{n}(T_n)\cup (\cup_{i\leq n} S_i)$. It follows from (1) and (2) that $\phibar_{n-1}(T_n)\cup D_{n-1}\subseteq 
\phibar_{n}(T_n)\cup D_n$. Clearly, $\phibar_{n}(T_n)\cup D_n \subseteq \phibar_{n}(T_{n+1})\cup D_n$. 

(ii) In Algorithm 3.1 we always work with stable colorings except during PEs, where only missing colors at supporting vertices are changed. So the desired
statement follows.
 
(iii) Let $f$ be an edge incident to $T_n$ with $\varphi_{n-1}(f) \in \phibar_{n-1}(T_n)\cup D_{n-1}$. 
If $\Theta_n=RE$, then $\varphi_n=\varphi_{n-1}$ by Algorithm 3.1, which implies $\phiv_n(f) = 
\phiv_{n-1}(f)$. So we may assume that $\Theta_n\ne RE$. Let $\pi_{n-1}$ be the $(T_n, D_{n-1}, \phiv_{n-1})$-stable 
coloring as specified in Algorithm 3.1. By the definition of stable colorings, we obtain $\pi_{n-1}(f) = 
\phiv_{n-1}(f)$. If $\Theta_n=SE$, then $\phiv_n(f)=  \pi_{n-1}(f)$ by Algorithm 3.1. Hence $\phiv_n(f) 
= \phiv_{n-1}(f)$. It remains to consider the case when $\Theta_n=PE$. Let $\pi_{n-1}'$ be the $(T_n, D_{n-1} 
\cup\{\delta_n\}, \pi_{n-1})$-stable coloring as specified in Algorithm 3.1. By Lemma \ref{sc2}, 
$\pi_{n-1}'$ is $(T_n, D_{n-1}, \phiv_{n-1})$-stable. Hence $\pi_{n-1}'(f) = \phiv_{n-1}(f)$. Since $\phiv_n = \pi'_{n-1}/P_{v_n}(\delta_n, \gamma_n, \pi_{n-1}')$ and $P_{v_n}(\delta_n, \gamma_n, \pi_{n-1}')$ contains only one edge $f_n$ incident to $T_n$ (see Algorithm 3.1), we have $\phiv_n(f) = \pi_{n-1}'(f)$, unless $f=f_n$. 
It follows that $\phiv_n(f) = \phiv_{n-1}(f)$, unless $f=f_n$; in this exceptional case, $\varphi_{n-1}(f)=\delta_n$ and $\varphi_{n}(f) = \gamma_n \in S_n$. Hence $\varphi_{n}(f) \in \phibar_{n-1}(T_n)\cup D_{n-1} \cup S_n \subseteq 
\phibar_{n}(T_n)\cup D_n \cup S_n =\phibar_{n}(T_n)\cup D_n$ by (i) and (2), as desired. 

(iv) Let us first prove the statement $\varphi_{n-1} \langle T_n \rangle \subseteq \phibar_{n-1}(T_n)\cup D_{n-1}$
by induction on $n$. As the statement holds trivially when $n=1$, we proceed to the induction step and assume
that the statement has been established for $n-1$; that is,  

(3) $\varphi_{n-2} \langle T_{n-1} \rangle \subseteq \phibar_{n-2}(T_{n-1})\cup D_{n-2}$.  

By (3) and (iii) (with $n-1$ in place of $n$), for each edge $f$ on $T_{n-1}$ we have $\varphi_{n-1}(f) \in  
\phibar_{n-1}(T_{n-1}) \cup D_{n-1} \subseteq \phibar_{n-1}(T_n)\cup D_{n-1}$.  For each edge $f \in T_n-T_{n-1}$, from 
Algorithm 3.1 and TAA we see that $\varphi_{n-1}(f) \in D_{n-1}$ if $f$ is a connecting edge and $\varphi_{n-1}(f) \in \phibar_{n-1}(T_n)$ otherwise. Combining these observations, we obtain $\varphi_{n-1}(f) \in \phibar_{n-1}(T_n)\cup D_{n-1}$. 
Hence $\varphi_{n-1} \langle T_n \rangle \subseteq \phibar_{n-1}(T_n)\cup D_{n-1}$, which together
with (iii) implies $\varphi_n \langle T_n \rangle \subseteq \phibar_n(T_n)\cup D_n$.  

It follows that for any edge $f$ on $T_n$, we have $\varphi_n(f) \in \phibar_n(T_n)\cup D_n$. Thus $\sigma_n(f)=\varphi_n(f)$ 
for any $(T_n, D_n, \varphi_n)$-stable coloring $\sigma_n$. 

(v) From the definitions of $\pi_{n-1}$, vertex $v_n$ (maximum defective vertex) and stable colorings, we see that edges 
in $\partial_{\pi_{n-1}, \delta_n}(T_n)$ are all incident to $V(T_n(v_n))$, and each color in  $\overline{\pi}_{n-1}(T_n)$ is 
closed in $T_n$ under $\pi_{n-1}$. So, by the definitions of $\pi_{n-1}'$ and stable colorings, edges in $\partial_{\pi_{n-1}', 
\delta_n}(T_n)$ are all incident to $V(T_n(v_n))$, and each color in $\overline{\pi}_{n-1}'(T_n)$ is closed in $T_n$ under $\pi_{n-1}'$. 
Thus the desired statements follow instantly from the definition of $\varphi_n$ in PE. \qed

\vskip 3mm

The following lemma generalizes Lemma \ref{hku}(iii) and ensures that colors on some edges incident to a 
tree remain intact if we grow it by using Algorithm 3.1.  

\begin{lemma}\label{samecolor}
For any $1\le i \le n$ and any edge $f$ incident to $T_i$, if
$\phiv_{i-1}(f)\in\phibar_{i-1}(T_i)\cup D_{i-1}$, then $\phiv_j(f) = \phiv_{i-1}(f)$ for any $j$ with $i \le j \le n$, unless 
$f=f_h \in F_h$ for some $h$ with $i\le h \le j$ and $\Theta_h=PE$. In particular, if $f$ is an edge in $G[T_i]$ with 
$\phiv_{i-1}(f)\in\phibar_{i-1}(T_i)\cup D_{i-1}$, then $\phiv_j(f) = \phiv_{i-1}(f)$ for any $j$ with $i\le j \le n$.
\end{lemma}

{\bf Proof.}  By Lemma \ref{hku}(i), we have $\phibar_{h-1}(T_h)\cup D_{h-1} \subseteq \phibar_h(T_{h+1})\cup D_h$ for 
all $h \ge 1$. So to establish the first half, it suffices to prove the statement for $j=i$, which is exactly the same as
Lemma \ref{hku}(iii).  

Note that if $f$ is an edge in $G[T_i]$, then $f\notin \partial(T_h)$ for any $h$ with $i\le h\le j$. Hence
$f\neq f_h\in F_h$ for any $h$ with $i\le h\le j$ and $\Theta_h=PE$. Thus the second half also holds. \qed

\vskip 3mm
The lemma below describes some interesting properties satisfied by a sequence of PEs with the same supporting vertex.

\begin{lemma}\label{uniquezang}
Let $u$ be a vertex of $T_n$ and let $B_{n}$ be the set of all iterations $j$ with $1\le j \le n$, such that $\Theta_j 
=PE$ and $v_j=u$. Suppose $B_{n}=\{i_1,i_2, \ldots, i_p\}$, where $1\le i_1<i_2< \ldots < i_p \le n$. 
Then the following statements hold:
\begin{itemize}
\vspace{-1mm}
\item[(i)] $\gamma_{i_2}=\delta_{i_1}, \, \gamma_{i_3}=\delta_{i_2},\, \ldots, \, \gamma_{i_p}=\delta_{i_{p-1}}$;
\vspace{-1mm}
\item[(ii)] $\phibar_{n}(u) \cap (\cup_{j \in B_{n}} S_j) =\phibar_{i_p}(u) \cap (\cup_{j \in B_{n}} S_j) = \{\delta_{i_p}\}$; and
\vspace{-1.5mm}
\item[(iii)] $\phibar_{i_1-1}(u)=(\phibar_{i_p}(u)-\{\delta_{i_p}\}) \cup \{\gamma_{i_1}\}$ and
$\phibar_{i_p}(u)=(\phibar_{i_1-1}(u)- \{\gamma_{i_1}\}) \cup \{\delta_{i_p}\}$. 
\end{itemize}
\end{lemma}

{\bf Proof.} From the definition of $B_{n}$, we see that for any $1\le j \le p-1$ and iteration $h$ with $i_j+1 \le h \le i_{j+1}-1$, 
if $v_h=u$, then $\Theta_h=RE$ or $SE$. By Lemma~\ref{hku}(ii), we have 

(1) $\phibar_{i_j}(u)=\phibar_{i_{j+1}-1}(u)$ for $1\le j \le p-1$. Similarly, $\phibar_{i_p}(u)=\phibar_{n}(u)$.

\noindent According to the choice of $\gamma_h$ in a general iteration $h$ involving PE,

(2) $\gamma_{i_{j+1}}=\delta_{i_j}$ for $1\le j \le p-1$, where $\gamma_{i_{j+1}} \in  \phibar_{i_{j+1}-1}(u)$ and $\delta_{i_j} \in \phibar_{i_j}(u)$. 
 
\noindent Thus (i) follows instantly from (2). Using (1) and (2), we obtain

(3) $\phibar_{i_j}(u)-\{\delta_{i_j}\}= \phibar_{i_{j+1}-1}(u)- \{\gamma_{i_{j+1}}\}$ for $1\le j \le p-1$. 

\noindent Since $\phibar_{i_j}(u)$ is obtained from $\phibar_{i_j-1}(u)$ by replacing $\gamma_{i_j}$ with $\delta_{i_j}$, 

(4) $\phibar_{i_j-1}(u)- \{\gamma_{i_j}\} = \phibar_{i_j}(u)-\{\delta_{i_j}\}$ for $1\le j \le p$. 
 
\noindent Combining (3) and (4), we deduce that 

(5) the $2p$ sets $\phibar_{i_j-1}(u)- \{\gamma_{i_j}\}$ and $\phibar_{i_j}(u)-\{\delta_{i_j}\}$ for $1\le j \le p$ are all equal. 

By (5), the set $\cup_{j=1}^p \{\gamma_{i_j}, \delta_{i_j}\}$ (and hence $\cup_{j=1}^p S_{i_j}$) is disjoint from all the $2p$ sets
displayed above. So $\phibar_{i_p}(u) \cap (\cup_{j=1}^p S_{i_j})= [(\phibar_{i_p}(u)- \{\delta_{i_p}\}) \cap  (\cup_{j=1}^p S_{i_j})]
\cup \{\delta_{i_p}\}=\{\delta_{i_p}\}$, which together with (1) yields (ii).

Again by (5), $\phibar_{i_1-1}(u)-\{\gamma_{i_1}\}=\phibar_{i_p}(u)-\{\delta_{i_p}\}$. As $\gamma_{i_1} 
\in \phibar_{i_1-1}(u)$ and $\delta_{i_p} \in \phibar_{i_p}(u)$, we get $\phibar_{i_1-1}(u)=(\phibar_{i_p}(u)-\{\delta_{i_p}\}) 
\cup \{\gamma_{i_1}\}$ and $\phibar_{i_p}(u)=(\phibar_{i_1-1}(u)- \{\gamma_{i_1}\}) \cup \{\delta_{i_p}\}$. Therefore (iii) also holds. \qed 

\begin{lemma}\label{Dnzang}
$|D_n|\leq n$.
\end{lemma}

{\bf Proof.}  Recall that $D_n=\cup_{i\leq n}S_i-\phibar_n(T_n)$ (so $D_0=\emptyset$). For $1 \le i \le n$, by Algorithm 3.1, we 
have $S_i=\{\delta_i\}$ if $\Theta_i=SE$ and $S_i=\{\delta_i, \gamma_i\}$ otherwise.

If $\Theta_n=RE$, then $\phiv_n = \phiv_{n-1}$ and $S_n = S_{n-1}$. 
So $D_n \subseteq D_{n-1}$. If $\Theta_n=SE$, then $S_n=\{\delta_n\}$ and $\phibar_n(T_n)=\phibar_{n-1}(T_n)$. 
It follows that $D_n \subseteq D_{n-1} \cup \{\delta_n\}$. It remains to consider the case when $\Theta_n=PE$. 
Now $\delta_n \notin \phibar_{n-1}(T_n)$ and $(\phibar_{n-1}(T_n)-\{\gamma_n\}) \cup  \{\delta_n \} \subseteq 
\phibar_{n}(T_n)$. So
\begin{equation*}
\begin{aligned}
D_n & =\cup_{i\leq n}S_i-\phibar_n(T_n)\\ 
& \subseteq \cup_{i\leq n-1} S_i \cup \{\delta_n,\gamma_n\} -[(\phibar_{n-1}(T_n)-\{\gamma_n\}) \cup \{ \delta_n \}]\\
& \subseteq \cup_{i\leq n-1} S_i \cup \{\gamma_n\} -(\phibar_{n-1}(T_n)-\{\gamma_n\})\\
& \subseteq [\cup_{i\leq n-1} S_i  - \phibar_{n-1}(T_n)] \cup \{\gamma_n\}\\
& \subseteq D_{n-1} \cup \{\gamma_n\}.
\end{aligned}
\end{equation*}

Combining the above three cases, we obtain $|D_n|\leq |D_{n-1}|+1$ for $n\ge 1$.  Hence $|D_n|\leq n$. \qed

\begin{lemma}\label{stablezang}
Suppose $\Theta_n=PE$. Let $\sigma_n$ be a $(T_n, D_n, \phiv_n)$-stable coloring and let 
$\sigma_{n-1}=\sigma_n/P_{v_{n}} (\gamma_{n},\delta_{n}, \sigma_n)$. If $P_{v_{n}}(\gamma_{n},\delta_{n},\sigma_n)
\cap T_{n}=\{v_{n}\}$, then $\sigma_{n-1}$ is  $(T_n, D_{n-1} \cup \{\delta_n\}, \pi_{n-1})$-stable
and hence is $(T_{n},D_{n-1}, \varphi_{n-1})$-stable. 
\end{lemma}

{\bf Proof.} Let $\pi_{n-1}'$ be as specified in Algorithm 3.1. Recall that 

(1) $\pi_{n-1}'$ is $(T_n, D_{n-1} \cup \{\delta_n\}, \pi_{n-1})$-stable. 

\noindent By definition, $\phiv_n = \pi_{n-1}'/P_{v_n}(\gamma_n, \delta_n, \pi'_{n-1})$. So 

(2) $\pi_{n-1}' =  \phiv_n /P_{v_n}( \gamma_n, \delta_n, \phiv_n)$.

We propose to show that 

(3) $\sigma_{n-1}$ is $(T_n, D_{n-1} \cup \{\delta_n\}, \pi_{n-1}')$-stable. 

By the definition of $\sigma_{n-1}$ and (2), we obtain

(4) $\sigma_n$ and $\sigma_{n-1}$ agree on every edge incident to $T_n$ except $f_n$, for which $\sigma_n(f_n)=\gamma_n$ 
and $\sigma_{n-1}(f_n)=\delta_n$; and

(5) $\varphi_n$ and $\pi_{n-1}'$ agree on every edge incident to $T_n$ except $f_n$, for which $\varphi_n(f_n)=\gamma_n$ 
and $\pi_{n-1}'(f_n)=\delta_n$. 

Since $D_n=D_{n-1}\cup \{\gamma_n\}$, $\delta_n\in\phibar_n(T_n)$ and $\sigma_n$ is $(T_n, D_n, \phiv_n)$-stable, (3) follows instantly from (4) and (5). 
Using (1), (3) and Lemma \ref{sc2}, we see that $\sigma_{n-1}$ is $(T_{n},D_{n-1} \cup \{\delta_n\}, \pi_{n-1})$-stable. So $\sigma_{n-1}$ is $(T_{n},
D_{n-1}, \pi_{n-1})$-stable. Since $\pi_{n-1}$ is $(T_{n},D_{n-1}, \varphi_{n-1})$-stable, from  Lemma \ref{sc2}
we conclude that $\sigma_{n-1}$ is $(T_{n},D_{n-1}, \varphi_{n-1})$-stable. \qed\\

Observe that an extended Tashkinov tree $T$ (see Definition \ref{wz1}) has a built-in ladder-like structure. So we propose to call the sequence 
$T_1\subset T_2 \subset \ldots \subset T_n \subset T$ the {\it ladder}\label{ladder} of $T$, and call $n$ the {\it rung number}\label{rung} of $T$ and 
denote it by $r(T)$. Moreover, we call $(\phiv_0, \phiv_1, \dots, \phiv_n)$ the {\it coloring sequence}\label{coloringsequence} of $T$, and call $\TT$ the 
Tashkinov series {\em corresponding} to $T$.

In our proof we shall frequently work with stable colorings; the following concept will be used to keep track of 
the structures of ETTs.   

\begin{definition} \label{wz2}
{\rm  Let $\TT=\{(T_i, \phiv_{i-1}, S_{i-1}, F_{i-1}, \Theta_{i-1}): 1\le i \le n+1\}$ be a Tashkinov series 
constructed from a $k$-triple $(G,e, \varphi)$ by using Algorithm 3.1. A coloring $\sigma_n\in  {\cal C}^k(G-e)$
is called $\varphi_n \bmod T_n$\label{modcoloring} if every tree-sequence $T^*\supset T_n$ obtained from $T_n+f_n$ 
(resp. $T_n$) by TAA under $\sigma_n$ when $\Theta_n=$ $RE$ or $SE$ (resp. when $\Theta_n=PE$) is an ETT under 
$\sigma_n$, with a corresponding Tashkinov series $\TT^*=\{(T_i^*, \sigma_{i-1}, S_{i-1}, F_{i-1}, \Theta_{i-1}): 
1\le i \le n+1\}$, satisfying the following conditions for all $i$ with $1\le i \le n$:

\vskip 1mm

$\bullet$ $T_i^*=T_i$ and

$\bullet$ $\sigma_i$ is a $(T_i, D_i, \varphi_i)$-stable coloring in ${\cal C}^k(G-e)$. 
\vskip 1mm

\noindent We call each $T^*$ an ETT {\em corresponding} to $(\sigma_n, T_n)$ (or simply {\em corresponding}\label{correspondingett} to 
$\sigma_n$ if no ambiguity arises).}
\end{definition}

\noindent {\bf Remark.} Comparing $\TT^*$ with $\TT$, we see that $T_{i+1}^*$ in $\TT^*$ is obtained from $T^*_i=T_i$ by 
using the same connecting edge, connecting color, and extension type as $T_{i+1}$ in $\TT$ for $1\le i \le n$. However, $T^*_{n+1}$ may be different from $T_{n+1}$.
Furthermore, $T_1\subset T_2 \subset \ldots \subset T_n \subset T^*$ is the ladder of $T^*$ and $r(T^*)=n$.  
Since $\sigma_i$ is a $(T_i, D_i, \varphi_i)$-stable coloring, by Lemma \ref{hku}(iv), we have $\sigma_i(f)=\varphi_i(f)$ 
for any edge $f$ on $T_i$ and $1\le i \le n$; this fact will be used repeatedly in our paper. 

\vskip 2mm

To ensure that the structures of ETTs are preserved under taking stable colorings, we impose some restrictions
on such trees. 

\begin{definition} \label{wz3}
{\rm Let $T$ be an ETT constructed from a $k$-triple $(G,e, \varphi)$ by using the Tashkinov series 
$\TT=\{(T_i, \phiv_{i-1}, S_{i-1}, F_{i-1}, \Theta_{i-1}): 1\le i \le n+1\}$. We say that $T$ has the {\em maximum property}\label{maximumproperty} 
(MP) under $(\varphi_0, \varphi_1, \ldots, \varphi_n)$ (or simply under $\varphi_n$ if no ambiguity arises), 
if $|T_1|$ is maximum among all Tashkinov trees $T_1'$ with respect to an edge $e'\in E$ and a coloring $\varphi'_0 \in  {\cal C}^k(G-e')$, 
and $|T_{i+1}|$ is maximum over all $(T_{i}, D_i, \phiv_{i})$-stable colorings for any $i$ with 
$1\le i \le n-1$; that is, $|T_{i+1}|$ is maximum over all tree-sequences $T_{i+1}'$, which is a closure of $T_i+f_i$ 
(resp. $T_i$) under a $(T_{i}, D_i, \phiv_{i})$-stable coloring $\phiv_i'$ if $\Theta_i=RE$ or $SE$
(resp. if $\Theta_i=PE$), where $f_i$ is the connecting edge in $F_i$. }
\end{definition}

Notice that in the above definition $|T_{n+1}|$ is not required to be maximum over all $(T_n, D_n, \varphi_n)$-stable colorings. This relaxation 
allows us to proceed by induction in our proofs.

As described before, the tree-sequence structure generated by Algorithm 3.1 might be very fragile, because RE does not allow change of coloring and PE requires 
the supporting vertex to be the exit of an exit-path. At this point, it is natural to ask whether there exists an ETT with MP and with an arbitrarily given rung number or
an arbitrarily given size. We shall demonstrate (see Corollary \ref{welldefined}) that the answer is in the affirmative. The statement below follows instantly from the 
above two definitions and Lemma \ref{sc2}.

\begin{lemma}\label{MP}
Let $T$ be an ETT constructed from a $k$-triple $(G,e, \varphi)$ by using the Tashkinov series $\TT=\{(T_i, 
\phiv_{i-1}, S_{i-1}, F_{i-1}, \Theta_{i-1}): 1\le i \le n+1\}$, let $\sigma_n$ be a $\varphi_n \bmod T_n$ coloring, 
and let $T^*$ be an ETT corresponding to $(\sigma_n, T_n)$ (see Definition \ref{wz2}). If $T$ satisfies MP under 
$\varphi_n$, then $T^*$ satisfies MP under $\sigma_n$.  \qed
\end{lemma}

Let us introduce one more notation and two more concepts before presenting our main theorem. For each $v\in V(T)$, we use $m(v)$ to denote the minimum 
subscript $i$ such that $v \in V(T_i)$. Let $\alpha$ and $\beta$ be two colors in $[k]$. We say that $\alpha$ and $\beta$ are
$T$-{\em interchangeable}\label{interchangable} under $\varphi_n$ if there is at most one $(\alpha, \beta)$-path with respect to $\varphi_n$
intersecting $T$. Note that in this situation we can easily find an $(\alpha, \beta)$-path disjoint from $T$ and then switch 
colors along it while keeping the resulting coloring stable. So this concept is very helpful for deriving elementary property satisfied 
by an ETT.  When $T$ is closed (that is, $T=T_{n+1}$), we also say that $T$ has the {\em interchangeability 
property}\label{interchangeproperty} with respect to $\varphi_{n}$ if under any $(T, D_n, \phiv_{n})$-stable coloring $\sigma_n$, any two 
colors $\alpha$ and $\beta$ are $T$-interchangeable, provided that $\overline{\sigma}_n(T) \cap \{\alpha, \beta\} 
\ne \emptyset$ (equivalently $\phibar_{n}(T) \cap \{\alpha, \beta\} \ne \emptyset$).    

We aim to show, by induction on the rung number, that every ETT satisfying MP is elementary. To carry out the induction step, we need several 
auxiliary results concerning ETTs with MP. Thus what we are going to prove is a stronger theorem (containing six statements) given below, 
in which the undefined symbols and notations can all be found in Algorithm 3.1. Together with Theorem \ref{egraph2}, statements (i) and 
(vi) will imply Theorem~\ref{ThmGS2}. Statements (ii)-(v) will be used in the proofs of (i) and (vi). Moreover, the proof of (iv) relies 
directly on MP and the design of RE, and the proofs of (iii) and (v) are based on the fact that the supporting and extension vertices 
involved in Algorithm 3.1 are maximum defective vertices over all stable colorings.

\begin{theorem} \label{thm:tech10}
	Let $T$ be an ETT constructed from a $k$-triple $(G,e, \varphi)$ by using the Tashkinov series 
	$\TT=\{(T_i, \phiv_{i-1}, S_{i-1}, F_{i-1}, \Theta_{i-1}): 1\le i \le n+1\}$. If $T$ has 
	MP under $\varphi_n$, then the following statements hold:
	\begin{itemize}
		\vspace{-1.2mm}
		\item[(i)] $V(T)$ is elementary with respect to $\phiv_n$.
		\vspace{-2mm}
		\item[(ii)] $T_{n+1}$ has the interchangeability property with respect to $\varphi_n$.
		\vspace{-2mm}
		\item[(iii)] For any $i \le n$,  if $v_i$ is a supporting vertex with $m(v_i) =j$, then every $(T_i, D_i, \varphi_i)$-stable 
		coloring $\sigma_i$ is $(T(v_i)-v_i, D_{j-1},\varphi_{j-1})$-stable, so $\sigma_i$ is $(T_{j-1},D_{j-1},\varphi_{j-1})$-stable. 
		Furthermore, for any two distinct supporting 		vertices $v_i$ and $v_j$ with $i,j \le n$,  if $m(v_i)=m(v_j)$, then $S_i \cap S_j=\emptyset$.
		\vspace{-2mm}
		\item[(iv)] If $\Theta_n = PE$, then $P_{v_n}(\gamma_n, \delta_n, \sigma_n)$ contains precisely one vertex,
		$v_n$, from $T_n$ for any $(T_n, D_n, \phiv_n)$-stable coloring $\sigma_n$.
		\vspace{-2mm}
		\item[(v)] For any $(T_n, D_n, \phiv_n)$-stable coloring $\sigma_n$ and any defective color $\delta$ of $T_n$ 
		with respect to $\sigma_n$, if $v$ is a vertex but not the smallest one (in the order $\prec$)  
		in $I[\partial_{\sigma_n, \delta}(T_n)]$, then $v \preceq v_i$ for any supporting or extension vertex  
		$v_i$ with $m(v) \le i$.
		\vspace{-2mm} 
		\item[(vi)] Every $(T_n,D_n,\varphi_n)$-stable coloring $\sigma_n$ is a $\phiv_n\bmod T_n$ coloring. (So
		every ETT corresponding to $(\sigma_n, T_n)$ (see Definition \ref{wz2}) satisfies 
		MP under $\sigma_n$ by Lemma~\ref{MP}.)
	\end{itemize}
\end{theorem}
\vskip 1mm

Let us show that Theorem \ref{ThmGS2} can be deduced easily from statement (i) and the following corollary (which relies on statement (vi)).  

\begin{corollary}\label{welldefined}
Let $\TT=\{(T_i, \phiv_{i-1}, S_{i-1}, F_{i-1}, \Theta_{i-1}): 1\le i \le n+1\}$ be a Tashkinov series  
constructed from a $k$-triple $(G,e, \varphi)$. Suppose $T_{n+1}$ has MP under $\varphi_n$.
Then there exists a Tashkinov series $\TT^*=\{(T_i^*, \sigma_{i-1}, S_{i-1}, F_{i-1}, \Theta_{i-1}): 1\le i \le n+1\}$, 
satisfying the following conditions for $1\le i \le n$:
\begin{itemize}
\vspace{-1mm}
\item[(i)] $T_i^*=T_i$;
\vspace{-2mm}
\item[(ii)] $\sigma_i$ is a $(T_i, D_i, \varphi_i)$-stable coloring in ${\cal C}^k(G-e)$; and 
\vspace{-2mm}
\item[(iii)] $|T_{i+1}^*|$ is maximum over all $(T_i, D_i, \sigma_i)$-stable colorings (note that Definition \ref{wz3} only requires 
this for $i\leq n-1$). 
\vspace{-1mm}
\end{itemize}
Furthermore, if $T_{n+1}^*$ is not strongly closed with respect to $\sigma_n$, then there exists a Tashkinov series 
$\{(T_i^*, \sigma_{i-1}, S_{i-1}, F_{i-1}, \Theta_{i-1}): 1\le i \le n+2\}$, such that $T_{n+1}^* \subset T_{n+2}^*$
and $T_{n+2}^*$ satisfies MP under $\sigma_{n+1}$. 
\end{corollary}	

{\bf Proof.} Let $\sigma_n$ be a $(T_n, D_n, \varphi_n)$-stable coloring such that a closure of $T_n$ (resp. of $T_n+f_n$) under $\sigma_n$, 
denoted by $T^*_{n+1}$, has maximum size over all $(T_n, D_n, \varphi_n)$-stable colorings if $\Theta_n=PE$ (resp. if $\Theta_n=RE$ or $SE$). 
By Lemma \ref{sc2}, every $(T_n, D_n, \sigma_n)$-stable coloring is a 
$(T_n, D_n, \varphi_n)$-stable coloring. So $|T_{n+1}^*|$ is also maximum over all $(T_n, D_n, \sigma_n)$-stable colorings.

Since $\sigma_n$ is $(T_n, D_n, \varphi_n)$-stable, it is $\phiv_n\bmod T_n$ by Theorem~\ref{thm:tech10}(vi). 
Thus Definition \ref{wz2} guarantees the existence of a 
Tashkinov series $\TT^*=\{(T_i^*, \sigma_{i-1}, S_{i-1}, F_{i-1}, \Theta_{i-1}): 1\le i \le n+1\}$ that satisfies conditions (i) 
and (ii) as described above. By Lemma \ref{MP}, $|T_{i+1}^*|$ is maximum over all $(T_i, D_i, \sigma_i)$-stable 
colorings as well for $1\le i \le n-1$.   

Suppose $T_{n+1}^*$ is not strongly closed with respect to $\sigma_n$. Then we can construct a new tuple  
$(T_{n+2}^*, \sigma_{n+1}, S_{n+1}, F_{n+1}, \Theta_{n+1})$ by using Algorithm 3.1. Clearly, $T_{n+1}^* \subset T_{n+2}^*$
and $T_{n+2}^*$ satisfies MP under $\sigma_{n+1}$. \qed 

\vskip 3mm
{\bf Proof of Theorem \ref{ThmGS2}.} Let $\TT=\{(T_i, \phiv_{i-1}, S_{i-1}, F_{i-1}, \Theta_{i-1}): 1\le i \le n+1\}$ be a 
Tashkinov series constructed from a $k$-triple $(G,e, \varphi)$, such that 

$(a)$ $T_{n+1}$ satisfies MP under $\varphi_n$;

$(b)$ subject to $(a)$, $|T_{n+1}|$ is maximum over all $(T_n, D_n, \varphi_n)$-stable colorings; and  

$(c)$ subject to $(a)$ and $(b)$, the integer $n$ is maximum. 

\noindent 
Since $G$ is finite, by Corollary \ref{welldefined}, such a Tashkinov series $\TT$ exists. Observe that $T_{n+1}$ is strongly closed, 
for otherwise, Corollary \ref{welldefined} would enable us to further extend $\TT$ to a longer Tashkinov series $\{(T_i^*, \sigma_{i-1}, S_{i-1}, 
F_{i-1}, \Theta_{i-1}): 1\le i \le n+2\}$ satisfying $(a)$ and $(b)$, contradicting $(c)$.  By Theorem \ref{thm:tech10}(i), $V(T_{n+1})$ is elementary with 
respect to $\phiv_n$. From Theorem \ref{egraph2}(i) and (iv), we thus deduce that $G$ is an elementary multigraph. \qed 

\vskip 3mm

The proof of Theorem \ref{thm:tech10} will take up the entire remainder of this paper.

\section{Auxiliary Results}

We prove Theorem~\ref{thm:tech10} by induction on the rung number $r(T)=n$. The present section is devoted to a proof 
of statement (ii) in Theorem~\ref{thm:tech10} in the base case and proofs of statements (iii)-(vi) in the general case.  

For $n=0$, statement (i) follows from Theorem \ref{TashTree},  statements (iii)-(vi) hold trivially, and statement (ii) 
is a corollary of the following more general lemma (because $T_1$ is also a closed Tashkinov tree with respect to $e$ and 
any $(T_1, D_0, \phiv_{0})$-stable coloring $\sigma_0$).  

\begin{lemma}\label{tchange}
Let $(G,e, \varphi)$ be a $k$-triple, let $T$ be a closed Tashkinov tree with respect to $e$ and $\varphi$,
and let $\alpha$ and $\beta$ be two colors in $[k]$ with $\overline{\varphi}(T) \cap \{\alpha,\beta\}\ne 
\emptyset$. Then there is at most one $(\alpha,\beta)$-path with respect to $\varphi$ intersecting $T$. 
\end{lemma}

{\bf Proof.} Assume the contrary: there are at least two $(\alpha,\beta)$-paths $Q_1$ and $Q_2$ with respect 
to $\varphi$ intersecting $T$. By Theorem \ref{TashTree}, $V(T)$ is elementary with respect to $\varphi$. So  
$T$ contains at most two vertices $v$ with $\overline{\varphi}(v) \cap \{\alpha, \beta\} \ne \emptyset$, 
which in turn implies that at least two ends of $Q_1$ and $Q_2$ are outside $T$. By hypothesis, $T$ is closed 
with respect to $\varphi$. Hence precisely one of $\alpha$ and $\beta$, say $\alpha$, is in $\overline{\varphi}(T)$. 
Thus we further deduce that at least three ends of $Q_1$ and $Q_2$ are outside $T$. Traversing $Q_1$ and $Q_2$ 
from these ends respectively, we can find at least three $(T, \varphi, \{\alpha, \beta\})$-exit paths 
$P_1,P_2,P_3$. We call the tuple $(\varphi,T, \alpha, \beta, P_1,P_2,P_3)$ a {\em counterexample} and use ${\cal K}$ 
to denote the set of all such counterexamples.

With a slight abuse of notation, let $(\varphi, T, \alpha, \beta, P_1,P_2,P_3)$ be a counterexample in ${\cal K}$ with 
the minimum $|P_1|+|P_2|+|P_3|$. For $i=1,2,3$, let $a_i$ and $b_i$ be the ends of $P_i$ with $b_i \in V(T)$, 
and $f_i$ be the edge of $P_i$ incident to $b_i$.  Renaming subscripts if necessary, we may assume that $b_1\prec b_2 
\prec b_3$.  Let $\gamma\in\phibar(b_3)$ and let $\sigma_1=\varphi/(G-T,\alpha,\gamma)$. Then $P_3$ is a $(\gamma,\beta)$-path under $\sigma_1$. Clearly, $\sigma_1 \in 
{\cal C}^k(G-e)$ and $T$ is also a Tashkinov tree with respect to $e$ and $\sigma_1$. Furthermore, $f_i$ is colored 
by $\beta$ under both $\varphi$ and $\sigma_1$ for $i=1,2,3$.  

Consider $\sigma_2=\sigma_1/P_3$. Note that $\beta \in \overline{\sigma}_2(b_3)$.
Let $T'$ be obtained from $T(b_3)$ by adding $f_1$ and $f_2$ and let $T''$ be a closure of $T'$ under $\sigma_2$.
Obviously, both $T'$ and $T''$ are Tashkinov trees with respect to $e$ and $\sigma_2$. By Theorem \ref{TashTree}, 
$V(T'')$ is elementary with respect to $\sigma_2$. 

Observe that none of $a_1,a_2,a_3$ is contained in $T''$, for otherwise, let $a_i \in V(T'')$ for some $i$
with $1\le i \le 3$. Since $\{\beta,\gamma\}\cap \overline{\sigma}_2(a_i) \ne \emptyset$ and $\beta \in \overline{\sigma}_2(b_3)$, 
we obtain $\gamma \in \overline{\sigma}_2(a_i)$. Hence from TAA we see that $P_1,P_2,P_3$ are all entirely contained in $G[T'']$, 
which in turn implies  $\gamma \in \overline{\sigma}_2(a_j)$ for $j=1,2,3$. So $V(T'')$ is not 
elementary with respect to $\sigma_2$, a contradiction. Each $P_i$ contains a subpath $L_i$, which is a $T''$-exit path 
with respect to $\sigma_2$. Since $f_1$ is not contained in $L_1$, we obtain $|L_1|+|L_2|+|L_3|<|P_1|+|P_2|+|P_3|$.
Thus the existence of the counterexample $(\sigma_2, T'', \gamma, \beta, L_1,L_2,L_3)$ violates the minimality 
assumption on $(\varphi, T, \alpha, \beta, P_1,P_2,P_3)$. \qed
 
\vskip 3mm

So Theorem~\ref{thm:tech10} is true in the base case. Suppose we have established that 

{\bf (4.1)} Theorem~\ref{thm:tech10} holds for all ETTs with at most $n-1$ rungs and satisfying MP, for some $n\ge 1$.

Let us proceed to the induction step. We postpone the proof of Theorem~\ref{thm:tech10}(i) and (ii) to 
Section 7, and present a proof of Theorem~\ref{thm:tech10}(iii)-(vi) in this section. In our proof of the $(i+2)$th 
statement in Theorem~\ref{thm:tech10} for $2\le i \le 4$, we further assume that

{\bf (4.} \hskip -2.8mm ${\bm i}${\bf )} the $j$th statement in \hskip 0.5mm Theorem~\ref{thm:tech10} \hskip 0.5mm holds for all ETTs with at most $n$ rungs and 
satisfying MP, for all $j$ with $3\le j \le i+1$.

Note that $(4. \hskip 0.4mmi)$ corresponds to $(4.2), (4.3)$ and $(4.4)$, respectively, for $i=2,3$ and $4$. 
For example, when we try to prove Theorem~\ref{thm:tech10}(v) (now $i=3$), we assume (4.3), which says that 
both Theorem~\ref{thm:tech10}(iii) and Theorem~\ref{thm:tech10}(iv) hold for all ETTs with at most $n$ rungs and satisfying MP. \\

We break the proof of the induction step into a series of lemmas. The following lemma derives some properties satisfied
by supporting vertices and connecting colors.

\begin{lemma}\label{extension rules}
(Assuming (4.1)) Theorem~\ref{thm:tech10}(iii) holds for all ETTs with $n$ rungs and satisfying MP; that is, 
for any $i \le n$,  if $v_i$ is a supporting vertex with $m(v_i) =j$, then every $(T_i, D_i, \varphi_i)$-stable 
coloring $\sigma_i$ is $(T(v_i)-v_i, D_{j-1},\varphi_{j-1})$-stable, so $\sigma_i$ is $(T_{j-1},D_{j-1},\varphi_{j-1})$-stable. 
Furthermore, for any two distinct supporting vertices $v_i$ and $v_j$ with $i,j \le n$, if $m(v_i)=m(v_j)$, then $S_i \cap S_j=\emptyset$.
\end{lemma}	

{\bf Proof}. By (4.1), Lemma~\ref{extension rules} holds for all ETTs with at most $n-1$ rungs and satisfying MP. So we may assume that
$T$ is an ETT with the corresponding Tashkinov series $\TT=\{(T_h, \phiv_{h-1}, S_{h-1}, F_{h-1}, \Theta_{h-1}): 1\le h \le n+1\}$ 
and satisfies MP under $\varphi_n$. Furthermore, $i = n$ throughout our proof.

In the first half of this lemma, $m(v_n) =j$ and $\sigma_n$ is a $(T_n, D_n, \varphi_n)$-stable coloring. Write $T^*=T(v_n)-v_n$
(so $T^*\subseteq T_j$). As $j\leq n$, repeated application of Lemma \ref{hku}(i) yields $\phibar_{j-1}(T_j)\cup D_{j-1}\subseteq 
\phibar_{n-1}(T_n)\cup D_{n-1} \subseteq \phibar_n(T_n) \cup D_n$. In particular, $D_{j-1}\subseteq \phibar_n(T_n) \cup D_n$. Hence $\sigma_n$ is a 
$(T^*, D_{j-1}, \varphi_n)$-stable coloring. By Lemma \ref{sc2}, to prove that $\sigma_n$ is $(T^*, D_{j-1},\varphi_{j-1})$-stable,
it suffices to show that $\phiv_n$ is $(T^*, D_{j-1}, \phiv_{j-1})$-stable. 

If $j=n$, then $v_n$ is the only supporting vertex contained inside $T_n$ but outside $T_{n-1}$. Recall that in Algorithm 3.1 the coloring
$\pi'_{n-1}$ is $(T_{n},D_{n-1}\cup\{\delta_n\},\pi_{n-1})$-stable and $\varphi_n=\pi'_{n-1}/P_{v_n}(\gamma_n,\delta_n,\pi'_{n-1})$, where 
$P_{v_n}(\gamma_n,\delta_n,\pi'_{n-1})\cap V(T_n)=\{v_n\}$. So $\phiv_n$ is a $(T^*, D_{n-1}, \pi_{n-1})$-stable
coloring. By Lemma \ref{sc2}, it is also a $(T^*, D_{n-1}, \phiv_{n-1})$-stable coloring, because $\pi_{n-1}$ is $(T_{n},D_{n-1},\varphi_{n-1})$-stable. 
Thus we assume hereafter that $j < n$. As $v_n$ is the largest vertex (in the order $\prec$) in $I[\partial_{\pi_{n-1}, \delta_{n}}(T_{n})]$
(see Algorithm 3.1), with $\delta_{n}=\pi_{n-1}(f_n)$, and $v_n$ is contained in $T_{j}\subseteq T_{n-1}$, we see that $\delta_{n}$ is a defective 
color of $T_{n-1}$ with respect to $\pi_{n-1}$, and $v_n$ is not the smallest vertex (in the order $\prec$) in $I[\partial_{\pi_{n-1}, 
\delta_{n}}(T_{n-1})]$. As $\pi_{n-1}$ is also a $(T_{n-1},D_{n-1},\varphi_{n-1})$-stable coloring, applying (4.1) and Theorem~\ref{thm:tech10}(v) 
to $v=v_n$ and $\pi_{n-1}$, we obtain $v_n \preceq v_h$ for any supporting vertex $v_h$ with $j \le h \le n-1$.  Thus $\phibar_{j-1}(v)=\phibar_{n}(v)$ 
for each vertex $v$ of $T^*$ by Lemma~\ref{hku}(ii). Furthermore, $\phiv_n(f) =\phiv_{j-1}(f)$ for each edge $f$ incident to $T^*$ with 
$\phiv_{j-1}(f) \in \phibar_{j-1}(T^*)\cup D_{j-1}$ by Lemma \ref{samecolor}. Hence $\phiv_n$ is $(T^*, D_{j-1}, \phiv_{j-1})$-stable, as desired. 

Let us proceed to the second half. Now $v_j$ is a supporting vertex with $j< n$ and $m(v_j)=m(v_n)$. 
To prove that $S_n \cap S_j = \emptyset$, we shall actually show that

(1) there are edges $f,g$ in $G[T_n]$ incident to $v_n$ with $\varphi_n(f)=\gamma_j$ and $\varphi_n(g)=\delta_j$.

\noindent Assuming (1), it follows instantly that $\gamma_j,\delta_j\notin S_n$, because $\varphi_n(f_n)=\gamma_n$, $f_n\in \partial(T_n)$ (so $f_n\notin G[T_n]$), 
and $\delta_n\in\phibar_{n}(v_n)$ (see Algorithm 3.1).  

To justify (1), recall that in Algorithm 3.1 coloring $\pi'_{j-1}$ is $(T_j, D_{j-1}\cup\{\delta_j\},\pi_{j-1})$-stable, 
$P=P_{v_j}(\gamma_j, \delta_j, \pi'_{j-1})$ contains only vertex $v_j$ from $T_j$, and $\phiv_j = \pi'_{j-1}/P$.
Write $Q=P_{v_n}(\gamma_j, \delta_j, \pi'_{j-1})$. Since $v_j\neq v_n$, $P$ and $Q$ are vertex-disjoint under $\pi'_{j-1}$.
For convenience, we still use $P$ and $Q$ to denote the corresponding paths under $\phiv_j$.  
By (4.1) and Theorem~\ref{thm:tech10}(ii), $T_{j+1}$ has the interchangeability property with 
respect to $\phiv_j$. So $P$ is the unique $(\gamma_j, \delta_j)$-path intersecting $T_{j+1}$ and $Q$ is a $(\gamma_j, \delta_j)$-cycle 
under $\phiv_j$.  Let $r>j$ be the smallest subscript with $\Theta_r \ne RE$. Since $\Theta_n = PE$, we have $r \le n$. From RE in Algorithm 3.1
we see that $Q$ is fully contained in $G[T_r]$.  Repeated application of Lemma \ref{hku}(i) yields  $\phibar_j(T_{j+1})\cup D_j \subseteq 
\phibar_{r-1}(T_r)\cup D_{r-1}$. Since  $\cup_{h=1}^j S_h \subseteq  \phibar_j(T_j) \cup D_j $, we have $S_j \subseteq \phibar_{r-1}(T_r)\cup D_{r-1}$. 
So $Q$ is also a $(\gamma_j, \delta_j)$-cycle containing $v_n$ under $\phiv_n$ by Lemma~\ref{samecolor} (with respect to $T_r$). Since 
$T_r \subseteq T_n$, we establish (1). \qed

\vskip 3mm
The following lemma asserts that parallel extensions (PEs) used in Algorithm 3.1 are preserved under taking stable 
colorings. Its proof is perhaps the most difficult part of the whole paper. After reading the 
proof, we may fully understand why RE is introduced in 
the algorithm.  

\begin{lemma}\label{extension base}
(Assuming (4.1) and (4.2)) Theorem~\ref{thm:tech10}(iv) holds for all ETTs with $n$ rungs and  satisfying MP; 
that is, if $\Theta_n = PE$, then $P_{v_n}(\gamma_n, \delta_n, \sigma_n)$ contains precisely one vertex,
$v_n$, from $T_n$ for any $(T_n, D_n, \phiv_n)$-stable coloring $\sigma_n$.
\end{lemma}

{\bf Proof}. Assume the contrary: $P_{v_n}(\gamma_n, \delta_n, \sigma_n)$ contains at least two vertices
from $T_n$ for some $(T_n, D_n, \phiv_n)$-stable coloring $\sigma_n$.  Let $j=m(v_n)$.  By applying a series
of Kempe changes to $\sigma_n$, we shall construct a certain $(T_j(v_n) - v_n, D_{j-1}, \phiv_{j-1})$-stable 
coloring $\mu$ and a certain ETT $T_j^{\mu}$ corresponding to $(\mu, T_{j-1})$ with ladder $T_1\subset T_2\subset \ldots \subset T_{j-1} 
\subset T_j^{\mu}$, such that either $|T_j^{\mu}| > |T_j|$ or $V(T_j^{\mu})$ is not elementary with respect to $\mu$, which
contradicts either the maximum property satisfied by $T$ or the induction hypothesis (4.1) on Theorem~\ref{thm:tech10}(i). 
We divide the proof into five parts; the assumption on the intersection of $P_{v_n}(\gamma_n, \delta_n, \sigma_n)$ 
and $T_n$ will only be used in the last part.

\vskip 2mm

{\bf (I)} In this part we exhibit some properties satisfied by supporting vertices $a$ with $m(a)=j$ and corresponding 
connecting colors in $T_j-T_{j-1}$, which allow us to restore missing color sets of these vertices except $v_n$ as 
under $\varphi_{j-1}$ later. 

Let $L$ be the set of all subscripts $s$ with $j\le s\le n$, such that $\Theta_s=PE$ and $m(v_s) \le j$, where
$v_s$ is the supporting vertex involved in iteration $s$.  

(1) For any $s, t\in L$ with $s< t$, we have $v_t \preceq v_s$. Consequently, $v_n\preceq  v_s$ and $m(v_s) =j$
for all $s\in L$. 

To justify this, let $\pi_{t-1}$, $S_t=\{\gamma_t, \delta_t\}$, and $f_t$ be the $(T_t,D_{t-1},\varphi_{t-1})$-stable coloring, 
the set of connecting colors, and the connecting edge, respectively, as specified in iteration $t$ of Algorithm 3.1, with 
$\Theta_t=PE$.  Recall that $\delta_{t}=\pi_{t-1}(f_t)$ is a defective color of $T_t$ with respect to $\pi_{t-1}$,
and $v_t$ is the largest vertex (in the order $\prec$) in $I[\partial_{\pi_{t-1}, \delta_{t}}(T_{t})]$.  Since $m(v_t) \le j \le s<t$,
we have $v_t\in V(T_{j}) \subseteq V(T_{t-1})$.  As $\pi_{t-1}$ is a $(T_{t-1}, D_{t-1}, \varphi_{t-1})$-stable coloring and $v_t$ 
is not the smallest vertex (in the order $\prec$) in $I[\partial_{\pi_{t-1}, \delta_t}(T_{t-1})]$, applying (4.1) and Theorem~\ref{thm:tech10}(v) 
to $\pi_{t-1}$, $T_{t-1}$, and $v=v_t$, we obtain $v_t \preceq v_s$. Hence (1) holds.  

(2)  For any $s, t\in L$ with $s\le t$, we have $\delta_{t}\notin \phibar_{s-1}(T_s)$ (so $\delta_t\ne \gamma_s$). Consequently,
$\delta_t \notin \phibar_{j-1}(T_j)$ for all $t\in L$. 

Assume the contrary: $\delta_{t}\in  \phibar_{s-1}(u)$ for some $u\in V(T_s)$. By Algorithm 3.1, $\delta_s \notin \phibar_{s-1}(T_s)$.
So $s<t$ and hence $V(T_s)$ is elementary with respect to $\varphi_{s-1}$ by (4.1) and Theorem~\ref{thm:tech10}(i). Let $v$ be an 
arbitrary vertex in $T_s-u$. Then $\delta_t \notin \phibar_{s-1}(v)$, so $v$ is incident to an edge $f$ with 
$\varphi_{s-1}(f)=\delta_t$. As described in Algorithm 3.1, $T_{s}$ is closed under $\varphi_{s-1}$ and thus $f$ is 
contained in $G[T_{s}]$. Hence $\varphi_{t-1}(f) = \phiv_{s -1}(f)=\delta_t$ by Lemma~\ref{samecolor} (for $\delta_{t}\in  \phibar_{s-1}(u)
\subseteq \phibar_{s-1}(T_s)$). From Lemma \ref{sc2} and the definitions of $\pi_{t-1}$ and $\pi_{t-1}'$ in Algorithm 3.1, we see that $\pi'_{t-1}$ is $(T_{t}, D_{t-1}, 
\phiv_{t-1})$-stable. By Lemma \ref{hku}(i), $\phibar_{s-1}(T_s)\cup D_{s-1}\subseteq \phibar_{t-1}(T_t)\cup D_{t-1}$. 
So $\pi'_{t-1}(f) =\phiv_{t-1}(f)=\delta_t$.  Since $f$ is contained in $G[T_s]$ and hence in $G[T_t]$, we have 
$v \not\in I[\partial_{\pi'_{t-1},  \delta_{t}}(T_{t})]$ for any vertex $v$ in $T_s-u$. In view of (1), $v_t \preceq v_s$, so
$T_t(v_t) \subseteq T_s(v_s) \subseteq T_s$. Therefore $v_t$ cannot be the supporting vertex of $T_{t}$ with respect to $\phiv_{t}$ and 
the connecting color $\delta_{t}$ (as $v_t$ is the maximum defective vertex of $T_t$ with corresponding defective color $\delta_t$ 
under $\pi'_{t-1}$ in Algorithm 3.1); this contradiction implies that $\delta_{t}\notin \phibar_{s-1}(T_s)$. Since $\gamma_s \in \phibar_{s-1}(T_s)$,
we conclude that $\delta_t\ne \gamma_s$. Finally, let $s$ be the smallest subscript in $L$. Then $\phibar_{j-1}(T_j)=\phibar_{s-1}(T_j)$
by Algorithm 3.1 (see Lemma~\ref{hku}(ii)). So $\delta_t \notin \phibar_{j-1}(T_j)$ and hence (2) is established. 

We partition $L$ into disjoint subsets $L_1, L_2, \dots, L_{\kappa}$, such that two subscripts $s, t\in L$ are in the same subset iff $v_s = v_t$.  
For $1\le i \le \kappa$, write $L_i=\{i_1,i_2, \ldots,i_{c(i)}\}$, where $i_1<i_2< \ldots <i_{c(i)}$, and let $w_i$ denote the 
common supporting vertex corresponding to $L_i$. For each $t \in L$, we have $v_t \notin V(T_{j-1})$ because $m(v_t)=j$ by (1). It follows that 
$w_i \notin V(T_{j-1})$ for $1\le i \le \kappa$.  Renaming subscripts if necessary, we may assume that $w_1\prec w_2\prec \ldots \prec w_{\kappa}$.  
By (1), we obtain

(3) $v_n = w_1$ (so $n=1_{c(1)}$) and $h_{c(h)} > h_{c(h)-1}> \ldots > h_1 > i_{c(i)} > i_{c(i)-1} > \ldots > i_1$ for any $1\le h<i \le \kappa$.

From (2) and Lemma \ref{uniquezang}(i) it is clear that  

(4) for any $1 \le i \le \kappa$, the colors in $\cup_{t\in L_i} S_t$ are 
\begin{center}
\vspace{-2mm}
$\gamma_{i_1}, \gamma_{i_2}=\delta_{i_1}, \gamma_{i_3}=\delta_{i_2}, \ldots, \gamma_{i_{c(i)}}=\delta_{i_{c(i)-1}}, \delta_{i_{c(i)}},$ 
\vspace{-2mm}
\end{center}
which are distinct. 

From (4.2), Theorem~\ref{thm:tech10}(iii), (1), and (4) we deduce that  

(5) for any $s,t\in L$ with $s<t$, the intersection $S_s\cap S_t\neq\emptyset$ iff $s$ and $t$ are two consecutive subscripts in the same $L_i$ for some 
$1\leq i\leq\kappa$; in this case, $S_s\cap S_t=\{\gamma_t\}=\{\delta_s\}$.

\vskip 2mm
{\bf (II)} In this part we derive some properties satisfied by $\sigma_{n}$ and establish a result on the so-called strong interchangeability 
property, which enable us to keep the rest of $T_j$ ``stable" while restoring missing color sets of supporting vertices in $T_j-
V(T_j(v_n))$ as under $\varphi_{j-1}$. 

For each $t$ with $1\le t \le n-1$ and $\Theta_t = PE$, let $\epsilon(t)$ be the smallest subscript $r>t$ such that 
$\Theta_r \ne RE$.  This $\epsilon(t)$ is well defined and $\epsilon(t) \le n$, as $\Theta_n = PE \ne RE$.  
Given a coloring $\phiv$ and two colors $\alpha$ and $\beta$, we say that  $\alpha$ and $\beta$ are $T_t$-{\em strongly interchangeable}\label{stronglyinterchange} 
($T_t$-SI) under $\varphi$ if for each vertex $v$ in $T_t-v_t$, the chain $P_v(\alpha, \beta, \phiv)$ is an $(\alpha, \beta)$-cycle avoiding 
$v_t$ and fully contained in $G[T_{\epsilon(t)}]$ (equivalently, $V(P_v(\alpha, \beta, \phiv)) \subseteq V(T_{\epsilon(t)})$). 

Recall that  $\alpha$ and $\beta$ are called $T_t$-interchangeable under $\varphi$ if there is at most one $(\alpha,\beta)$-path with respect 
to $\varphi$ intersecting $T_t$; that is, all $(\alpha, \beta)$-chains intersecting $T_t$ are $(\alpha, \beta)$-cycles, with possibly one
exception. Therefore, if $\alpha$ and $\beta$ are $T_t$-SI under $\varphi$, then they are $T_t$-interchangeable 
under $\varphi$.

The following observations reveal some connections between colorings $\sigma_n$ and $\varphi_{j-1}$.  

\vskip 2mm
\noindent {\bf Claim 4.1}. {\em The coloring $\sigma_n$ satisfies the following properties: 
\begin{itemize}
\vspace{-2mm}
\item[(a1)] $\sigma_n$ is $(T_j(v_n)-v_n,D_{j-1},\varphi_{j-1})$-stable;
\vspace{-2mm}
\item [(a2)]  $\sigma_n(f) = \phiv_{j-1}(f)$ for all edges $f$ in $G[T_{j}]$ with $\phiv_{j-1}(f) \in \phibar_{j-1}(T_j)
\cup D_{j-1}$; in particular, this equality holds for all edges on $T_j$;  
\vspace{-2mm}
\item [(a3)] $\overline{\sigma}_n(v) = \overline{\phiv}_{j-1}(v)$ for all $v\in V(T_j) - \{w_1, w_2, \ldots, w_{\kappa}\}$; 
\vspace{-2mm} 
\item [(a4)] $\sigmabar_{n}(w_i) \cap (\cup_{t \in L_{i}} S_t) =\{\delta_{i_{c(i)}}\}$ and $\phibar_{j-1}(w_i) 
= (\sigmabar_n(w_i) - \{\delta_{i_{c(i)}} \} )\cup \{\gamma_{i_1}\}$ for each $i=1, 2, \dots, \kappa$; 
\vspace{-2mm} 
\item [(a5)] for any $t\in L-\{n\}$, the colors $\gamma_t$ and $\delta_t$ are $T_t$-SI under $\sigma_n$.
\end{itemize}}

To justify this claim, observe that $\phibar_{j-1}(T_j) \cup D_{j-1}\subseteq \phibar_{n-1}(T_n) \cup D_{n-1}\subseteq  \phibar_{n}(T_n) \cup D_{n}$ 
by Lemma \ref{hku}(i) and that $\cup_{t \in L} S_t \subseteq \cup_{t\le n} S_t \subseteq \phibar_{n}(T_n) \cup D_{n}$. 
Since $\sigma_n$ is a $(T_n, D_n, \varphi_n)$-stable coloring,  it suffices to prove $(a1)$-$(a5)$ for $\phiv_n$ (instead of $\sigma_n$).  

Clearly, $(a1)$ follows from (4.2) and Theorem \ref{thm:tech10}(iii), and $(a3)$ follows from Lemma \ref{hku}(ii). 

$(a2)$ By Lemma \ref{samecolor}, we have $\phiv_n(f) = \phiv_{j-1}(f)$ for all edges $f\in G[T_j]$ with $\phiv_{j-1}(f) \in 
\phibar_{j-1}(T_j)\cup D_{j-1}$. By Lemma \ref{hku}(iv), each edge $f$ on $T_j$ satisfies $\phiv_{j-1}(f) \in \phibar_{j-1}(T_j) 
\cup D_{j-1}$, so the equality $\varphi_n(f) = \phiv_{j-1}(f)$ holds for all edges $f$ on $T_j$.

$(a4)$ By Lemma~\ref{uniquezang}(ii), we obtain $\phibar_{n}(w_i) \cap (\cup_{t \in L_{i}} S_t) =\{\delta_{i_{c(i)}}\}$ (with $w_i$ in place of $u$).
Since $L_i$ consists of all subscripts $t$ with $j\le t \le n$, such that $v_t=w_i$ and $\Theta_t=PE$, there hold $\phibar_{j-1}(w_i)=\phibar_{i_1-1}(w_i)$ 
and  $\phibar_n(w_i)=\phibar_{i_{c(i)}} (w_i)$ by Lemma~\ref{hku}(ii). Furthermore, $\phibar_{i_1-1}(w_i)=(\phibar_{i_{c(i)}}(w_i)-\{\delta_{i_{c(i)}}\}) 
\cup \{\gamma_{i_1}\}$ by Lemma \ref{uniquezang}(iii) (with $w_i$ in place of $u$). So $\phibar_{j-1}(w_i)=(\phibar_n(w_i)-\{\delta_{i_{c(i)}}\}) \cup \{\gamma_{i_1}\}$ for 
$1\le i \le \kappa$.  

$(a5)$ Let $t\in L-\{n\}$. Then $t<n$. By the induction hypothesis (4.1) on Theorem~\ref{thm:tech10}(ii), $\gamma_t$ and $\delta_t$ are $T_{t+1}$- and hence
$T_t$-interchangeable under $\phiv_t$. So all but at most one $(\gamma_t, \delta_t)$-chains intersecting $T_t$ under $\phiv_t$ are $(\gamma_t, \delta_t)$-cycles.
According to Algorithm 3.1,  $P_{v_t}(\gamma_t, \delta_t, \phiv_t)$ is a path containing only one vertex $v_t$ from $T_t$.
Hence, for each vertex $v$ in $T_t-v_t$, $P_v(\gamma_t, \delta_t, \phiv_t)$ is a $(\gamma_t, \delta_t)$-cycle avoiding $v_t$. Since  $RE$ has 
priority over $PE$ and $SE$ in Algorithm 3.1, $P_v(\gamma_t, \delta_t, \phiv_t)$ is fully contained in 
$G[T_{\epsilon(t)}]$, for otherwise, we would have $\Theta_{\epsilon(t)}=RE$, contradicting the definition of $\epsilon(t)$.
It follows that $\gamma_t$ and $\delta_t$ are $T_t$-SI under $\phiv_t$. By Lemma~\ref{samecolor} (with respect to $T_{\epsilon(t)}$),  
we obtain $\phiv_t (f) = \phiv_n(f)$ for each edge $f$ on $P_{v}(\gamma_t, \delta_t, \phiv_t)$, because $\{\gamma_t, \delta_t\}
\subseteq \cup_{i\le r-1} S_i \subseteq \phibar_{r-1}(T_{r-1}) \cup D_{r-1} \subseteq \phibar_{r-1}(T_r) \cup D_{r-1}$, where 
$r=\epsilon(t)$. Therefore $\gamma_t$ and $\delta_t$ are $T_t$-SI under $\phiv_n$ as well. This establishes Claim 4.1. 

\vskip 2mm
The following technical statement will be used repeatedly in our proof. 

\vskip 2mm
\noindent {\bf Claim 4.2}. {\em Let $t\in L_i$ for some $1\le i \le \kappa$ and let $P$ be an arbitrary $(\gamma_t, \delta_t)$-path. If 
the connecting colors $\gamma_t, \delta_t$ are $T_t$-SI under a coloring $\phiv\in \CC(G-e)$ then, for any $s\in L_h$ with  $h \ne i$ 
or $s<  t$, the colors $\gamma_s, \delta_s$ are $T_s$-SI under $\phiv^* = \phiv/P$, provided that $\gamma_s, \delta_s$ are $T_s$-SI 
under $\phiv$. }

\vskip 2mm
To justify this, we assume that $\gamma_s, \delta_s$ are $T_s$-SI under coloring $\phiv$. For each $v\in V(T_s-v_s)$, we propose to show that 
$P_v(\gamma_s, \delta_s, \phiv^*)=P_v(\gamma_s, \delta_s, \phiv)$, which is a $(\gamma_s, \delta_s)$-cycle 
avoiding $v_s$ and fully contained in $G[T_{\epsilon(s)}]$. Consequently, $\gamma_s, \delta_s$ are also $T_s$-SI under coloring $\phiv^*$.  

If $h\ne i$, then $\{\gamma_s, \delta_ s\} \cap \{\gamma_t, \delta_t\} = \emptyset$ by (5). In this case, clearly $P_v(\gamma_s, \delta_s, \phiv^*) = P_v(\gamma_s, \delta_s, \phiv)$. 
So we assume that $h =i$ and $s< t$. Now $v_s=v_t=w_i$.  Observe that $P$ contains at most one vertex $v_t$ from $T_t$, because $\gamma_t, \delta_t$ 
are $T_t$-SI under $\phiv$. Furthermore, $\epsilon(s) \le t$, because $s <t$ and $\Theta_t =PE \ne RE$.  As $\gamma_s, \delta_s$ are $T_s$-SI under coloring $\phiv$, 
the chain $P_v(\gamma_s, \delta_s, \phiv)$ is a cycle avoiding $v_s$ and fully contained in $G[T_{\epsilon(s)}] \subseteq G[T_t]$, so it is disjoint 
from $P$. Thus $P_v(\gamma_s, \delta_s, \phiv)$ is still a $(\gamma_s, \delta_s)$-cycle avoiding $v_s$ and fully contained in $G[T_{\epsilon(s)}] 
\subseteq G[T_t]$ under $\phiv^*$ and hence $P_v(\gamma_s, \delta_s, \phiv^*) = P_v(\gamma_s, \delta_s, \phiv)$, as desired.  

\vskip 2mm

{\bf (III)} With the preparations made in the first two parts, now we can move on to the aforementioned restoration of missing color sets at
certain supporting vertices. 

Write $L^* = L-L_1$.  By (3), the subscripts in $L^*$ satisfies $h_{c(h)} > h_{c(h)-1}> \ldots > h_1 > i_{c(i)} > i_{c(i)-1} > \ldots > i_1$ 
for any $2\le h<i \le \kappa$. So $2_{c(2)}$ (resp. $\kappa_1$) is the largest (resp. smallest) subscripts in $L^*$. Starting from $\sigma_n$
and following the decreasing order of subscripts $t$ in $L^*$, we perform a sequence of $(\gamma_t, \delta_t)$-Kemple changes at $v_t$ for 
all $t\in L^*$ and get a new coloring in $\CC(G-e)$, under which each $w_i$, for $i \ge 2$, has the same set of missing colors as under 
$\phiv_{j-1}$. A detailed description of the algorithm is given below.

\vskip 1mm

{\flushleft \bf (A)} \hskip 1mm Let $I=\emptyset$ and $\sigma= \sigma_n$. While $I \ne L^*$, do: let $t$ be the largest member of $L^* -I$ and set 
\[
{\bf A}(t) : \qquad  \sigma = \sigma/P_{v_t}(\gamma_t, \delta_t, \sigma) \quad \mbox{ and }\quad  I = I\cup\{t\}. 
\]

Let us make some observations about this algorithm.

\vskip 1mm
(6) Let $I, t, \sigma$ be as specified in Algorithm $(A)$ before performing the iteration $A(t)$. Then 
$P_{v_t}(\gamma_t, \delta_t, \sigma)$ is a path containing precisely one vertex $v_t$ from $T_t$,
with $\delta_t \in \overline{\sigma}(v_t)$. Furthermore, let $\sigma'= \sigma / P_{v_t} (\gamma_t, \delta_t, \sigma)$
and $I'= I\cup \{t\}$ denote the objects generated in the iteration $A(t)$. Then for any $s\in L-\{n\}-I'$, 
the colors $\gamma_s$ and $\delta_s$ are $T_s$-SI under the coloring $\sigma'$. 

To justify this, recall that 

(7) $\delta_{i_{c(i)}}\in\overline{\sigma}_n(w_i)$ for each $2\leq i\leq \kappa$ by $(a4)$ in Claim 4.1 and 

(8) for any $s\in L-\{n\}$, the colors $\gamma_s$ and $\delta_s$ are $T_s$-SI under $\sigma_n$ by $(a5)$ in Claim 4.1. 

\noindent In particular, (8) holds for $t=2_{c(2)}$, the largest subscript in $L^*$, which implies that now
$P_{v_t}(\gamma_t, \delta_t, \sigma_n)$ is a path containing precisely one vertex $v_t=w_2$ from $T_t$, 
with $\delta_t \in \overline{\sigma}_n(v_t)$ by (7). Keep in mind that this $P_{v_t}(\gamma_t, \delta_t, \sigma_n)$
is the first path employed in Algorithm (A).

As the algorithm proceeds in the decreasing order of subscripts in $L^*$, using (4), (5), (7), (8)
and applying Claim 4.2 repeatedly, we see that (6) is true.  

\vskip 3mm
\noindent {\bf Claim 4.3.} {\em
Let $\varrho_1$ denote the coloring $\sigma$ output by Algorithm $(A)$.  Then the following statements hold:
\begin{itemize}
	\vspace{-2mm}
	\item[(b1)] $\varrho_1$ is $(T_j(v_n)-v_n,D_{j-1},\varphi_{j-1})$-stable;
	\vspace{-2mm}
	\item[(b2)] $\overline{\varrho}_1 (v)=\phibar_{j-1}(v)$ for all $v\in V(T_j-v_n)$,
	$\overline{\varrho}_1(v_n)=\sigmabar_n(v_n)$, and $\varrho_1(f)=\sigma_n(f) =\varphi_{j-1}(f)$ for all edges $f$ on $T_j$;
	\vspace{-2mm}
	\item[(b3)] for any edge $f\in E(G)$, if $\varrho_1(f)\ne \sigma_n(f)$, then $f$ is not contained in $G[T_j]$ and 
	$\{\sigma_n(f), \varrho_1(f)\} \\ \subseteq \cup_{t \in L^*} S_t$;  and
	\vspace{-2mm}
	\item[(b4)] for any $i\in L_1-\{n\}$ (so $v_i=v_n$), the colors $\gamma_i$ and $\delta_i$ are
	$T_i$-SI under $\varrho_1$. 
\end{itemize}}		

To justify this claim, recall from (6) that 

(9) at each iteration $A(t)$ of Algorithm $(A)$,  the chain $P_{v_t}(\gamma_t, \delta_t, \sigma)$ is a path 
containing precisely one vertex $v_t$ from $T_t$, with $\delta_t \in \overline{\sigma}(v_t)$. 

By (3) and the definitions of $L$ and $w_i$'s, we have

(10) $v_n =w_1 \prec w_i$ for all $i\ge 2$. Besides, $v_n\prec v_t$ and $T_j \subset T_t$ for each iteration $A(t)$ of Algorithm $(A)$.  

It follows from (9) and (10) that $\overline{\sigma}(v) = \overline{\sigma}_n(v)$ for each $v\in V(T_j(v_n)-v_n)$ 
and $\sigma(f) = \sigma_n(f)$ for all edges $f$ incident to $T_j(v_n) - v_n$ during each iteration of Algorithm (A). So $\sigma$
and hence $\varrho_1$ is a $(T_j(v_n) - v_n, D_{j-1}, \sigma_n)$-stable coloring.  By (4.2) and Theorem \ref{thm:tech10}(iii), 
$\sigma_n$ is $(T_j(v_n)-v_n,D_{j-1},\varphi_{j-1})$-stable. From Lemma \ref{sc2} we deduce that $\varrho_1$ is 
$(T_j(v_n)-v_n,D_{j-1},\varphi_{j-1})$-stable. So $(b1)$ holds. 

By $(a4)$ in Claim 4.1, we have

(11) $\phibar_{j-1}(w_i)= (\overline{\sigma}_n(w_i)- \{\delta_{i_{c(i)}}\}) \cup \{\gamma_{i_1}\}$ for each vertex $w_i$ with $i\ge 2$.

Recall that $S_p \cap S_q =\emptyset$ whenever $p$ and $q$ are contained in different $L_i$'s by (5). After executing Algorithm $(A)$, 
using (4) and Lemma \ref{uniquezang}(iii) (more precisely, the same argument), we obtain $\overline{\varrho}_1(w_i)=(\overline{\sigma}_n(w_i)-\{\delta_{i_{c(i)}}\}) 
\cup  \{\gamma_{i_1}\}$, so $\overline{\varrho}_1(w_i)=\phibar_{j-1}(w_i)$ for $i\ge 2$ by (11). Combining this with $(a3)$ in 
Claim 4.1, we see that $\overline{\varrho}_1 (v)=\phibar_{j-1}(v)$ for all $v\in V(T_j-v_n)$. By (6), the 
path $P_{v_t}(\gamma_t, \delta_t, \sigma)$ involved in each iteration $A(t)$ of Algorithm (A) is disjoint from
$v_n=w_1$. So $\overline{\varrho}_1(v_n)= \overline{\sigma}_n(v_n)=\phibar_n(v_n)$. In view of (9) and (10), we get 
$\sigma(f) = \sigma_n(f)$ for all edges $f$ on $T_j$ at each iteration $A(t)$ of Algorithm $(A)$.  Hence $\varrho_1(f)
=\sigma_n(f) =\varphi_{j-1}(f)$ for all edges $f$ on $T_j$, where the second equality follows from $(a2)$ in Claim 4.1.
Thus $(b2)$ is established.

Since the Kempe changes performed in Algorithm $(A)$ only involve edges outside $G[T_j]$ and colors in $\cup_{t\in L^*} S_t$ by the first half
of (6), we immediately get $(b3)$. Clearly, $(b4)$ follows from the second half of $(6)$. This proves Claim 4.3. 

\vskip 2mm
By analyzing two cases in the last two parts, we now demonstrate that the desired coloring can indeed be obtained by making $\gamma_{1_1}$ missing 
at a certain vertex $u$ (to be introduced) outside $T_j$ but inside a closure of $T_{j}(v_n)$. 

Consider the coloring $\varrho_1 \in {\cal C}^k(G-e)$ described in Claim 4.3.  By $(b1)$, $\varrho_1$ is $(T_j(v_n)-v_n,D_{j-1},
\varphi_{j-1})$-stable, so it is a $(T_{j-1},D_{j-1},\varphi_{j-1})$-stable coloring and hence is a $\phiv_{j-1}\bmod T_{j-1}$ coloring by
(4.1) and Theorem \ref{thm:tech10}(vi), which implies that every ETT corresponding to $(\varrho_1, T_{j-1})$ satisfies MP.  By $(b2)$, we have 
$\overline{\varrho}_1(v)=\phibar_{j-1}(v)$ for each $v\in V(T_j(v_n)-v_n)$ and $\varrho_1(f)=\varphi_{j-1}(f)$ 
for any edge $f$ on $T_j(v_n)$. Thus $T_j(v_n)$ is an ETT satisfying MP under $\varrho_1$. Let $T_j'$ be a closure of $T_{j}(v_n)$ under $\varrho_1$. 
(We point out that the first edge added to $T_j'-T_{j}(v_n)$ by TAA is incident to $V(T_n(v_n)-v_n)$ and colored with 
$\delta_n$ under $\varrho_1$ by Lemma \ref{hku}(v), $(b2)$ and $(b3)$, though we do not need this in our proof.) 
Then 

(12) $T_j'$ is an ETT satisfying MP under $\varrho_1$. Hence $V(T_j')$ is elementary with respect to $\varrho_1$ by (4.1) and Theorem 
\ref{thm:tech10}(i) (as $j \le n$). 

\vskip 2mm
Depending on the intersection of $\overline{\varrho}_1(T_j'-v_n)$ and $\cup_{i\in L_1}S_i$, we consider two cases.
\vskip 1mm

{\bf (IV)} This part is devoted to the study of the situation when the intersection is nonempty. 

\vskip 1mm
{\bf Case 1.} $\overline{\varrho}_1(T_j'-v_n)\cap (\cup_{i\in L_1}S_i)\ne \emptyset$.

\vskip 1mm
Let $u$ be the smallest vertex (in the order $\prec$) in $T_j'-v_n$ (so $u \ne v_n$), such that $\overline{\varrho}_1(u)
\cap (\cup_{i\in L_1}S_i)\neq \emptyset$. By (12), $V(T_j')$ is elementary with respect to
$\varrho_1$. Since $\delta_n\in \overline{\varphi}_n(v_n) =\overline{\varrho}_1(v_n)$ by $(b2)$, we obtain
$\delta_n\notin \overline{\varrho}_1(T_j'-v_n)$; in particular, $\delta_n\notin \overline{\varrho}_1(u)$. Hence,
by (4) and the definition of $u$, there exists a minimum member $r$ (as an integer) of $L_1$, such that
$\gamma_r \in  \overline{\varrho}_1(u) \cap (\cup_{i\in L_1}S_i)$. Since $m(v_r)=j$ by $(1)$, there holds
$r \geq j$. We propose to show, by using $\gamma_r$, that

(13) $u\in V(T_j')-V(T_j)$.

Indeed, if $r=1_1$, then $\gamma_r\in \phibar_{r-1}(v_n)=\phibar_{j-1}(v_n)$ by Algorithm 3.1. Since $V(T_j)$ is elementary with respect to 
$\phiv_{j-1}$ by (4.1) and Theorem~\ref{thm:tech10}(i) (for $j\leq n$),  we have $\gamma_r\notin \phibar_{j-1}(T_j - v_n)$. If $r > 1_1$, 
then $\gamma_r = \delta_t$ for some $t\in L_1$ by (4). Note that $\delta_t\notin \phibar_{j-1}(T_j)$ by (2). So we also have 
$\gamma_r\notin \phibar_{j-1}(T_j - v_n)$. It follows from $(b2)$ that $\gamma_r\notin \rhobar_1(T_j-v_n)$ in either subcase.  As 
$u\neq v_n$ and $\gamma_r\in \rhobar_1(u)$, we obtain $u \notin V(T_j)$. This proves (13).

(14) $\rhobar_1(T'_j(u)-u)\cap (\cup_{i \in L_1} S_i-\{\delta_n\}) = \emptyset$. 

By the minimality assumption on $u$, we have  $\overline{\varrho}_1(T_j'(u)-\{v_n,u\}) \cap (\cup_{i\in L_1}S_i)=\emptyset$. 
Using Lemma \ref{uniquezang}(ii), we obtain $\phibar_n(v_n) \cap (\cup_{i \in L_1} S_i)=\{\delta_n\}$. It follows
from $(b2)$ in Claim 4.3 that $\overline{\varrho}_1(v_n) \cap (\cup_{i \in L_1} S_i)=\{\delta_n\}$. Thus
(14) holds. 

Let $r$ be the subscript as defined above (13). Then $r=1_p$ for some $1\le p \le c(1)$. By (4), we have 
$\gamma_r=\gamma_{1_p}=\delta_{1_{p-1}}$ if $p \ge 2$. Let $L_1^* =\{1_1, 1_2, \ldots, 1_{p-1}\}$ (so $L_1^*=\emptyset$ if $p=1$). 
Since $1_{p-1}<1_p=r\le n$, we have $n \notin L_1^*$. Observe that

(15) $\delta_n \notin \cup_{i\in L_1^*} S_i$ and $\overline{\varrho}_1(v_n) \cap (\cup_{i\in L_1^*} S_i) = \emptyset$.

Indeed, by $(b2)$ in Claim 4.3 and Lemma \ref{uniquezang}(ii), we obtain $\overline{\varrho}_1(v_n)=
\phibar_n(v_n)$ and $\phibar_n(v_n) \cap (\cup_{i \in L_1} S_i)=\{\delta_n\}$.  As $n \notin L_1^*$, from (4)
we see that $\delta_n \notin \cup_{i\in L_1^*} S_i$. So $\phibar_n(v_n) \cap (\cup_{i \in L_1^*} S_i)=\emptyset$.
Hence (15) follows.

\vskip 2mm
We construct a new coloring from $\varrho_1$ by using the following algorithm.

{\flushleft \bf (B)} Let $I = \emptyset$ and $\varrho= \varrho_1$. While $I \ne L_1^*$, do: let $t$ be the largest 
member of $L_1^* -I$ and set 
\[
{\bf B}(t):  \qquad \varrho  = \varrho/P_u(\gamma_t, \delta_t, \varrho) \quad \mbox{  and  } \quad I = I\cup \{t\}. 
\] 

Let us exhibit some properties satisfied by this algorithm.  

(16) Let $I, t, \varrho$ be as specified in Algorithm $(B)$ before performing the iteration $B(t)$. Then
$\delta_t \in \overline{\varrho}(u)$, and $P_u(\gamma_t, \delta_t, \varrho)$ is a path containing at most one vertex 
$v_n$ from $T_t$, but $v_n$ is not an end of $P_u(\gamma_t, \delta_t, \varrho)$.  Furthermore, let $\varrho'  = 
\varrho/P_u(\gamma_t, \delta_t, \varrho)$ and $I'= I\cup \{t\}$ denote the objects generated in the iteration 
$B(t)$. Then for any $s\in L_1^*-I'$, the colors $\gamma_s$ and $\delta_s$ are $T_s$-SI under the 
coloring $\varrho'$. 

To justify this, recall that $\delta_{1_{p-1}}=\gamma_r \in\rhobar_1(u)$ and 

(17) for any $i\in L_1-\{n\}$ (so $v_i=v_n$), the colors $\gamma_i$ and $\delta_i$ are $T_i$-SI under $\varrho_1$ 
by $(b4)$. 

\noindent In particular, (17) holds for $t=1_{p-1}$, the largest subscript in $L^*_1$, which implies that now $P_u(\gamma_t, \delta_t, \varrho_1)$ 
is a path containing at most one vertex $v_n$ from $T_t$, but $v_n$ is not an end of $P_u(\gamma_t, \delta_t, \varrho_1)$ by (15). 
Keep in mind that this $P_u(\gamma_t, \delta_t, \varrho_1)$ is the first path employed in Algorithm (B).

Since the algorithm proceeds in the decreasing order of subscripts in $L^*_1$, using (4), (5), (15), (17),
and applying Claim 4.2 repeatedly, we see that (16) is true.  

\vskip 2mm
\noindent {\bf Claim 4.4.} {\em 
Let $\varrho_2$ denote the coloring $\varrho$ output by Algorithm $(B)$. Then the following statements hold:
\begin{itemize}
\vspace{-2mm}
\item[(c1)] $\varrho_2 $ is $(T_j(v_n)-v_n,D_{j-1},\varphi_{j-1})$-stable;
\vspace{-2mm}
\item[(c2)] $\overline{\varrho}_2 (v)=\overline{\varrho}_1(v)$ for all $v\in V(T_j\cup T_j'(u)-u)$ and $\varrho_2 (f)=\varrho_1(f)$ for all $f\in E(T_j\cup T_j'(u))$;
\vspace{-2mm}
\item[(c3)] $\gamma_{1_1}\in \overline{\varrho}_2 (u)$.
\end{itemize}	}

To justify this claim, recall from (16) that 

(18) at each iteration $B(t)$, the path $P_u(\gamma_t, \delta_t, \varrho)$ contains at most one vertex $v_n$ from 
$T_t$, but $v_n$ is not an end of $P_u(\gamma_t, \delta_t, \varrho)$. 

Since $T_j \subseteq T_t$, we have  $\rhobar(v) = \rhobar_1(v)$ for each $v\in V(T_j(v_n)-v_n)$  and 
$\varrho(f) = \varrho_1(f)$ for each edge $f$ incident to $T_j(v_n) - v_n$ during each iteration of Algorithm (B) 
by (18). It follows that $\rho$ and hence $\varrho_2$ is a $(T_j(v_n) - v_n, D_{j-1}, \varrho_1)$-stable coloring. 
By $(b1)$ in Claim 4.3, $\varrho_1$ is a $(T_j(v_n) -v_n, D_{j-1}, \phiv_{j-1})$-stable coloring. From 
Lemma \ref{sc2} we see that $(c1)$ holds. 

Similarly, from (18) we deduce that $\rhobar_2(v) = \rhobar_1(v)$ for all $v\in V(T_j)$ and $\varrho_2(f) = \varrho_1(f)$ 
for all $f\in E(T_j)$. By (14) and (15), we also have $\rhobar_1(T'_j(u)-u)\cap (\cup_{i \in L_1^*} S_i) = \emptyset$. 
So $T_j'(u)$ does not contain the other end of $P_u(\gamma_t, \delta_t, \varrho)$ at each iteration $B(t)$, and hence
$\rhobar_2(v) = \rhobar_1(v)$ for each $v\in  V(T_j'(u) -u)$. Since $T_{j}'$ is a 
closure of $T_j(v_n)$ under $\varrho_1$, from TAA we deduce that  $\varrho_1 \langle T_j'(u) - T_j(v_n) \rangle 
\cap (\cup_{i\in L_1^*}S_i) = \emptyset$. It follows that $\varrho(f) = \varrho_1(f)$ for all edges $f$ in $T_j'(u) - T_j(v_n)$ 
at each iteration $B(t)$. So $\varrho_2(f) = \varrho_1(f)$ for all edges $f$ in $T_j'(u) - T_j(v_n)$
and hence $(c2)$ holds.    

By (16), we have $\delta_t \in \overline{\varrho}(u)$ before each iteration $B(t)$. So $\gamma_t$ becomes a 
missing color at $u$ after performing iteration $B(t)$. It follows that $\gamma_{1_1}\in \overline{\varrho}_2 (u)$ (see (4)).
Hence $(c3)$ and therefore Claim 4.4 is established.  

\vskip 2mm

By $(c1)$ in Claim 4.4, $\varrho_2$ is $(T_j(v_n)-v_n,D_{j-1}, \varphi_{j-1})$-stable. So it is a $(T_{j-1},D_{j-1},\varphi_{j-1})$-stable coloring and hence is a
$\phiv_{j-1}\bmod T_{j-1}$ coloring by (4.1) and Theorem \ref{thm:tech10}(vi), which implies that every ETT corresponding to $(\varrho_2, T_{j-1})$ satisfies MP. 
By $(b2)$ and $(c2)$, we have $\varrho_2(f)=\varphi_{j-1}(f)$ for each edge $f$ on $T_j(v_n)$. So $T_j(v_n)$ an ETT satisfying MP under $\varrho_2$. Since 
$T'_j(u)$ is obtained from $T_j(v_n)$ by TAA under $\varrho_1$, it can also be obtained from $T_j(v_n)$ by TAA under $\varrho_2$ by $(c2)$. Thus $T'_j(u)$ 
is an ETT satisfying MP under $\varrho_2$ as well. 

In view of $(b2)$ and $(c2)$, we have $\overline{\varrho}_2(v)=\phibar_{j-1}(v)$ for all $v\in V(T_j-v_n)$, $\overline{\varrho}_2(v_n)=
\overline{\varphi}_{n}(v_n)$, and $\varrho_2 (f)=\varphi_{j-1}(f)$ for all $f\in E(T_j)$. Moreover, by Lemma \ref{uniquezang}(iii) and $(c3)$, 
we obtain $\overline{\varphi}_{j-1}(v_n) =\overline{\varphi}_{1_1-1}(v_n)  \subseteq \overline{\varphi}_{1_{c(1)}}(v_n) \cup \{\gamma_{1_1}\} = 
\overline{\varphi}_{n}(v_n) \cup \{\gamma_{1_1}\}=\overline{\varrho}_2(v_n) \cup \{\gamma_{1_1}\} \subseteq \overline{\varrho}_2(v_n) \cup 
\overline{\varrho}_2(u)$. Therefore we can further grow $T'_j(u)$ by adding all edges on $T_j$ but outside $G[T'_j(u)]$ using TAA under 
$\varrho_2$; let $T_j^1$ denote the resulting tree-sequence. Clearly, $T_j^1$ is an ETT satisfying MP under $\varrho_2$ and $V(T_j \cup T_j'(u))
\subseteq V(T_j^1)$, which contradicts MP satisfied by $T$ under $\varphi_n$, because $u \notin V(T_j)$ by (13).

\vskip 2mm

{\bf (V)} Let us give an analysis of the opposite situation, which is the last part of this long proof.
 
\vskip 1mm
{\bf Case 2.} $\overline{\varrho}_1(T_j'-v_n)\cap (\cup_{i\in L_1}S_i)= \emptyset$.

\vskip 1mm

Recall that $L_1=\{1_1,1_2,...,1_{c(1)}\}$. Set $S'=\cup_{i\in L_1}S_i$. Let us make some simple observations about
$T_j'$, $T_j$ and $T_n$.

\vskip 1mm
(19) $\overline{\varrho}_1 (T_j') \cap S'=\overline{\varrho}_1(v_n) \cap S'=\{\delta_n\}$ and $\varrho_1 \langle T_j' -T_j(v_n) \rangle \cap 
(S'-\{\delta_n\})=\emptyset$.

To justify this, note that $V(T_j')$ is elementary with respect to $\varrho_1$ by (12) and that
$\overline{\varrho}_1(v_n)=\phibar_n(v_n)$ by $(b2)$. By Lemma \ref{uniquezang}(ii), we have $\phibar_n(v_n) \cap  S'=
\{\delta_n\}$. So $\overline{\varrho}_1(v_n) \cap  S'=\{\delta_n\}$ and hence $\delta_n \notin \overline{\varrho}_1
(T_j'-v_n)$. By the hypothesis of the present case, we obtain $\overline{\varrho}_1 (T_j') \cap S'=\overline{\varrho}_1(v_n) \cap S'=\{\delta_n\}$. 
Since $T_j'$ is a closure of $T_{j}(v_n)$ under $\varrho_1$, from TAA we see that $\varrho_1 \langle T_j' -T_j(v_n) 
\rangle \cap (S'-\{\delta_n\})=\emptyset$. Hence (19) holds.

(20) $\partial_{\varrho_1, \gamma_n}(T_n)=\{f_n\}$ and $\partial_{\varrho_1, \delta_n}(T_n \cup T_j') =\emptyset$.

To justify this, note from Lemma \ref{hku}(v) that $\partial_{\varphi_n, \gamma_n}(T_n)=\{f_n\}$ and edges in $\partial_{\varphi_n, \delta_n}(T_n)$
are all incident to $V(T_n(v_n)-v_n)$. Since $\sigma_n$ is $(T_n, D_n, \phiv_n)$-stable, by $(b3)$ in Claim 4.3 and (5), we obtain 
$\partial_{\varrho_1, \alpha}(T_n)=\partial_{\sigma_n, \alpha}(T_n)=\partial_{\varphi_n, \alpha}(T_n)$ for $\alpha=\gamma_n, 
\delta_n$. In particular, $\partial_{\varrho_1, \gamma_n}(T_n)=\{f_n\}$. Since $T_j'$ is a closure of $T_j(v_n)$ under $\varrho_1$ and 
$\delta_n \in \phibar_n(v_n)=\overline{\varrho}_1(v_n)$ by $(b2)$, from TAA we see that $\partial_{\varrho_1, \delta_n}(T_n \cup T_j') =\emptyset$. 
Hence (20) is true.

\vskip 3mm
Consider the path $P_{v_n}(\gamma_n, \delta_n, \sigma_n)$ specified in the present lemma. By Algorithm 3.1, (4.1) and Theorem~\ref{thm:tech10}(i), 
we have $\delta_n, \gamma_n \notin \overline{\varphi}_{n-1}(T_n-v_n)$. So $\delta_n, \gamma_n \notin \overline{\varphi}_n(T_n-v_n)$ and hence
$\delta_n, \gamma_n \notin \overline{\sigma}_n(T_n-v_n)$.  It follows that the other end $x$ of $P_{v_n}(\gamma_n, \delta_n, \sigma_n)$ 
is outside $T_n$. Let $P$ denote $P_{v_n}(\gamma_n, \delta_n, \varrho_1)$. Then $P=P_{v_n}(\gamma_n, \delta_n, \sigma_n)$
by $(b3)$ and $(5)$. From the hypothesis of the present case, we deduce that $x$ is outside $T_j'-v_n$. Combining these two observations,
we see that $x$ is outside $T_n \cup T_j'$. Let $u$ be the vertex of $P$ such that the subpath $P[u,x]$ is a
$(T_n \cup T_j')$-exit path with respect to $\varrho_1$. At the beginning of our proof, we assume that $P_{v_n}(\gamma_n, \delta_n, \sigma_n)$
(and hence $P$) contains at least two vertices from $T_n$. So $u\ne v_n$. By (20), all edges in $E(P) \cap \partial(T_n \cup T_j')$ are colored 
by $\gamma_n$ under $\varrho_1$ and $f_n$ is the only edge in $\partial_{\varrho_1, \gamma_n}(T_n)$. Therefore $u$ is not incident to $f_n$ and furthermore

(21) $u \in V(T_j')-V(T_n)$.

Figure 3 gives an illustration of $P$ under $\varrho_1$.

\vskip 2mm
\begin{figure}[htpb]
	\vspace{-1mm}
	\centerline{\includegraphics[width=6cm]{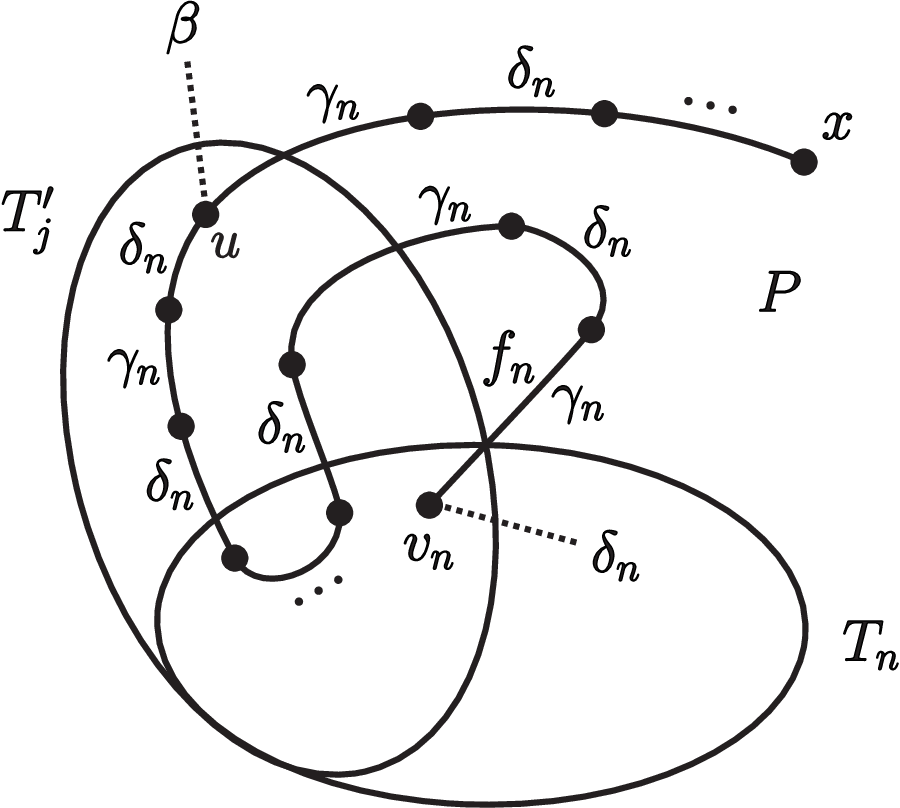}}
	\vspace{-1mm}\caption{The path $P$ under $\varrho_1$}
\end{figure}
\vskip 2mm

Let $\beta \in \overline{\varrho}_1(u)$. By the hypothesis of the present case, we have

(22) $\beta\notin S'$.

\noindent If $\beta \in \overline{\varrho}_1(T_j-V(T_j(v_n)))$, let $z$ be the smallest vertex in $T_j-V(T_j(v_n))$ in the order $\prec$ 
such that $\beta \in \overline{\varrho}_1(z)$; otherwise, let $z$ be the largest vertex of $T_j$ in the order $\prec$ (now $T_j(z)=T_j$).

(23) $\beta\notin \overline{\varrho}_1(T_j(z)-z)$ and $\beta\notin \varrho_1\langle T_j(z)-T_j(v_n) \rangle $.

By the definition of $z$, we have $\beta\notin \overline{\varrho}_1(T_j(z)-V(T_j(v_n))-z)$. Since $\beta \in \overline{\varrho}_1(u)$, by (12) and (21) 
we obtain $\beta\notin \overline{\varrho}_1(T_j(v_n))$. So $\beta\notin \overline{\varrho}_1(T_j(z)-z)$.  From $(b2)$ it follows that 
$\beta \notin \overline{\varphi}_{j-1}(T_j(z)-z-v_n)$ and $\beta\notin \overline{\sigma}_n(v_n)=\phibar_n(v_n)$. By Lemma \ref{uniquezang}(iii),  
$\phibar_{j-1}(v_n)= \phibar_{1_1-1}(v_n) \subseteq \phibar_{1_{c(1)}}(v_n) \cup \{\gamma_{1_1}\} = \phibar_n(v_n) \cup \{\gamma_{1_1}\}$.
Since $\beta\ne \gamma_{1_1}$ by (22), we get $\beta\notin\phibar_{j-1}(v_n)$. Hence $\beta \notin \overline{\varphi}_{j-1}(T_j(z)-z)$.
As $T_j(z)$ is obtained from $T_j(v_n)$ by TAA under $\varphi_{j-1}$, $\beta\notin \varphi_{j-1} \langle T_j(z)-T_j(v_n) \rangle $. 
Therefore $\beta\notin \varrho_1 \langle T_j(z)-T_j(v_n) \rangle $ by $(b2)$. This justifies (23). 

\vskip 2mm
\noindent {\bf Claim 4.5.} {\em
There exists a coloring $\varrho_3 \in {\cal C}^k(G-e)$ with the following properties:
\begin{itemize}
\vspace{-2mm}
\item[(d1)] $\varrho_3$ is $(T_j(v_n)-v_n,D_{j-1},\varphi_{j-1})$-stable;
\vspace{-2mm}
\item[(d2)] $\overline{\varrho}_3(v)=\overline{\varrho}_1(v)$ for all $v\in V(T_j(z)\cup T_j'(u))-\{u,z\}$ and $\varrho_3(f)={\varrho}_1(f)$ 
for all $f\in E(T_j(z)\cup T_j'(u))$. Furthermore, $\delta_n\in \overline{\varrho}_3(z)$ if $\beta \in \overline{\varrho}_1(z)$; and
\vspace{-2mm}		
\item[(d3)] $\gamma_{1_1}\in \overline{\varrho}_3(u)$.
\end{itemize}	}

{\bf (Assuming Claim 4.5)}  By $(d1)$ in Claim 4.5, $\varrho_3$ is a $(T_j(v_n)-v_n,D_{j-1},\varphi_{j-1})$-stable
coloring. So it is a $(T_{j-1},D_{j-1},\varphi_{j-1})$-stable coloring and hence is a $\phiv_{j-1}\bmod T_{j-1}$ coloring by
(4.1) and Theorem \ref{thm:tech10}(vi), which implies that every ETT corresponding to $(\varrho_3, T_{j-1})$ satisfies MP. By $(b2)$ and $(d2)$, 
we have $\varrho_3(f)=\varrho_1(f)=\varphi_{j-1}(f)$ for each edge $f$ on $T_j(v_n)$. So $T_j(v_n)$ is an ETT satisfying MP under 
$\varrho_3$. Since $T'_j(u)$ is obtained from $T_j(v_n)$ by TAA under $\varrho_1$, it can also be obtained from $T_j(v_n)$ by 
TAA under $\varrho_3$ by $(d2)$. Thus $T'_j(u)$ is an ETT satisfying MP under $\varrho_3$ as well. 

In view of $(b2)$ and $(d2)$, we have $\overline{\varrho}_3(v)=\phibar_{j-1}(v)$ for all $v\in V(T_j(z)-\{v_n,z\})$, $\overline{\varrho}_3(v_n)=
\overline{\varphi}_{n}(v_n)$, and $\varrho_3 (f)=\varphi_{j-1}(f)$ for all $f\in E(T_j(z))$. Moreover, by Lemma \ref{uniquezang}(iii) and $(d3)$, 
we obtain $\overline{\varphi}_{j-1}(v_n) =\overline{\varphi}_{1_1-1}(v_n)  \subseteq \overline{\varphi}_{1_{c(1)}}(v_n) \cup \{\gamma_{1_1}\} = 
\overline{\varphi}_{n}(v_n) \cup \{\gamma_{1_1}\}=\overline{\varrho}_3(v_n) \cup \{\gamma_{1_1}\} \subseteq \overline{\varrho}_3(v_n) \cup 
\overline{\varrho}_3(u)$. Therefore we can further grow $T'_j(u)$ by adding all edges on $T_j(z)$ but outside $G[T'_j(u)]$ using TAA under 
$\varrho_3$; let $T_j^2$ denote the resulting tree-sequence. Clearly, $T_j^2$ is an ETT satisfying MP. So $V(T_j^2)$ is elementary with respect to
$\varrho_3$ by (4.1) and Theorem~\ref{thm:tech10}(i). If $z$ is the largest vertex of $T_j$ in the order $\prec$, then
$V(T_j \cup T_j'(u))= V(T_j(z) \cup T_j'(u)) \subseteq V(T_j^2)$, which contradicts MP satisfied by $T$, as $u \notin V(T_j)$ by (21); otherwise,
$\delta_n\in \overline{\varrho}_3(z) \cap \overline{\varrho}_3(v_n)$ by $(d2)$ and $(19)$, which contradicts the elementary property
satisfied by $V(T_j^2)$ under $\varrho_3$. 

\vskip 2mm
To prove Claim 4.5, we consider the coloring $\varrho_0=\varrho_1/(G-T_j',\beta,\delta_n)$. Since $T_j'$ is closed with
respect to $\varrho_1$ and $\{v_n, u\} \subseteq V(T_j')$, no boundary edge of $T_j'$ is colored by $\beta$ or $\delta_n$ under $\varrho_1$ (see (19)). 
So $\varrho_0$ is $(T_j', D_{j-1}, \varrho_1)$-stable and hence is $(T_j(v_n)-v_n, D_{j-1}, \varrho_1)$-stable. Clearly, 
$P_u(\gamma_n, \beta, \varrho_0) = P_u(\gamma_n, \delta_n, \varrho_1)$.  Thus $u$ is the only vertex shared by $P_u(\gamma_n, \beta, \varrho_0)$
and $T_n \cup T_j'$. Define $\mu_0 = \varrho_0/P_u(\gamma_n, \beta, \varrho_0)$.  

\vskip 2mm
\noindent {\bf Claim 4.6}. {\em The coloring $\mu_0$ satisfies the following properties:  
\begin{itemize}
\vspace{-2mm}
\item[(e1)] $\mu_0$ is a $(T_j(v_n)-v_n,D_{j-1},\varphi_{j-1})$-stable coloring;
\vspace{-2mm}
\item[(e2)] $\mubar_0(v)=\overline{\varrho}_1(v)$ for all $v\in V(T_j(z)\cup T_j'(u))-\{u,z\}$ and 
$\mu_0(f)={\varrho}_1(f)$ for all $f\in E(T_j(z)\cup T_j'(u))$. Furthermore, $\delta_n\in \overline{\mu}_0(z)$
if $\beta \in \overline{\varrho}_1(z)$; 
\vspace{-2mm}		
\item[(e3)] $\gamma_n = \delta_{1_{c(1)-1}}\in \mubar_0(u)$ and $\beta\notin \mubar_0(u)$; 
\vspace{-2mm}		
\item[(e4)] for any $t\in L_1-\{n\}$, the colors $\gamma_t$ and $\delta_t$ are $T_t$-SI under $\mu_0$; and
\vspace{-2mm}		
\item[(e5)] $\overline{\mu}_0 (T_j'-u) \cap S'= \overline{\mu}_0 (v_n) \cap S' = \{\delta_n\}$ and 
$\mu_0 \langle T_j' -T_j(v_n) \rangle \cap (S'-\{\delta_n\})=\emptyset$.
\end{itemize}}	

To justify this, recall that $\varrho_1$ is $(T_{j}(v_n)-v_n,D_{j-1},\varphi_{j-1})$-stable by $(b1)$. By the definitions
of $\varrho_0$ and $\mu_0$, the transformation from $\varrho_1$ to $\mu_0$ only changes colors on some edges disjoint 
from $V(T_j(v_n))$. So $(e1)$ holds.  Statement $(e3)$ follows instantly from the definition of $\mu_0$. Note that
$\delta_n, \beta \notin \cup_{t\in L_1-\{n\}}S_t$ by (4), (5) and (22), and that $T_{\epsilon(t)} \subseteq T_n$ for each 
$t\in L_1 - \{n\}$.  Furthermore, $P_u(\gamma_n, \beta, \varrho_0)$ is disjoint from $V(T_n)$. So $(e4)$ can be deduced 
from $(b4)$ immediately. Using (19) and the definitions of $\varrho_0$ and $\mu_0$, we obtain $(e5)$. 

It remains to prove $(e2)$. Recall from (23) that $\beta\notin \overline{\varrho}_1(T_j(z)-z)$ and $\beta\notin \varrho_1\langle 
T_j(z)-T_j(v_n) \rangle $. By (2), we obtain $\delta_n\notin\phibar_{j-1}(T_j)$ and hence $\delta_n\notin 
\varphi_{j-1} \langle T_j(z)-T_j(v_n) \rangle $ by TAA. From $(b2)$ we deduce that $\delta_n\notin \overline{\varrho}_1(T_j(z)-v_n)$
and $\delta_n\notin \varrho_1 \langle T_j(z)-T_j(v_n) \rangle$.  From the definition of $\varrho_0$ and $\mu_0$, we see that $(e2)$ holds. 
So Claim 4.6 is established.

\vskip 2mm

Let $L_1^*= L_1 -\{n\}$. We construct a new coloring from $\mu_0$ by using the following algorithm.

{\flushleft \bf (C)} Let $I= \emptyset$ and $\mu  = \mu_0$. While $I \ne L_1^*$, do: let $t$ be the largest member in 
$L_1^* -I$ and set 
\[
\mbox{ {\bf C}(t):} \qquad \mu = \mu/P_u(\gamma_t, \delta_t, \mu) \quad \mbox{ and } \quad 
I = I\cup \{t\}.
\]

Let $\varrho_3$ denote the coloring $\mu$ output by Algorithm $(C)$. We aim to show that $\varrho_3$ is as described in
Claim 4.5; our proof is based on the following statement. 

(24) Let $I, t, \mu$ be as specified in Algorithm $(C)$ before performing the iteration $C(t)$. Then
$\delta_t \in \overline{\mu}(u)$, and $P_u(\gamma_t, \delta_t, \mu)$ is a path containing at most one vertex 
$v_n$ from $T_t$, but $v_n$ is not an end of $P_u(\gamma_t, \delta_t, \mu)$.  Furthermore, let $\mu'  = 
\mu/P_u(\gamma_t, \delta_t, \mu)$ and $I'= I\cup \{t\}$ denote the objects generated in the iteration 
$C(t)$. Then for any $s\in L_1^*-I'$, the colors $\gamma_s$ and $\delta_s$ are $T_s$-SI under the 
coloring $\mu'$. 

To justify this, observe that 

(25) $\overline{\mu}_0 (v_n) \cap (\cup_{i\in L_1^*} S_i)= \emptyset$ by (4), (5) and $(e5)$.

\noindent Furthermore, 

(26) for any $s\in L_1^*$, the colors $\gamma_s$ and $\delta_s$ are $T_s$-SI under $\mu_0$ by $(e4)$. 

\noindent In particular, (26) holds for $t=1_{c(1)-1}$, the largest subscript in $L^*_1$, which implies
that now $P_u(\gamma_t, \delta_t, \mu_0)$ is a path containing at most one vertex $v_t=v_n$ from $T_t$, 
but $v_n$ is not an end of $P_u(\gamma_t, \delta_t, \mu_0)$ by (25). In view of $(e3)$, we have $\delta_t \in 
\overline{\mu}_0(u)$. Keep in mind that this $P_u(\gamma_t, \delta_t, \mu_0)$ is the first path employed in
Algorithm (C).

As the algorithm proceeds in the decreasing order of subscripts in $L^*_1$, using (4), (5), (25), (26)
and applying Claim 4.2 repeatedly, we see that (24) is true.  

\vskip 2mm

To justify Claim 4.5, recall from (24) that 

(27) at each iteration $C(t)$, the path $P_u(\gamma_t, \delta_t, \mu)$ contains at most one vertex $v_n=v_t$ from 
$T_t$, but $v_n$ is not an end of $P_u(\gamma_t, \delta_t, \mu)$. 

Since $T_j \subseteq T_t$, we have $\overline{\mu}(v) = \overline{\mu}_0(v)$ for each $v\in V(T_j(v_n)-v_n)$ 
and $\mu(f) = \mu_0(f)$ for each edge $f$ incident to $T_j(v_n) - v_n$ during each iteration of Algorithm (C) 
by (27). It follows that $\mu$ and hence $\varrho_3$ is a $(T_j(v_n) - v_n, D_{j-1}, \mu_0)$-stable coloring. 
By $(e1)$ in Claim 4.6, $\mu_0$ is a $(T_j(v_n) -v_n, D_{j-1}, \phiv_{j-1})$-stable coloring. From 
Lemma \ref{sc2} we see that $(d1)$ holds. 

Since $T_j(z) \subseteq T_t$, from $(e2)$ and (27) we deduce that $\rhobar_3(v) = \overline{\mu}_0(v)$ for all $v\in V(T_j(z)-z)$ 
and $\varrho_3(f) = \mu_0(f)$ for all $f\in E(T_j(z))$.  By $(e5)$, we have $\overline{\mu}_0 (T_j'-u) \cap S'= \overline{\mu}_0 (v_n) \cap S' = 
\{\delta_n\}$ and $\mu_0 \langle T_j' -T_j(v_n) \rangle \cap (S'-\{\delta_n\})=\emptyset$. By (4) and (5), 
we obtain $\delta_n \notin \cup_{i \in L_1^*} S_i$. So at each iteration $C(t)$ the path  
$P_u(\gamma_t, \delta_t, \mu)$ neither contains any edge from $T_j'(u)$ nor terminate at a 
vertex in $T_j'(u)-u$. It follows that $\rhobar_3(v) = \overline{\mu}_0(v)$ for all $v\in V(T_j'(u)-u)$
and $\varrho_3(f) = \mu_0(f)$ for all edges $f$ in  $T_j'(u) -T_j(v_n)$. Hence $\overline{\varrho}_3(v)=
\overline{\mu}_0(v)$ for all $v\in V(T_j(z)\cup T_j'(u))-\{u,z\}$ and $\varrho_3 (f)={\mu}_0(f)$ for all $f\in 
E(T_j(z)\cup T_j'(u))$. Combining this with $(e2)$, we see that $(d2)$ holds. 

By (24), we have $\delta_t \in \overline{\mu}(u)$ before each iteration $C(t)$. So $\gamma_t$ becomes a 
missing color at $u$ after performing iteration $C(t)$. It follows that $\gamma_{1_1}\in \overline{\varrho}_3 (u)$ (see (4)).
Therefore $(d3)$ is established. This completes the proof of Claim 4.5 and hence of Lemma \ref{extension base}. \qed

\vskip 3mm

The following lemma asserts that supporting and extension vertices are subject to some order.

\begin{lemma}\label{basic2}
(Assuming (4.1) and (4.3)) Theorem~\ref{thm:tech10}(v) holds for all ETTs with $n$ rungs and satisfying MP; 
that is, for any $(T_n, D_n, \phiv_n)$-stable coloring $\sigma_n$ and any defective color $\delta$ of $T_n$ 
with respect to $\sigma_n$, if $v$ is a vertex but not the smallest one (in the order $\prec$)  
in $I[\partial_{\sigma_n, \delta}(T_n)]$, then $v \preceq v_i$ for any supporting or extension vertex  
$v_i$ with $m(v) \le i$.
 
\end{lemma}

{\bf Proof.} By the hypothesis of Theorem~\ref{thm:tech10}, $T$ is an ETT with the corresponding Tashkinov series 
$\TT=\{(T_i, \phiv_{i-1}, S_{i-1}, F_{i-1}, \Theta_{i-1}): 1\le i \le n+1\}$, and $T$ satisfies MP under $\varphi_n$. 
Depending on the extension type $\Theta_n$, we consider two cases.

{\bf Case 1.} $\Theta_n=PE$.  In this case, $\partial_{\varphi_n, \gamma_n}(T_n)=\{f_n\}$ by Lemma \ref{hku}(v).
Since $\sigma_n$ is $(T_n, D_n, \phiv_n)$-stable and $\gamma_n\in S_n\subseteq \phibar_n(T_n)\cup D_n$, we have
$\partial_{\sigma_n, \gamma_n}(T_n)=\{f_n\}$. So $\delta\ne \gamma_n$.  

By Theorem~\ref{thm:tech10}(iv), $P_{v_{n}}(\gamma_{n},\delta_{n},\sigma_n)\cap T_{n}=\{v_{n}\}$. Define $\sigma_{n-1}=\sigma_n/P_{v_{n}} (\gamma_{n},\delta_{n}, \sigma_n)$.  Then

(1) $\sigma_{n-1}$ is $(T_{n},D_{n-1},\varphi_{n-1})$-stable by Lemma \ref{stablezang} and hence it
is also $(T_{n-1},D_{n-1},\varphi_{n-1})$-stable. Furthermore, $\partial_{\sigma_n, \delta}(T_n) \subseteq 
\partial_{\sigma_{n-1}, \delta}(T_n)$ (because $\delta\ne \gamma_n$ and possibly $\delta=\delta_n$).

If $i=n$, then $v\preceq v_n$ by (1), as $v_n$ is the maximum defective vertex over all $(T_n, D_{n-1}, \phiv_{n-1})$-stable colorings. So we 
assume that $i<n$. Then $v\in T_{n-1}$ because $m(v)\leq i< n$. Since $v$ is not the smallest vertex in $I[\partial_{\sigma_n, 
\delta}(T_n)]$ and $v\in T_{n-1}$, from (1) it can be seen that  $\delta$ is a defective 
color of $T_{n-1}$ with respect to $\sigma_{n-1}$, and $v$ is not the smallest vertex in $I[\partial_{\sigma_{n-1}, 
	\delta}(T_{n-1})]$. Applying (4.3) and Theorem~\ref{thm:tech10}(v) to $T_{n-1}$ and $\sigma_{n-1}$ 
(see (1)), we obtain $v\preceq v_i$. 

{\bf Case 2.} $\Theta_n=RE$ or $SE$.  In this case, $\varphi_n$ is $(T_{n},D_{n-1},\varphi_{n-1})$-stable (see Algorithm 3.1).
Since $\sigma_n$ is $(T_n, D_n, \phiv_n)$-stable and $\phibar_{n-1}(T_n)\cup D_{n-1}\subseteq \phibar_{n}(T_n)\cup D_n$
by Lemma \ref{hku}(i), $\sigma_n$ is $(T_n,D_{n-1},\varphi_{n-1})$-stable and hence is also $(T_{n-1},D_{n-1},\varphi_{n-1})$-stable.  
If $i=n$, then $v\preceq v_n$, because $v_n$ the maximum defective vertex over all $(T_n, D_{n-1}, \phiv_{n-1})$-stable colorings. 
So we assume that $i<n$. Then $m(v)\le i <n$. Since $v\in T_{n-1}$ and $v$ is not the smallest vertex in $I[\partial_{\sigma_n, \delta}(T_n)]$, 
$\delta$ is a defective color of $T_{n-1}$ with respect to $\sigma_n$, and $v$ is not the smallest vertex in $I[\partial_{\sigma_n, 
\delta}(T_{n-1})]$. Since $\sigma_n$ is $(T_{n-1},D_{n-1},\varphi_{n-1})$-stable, from (4.3) and Theorem~\ref{thm:tech10}(v)
we conclude that $v\preceq v_i$.  \qed

\vskip 3mm

By Definition \ref{wz2}, every $\phiv_n\bmod T_n$ coloring is a $(T_n,D_n,\varphi_n)$-stable coloring. The lemma below says that 
the converse also holds when MP is satisfied, so these two concepts are equivalent in this case. 

\begin{lemma}\label{td}
(Assuming (4.1) and (4.4)) Theorem~\ref{thm:tech10}(vi) holds for all ETTs with $n$ rungs and  
satisfying MP; that is, every $(T_n,D_n,\varphi_n)$-stable coloring $\sigma_n$ is a $\phiv_n\bmod T_n$ coloring.
(Thus every ETT  corresponding to $(\sigma_n, T_n)$ satisfies MP under $\sigma_n$ by Lemma~\ref{MP}.)
\end{lemma}

{\bf Proof.}  By the hypothesis of Theorem~\ref{thm:tech10}, $T$ is an ETT with the corresponding Tashkinov series 
$\TT=\{(T_i, \phiv_{i-1}, S_{i-1}, F_{i-1}, \Theta_{i-1}): 1\le i \le n+1\}$, and $T$ satisfies MP under $\varphi_n$. 
Clearly, every tree-sequence $T^*$ obtained from $T_n$ (resp. $T_n+f_n$) by TAA under $\sigma_n$ if $\Theta_n=PE$ 
(resp. if $\Theta_n=$ $RE$ or $SE$) is a sub-sequence of some closure of $T_n$ (resp. $T_n+f_n$) under $\sigma_n$. 
So to prove that $\sigma_n$ is $\phiv_n\bmod T_n$, it suffices to show that, for an arbitrary closure $T^*_{n+1}$ of $T_n$ 
(resp. $T_n+f_n$) under $\sigma_n$, there exists a Tashkinov series $\TT^*=\{(T_i^*, \sigma_{i-1}, 
S_{i-1}, F_{i-1}, \Theta_{i-1}): 1\le i \le n+1\}$, satisfying the following conditions for all $i$ with $1\le i \le n$:

\vskip 1mm

(1) $T_i^*=T_i$ and

(2) $\sigma_i$ is a $(T_i, D_i, \varphi_i)$-stable coloring in ${\cal C}^k(G-e)$.

For this purpose,  we shall define a coloring $\sigma_{n-1}$ based on $\sigma_n$, such that

(3) $\sigma_{n-1}$ is $(T_{n},D_{n-1},\varphi_{n-1})$-stable and hence is also $(T_{n-1},D_{n-1},\varphi_{n-1})$-stable.

By Lemma~\ref{hku}(iv), we have $\varphi_{n-1} \langle T_n \rangle \subseteq \phibar_{n-1}(T_n)\cup D_{n-1}$, which together with (3) implies that  
$\sigma_{n-1}(f)=\varphi_{n-1}(f)$ for every edge $f$ on $T_n$ and $\overline{\sigma}_{n-1}(v)=\phibar_{n-1}(v)$ for every 
vertex $v$ in $T_n$. Hence $T_n$ can be obtained by TAA from $T_{n-1}$ (resp. $T_{n-1}+f_{n-1}$) under $\sigma_{n-1}$ if $\Theta_{n-1}=PE$ 
(resp. if $\Theta_{n-1}=RE$ or $SE$) in the same way as it under $\varphi_{n-1}$. Moreover, since $T_n$ is closed under $\varphi_{n-1}$, 
it is also a closure of $T_{n-1}$ (resp. $T_{n-1}+f_{n-1}$) under $\sigma_{n-1}$ if $\Theta_{n-1}=PE$ (resp. if $\Theta_{n-1}=RE$ or $SE$). 
By (3), (4.1) and Theorem~\ref{thm:tech10}(vi), $\sigma_{n-1}$ is a $T_{n-1}\bmod\varphi_{n-1}$ coloring. Therefore  

(4) there exists a Tashkinov series $\TT'=\{(T_i^*, \sigma_{i-1}, S_{i-1}, F_{i-1}, \Theta_{i-1}): 1\le i \le n\}$, which together with $\sigma_n$ satisfies 
(1) and (2) for $1\le i \le n$.

We shall then show, using Algorithm 3.1, that the desired Tashkinov series $\TT^*$ can be built from $\TT'$ by adding the tuple $(T^*_{n+1}, \sigma_{n}, S_{n}, F_{n}, \Theta_{n})$. 

Let us now give detailed descriptions. Depending on the extension type, we distinguish between two cases.

{\bf Case 1.} $\Theta_n=RE$. In this case, define $\sigma_{n-1}=\sigma_n$. Since $\sigma_n$ is a $(T_n,D_n,\varphi_n)$-stable coloring, so is $\sigma_{n-1}$. 
Recall that $\varphi_n=\varphi_{n-1}$ by RE in Algorithm 3.1 and that $\phibar_{n-1}(T_n)\cup D_{n-1}\subseteq \phibar_{n}(T_n)\cup D_n$ by Lemma \ref{hku}(i).
So $\sigma_{n-1}$ satisfies (3) and hence (4) holds.

According to Algorithm 3.1, there is a subscript $h \le n-1$ with $\Theta_h=PE$ and $S_h=\{\delta_h, \gamma_h\}$, such that $\Theta_i=RE$ for all $i$ 
with $h+1\le i \le n-1$, if any, and that some $(\gamma_h,\delta_h)$-cycle $O$ with respect to $\varphi_{n-1}$ contains a sub-path $L$ 
with $V(L)\subseteq V(T_n)$ connecting the edge $f_n$ and $V(T_h)$. Note that $\varphi_h=\varphi_{h+1}= \ldots =\varphi_{n-1}$.  
Since $v_h$ is an end of the exit-path $P_{v_h}(\gamma_h,\delta_h,\varphi_h)=P_{v_h}(\gamma_h,\delta_h,
\varphi_{n-1})$, it is outside $O$. Take $w$ in $V(L)\cap V(T_h)$. Then $w\neq v_h$.
As $\sigma_{n-1}$ is $(T_{n-1},D_{n-1},\varphi_{n-1})$-stable by (3) and $\{\delta_h,\gamma_h\} \subseteq \cup_{i \le n-1} S_i \subseteq
\phibar_{n-1}(T_{n-1})\cup D_{n-1}$, every edge of $L$ is colored the same under $\sigma_{n-1}$ as under $\varphi_{n-1}$. 

Let $O^*$ be the $(\gamma_h,\delta_h)$-chain containing $L$ under $\sigma_{n-1}$. Then $O^*$ intersects $T_h$. By (4), $\TT'$ is a 
Tashkinov series, $\Theta_h=PE$, and $\Theta_i=RE$ for all $i$ with $h+1\le i \le n-1$. From Algorithm 3.1 it follows that 
$\sigma_h=\sigma_{n-1}$, $\delta_h\in\overline{\sigma}_h(v_h)$, and $P_{v_h}(\gamma_h,\delta_h,\sigma_h)\cap V(T_h)=\{v_h\}$.
Hence $\delta_h\in\overline{\sigma}_{n-1}(v_h)$, $P_{v_h}(\gamma_h,\delta_h,\sigma_{n-1})=P_{v_h}(\gamma_h,\delta_h,\sigma_h)$,
and $O^*$ is disjoint from the $(\gamma_h,\delta_h)$-path $P_{v_h}(\gamma_h,\delta_h,\sigma_{n-1})$ (because $w\in V(L)\subseteq V(O^*)$). Applying (4.1) and 
Theorem~\ref{thm:tech10}(ii) to $T_n$ under $\sigma_{n-1}$, we see that there is at most one $(\gamma_h,\delta_h)$-path 
intersecting $T_n$. So $O^*$ must be a $(\gamma_h,\delta_h)$-cycle containing $L$ as a sub-path under $\sigma_{n-1}$. 
Therefore $f_n$ can be chosen as an RE connecting edge for $T_n$ under $\sigma_{n-1}$, and $\TT^*$ can thus be built from $\TT'$ 
by adding the tuple $(T^*_{n+1}, \sigma_{n}, S_{n}, F_{n}, \Theta_{n})$ using RE of Algorithm 3.1.

{\bf Case 2.} $\Theta_n=SE$ or $PE$. In this case, define $\sigma_{n-1}=\sigma_n$ if $\Theta_n=SE$ and $\sigma_{n-1}=\sigma_n/P_{v_{n}} (\gamma_{n},
\delta_{n}, \sigma_n)$ if $\Theta_n=PE$. By (4.4) and Theorem~\ref{thm:tech10}(iv), we have 

(5) $P_{v_n}(\gamma_n,\delta_n,\sigma_{n})\cap V(T_n)=\{v_n\}$ when $\Theta_n=PE$, because $\sigma_n$ is $(T_n,D_n,\varphi_n)$-stable. 

According to Algorithm 3.1, $\pi_{n-1}$ is a $(T_n, D_{n-1}, \phiv_{n-1})$-stable coloring whose largest defective vertex $v_n$
is maximum over all $(T_n, D_{n-1}, \phiv_{n-1})$-stable colorings, and $f_n$ is colored by $\delta_n$ under $\pi_{n-1}$.
Observe that

(6) $\sigma_{n-1}$ is $(T_n,D_{n-1}\cup\{\delta_n\},\pi_{n-1})$-stable and hence is $(T_n,D_{n-1},\varphi_{n-1})$-stable. 

Indeed, if $\Theta_n=SE$, then $\phiv_n = \pi_{n-1}$ by SE in Algorithm 3.1. Since $\sigma_{n-1}=\sigma_n$ is 
$(T_n,D_n,\varphi_n)$-stable and $\delta_n\in S_n\subseteq \phibar_n(T_n)\cup D_n$, the desired statement (6) holds. 
If $\Theta_n=PE$, then (6) follows instantly from (5) and Lemma~\ref{stablezang}. 

From (6) we see that both (3) and (4) hold true. Furthermore,  

(7) $\partial_{\sigma_{n-1}, \delta_n}(T_n)= \partial_{\pi_{n-1}, \delta_n}(T_n)$.

By (7), we obtain $\sigma_{n-1}(f_n)=\pi_{n-1}(f_n)=\delta_n$. By (6) and Lemma~\ref{sc2}, every $(T_n, D_{n-1}, \sigma_{n-1})$-stable coloring is
also $(T_n,D_{n-1},\varphi_{n-1})$-stable. So $\sigma_{n-1}$ is a $(T_n, D_{n-1}, \sigma_{n-1})$-stable coloring whose largest 
defective vertex $v_n$ is maximum over all $(T_n, D_{n-1}, \sigma_{n-1})$-stable colorings. 

Again by (6) and Lemma~\ref{sc2}, a coloring is $(T_n,D_{n-1}\cup\{\delta_n\},\pi_{n-1})$-stable 
iff it is $(T_n,D_{n-1}\cup\{\delta_n\},\sigma_{n-1})$-stable. So the equality $\overline{\sigma}(u_n) \cap \overline{\sigma}(T_n)=\emptyset$ holds for every 
$(T_n, D_{n-1}\cup\{\delta_n\}, \sigma_{n-1})$-stable coloring $\sigma$ iff the equality $\overline{\pi}(u_n) \cap \overline{\pi}(T_n)=\emptyset$ holds for every 
$(T_n, D_{n-1}\cup\{\delta_n\}, \pi_{n-1})$-stable colorings $\pi$, where $u_n$ is the end of $f_n$ outside $T_n$ (see Algorithm 3.1). Moreover,
if $\Theta_n=PE$, then $v_n$ is also a $(T_n, \sigma_{n-1}, \{\gamma_n, \delta_n\})$-exit by  the definition of $\sigma_{n-1}$ and (5). 
From Algorithm 3.1 we thus deduce that if RE does not apply to the coloring $\sigma_{n-1}$, then we can construct $\TT^*$ from $\TT'$ (see (4))
by adding the tuple $(T^*_{n+1}, \sigma_{n}, S_{n}, F_{n}, \Theta_{n})$ under $\sigma_n$, using the same extension type, SE or PE, as specified in $\Theta_n$.

It remains to verify that indeed RE does not apply to the coloring $\sigma_{n-1}$. (Recall that RE has priority over both SE and PE in the construction of a
Tashkinov series using Algorithm 3.1 (see (3.3)). That is why we need to check this.)

Assume the contrary: under $\sigma_{n-1}$, there exist an edge $f\in\partial_{\sigma_{n-1},\gamma_h}(T_n)$ and a $(\gamma_h,\delta_h)$-cycle $O$ 
containing a sub-path $L$ with $V(L)\subseteq V(T_n)$ connecting the edge $f$ and $V(T_h)$, where $\Theta_h=PE$, $S_h=\{\delta_h,\gamma_h\}$,
and $\Theta_i=RE$ for all $i$ with $h+1\le i \le n-1$. Then $\varphi_{h}=\varphi_{h+1}=\ldots =\varphi_{n-1}$ by (4). Since $S_h\subseteq\phibar_{n-1}
(T_{n-1})\cup D_{n-1}$, from (5) we see that $\sigma_{n-1}(f)=\varphi_{n-1}(f)$ and that every edge of $L$ is colored the same under $\sigma_{n-1}$ 
as under $\varphi_{n-1}=\varphi_{h}$.  

Let $O^*$ be the $(\gamma_h,\delta_h)$-chain containing $L$ under $\varphi_{n-1}$. By an argument parallel to that used for Case 1 (see the paragraph right above
Case 2), we can ensure that $O^*$ is a $(\gamma_h,\delta_h)$-cycle under $\varphi_{n-1}$ containing the sub-path $L$ connecting the edge $f$ and 
$V(T_h)$. Therefore $\Theta_n=RE$ with respect to $\varphi_{n-1}$ (see Algorithm 3.1), contradicting the hypothesis of the present case. \qed

\section{Good Hierarchies}

It is well known that Kempe changes play a fundamental role in edge-coloring theory. To ensure that an ETT under a coloring remains to
be an ETT under a new coloring arising from Kempe changes, in this section we develop an effective control mechanism over such operations, 
the so-called good hierarchy of an ETT, which will serve as a powerful tool in the proof of Theorem~\ref{thm:tech10}(i). As stated before, a prototype
of this mechanism can be found in Chen and Jing in \cite{GG} (see Condition R2 therein). Throughout this 
section, we assume that

{\bf (5.1)} Theorem~\ref{thm:tech10}(i) and (ii) hold for all ETTs with at most $n-1$ rungs and satisfying MP, and 
Theorem~\ref{thm:tech10}(iii)-(vi) hold for all ETTs with at most $n$ rungs and satisfying MP.

In the case of $\Theta_n=PE$, let $J_n$ be a closure of $T_{n}(v_n)$ under a $(T_n,D_n,\varphi_n)$-stable coloring $\sigma_n$. By
Algorithm 3.1, $\delta_n\in \phibar_n(v_n)$ and $|\partial_{\varphi_n, \delta_n}(T_n)|\ge 2$ (see (3.2)). By Lemma \ref{hku}(v), edges in 
$\partial_{\varphi_n, \delta_n}(T_n)$ are all incident to $V(T_n(v_n)-v_n)$. Since $\bar{\sigma}_n(v_n)=\phibar_n(v_n)$
and $\partial_{\sigma_n, \delta_n}(T_n)=\partial_{\varphi_n, \delta_n}(T_n)$, there holds $V(J_n)-V(T_n)\neq\emptyset$. 
We use $T_n\vee J_n$ to denote the tree-sequence obtained from $T_n$ by adding all vertices in $V(J_n)-V(T_n)$ to $T_n$ one by one, 
following the linear order $\prec$ in $J_n$, and using edges in $J_n$.

\begin{lemma}\label{elementary}
(Assuming (5.1)) Let $T$ be an ETT constructed from a $k$-triple $(G,e, \varphi)$ by using the Tashkinov series 
$\TT=\{(T_i, \phiv_{i-1}, S_{i-1}, F_{i-1}, \Theta_{i-1}): 1\le i \le n+1\}$. Suppose $\Theta_n=PE$ and $T$ satisfies 
MP under $\varphi_n$. If $J_n$ is a closure of $T_{n}(v_n)$ under a $(T_n,D_n,\varphi_n)$-stable coloring
$\sigma_n$, then $V(T_n\vee J_n)$ is elementary with respect to $\sigma_n$.
\end{lemma}

{\bf Proof.} 	Clearly, $T_n$ is an ETT with corresponding Tashkinov series $\TT=\{(T_i, \phiv_{i-1}, S_{i-1}, F_{i-1}, \\
\Theta_{i-1}): 1\le i \le n\}$ and satisfies MP under $\varphi_{n-1}$. Since $r(T_n)=n-1$, by (5.1) and Theorem~\ref{thm:tech10}(i), 
$V(T_n)$ is elementary with respect to $\varphi_{n-1}$. Let $\pi_{n-1}$ and $\pi_{n-1}'$ be as specified in Algorithm 3.1. Since 
$\pi_{n-1}$ is a $(T_n,D_{n-1}, \varphi_{n-1})$-stable coloring and $\pi_{n-1}'$ is a $(T_n,D_{n-1} 
\cup \{\delta_n\}, \pi_{n-1})$-stable coloring, by definition $V(T_n)$ is also elementary with respect to $\pi'_{n-1}$. As 
$\varphi_n=\pi'_{n-1}/P_{v_n}(\delta_n,\gamma_n, \pi'_{n-1})$ and $\delta_n\notin\overline{\pi}'_{n-1}(T_n)$, we further obtain

(1) $V(T_n)$ is elementary with respect to $\varphi_{n}$ and hence elementary with respect 
to $\sigma_n$. 

Since $\sigma_n$ is a $(T_n,D_n,\varphi_n)$-stable coloring, it follows from (5.1) and Theorem~\ref{thm:tech10}(iii) that 
$\sigma_n$ is $(T_{j}(v_n)-v_n, D_{j-1},\varphi_{j-1})$-stable and hence is $(T_{j-1}, D_{j-1}, \varphi_{j-1})$-stable, where  
$j=m(v_n)$. By Theorem \ref{thm:tech10}(vi), $\sigma_n$ is a $\phiv_{j-1}\bmod T_{j-1}$ coloring, so every ETT corresponding 
to $(\sigma_n, T_{j-1})$ satisfies MP.  Using Lemma \ref{hku}(iv) and Lemma \ref{samecolor}, we obtain $\sigma_n(f)=\varphi_n(f)=
\varphi_{j-1}(f)$ for each edge $f$ on $T_j$. Hence $J_n$ is a closure of $T_n(v_n)=T_j(v_n)$ under $\sigma_n$. Consequently,
$J_n$ is an ETT corresponding to $(\sigma_n, T_{j-1})$ and satisfies MP. Since $r(J_n)=j-1\le n-1$,

(2) $V(J_n)$ is elementary and closed with respect to $\sigma_{n}$ by (5.1) and Theorem~\ref{thm:tech10}(i).

Suppose on the contrary that $V(T_n\vee J_n)$ is not elementary with respect to $\sigma_n$. Then $T_n\vee J_n$ contains 
two distinct vertices $u$ and $v$ such that $\overline{\sigma}_n(u) \cap \overline{\sigma}_n(v) \ne \emptyset$.  By (1) 
and (2), we may assume that $u \in V(T_n)-V(J_n)$ and $v \in V(J_n)-V(T_n)$. So $v\neq v_n$. Let $\alpha\in \overline{\sigma}_n(u) \cap \overline{\sigma}_n(v)$. Then $\alpha\neq \delta_n$ by (2), because $\delta_n \in \overline{\varphi}_n(v_n)=
\overline{\sigma}_n(v_n)$. Moreover, since $\gamma_n\in\phibar_{n-1}(v_n)$ and $V(T_n)$ is elementary with respect to 
$\varphi_{n-1}$, from $PE$ of Algorithm 3.1 and the definition of stable colorings, we deduce that $\gamma_n\notin
\phibar_n(T_n)$ and hence $\gamma_n\notin \overline{\sigma}_n(T_n)$. So $\alpha\neq\gamma_n$. Consequently, 

(3) $\alpha\notin S_n$. 

Since $v_n$ is a maximum defective vertex according to Algorithm 3.1, $T_{n}(v_n)$ contains a vertex $w \ne v_n$. Note that $w$ is contained in
both $T_n$ and $J_n$. Let $\beta \in \overline{\sigma}_n(w)$. Since  $\delta_n \in \overline{\sigma}_n(v_n)$
and $\gamma_n\notin \overline{\sigma}_n(T_n)$, by (2)  we obtain

(4) $\beta \notin S_n$ and the other end of $P_v(\alpha, \beta, \sigma_n)$ is $w$.

From (3), (4), and Algorithm 3.1, we see that $\partial(T_n)$ contains no edge colored by $\alpha$ or 
$\beta$ under $\varphi_{n}$ and hence under $\sigma_n$, because $\sigma_n$ is $(T_n,D_n,\varphi_n)$-stable. 
Combining this with (1), we conclude that the other end of $P_u(\alpha, \beta, \sigma_n)$ is 
also $w$. Thus $P_w(\alpha, \beta, \sigma_n)$ terminates at both $u$ and $v$, a contradiction. \qed 

\vskip 3mm

Let $T$ be an ETT as specified in Theorem \ref{thm:tech10}; that is, $T$ is constructed from a $k$-triple $(G,e, \varphi)$ 
by using the Tashkinov series $\TT=\{(T_i, \phiv_{i-1}, S_{i-1}, F_{i-1}, \Theta_{i-1}): 1\le i \le n+1\}$. To prove that
$V(T)$ is elementary with respect to $\phiv_n$, we shall turn to considering a restricted ETT $T'$ with 
ladder $T_1\subset T_2 \subset \ldots \subset T_n \subset T'$ and $V(T')=V(T_{n+1})$, and then show that 
$V(T')$ is elementary with respect to $\phiv_n$. For convenience, we may simply view $T'$ as $T_{n+1}$.   

\vskip 2mm
In the remainder of this paper, we reserve the symbol $R_n$ for a fixed closure of $T_n(v_n)$ under $\phiv_n$, if
$\Theta_n=PE$. Let $T_n\vee R_n$ be the tree-sequence as defined above Lemma \ref{elementary}. We assume hereafter that

{\bf (5.2)} $T_{n+1}$ is a closure of $T_n\vee R_n$ under $\phiv_n$, which is a special closure of $T_n$ under $\phiv_n$
(see PE in Algorithm 3.1), when $\Theta_n=PE$.

By Lemma~\ref{elementary}, $V(T_n\vee R_n)$ is elementary with respect to $\varphi_n$, so we may further assume that

{\bf (5.3)} $T \ne T_n\vee R_n$ if $\Theta_n=PE$, which together with (5.2) implies that $T_n\vee R_n$ is not
closed with respect to $\varphi_n$.  

{\bf (5.4)} If $\Theta_n = PE$, then each color in $\overline{\varphi}_n(T_n) \cap \overline{\varphi}_n(R_n)$ is 
closed in $T_n \vee R_n$ with respect to $\varphi_n$. 

To justify this, note that each color in $\overline{\varphi}_n(R_n)$ is closed in $R_n$ under $\varphi_n$ because
$R_n$ is a closure. By Lemma \ref{hku}(v), each color in  $\overline{\varphi}_n(T_n)-\{\delta_n\}$ is closed in 
$T_n$ under $\varphi_n$. Hence each color in $\overline{\varphi}_n(T_n) \cap \overline{\varphi}_n(R_n)-\{\delta_n\}$ 
is closed in $T_n \vee R_n$ with respect to $\varphi_n$. Lemma \ref{hku}(v) also asserts that edges in $\partial_{\varphi_n, \delta_n}(T_n)$ are all incident to $V(T_n(v_n)-v_n)$. So $\delta_n$ is closed in $T_n \vee R_n$ as well, because it is closed in $R_n$. Hence (5.4) follows.  

\vskip 2mm
To prove Theorem \ref{thm:tech10}(i), we shall appeal to a {\em hierarchy}\label{hierar} of $T$ of the form

{\bf (5.5)} $T_n=T_{n,0} \subset T_{n,1} \subset \ldots \subset T_{n,q} \subset T_{n,q+1}=T$, such that $T_n\vee R_n 
\subset T_{n,1}$ if $\Theta_n=PE$ and that $T_{n,i}=T(u_i)$ for $1\le i\leq q$, where $u_1\prec u_2\prec \ldots 
\prec u_{q}$ are some vertices in $T-V(T_n)$, called {\it dividers} of $T$. (So $T$ has $q$ dividers in total.)

As introduced before, $D_n=\cup_{h\leq n}S_h-\phibar_n(T_n)$, where $S_h=\{\delta_h\}$ 
if $\Theta_h=SE$ and $S_h=\{\delta_h,\gamma_h\}$ otherwise. By Lemma \ref{Dnzang}, we have

{\bf (5.6)} $|D_n|\leq n$. 

\noindent Write $D_n=\{\eta_1, \eta_2, \ldots ,\eta_{n'}\}$. In Definition \ref{R2} given below and the remainder of
this paper,

\vskip 1mm
$\bullet$ $T_{n,0}^*= T_{n}\vee R_n$ if $\Theta_n=PE$ and $T_{n,0}^*= T_{n}$ otherwise, and $T_{n,j}^*=T_{n,j}$ if $j \ge 1$;

$\bullet$ $D_{n,j}=\cup_{h\leq n}S_h-\phibar_n(T_{n,j}^*)$ for $0 \le j \le q$ (so $D_{n,j}\subseteq D_n$);  

$\bullet$ $v_{\eta_h}$, for $\eta_h \in D_n$, is defined to be the first vertex $u$ of $T$ in the order $\prec$ with  $\eta_h \in 
\phibar_n (u)$, if  

\hskip 3mm any, and defined to be the last vertex of $T$ in the order $\prec$ otherwise;

$\bullet$ $T_{n,j}(v_{\eta_h})=T_{n,j}$ if $v_{\eta_h}$ is outside $T_{n,j}$ for $1 \le j \le q$ and $\eta_h \in D_n$; and

$\bullet$ $\Gamma^{j}=\cup_{\eta_h\in D_{n,j}}\Gamma^{j}_h$ for $0 \le j \le q$.

\vskip 1mm

Let $H$ be a subgraph of $G$ and let $C$ be a subset of $[k]$. We say that $H$ is $C$-{\em closed}\label{cclosed} with respect to 
$\varphi_n$ if $\partial_{\varphi_n, \alpha}(H)=\emptyset$ for any $\alpha \in C$, and say that $H$ is $C^-$-{\em closed}
with respect to $\varphi_n$ if it is $(\phibar_n(H)-C)$-closed with respect to $\varphi_n$.    

\begin{definition} \label{R2}
{\rm Hierarchy (5.5) of $T$ is called {\em good}\label{goodhierarchy} with respect to $\varphi_n$ if for any $j$ with $0 \le j \le q$ and 
any $\eta_h\in D_{n,j}$, there exists a $2$-color subset $\Gamma^{j}_h=\{\gamma^{j}_{h_1},\gamma^{j}_{h_2}\}\subseteq 
[k]$, such that
\begin{itemize}
\vspace{-2mm}
\item[(i)] $\Gamma^0_h \subseteq \phibar_n(T_n)-\varphi_n \langle T_{n,1}(v_{\eta_h})-T_{n,0}^* \rangle$ 
and $\Gamma^j_h \subseteq \phibar_n(T_{n,j})-\varphi_n \langle T_{n,j+1}(v_{\eta_h})-T_{n,j} \rangle$ for $1 \le j \le q$
(so neither color in $\Gamma^{j}_h$ can be used by edges on $T_{n,j+1}-T_{n,j}^*$ until after $\eta_h$ becomes 
missing at the vertex $v_{\eta_h}$ in $T_{n,j+1}$ for $0 \le j \le q$); 
\vspace{-2mm}
\item[(ii)] $\Gamma^{j}_{g} \cap \Gamma^{j}_{h}=\emptyset$ whenever $\eta_{g}$ and $\eta_{h}$ are two distinct colors in $D_{n,j}$;
\vspace{-2mm}
\item[(iii)] for any $j$ with $1\le j \le q$, there exists precisely one color $\eta_g\in D_{n,j}$, such that  
$\Gamma^{j}_{g} \subseteq \phibar_n(T_{n,j}-V(T_{n,j-1}^*))$ (so $\Gamma^{j}_{g} \cap \Gamma^{j-1}_{g} =\emptyset$)
and $\Gamma^{j}_{h}=\Gamma^{j-1}_{h}$ for all $\eta_h\in D_{n,j}-\{\eta_g\}$; 
\vspace{-2mm}
\item[(iv)] if $\Theta_n=PE$, then $T_{n}\vee R_n$ is not $(\Gamma^0)^-$-closed with respect to $\varphi_n$ and, subject 
to this, $|\phibar_n(T_n) \cap \phibar_n(R_n)-\Gamma^0|\ge 4$; and
\vspace{-2mm}
\item[(v)] $T_{n,j}$ is $(\cup_{\eta_h\in D_{n,j}}\Gamma^{j-1}_h)^-$-closed with respect to $\varphi_n$ 
for all $j$ with $1\le j \le q$.
\vspace{-2mm}
\end{itemize}
\noindent The sets $\Gamma^{j}_h$ are referred to as $\Gamma$-{\em sets}\label{gammasets} of the hierarchy (or of $T$) under $\varphi_n$.} 
\end{definition}

At first glance, the concept of good hierarchies is very complicated. After reading the constructive proof of Theorem \ref{good} shortly,
one may realize that it is, nevertheless, fairly easy to understand. The following remarks may foster a better grasp of 
this concept. 

{\bf (5.7)} For $0 \le j\le q$ and $\eta_h\in D_{n,j}$, we have $\Gamma^{j}_h \subseteq \phibar_n(T_{n,j}^*)$ by Condition (i).
So  $\Gamma^{j}_h \cap D_{n,j}=\emptyset$ and hence $\Gamma^{j} \cap D_{n,j}=\emptyset$.

{\bf (5.8)} Condition (iv) implies that $T_{n,1}\neq T_{n}\vee R_n$ if $\Theta_n=PE$.

{\bf (5.9)} When $\Theta_n=RE$ or $SE$, the first edge added to $T_{n,1}-T_{n,0}$ is $f_n$ (see (5.5) and Algorithm 3.1). 
For $1\le j\le q$, by definitions, $D_{n,j} \subseteq D_{n,j-1}$, so $\Gamma^{j-1}_h$ is well defined for any
$\eta_h\in D_{n,j}$ and $\cup_{\eta_h\in D_{n,j}}\Gamma^{j-1}_h \subseteq \Gamma^{j-1}$. In view of Condition 
(v), the first edge added to $T_{n,j+1}- T_{n,j}$ is colored by a color $\alpha$ in $\Gamma^{j-1}_g$
for some $g$ with $\eta_g\in D_{n,j}$. From Condition (i)  we see that $\alpha\notin \Gamma^{j}_g$. So 
$ \Gamma^{j}_g \ne \Gamma^{j-1}_g$. According to Condition (iii), now $\Gamma^{j}_g$ consists of two colors in 
$\phibar_n(T_{n,j} -V(T_{n,j-1}^*))$. Thus $\Gamma^{j-1}_g \cap \Gamma^{j}_g=\emptyset$ and hence $\alpha \notin 
\Gamma^j$.

{{\bf (5.10)} If a color $\alpha \in \phibar_n(T_{n,j}-V(T_{n,j-1}^*))$ for some $j$ with $1 \le j \le q$, then
$\alpha \notin \Gamma^{j-1}$ by Condition (i), and hence $\alpha$ is closed in $T_{n,j}$ with respect to $\varphi_n$ 
by Condition (v). This simple observation will be used repeatedly in subsequent proofs.

{{\bf (5.11)} Note that not every ETT admits a good hierarchy. Suppose $T$ does have such a hierarchy. To prove that
$V(T)$ is elementary with respect to $\varphi_n$, as usual, we shall perform a sequence of Kempe changes to reduce a minimum counterexample 
to an even smaller one, thereby reaching a contradiction. (The adjectives {\em minimum} and {\em smaller} used here are not meant 
with respect to the number of vertices. The rigorous definition of a minimum counterexample will be given in the next section; see (6.2)-(6.5).) 
Since interchanging with a color in $D_{n,j}$ by a Kempe change often results in 
a coloring which is not stable, in our proof we shall use colors in $\Gamma^{j}_h$ as stepping stones to interchange with the color 
$\eta_h$ in $D_{n,j}$ while maintaining stable colorings in subsequent proofs (such an interchange property indeed holds, as we shall see). 
So we may think of $\Gamma^{j}_h$ as a color set exclusively reserved for $\eta_h$ (see Condition (ii)) and think of a good hierarchy as a 
control mechanism over Kempe changes. We point out that Condition (i) can be used to preserve colors on edges of $T_{n,j}(v_{\eta_h})-T_{n,j-1}^*$ 
under Kempe changes for $\eta_h$ and a color in $\Gamma^{j}_h$. Condition (v) ensures that the aforementioned 
interchange property is satisfied by colors closed in $T_{n,j}$. Moreover, extending $T$ by TAA while keeping condition 
(i) for $j=q$ leads to Condition (v). Unless $T$ is already closed, Condition (iii) allows us to further extend $T$ by TAA
while keeping the good hierarchy property, provided that Condition (v) holds for $T=T_{n,q+1}$.  

\vskip 3mm

We break the proof of Theorem \ref{thm:tech10}(i) into the following two theorems. Although the first theorem appears to
be weaker than Theorem \ref{thm:tech10}(i), the second one implies that they are actually equivalent. 
We only present a proof of the second theorem in this section, and will give a proof of the first one 
in the next two sections.

\begin{theorem}\label{hierarchy}
(Assuming (5.1)) Let $T$ be an ETT constructed from a $k$-triple $(G,e, \varphi)$ by using the Tashkinov series $\TT=\{(T_i, \phiv_{i-1}, S_{i-1}, F_{i-1}, \Theta_{i-1}): 1\le i \le n+1\}$. Suppose $T$ admits a good hierarchy and satisfies
MP with respect to $\varphi_n$. Then $V(T)$ is elementary with respect to $\varphi_n$. 
\end{theorem}

\begin{theorem}\label{good}
(Assuming (5.1)) Let $T$ be an ETT constructed from a $k$-triple $(G,e, \varphi)$ by 
using the Tashkinov series $\TT=\{(T_i, \phiv_{i-1}, S_{i-1}, F_{i-1}, \Theta_{i-1}): 1\le i \le n+1\}$. If $T$ 
satisfies MP under $\varphi_n$, then there exists a closed ETT $T'$ corresponding to $(\varphi_n,T_n)$ with $V(T')=V(T_{n+1})$, 
such that $T'$ admits a good hierarchy and satisfies MP with respect to $\varphi_n$. 
\end{theorem}

\noindent {\bf Remark}. Our proof of Theorem \ref{good} is based on Theorem \ref{hierarchy},
while the proof of Theorem \ref{hierarchy} is completely independent of Theorem \ref{good}.

\vskip 2mm

{\bf Proof of Theorem \ref{good}.} By (5.1) and Theorem~\ref{thm:tech10}(i), $V(T_i)$ is elementary with respect to 
$\varphi_{i-1}$ for $1 \le i\leq n$. So each $|T_i|$ is an odd number. Thus $|T_{i}|-|T_{i-1}|\geq 2$ for each $1\le 
i\leq n$. By Theorem \ref{ThmScheide}, if $|T_1|\le 10$, then $G$ is an elementary multigraph, thereby proving Theorem 
\ref{ThmGS2} in this case. So we may assume that $|T_1|\ge 11$. Hence

(1) $|T_i| \ge 2i+9$ for $1\le i \le n$. 

We shall actually construct an ETT $T'$ from $T_n$ by using the same connecting edge, connecting 
color, and extension type as $T$, which has a good hierarchy:

(2) $T_n=T_{n,0} \subset T_{n,1} \subset \ldots \subset T_{n,q+1}=T'$, such that $T_n\vee R_n \subset T_{n,1}$ 
if $\Theta_n=PE$ and such that $V(T')=V(T_{n+1})$. 

Since $V(T_n)$ is elementary with respect to $\varphi_{n-1}$, by (1) we have $|\phibar_{n-1}(T_{n})|\ge 2n+11$
 (as $e$ is uncolored).  From Algorithm 3.1 we see that $|\phibar_{n-1}(T_{n})|=|\phibar_n(T_{n})|$. So 

(3) $|\phibar_n(T_{n})|\ge 2n+11$. Moreover, $|D_{n,0}|\le |D_n|\le n$ by (5.6). 

(4) If $\Theta_n=PE$, then we can find a $2$-color set $\Gamma^{0}_h=\{\gamma^{0}_{h_1},\gamma^{0}_{h_2}\}\subseteq 
\phibar_n(T_{n})$ for each $\eta_h\in D_{n,0}=\cup_{h\leq n}S_h-\phibar_n(T_n\vee R_n)$, such that $\Gamma^{0}_{g} \cap \Gamma^{0}_{h}=\emptyset$ whenever $\eta_{g}$ and $\eta_{h}$ are two distinct colors in $D_{n,0}$, and such that $T_{n}\vee R_n$ 
is not $(\Gamma^0)^-$-closed with respect to $\varphi_n$, where $\Gamma^{0}=\cup_{\eta_h\in D_{n,0}}\Gamma^{0}_h$.

To justify this, let $\alpha$ be a color in $\overline{\varphi}_n(T_n \vee R_n)$ that is not closed in $T_n \vee R_n$
under $\varphi_n$; such a color exists by (5.3). In view of (3), $\overline{\varphi}_n(T_n)-\{\alpha\}$ 
contains at least $2n+10$ colors. So (4) follows if we pick all colors in $\Gamma^{0}$ from $\overline{\varphi}_n(T_n)-
\{\alpha\}$.  

(5) If $\Theta_n=PE$, then there exists a $2$-color set $\Gamma^{0}_h=\{\gamma^{0}_{h_1},\gamma^{0}_{h_2}\}\subseteq 
\phibar_n(T_{n})$ for each $\eta_h\in D_{n,0}$ as described in (4), such that $|\phibar_n(T_n) \cap \phibar_n(R_n)-\Gamma^0|\ge 4$. 

To justify this, let $\alpha$ be as specified in the proof of (4). Then $\alpha \notin \overline{\varphi}_n(T_n) \cap \overline{\varphi}_n(R_n)$ by (5.4). Since $v_n$ is a maximum defective vertex and $v_n\in T_n\cap R_n$, the ends of the uncolored edge $e$ are contained in both 
$T_n$ and $R_n$. So $|\overline{\varphi}_n(T_n) \cap \overline{\varphi}_n(R_n)|\geq 4$.  If we pick all colors in $\Gamma^{0}$ from $\overline{\varphi}_n(T_n)-\{\alpha\}$, 
with priority given to those in $\overline{\varphi}_n(T_n)-\overline{\varphi}_n(R_n)$, then $|\overline{\varphi}_n(T_n) \cap \overline{\varphi}_n(R_n) -\Gamma^0| \ge 4$ 
by (3), thereby establishing (5).

Thus Definition \ref{R2}(iv) is satisfied by these sets $\Gamma^{0}_h$. Using (3), we can similarly get 
the following statement. 

(6) If $\Theta_n \ne PE$, then we can find a $2$-color set $\Gamma^{0}_h=\{\gamma^{0}_{h_1},\gamma^{0}_{h_2}\}\subseteq 
\phibar_n(T_{n})$ for each $\eta_h\in D_{n,0}=D_n$, such that $\Gamma^{0}_{g} \cap \Gamma^{0}_{h}=\emptyset$ whenever 
$\eta_{g}$ and $\eta_{h}$ are two distinct colors in $D_{n,0}$.

So Definition \ref{R2}(ii) is also satisfied by these sets $\Gamma^{0}_h$.
Let us construct $T'$ by the following Algorithms 5.5 and 5.6. Recall that   
$v_{\eta_h}$ is defined to be the first vertex of $T'$ in the order $\prec'$ for which $\eta_h \in 
\phibar_n (v_{\eta_h})$, if any, and defined to be the last vertex of $T'$ in the order $\prec'$ otherwise;
and $T_{n,j+1}(v_{\eta_h})=T_{n,j+1}$ if $v_{\eta_h}$ is not contained in $T_{n,j+1}$ for $0\le j \le q$.

Given  $\{\Gamma^0_h:\, \eta_h\in D_{n,0}\}$, let us construct $T_{n,1}$ using the following procedure.

\vskip 3mm
\noindent {\bf Algorithm 5.5}
\vskip 2mm
\noindent {\bf Step 0.} Set $T_{n,1}= T_{n}\vee R_n$ if $\Theta_n=PE$ and $T_{n,1}=T_n+f_n$ otherwise, 
where $f_n$ is the connecting edge used in Algorithm 3.1, depending on $\Theta_n$. 

\vskip 2mm
\noindent {\bf Step 1.} While there exists $f\in\partial(T_{n,1})$ with $\varphi_n(f)\in\phibar_n(T_{n,1})$, 
do: set $T_{n, 1}=T_{n,1}+f$ if the resulting $T_{n, 1}$ satisfies $\Gamma^{0}_h \cap \varphi_n \langle T_{n,1}(v_{\eta_h})
-T_{n,0}^*  \rangle =\emptyset$ for all $\eta_h\in D_{n,0}$, where $T_{n,0}^*= T_{n}\vee R_n$ if $\Theta_n=PE$ and
$T_{n,0}^*= T_{n}$ otherwise.

\vskip 2mm
\noindent {\bf Step 2.} Return $T_{n, 1}$. 

\vskip 3mm
Note that if $\Theta_n=PE$, then $T_{n}\vee R_n$ is not $(\Gamma^0)^-$-closed with respect to $\varphi_n$ by (4) and (5). 
So Step 1 is applicable to $T_{n}\vee R_n$, and hence $T_{n,1}\ne  T_{n}\vee R_n$. If $\Theta_n=RE$ or $SE$, then 
$T_{n,1}\ne  T_{n}$ by the algorithm.  For each $\eta_h\in D_{n,0}$, it follows from (5), (6), and Step 1 that $\Gamma^{0}_h
\subseteq \phibar_n(T_n)- \varphi_n \langle T_{n,1}(v_{\eta_h})- T_{n,0}^*  \rangle $. So  
Definition \ref{R2}(i) is satisfied. Moreover, $T_{n,1}$ is $(\cup_{\eta_h\in D_{n,1}}\Gamma^{0}_h)^-$-closed with 
respect to $\varphi_n$, as stated in Definition \ref{R2}(v). To justify this, assume the contrary: there exists 
$f\in\partial(T_{n,1})$ with $\varphi_n(f) \in \phibar_n(T_{n,1})- (\cup_{\eta_h\in D_{n,1}} \Gamma^{0}_h)$. Then either 
$\varphi_n(f) \in \phibar_n(T_{n,1})- (\cup_{\eta_h\in D_{n,0}} \Gamma^{0}_h)$ or $\varphi_n(f) \in \Gamma^{0}_h$ for some 
$\eta_h\in D_{n,0}$ but $\eta_h\notin D_{n,1}$; in the latter case, ${\eta_h}$ has become a missing color at the vertex $v_{\eta_h}$ in $T_{n,1}$. 
Thus we can further grow $T_{n,1}$ by using $f$ and Step 1 in either case, a contradiction. Since Definition \ref{R2}(iii) starts with $j=1$, 
$\{\Gamma^{0}_h:\, \eta_h\in D_{n,0}\}$ and $T_{n,1}$ satisfy all the conditions specified in Definition \ref{R2}. 

\vskip 2mm
Suppose we have constructed  $\{\Gamma^{i-1}_h:\, \eta_h\in D_{n,i-1}\}$ and $T_{n,i}$ for all $i$ with 
$1 \le i \le j$, which are as described in Definition \ref{R2}. If $T_{n,j}$ is closed with respect to
$\varphi_n$ (equivalently $V(T_{n,j})=V(T_{n+1}))$, set $T'=T_{n,j}$. Otherwise, we proceed to the construction 
of $\{\Gamma^{j}_h:\, \eta_h\in D_{n,j}\}$ and $T_{n,j+1}$ using the following procedure. 

\vskip 3mm
\noindent {\bf Algorithm 5.6}
\vskip 2mm
\noindent {\bf Step 0.} Set $\Gamma^{j}_h =\Gamma^{j-1}_h$ for each $\eta_h\in D_{n, j}$.

\vskip 2mm
\noindent {\bf Step 1.}  Let $f$ be an edge in $\partial(T_{n,j})$ with $\varphi_n(f) \in \Gamma^{j-1}_h$ for some 
$\eta_h\in D_{n,j}$, let $T_{n,j+1}=T_{n,j}+f$, and let $\{\gamma^{j}_{h_1},\gamma^{j}_{h_2}\}$ be a 2-subset of $\phibar_n(T_{n,j}-V(T_{n,j-1}))$. Replace $\Gamma^{j}_h$ by $\{\gamma^{j}_{h_1},\gamma^{j}_{h_2}\}$.

\vskip 2mm
\noindent {\bf Step 2.} While there exists $f\in\partial(T_{n,j+1})$ with $\varphi_n(f)\in\phibar_n(T_{n,j+1})$, 
do: set $T_{n,j+1}=T_{n,j+1}+f$ if the resulting $T_{n,j+1}$ satisfies $\Gamma^{j}_h \cap \varphi_n \langle T_{n,j+1}(v_{\eta_h})-T_{n,j}  \rangle =\emptyset$ for all $\eta_h\in D_{n,j}$. 

\vskip 2mm
\noindent {\bf Step 3.} Return $\{\Gamma^{j}_h:\, \eta_h\in D_{n,j}\}$ and $T_{n,j+1}$.

\vskip 2mm
Let us make some observations about this algorithm and its output. 

As $T_{n,j}$ is not closed with respect to $\varphi_n$, $V(T_{n,j})$ is a proper subset of $V(T_{n+1})$. By  Definition 
\ref{R2}(v), $T_{n,j}$ is $(\cup_{\eta_h\in D_{n,j}}\Gamma^{j-1}_h)^-$-closed with respect to $\varphi_n$.
So there exists a color $\beta \in \cup_{\eta_h\in D_{n,j}}\Gamma^{j-1}_h$, such that $\partial_{\varphi_n, \beta}(T_{n,j}) \ne \emptyset$. Hence the edge $f$ specified in Step 1 is available. 

For $1\le i \le j$, we have $|\phibar_n(T_{n,i})|\ge |\phibar_n(T_{n})|\ge 2n+11$ and $|D_{n,i}| 
\le |D_{n,0}|\le |D_n|\le n$  by (3). So $\phibar_n(T_{n,i})-(\cup_{\eta_h\in D_{n,i}}\Gamma^{i-1}_h) \ne \emptyset$; 
let $\alpha$ be a color in this set. By Theorem \ref{hierarchy} (see the remark right above the proof of this theorem), 
$V(T_{n,i})$ is elementary with respect to $\varphi_n$, which implies that $|T_{n,i}|$ is odd, because $\alpha$ is closed 
in $T_{n,i}$ under $\varphi_n$ by Definition \ref{R2}(v). It follows that $|T_{n,j}|- |T_{n,j-1}|\ge 2$.  So $\phibar_n(T_{n,j}-V(T_{n,j-1}))$ contains at least two distinct colors, and hence the $2$-subset $\{\gamma^{j}_{h_1},\gamma^{j}_{h_2}\}$ involved in Step 1 exists. Thus Definition \ref{R2}(iii) is satisfied. Since $T_{n,j}$ is elementary by Theorem \ref{hierarchy} and Definition \ref{R2}(ii) is satisfied by $T_{n,j-1}$, from Step 0 and Step 1 we see that Definition \ref{R2}(ii) holds for $T_{n,j}$.

Note that each color in $\phibar_n(T_{n,j+1})-(\cup_{\eta_h\in D_{n,j+1}}\Gamma^{j}_h)$ is closed in $T_{n,j+1}$ with respect to $\varphi_n$, for otherwise, 
$T_{n,j+1}$ can be augmented further using Step 2 (see the paragraph succeeding Algorithm 5.5 for a proof). Thus $T_{n,j+1}$ is $(\cup_{\eta_h\in D_{n,j+1}}\Gamma^{j}_h)^-$-closed with respect to $\varphi_n$, and hence Definition \ref{R2}(v) holds. From the algorithm it follows that $\Gamma^{j}_h \subseteq \phibar(T_{n,j})-\varphi_n \langle T_{n,j+1}(v_{\eta_h})-T_{n,j}^* \rangle$ for all $\eta_h\in D_{n,j}$, so Definition \ref{R2}(i) holds. Thus 
$\{\Gamma^{j}_h:\, \eta_h\in D_{n,j}\}$ and $T_{n,j+1}$ satisfy all the conditions in Definition \ref{R2} and hence are as desired.

Repeating the process, we can eventually get a closed ETT $T'$, with $V(T')=V(T_{n+1})$, that admits a good hierarchy
with respect to $\varphi_n$. Clearly, $T'$ also satisfies MP under $\varphi_n$.  \qed \\

Consider the case when $\Theta_n=PE$. By the definition of hierarchy (see (5.5)), $T_n\vee R_n$ is fully contained
in $T_{n,1}$. To maintain the structure of $T_n\vee R_n$ under Kempe changes, we need the following concept in subsequent proofs.  
A coloring $\sigma \in {\cal C}^k(G-e)$ is called a $(T_n \oplus R_n, D_n,\varphi_n)$-{\em stable coloring}\label{trstable} if it is both $(T_n,D_n,\varphi_n)$-stable and $(R_n,\emptyset,\varphi_n)$-stable; that is, the following conditions are satisfied:

$\bullet$ $\sigma(f)=\varphi_n(f)$ for any edge $f$ incident to $T_n$ with $\varphi_n(f)\in \phibar_n(T_n)\cup D_n$;

$\bullet$ $\sigma(f)=\varphi_n(f)$ for any edge $f$ incident to $R_n$ with $\varphi_n(f)\in\phibar_n(R_n)$; and

$\bullet$ $\overline{\sigma}(v)=\phibar_n(v)$ for any $v\in V(T_n\cup R_n)$.

\vskip 2mm

{{\bf (5.12)} If $\sigma$ is a $(T_n \oplus R_n, D_n,\varphi_n)$-stable coloring, then $\sigma(f)=\varphi_n(f)$ for any 
edge $f$ on $T_n \cup R_n$, and $R_n$ is also a closure of $T_n(v_n)$ under $\sigma$. To justify this, note that, for any edge $f$ on $T_n$, 
this equality holds by Lemma \ref{hku}(iv). For any edge $f$ in $R_n -T_n$, we have $\varphi_n(f)\in\phibar_n(R_n)$ by the definition of 
$R_n$ and TAA. It follows from the above definition that $\sigma(f)=\varphi_n(f)$.  Since $\sigma$ is $(R_n,\emptyset,\varphi_n)$-stable, 
$R_n$ is a closure of $T_n(v_n)$ under $\sigma$ as well.

From Lemma \ref{sc2} it is clear that being $(T_n \oplus R_n, D_n,\cdot)$-stable is also an equivalence relation 
on ${\cal C}^k(G-e)$. Moreover, every $(T_n\vee R_n, D_n,\varphi_n)$-stable coloring is $(T_n\oplus R_n,D_n,\varphi_n)$-stable,
but the converse need not hold. 

Observe that, in the case of $\Theta_n=PE$, even when $T_n=T_{n,0}\subset T_{n,1} \subset 	\ldots \subset T_{n,q}\subset T_{n,q+1}=T$ is a 
hierarchy of $T$ (see (5.5)) under $\varphi_n$, and $T$ remains an ETT under a $(T_n, D_n,\varphi_n)$-stable coloring $\sigma_n$, there is
no guarantee that $T_n=T_{n,0}\subset T_{n,1} \subset 	\ldots \subset T_{n,q}\subset T_{n,q+1}=T$ is a hierarchy of $T$ under $\sigma_n$, 
because $R_n$ may not be a closure of $T_n(v_n)$ under $\sigma_n$. Nonetheless, we can establish the following statement.

\setcounter{theorem}{6}

\begin{lemma}\label{splitter}
Let $T$ be an ETT constructed from a $k$-triple $(G,e, \varphi)$ by using the Tashkinov series 
$\TT=\{(T_i, \phiv_{i-1}, S_{i-1}, F_{i-1}, \Theta_{i-1}): 1\le i \le n+1\}$. Suppose $\Theta_n=PE$ and
$T$ satisfies MP under $\varphi_n$. Let $T_n=T_{n,0}\subset T_{n,1} \subset 	\ldots \subset T_{n,q}\subset T_{n,q+1}=T$ be a 
hierarchy of $T$ under $\varphi_n$, and let $\sigma_n$ be a $(T_n \oplus R_n, D_n,\varphi_n)$-stable coloring. If $T$ can be built from 
$T_n$ by TAA under $\sigma_n$, then $T$ is also an ETT satisfying MP with respect to $\sigma_n$, 
and $T_n=T_{n,0} \subset T_{n,1}\subset \ldots \subset T_{n,q}\subset T_{n,q+1}=T$ remains to be a hierarchy of $T$ 
under $\sigma_n$.
\end{lemma}

{\bf Proof.}  By hypothesis, $\sigma_n$ is a $(T_n \oplus R_n, D_n,\varphi_n)$-stable coloring. So $\sigma(f)=\varphi_n(f)$ for any edge $f$ 
on $T_n \vee R_n$ and $R_n$ is also a closure of $T_n(v_n)$ under $\sigma_n$ by (5.12). Furthermore, $\sigma_n$ is a $(T_n,D_n,\varphi_n)$-stable 
coloring and hence is a $\phiv_n\bmod T_n$ coloring by (5.1) and Theorem~\ref{thm:tech10}(vi).  As $T$ can be built from $T_n$ by TAA under $\sigma_n$,
it is an ETT corresponding to $\sigma_n$ and satisfies MP under $\sigma_n$ by Theorem~\ref{thm:tech10}(vi). In view of the hierarchy of $T$ under 
$\varphi_n$, we obtain $T_n\vee R_n\subset T_{n,1}$. Hence $T_n=T_{n,0} \subset T_{n,1}\subset \ldots \subset T_{n,q}\subset T_{n,q+1}=T$ 
remains to be a hierarchy of $T$ under $\sigma_n$.  \qed

\vskip 2mm
From the above lemma we see that if $\Theta_n=PE$, $\sigma_n$ is a $(T_n \oplus R_n, D_n,\varphi_n)$-stable coloring, 
and $T$ is also an ETT under $\sigma_n$, then each hierarchy of $T$ under $\varphi_n$ is also a hierarchy under $\sigma_n$. Thus, 
to check whether a good hierarchy of $T$ remains good under a $(T_n \oplus R_n, D_n,\varphi_n)$-stable coloring in subsequent proofs, 
we shall only check whether it satisfies Definition~\ref{R2}, without even stating that it is a hierarchy by Lemma~\ref{splitter}.

We define one more term before proceeding. Let $T$ be a tree-sequence with respect to $G$ and $e$. A coloring 
$\pi \in {\cal C}^k(G-e)$ is called $(T, \varphi_n)$-{\em invariant}\label{tinvariant} 
if $\pi(f) = \varphi_n(f)$ for any $f\in E(T-e)$ and $\overline{\pi} (v) = \overline{\varphi}_n(v)$ for any $v\in V(T)$. 
Clearly, being $(T, \cdot)$-invariant is also an equivalence relation on ${\cal C}^k(G-e)$. Note that for any subset 
$C$ of $[k]$, a $(T, C, \varphi_n)$-stable coloring $\pi$ is also $(T, \varphi_n)$-invariant, provided that ${\pi}\langle 
T \rangle \subseteq \overline{\varphi}_n(T)\cup C$. Thus, if a coloring $\sigma_n$ is both $(T, \varphi_n)$-invariant and 
$(T_n \oplus R_n, D_n,\varphi_n)$-stable, then each hierarchy of $T$ under $\varphi_n$ is also a hierarchy under $\sigma_n$.

\vskip 3mm
\begin{lemma}\label{LEM:Stable} 
(Assuming (5.1)) Let $T$ be an ETT constructed from a $k$-triple $(G,e, \varphi)$ by using the Tashkinov series 
$\TT=\{(T_i, \phiv_{i-1}, S_{i-1}, F_{i-1}, \Theta_{i-1}): 1\le i \le n+1\}$. Suppose $T$ satisfies MP under 
$\varphi_n$. Let $\sigma_n$ be obtained from $\phiv_n$ by recoloring some $(\alpha, \beta)$-chains fully contained in 
$G-V(T)$. Then the following statements hold:
\begin{itemize}
\vspace{-1.5mm}
\item[(i)] $\sigma_n$ is $(T,D_n,\varphi_n)$-stable. In particular, $\sigma_n$ is $(T,\varphi_n)$-invariant.
Furthermore, if $\Theta_n=PE$ and $T_n\vee R_n\subseteq T$, then $\sigma_n$ is $(T_n\oplus R_n,D_n,\varphi_n)$-stable.  
\vspace{-2mm} 
\item[(ii)] $T$ is an ETT satisfying MP with respect to $\sigma_n$.
\vspace{-2mm}
\item[(iii)] If $T$ admits a good hierarchy $T_n=T_{n,0} \subset T_{n,1} \subset \ldots \subset 
T_{n,q+1}=T$ under $\varphi_n$, then this hierarchy of $T$ remains good under $\sigma_n$, with the same $\Gamma$-sets  (see Definition \ref{R2}). Furthermore, if $T$ is $(\cup_{\eta_h\in D_{n,q+1}}\Gamma^{q}_h)^-$-closed with respect to $\varphi_n$, 
then $T$ is also $(\cup_{\eta_h\in D_{n,q+1}}\Gamma^{q}_h)^-$-closed with respect to $\sigma_n$.  
\end{itemize}
\end{lemma}

{\bf Proof.} Since the recolored $(\alpha, \beta)$-chains are fully contained in $G-V(T)$, we have

(1) $\sigma_n(f) = \varphi_n(f)$ for each edge $f$ incident to $V(T)$ and $\phibar_n(v) =
\overline{\sigma}_n(v)$ for each $v\in V(T)$.

\noindent Our proof relies heavily on this observation.

(i) From (1) and definitions, it is clear that $\sigma_n$ is a $(T,D_n,\varphi_n)$-stable. In particular,
$\sigma_n$ is $(T,\varphi_n)$-invariant. Furthermore,  if $\Theta_n=PE$ and $T_n\vee R_n\subseteq T$, then $\sigma_n$ is $(T_n\vee R_n,D_n,
\varphi_n)$-stable, which implies that $\sigma_n$ is $(T_n\oplus R_n,D_n,\varphi_n)$-stable.

(ii) In view of (1), $T$ can also be obtained by TAA from $T_n$ (resp. $T_n+f_n$) under $\sigma_n$ when $\Theta_n=PE$ (resp. $\Theta_n$= $RE$ or $SE$).
Besides, $\sigma_n$ is a $(T_n,D_n,\varphi_n)$-stable coloring. Hence, by Theorem \ref{thm:tech10}(vi), 
$T$ remains to be an ETT and satisfies MP under $\sigma_n$. 

(iii) 
By (ii), $T$ is also an ETT under $\sigma_n$. By hypothesis, $T_n=T_{n,0} \subset T_{n,1} \subset \ldots  \subset T_{n,q+1}=T$ is a 
good hierarchy of $T$ under $\varphi_n$. Consider the $\Gamma$-sets specified in Definition \ref{R2} with respect 
to $\varphi_n$.  Using (1) it is routine to check that these $\Gamma$-sets satisfy all the conditions in Definition 
\ref{R2} with respect to $\sigma_n$. So the given hierarchy of $T$ remains good under $\sigma_n$, with the same 
$\Gamma$-sets.  Furthermore, if $T$ is $(\cup_{\eta_h\in D_{n,q+1}}\Gamma^{q}_h)^-$-closed with respect to $\varphi_n$, 
then $T$ is also $(\cup_{\eta_h\in D_{n,q+1}}\Gamma^{q}_h)^-$-closed with respect to $\sigma_n$. \qed

\vskip 2mm

In subsequent proofs, if we say that a hierarchy of an ETT under one coloring remains good under another coloring without giving 
the $\Gamma$-sets, we mean that it is a good hierarchy with the same $\Gamma$-sets.  

\section{Basic Properties}

As we have seen, Theorem \ref{thm:tech10}(i) follows from Theorems \ref{hierarchy} and \ref{good}. 
In the preceding section we have proved Theorem \ref{good}. The remainder of this paper is devoted to a
proof of Theorem \ref{hierarchy}. In this section we make some technical preparations; the reader is referred
to \cite{GG} for prototypes of some lemmas to be established herein. 

Let $T$ is an ETT that admits a good hierarchy $T_n=T_{n,0} \subset T_{n,1} \subset \ldots \subset
T_{n,q} \subset T_{n,q+1}=T$ and satisfies MP with respect to the generating coloring $\varphi_n$. To prove
Theorem \ref{hierarchy} (that is, $V(T)$ is elementary with respect to $\varphi_n$), we apply induction on
$q$; the induction base is Theorem \ref{thm:tech10}(i) for $T_n$. For convenience, we view $T_{n,0}$ as 
an ETT with $-1$ divider and $n$ rungs in the following assumption. Throughout this section we assume that 

{\bf (6.1)} In addition to (5.1), Theorem \ref{hierarchy} holds for every ETT that admits a good hierarchy
and satisfies MP, with $n$ rungs and at most $q-1$ dividers, where $q\ge 0$.  

Let us first prove two technical lemmas that will be used in the proof of Theorem \ref{hierarchy}. 

\begin{lemma}\label{interchange}

(Assuming (5.1)) Let $T$ be an ETT constructed from a $k$-triple $(G,e, \varphi)$ by using the Tashkinov series 
$\TT=\{(T_i, \phiv_{i-1}, S_{i-1}, F_{i-1}, \Theta_{i-1}): 1\le i \le n+1\}$. Suppose $\Theta_n=PE$ and $T$ satisfies
MP under $\varphi_n$. Let $\sigma_n$ be a $(T_n \oplus R_n, D_n,\varphi_n)$-stable coloring and let $\alpha$ and $\beta$
be two colors in $[k]$. Then the following statements hold:
\begin{itemize}
\vspace{-2mm}
\item[(i)] $\alpha$ and $\beta$ are $R_n$-interchangeable under $\sigma_n$ if $\alpha\in \overline{\sigma}_n(R_n)$;
\vspace{-2mm}
\item[(ii)] $\alpha$ and $\beta$ are $T_n$-interchangeable under $\sigma_n$ if $\alpha\in \overline{\sigma}_n(T_n)$;
\vspace{-2mm}
\item[(iii)] $\alpha$ and $\beta$ are $T_n\vee R_n$-interchangeable under $\sigma_n$ if $\alpha \in \overline{\sigma}_n
(T_n \vee R_n)$ is closed in $T_n\vee R_n$ under $\sigma_n$; and
\vspace{-2mm}
\item[(iv)] $\alpha$ and $\beta$ are $T_n\vee R_n$-interchangeable under $\sigma_n$ if $\alpha\in
\overline{\sigma}_n(T_n)$ and $\beta\in \overline{\sigma}_n(R_n)$.
\end{itemize}	
\end{lemma}	

{\bf Proof.} Since $\sigma_n$ is a $(T_n \oplus R_n, D_n,\varphi_n)$-stable coloring, it is  $(T_n, D_{n},\varphi_{n})$-stable 
by definition. Let $j=m(v_n)$. It follows from (5.1) and Theorem \ref{thm:tech10}(iii) that $\sigma_n$ is a
$(T_{j}(v_n)-v_n, D_{j-1},\varphi_{j-1})$-stable coloring. So $\sigma_n$ is $(T_{j-1}, D_{j-1}, \varphi_{j-1})$-stable and hence, 
by (5.1) and Theorem \ref{thm:tech10}(vi), it is a $\phiv_{j-1}\bmod T_{j-1}$ coloring, and every ETT corresponding to 
$(\sigma_n,T_{j-1})$ satisfies MP. Furthermore, $\sigma(f)=\varphi_n(f)$ 
for any edge $f$ in $T_n \cup R_n$ by (5.12) and $\overline{\sigma}_n(v)=\phibar_n(v)$ for all $v\in V(T_n \cup R_n)$. 

(i) Since $R_n$ is a closure of $T_n(v_n)$ under $\varphi_n$ and $\sigma_n$ is $(R_n, \emptyset, \varphi_n)$-stable,
$R_n$ is also a closure of $T_n(v_n)$ under $\sigma_n$. Since $\sigma_n$ is $\phiv_{j-1}\bmod T_{j-1}$,
$R_n$ is an ETT corresponding to $(\sigma_n,T_{j-1})$ and satisfies MP under $\sigma_n$. Let $\alpha$ and $\beta$ be as 
specified in the lemma. As $r(R_n)=j-1<n$, by (5.1) and Theorem \ref{thm:tech10}(ii), there is at most one $(\alpha, 
\beta)$-path with respect to $\sigma_n$ intersecting $R_n$. Hence $\alpha$ and $\beta$ are $R_n$-interchangeable under $\sigma_n$. 

Let us make some observations before proving statements (ii) and (iii). By (5.4), each color in 
$\overline{\varphi}_n(T_n) \cap \overline{\varphi}_n(R_n)$ is closed in $T_n\vee R_n$ with respect to $\varphi_n$.
Since $\sigma_n$ is a $(T_n \oplus R_n, D_n,\varphi_n)$-stable coloring, by definition we obtain

(1) each color in $\overline{\sigma}_n(T_n) \cap \overline{\sigma}_n(R_n)$ is closed in $T_n\vee R_n$ under $\sigma_n$.

(2) $\alpha$ and $\beta$ are $T_n$-interchangeable under $\sigma_n$ if $\alpha\in \overline{\sigma}_n(T_n)$, 
$\alpha \ne \delta_n$, and $\beta \ne \delta_n$. 

To justify (2), note that $\alpha\neq \gamma_n$, because $\gamma_n\notin\phibar_n(T_n)=\overline{\sigma}_n(T_n)$. So
$\alpha \notin S_n$. Nevertheless, the case $\beta=\gamma_n$ may occur.

Let us first consider the case when $\beta \ne \gamma_n$. Since $\sigma_n$ is $(T_n, D_{n},\varphi_{n})$-stable, $P_{v_n}(\gamma_n,\delta_n,\sigma_n)\cap T_n
=\{v_n\}$ by (5.1) and Theorem \ref{thm:tech10}(iv). Define $\sigma_n'=\sigma_n/P_{v_n}(\gamma_n,\delta_n,\sigma_n)$. By Lemma \ref{stablezang},
$\sigma_n'$ is $(T_n,D_{n-1}, \varphi_{n-1})$-stable.  From (5.1) and Theorem \ref{thm:tech10}(ii) 
we deduce that $\alpha$ and $\beta$ are $T_n$-interchangeable under $\sigma_n'$. So they are $T_n$-interchangeable 
under $\sigma_n$ because $\{\alpha,\beta\}\cap S_n=\emptyset$.

It remains to consider the case when $\beta=\gamma_n$. In this case, $f_n$ is the only edge in
$\partial_{\sigma_n, \gamma_n}(T_n)=\partial_{\varphi_n, \gamma_n}(T_n)$ by Lemma \ref{hku}(v). Since $V(T_n)$ is elementary
with respect to $\varphi_n$, it is also elementary with respect to $\sigma_n$. Together with $\partial_{\sigma_n, \alpha}(T_n)
=\emptyset$, we see that there is at most one $(\alpha,\gamma_n)$-path with respect to $\sigma_n$ intersecting $T_n$. So 
$\alpha$ and $\beta$ are $T_n$-interchangeable under $\sigma_n$. Thus (2) is established.

By (1), $\delta_n$ is closed in $T_n\vee R_n$ with respect to $\sigma_n$. So statement (ii) follows instantly from
(2) and statement (iii).

(iii) Assume the contrary: there are at least two $(\alpha,\beta)$-paths $P_1$ and $P_2$ with respect to 
$\sigma_n$ intersecting $T_n \vee R_n$. We may assume that

(3)  $\alpha\in \overline{\sigma}_n(T_n) \cap \overline{\sigma}_n(R_n)$. 

To justify this, let $A$ be the set of four ends of $P_1$ and $P_2$. Then at least two vertices from $A$ are 
outside $T_n\vee R_n$ because, by Lemma \ref{elementary}, $V(T_n\vee R_n)$ is elementary with respect to $\sigma_n$. 
Thus $P_1 \cup P_2$ contains two vertex-disjoint subpaths $Q_1$ and $Q_2$, which are two $T_n\vee R_n$-exit 
paths with respect to $\sigma_n$. Let $u \in V(T_n) \cap V(R_n)$, let $\eta \in \overline{\sigma}_n(u)$, and let $\sigma_n'=
\sigma_n/(G-T_n\vee R_n,\alpha, \eta)$. By (1), $\eta$ is closed in $T_n\vee R_n$ with respect to $\sigma_n$; so is $\alpha$ by hypothesis. Hence $\sigma_n'$ is a $(T_n \oplus R_n, D_n,\varphi_n)$-stable coloring, and $Q_1$ and $Q_2$ are two $T_n\vee R_n$-exit paths with respect to $\sigma_n'$. Since $P_u(\eta, \beta, \sigma_n')$ contains at most one of $Q_1$ and $Q_2$, replacing
$\sigma_n$ and $\alpha$ by $\sigma_n'$ and $\eta$, respectively, we obtain (3). 

Let $v$ be a vertex in $V(T_n)\cap V(R_n)$ with $\alpha\in \overline{\sigma}_n(v)$. Clearly, we may assume that $P_1=
P_v(\alpha,\beta,\sigma_n)$. By (i), we may further assume that $P_2$ is disjoint from $R_n$. So $P_2$ intersects $T_n$.
Therefore $\alpha$ and $\beta$ are not $T_n$-interchangeable under $\sigma_n$. Since $\gamma_n\notin\phibar_n(T_n)=
\overline{\sigma}_n(T_n)$, we have $\alpha\neq\gamma_n$. By (2), we may assume that $\alpha = \delta_n$ or $\beta=\delta_n$.

Suppose $\beta=\delta_n$. By Lemma \ref{hku}(v) and the definition of stable colorings, edges in $\partial_{\sigma_n, \delta_n}(T_n)$ 
are all incident to $V(T_n)\cap V(R_n)$. Thus both $P_1$ and $P_2$ intersect $V(T_n)\cap V(R_n)$, contradicting statement (i).

Suppose $\alpha=\delta_n$. By (1), $\delta_n$ is closed in $T_n\vee R_n$ under $\sigma_n$. Since $v_n$ is a maximum defective vertex, $V(T_n) \cap V(R_n)$ 
contains both ends of the uncolored edge $e$, so there exists a color $\theta \in \overline{\sigma}_n(T_n) \cap \overline{\sigma}_n(R_n)-\{\delta_n\}$. Let $\sigma_n''= \sigma_n/(G-T_n\vee R_n,\delta_n, \theta)$. Then $\sigma_n''$ is also 
$(T_n\oplus R_n, D_n,\varphi_n)$-stable. From the existence of $P_1$ and $P_2$, we see that $\theta$ and $\beta$ are not $T_n\vee R_n$-interchangeable under $\sigma_n''$, contradicting our observation (2) above the case when $\alpha \ne \delta_n$ and
$\beta \ne \delta_n$. 

(iv) Assume the contrary: there are at least two $(\alpha,\beta)$-paths $P_1$ and $P_2$ with respect to $\sigma_n$ 
intersecting $T_n \vee R_n$. Let $u$ be a vertex in $T_n$ with $\alpha \in \overline{\sigma}_n (u)$ and let $v$ be 
a vertex in $R_n$ with $\beta \in \overline{\sigma}_n(v)$. By (ii) (resp. (i)), $P_u(\alpha,\beta,\sigma_n)$
(resp.  $P_v(\alpha,\beta,\sigma_n)$) is the only $(\alpha,\beta)$-path with respect to $\sigma_n$ intersecting $T_n$ 
(resp. $R_n$). Hence we may assume that $P_1=P_u(\alpha,\beta,\sigma_n)$, $P_2=P_v(\alpha, \beta,\sigma_n)$ (rename subscripts if 
necessary), and $P_u(\alpha,\beta,\sigma_n) \ne P_v(\alpha,\beta,\sigma_n)$. Moreover, neither $P_u(\alpha,\beta,\sigma_n)$ 
nor $P_v(\alpha, \beta,\sigma_n)$ has an end in $V(T_n) \cap V(R_n)$, which in turn implies that    

(4) $u \in V(T_n)-V(R_n)$ and $v \in V(R_n)-V(T_n)$.
 
By (4) and statement (ii), $P_v(\alpha,\beta,\sigma_n)$ is disjoint from $T_n$. Let $\sigma_n'=\sigma_n/P_v(\alpha,
\beta,\sigma_n)$. By Lemma \ref{LEM:Stable}, $\sigma_n'$ is a $(T_n,D_n,\varphi_n)$-stable coloring.  By Lemma~\ref{elementary},  
$V(T_n\vee R_n)$ is elementary with respect to $\sigma_n$. Since $\alpha \in \overline{\sigma}_n(u)$ and $\beta\in \overline{\sigma}_n(v)$, 
from TAA we see that no edge in $R_n(v)-T_n(v_n)$ is colored by $\alpha$ or $\beta$
under both $\varphi_n$ and $\sigma_n$. Thus edges in $R_n(v)-T_n(v_n)$ are colored exactly the same under 
$\sigma_n'$ as under $\sigma_n$ and $\overline{\sigma}_n(x)= \overline{\sigma}_n'(x)$ for any $x \in  
V(R_n(v)-v)) \cup V(T_n)$.  Let $R_n'$ be a closure of $T_n(v_n)$ under $\sigma_n'$. Then $v \in V(R_n')$. In view of 
Lemma~\ref{elementary},  $V(T_n\vee R_n')$ is elementary with respect to $\sigma_n'$. However,  
$\alpha \in \overline{\sigma}'_n(u) \cap \overline{\sigma}'_n(v)$, a contradiction. \qed
\vskip 3mm




%

\begin{lemma}\label{step1}

(Assuming (6.1)) Let $T$ be an ETT satisfying MP constructed from a $k$-triple $(G,e, \varphi)$ by using the Tashkinov series 
$\TT=\{(T_i, \phiv_{i-1}, S_{i-1}, F_{i-1}, \Theta_{i-1}): 1\le i \le n+1\}$. Suppose $T$ has a good hierarchy 
$T_n=T_{n,0} \subset T_{n,1} \subset \ldots \subset T_{n,q} \subset T_{n,q+1}=T$. Let $p$ be a subscript with
$1\le p \le q$, and let $\alpha \in \overline{\varphi}_n(T_{n,p})$ and $\beta \in [k]-\{\alpha\}$.  If $\alpha$ is 
closed in $T_{n,p}$ under $\varphi_n$, then $\alpha$ and $\beta$ are $T_{n,p}$-interchangeable under $\varphi_n$.
\end{lemma}

{\bf Proof.} 	Assume the contrary: Let $p$ be the smallest index such that there exist two $(\alpha,\beta)$-paths $Q_1$ and $Q_2$ 
with respect to $\varphi_n$ intersecting $T_{n,p}$.  Let us make some simple observations about $T_{n,p}$ before proceeding. 
Since $T_{n,p}$ satisfies MP under $\varphi_n$ and $p\leq q$,

(1) $V(T_{n,p})$ is elementary with respect to $\varphi_n$
by (6.1) and  Theorem \ref{hierarchy}.

\noindent By hypthesis, $\alpha\in\phibar_n(T_{n,p})$ is closed in $T_{n,p}$ with respect to $\varphi_n$, which together with (1) yields

(2) $|T_{n,p}|$ is odd.

\noindent As $T_n=T_{n,0} \subset T_{n,1} \subset \ldots \subset T_{n,q} \subset T_{n,q+1}=T$ is a good hierarchy, 

(3) $T_{n,p}$ is $(\cup_{\eta_h\in D_{n,p}}\Gamma^{p-1}_h)^-$-closed with respect to $\varphi_n$ by Definition~\ref{R2}(v).

Depending on whether $\beta$ is contained in $\overline{\varphi}_n(T_{n,p})$, we consider two cases.

{\bf Case 1}. $\beta\in \overline{\varphi}_n(T_{n,p})$.
	
In this case, by (1) and (2), $|\partial_{\varphi_n, \beta}(T_{n,p})|$ is even. From the existence of $Q_1$ and $Q_2$, we see that
$G$ contains two vertex-disjoint $(T_{n,p}, \varphi_n, \{\alpha,\beta\})$-exit paths $P_1$ and $P_2$. 
For $i=1,2$, let $a_i$ and $b_i$ be the ends of $P_i$ with $b_i \in V(T_{n,p})$. Renaming subscripts if necessary, we may assume 
that $b_1 \prec b_2$.  We distinguish between two subcases according to the location of $b_2$.

{\bf Subcase 1.1}.  $b_2\in V(T_{n,p})-V(T_{n,p-1}^*)$.	
	
Since the edge on $P_2$ incident to $b_2$ is a boundary edge of $T_{n,p}$ and is colored by $\beta$, we have  
$\beta \in \Gamma^{p-1}_h$ for some $h$ with $\eta_h\in D_{n,p}$ by (3), which together with Definition \ref{R2}(i)
implies that $\beta \in \phibar(T_{n,p-1})$. Let $\gamma\in \overline{\varphi}_n(b_2)$. By the assumption of
the present subcase and Definition \ref{R2}(i), we have $\gamma\notin \Gamma^{p-1}$. Hence $\gamma$ is closed 
with respect to $\varphi_n$ in $T_{n,p}$ by (3) (see (5.10) for details).  So 

(4) both $\alpha$ and $\gamma$ are closed in $T_{n,p}$ under $\varphi_n$. 

Let $\mu_1=\varphi_n/(G-T_{n,p},\alpha,\gamma)$. By Lemma~\ref{LEM:Stable}, 

(5) the given hierarchy of $T_{n,p}$ remains good under $\mu_1$, with the same $\Gamma$-sets as those under $\varphi_n$  
(see Definition \ref{R2}). Furthermore, $T_{n,p}$ is $(\cup_{\eta_h\in D_{n,p}}\Gamma^{p-1}_h)^-$-closed under $\mu_1$ and $\beta \in \overline{\mu}_1(T_{n,p-1})$. 

Note that $P_1$ and $P_2$ are two $(T_{n,p}, \mu_1, \{\gamma, \beta\})$-exit paths.
Let $\mu_2=\mu_1/P_{b_2}(\gamma, \beta, \mu_1)$. Since $P_{b_2}(\gamma, \beta, \mu_1) \cap T_{n,p}=\{b_2\}$,  all edges incident to 
$V(T_{n,p}(b_2)-b_2)$ are colored the same under $\mu_2$ as under $\mu_1$. So $\beta \in \overline{\mu}_2(T_{n,p-1})$. By (5) and Lemma~\ref{LEM:Stable},
$T_n=T_{n,0} \subset T_{n,1}  \subset \ldots \subset T_{n,p-1} \subseteq T_{n,p}(b_2)-b_2$ is a good hierarchy 
of the ETT $T_{n,p}(b_2)-b_2$ (satisfying MP) under $\mu_2$, with the same $\Gamma$-sets as $T_{n,p}$ under $\varphi_n$. So

(6) $T_n=T_{n,0} \subset T_{n,1}  \subset \ldots \subset T_{n,p-1} \subset T_{n,p}(b_2)$ is a good hierarchy 
of the ETT $T_{n,p}(b_2)$ (satisfying MP) under $\mu_2$, with the same $\Gamma$-sets as $T_{n,p}$ under $\varphi_n$. 

Thus from (6) and (6.1) on Theorem \ref{hierarchy}, we conclude that $V(T_{n,p}(b_2))$ is elementary with respect to 
$\mu_2$. However, $\beta\in \overline{\mu}_2(T_{n,p-1}) \cap \overline{\mu}_2(b_2)$, a contradiction.
	
{\bf Subcase 1.2}. $b_2\in V(T_{n,p-1}^*)$.			
	
We propose to show that

(7) there exists a color $\theta \in\overline{\varphi}_n(T_n)$ that is closed in both $T_{n,0}^*$ and $T_{n,1}$ under 
$\varphi_n$ if $p=1$, and a color $\theta \in\overline{\varphi}_n(T_{n,p-1})$ that is closed in both $T_{n,p-1}$ 
and $T_{n,p}$ under $\varphi_n$ if $p\ge 2$.

Our proof is based on the following simple observation (see (3) in the proof of Theorem \ref{good}).

(8) $|\overline{\varphi}_n(T_n)|\ge 2n+11$ and $|D_{n,i}|\le |D_n|\le n$ for $0\le i \le q$. 

Let us first assume that $p=1$.  When $\Theta_n \ne PE$, let $\theta$ be a color in 
$\overline{\varphi}_n(T_n)-(\cup_{\eta_h\in D_{n,1}}\Gamma^{0}_h)$; such a color exists by (8). 
From Algorithm 3.1 we see that $T_n$ is closed under $\varphi_n$. 
By (3), $T_{n,1}$ is $(\cup_{\eta_h\in D_{n,1}}
\Gamma^{0}_h)^-$-closed under $\varphi_n$. So $\theta$ is as desired. When $\Theta_n=PE$, we have
$|\overline{\varphi}_n(T_n) \cap \overline{\varphi}_n(R_n)-\Gamma^0|\ge 4$ by Definition \ref{R2}(iv).
Let $\theta \in \overline{\varphi}_n(T_n) \cap \overline{\varphi}_n(R_n)-\Gamma^0-\{\delta_n\}$. 
It follows from (5.4) that $\theta$ is closed in $T_n \vee R_n$ under $\varphi_n$.  Since 
$T_{n,1}$ is $(\cup_{\eta_h\in D_{n,1}} \Gamma^{0}_h)^-$-closed with respect to $\varphi_n$,  $\theta$ also closed in $T_{n,1}$ under $\varphi_n$ as desired.
 
Next we assume that $p\geq 2$. By (8), we have $|\overline{\varphi}_n(T_{n,p-2})|\ge |\overline{\varphi}_n(T_{n})|
\ge 2n+11$ and $|D_{n,p-1}|\le |D_n|\le n$. So there exists a color $\theta$ in 
$\overline{\varphi}_n(T_{n,p-2})-(\cup_{\eta_h\in D_{n,p-1}}\Gamma^{p-2}_h)$. Since $\overline{\varphi}_n(T_{n,p-2}) 
\subseteq \overline{\varphi}_n(T_{n,p-1})$, we get $\theta \in \overline{\varphi}_n(T_{n,p-1})-(\cup_{\eta_h\in D_{n,p-1}}\Gamma^{p-2}_h)$. 
By Definition \ref{R2}(v), $\theta$ is closed in $T_{n, p-1}$ under 
$\varphi_n$.  From the definition of $\theta$ and Definition \ref{R2}(iii), it follows that 
$\theta \notin \Gamma^{p-1}$. So $\theta \in \overline{\varphi}_n(T_{n,p})- \Gamma^{p-1} \subseteq  \overline{\varphi}_n(T_{n,p}) - 
(\cup_{\eta_h\in D_{n,p}}\Gamma^{p-1}_h)$. By (3), $\theta$ is closed in $T_{n, p}$ under $\varphi_n$. Hence (7) is 
established.

Let $\mu_3=\varphi_n/(G-T_{n,p}, \alpha,\theta)$. Since both $\alpha$ and $\theta$ are closed in $T_{n,p}$ with
respect to $\varphi_n$, by Lemma~\ref{LEM:Stable}, $T_{n,p}$ admits a good hierarchy and satisfies
MP with respect to $\mu_3$. Thus $T_{n,p-1}$ also admits a good hierarchy and satisfies
MP with respect to $\mu_3$ if $p\ge 2$. By (7), $\theta$ is closed in $T_{n,0}^*$ if $p=1$ and closed in  $T_{n,p-1}$
if $p \ge 2$ under $\mu_3$.  Note that both $P_1$ and $P_2$ are  $(T_{n,p-1}^*, \mu_3, \{\theta, \beta\})$-exit 
paths. So $\theta$ and $\beta$ are not $T_{n,0}^*$-interchangeable under $\mu_3$ if $p=1$ 
and not $T_{n,p-1}$-interchangeable under $\mu_3$ if $p\ge 2$, which contradicts 
Lemma \ref{interchange}(iii) or the interchangeability property of $T_n$ when $p=1$, and the minimality assumption on $p$ when $p\ge 2$. 

\vskip 2mm
{\bf Case 2}. $\beta\notin \overline{\varphi}_n(T_{n,p})$.

In this case, $|\partial_{\varphi_n, \beta}(T_{n,p})|$ is odd and at least three by (1) and (2). From the existence of 
$Q_1$ and $Q_2$, we see that $G$ contains at least three $(T, \varphi_n, \{\alpha, \beta\})$-exit paths $P_1,P_2,P_3$. 
For $i=1,2,3$, let $a_i$ and $b_i$ be the ends of $P_i$ with $b_i \in V(T)$, and 
$f_i$ be the edge of $P_i$ incident to $b_i$.  Renaming subscripts if necessary, we may assume that $b_1\prec b_2 \prec b_3$. 

{\bf Subcase 2.1}. $b_3\in V(T_{n,p})-V(T_{n,p-1}^*)$.	

For convenience, we call the tuple $(\varphi_n, T_{n,p}, \alpha, \beta, P_1,P_2,P_3)$ a {\em counterexample} and use ${\cal K}$ 
to denote the set of all such counterexamples. With a slight abuse of notation, we still use $(\varphi_n, T_{n,p}, \alpha, \beta, P_1,
P_2,P_3)$ to denote a counterexample in ${\cal K}$ with the minimum $|P_1|+|P_2|+|P_3|$.  Let $\gamma\in\phibar(b_3)$.
By the hypothesis of the present subcase and Definition \ref{R2}(i), we have $\gamma\notin \Gamma^{p-1}$. So $\gamma$ 
is closed in $T_{n,p}$ under $\varphi_n$ by (3). Note that $\gamma$ might be some $\eta_h \in D_n$. 

Let $\mu_4=\varphi_n/(G-T_{n,p}, \alpha,\gamma)$. By Lemma~\ref{LEM:Stable},  $T_{n,p}$ admits a good hierarchy and satisfies MP under $\mu_4$. 
Note that $P_1,P_2, P_3$ are three $(T_{n,p}, \mu_4, \{\gamma, \beta\})$-exit paths.

Consider $\mu_5=\mu_4/P_{b_3}(\gamma, \beta, \mu_4)$. Clearly, $\beta \in \overline{\mu}_5(b_3)$ and $\beta\notin 
\Gamma^{p-1}$. Since $P_{b_3}(\gamma, \beta, \mu_4) \cap T_{n,p}=\{b_3\}$, it is easy to see that all edges incident to 
$V(T_{n,p}(b_3)-b_3)$ are colored the same under $\mu_5$ as under $\mu_4$. By Lemma \ref{LEM:Stable}, $T_n=T_{n,0} \subset T_{n,1}  
\subset \ldots \subset T_{n,p-1} \subset T_{n,p}(b_3)-b_3$ is a good hierarchy of the ETT $T_{n,p}(b_3)-b_3$ satisfying MP under $\mu_5$, with 
the same $\Gamma$-sets as $T_{n,p}$ under $\varphi_n$. So 

(9) $T_n=T_{n,0} \subset T_{n,1}  \subset \ldots \subset T_{n,p-1} \subset T_{n,p}(b_3)$ is a good hierarchy 
of $T_{n,p}(b_3)$ under $\mu_5$, with the same $\Gamma$-sets as $T_{n,p}$ under $\varphi_n$. 

Let $H$ be obtained from $T_{n,p}(b_3)$ by adding $f_1$ and $f_2$. Since $\beta\notin\Gamma^{p-1}$, it can be seen from (9) that

(10) $T_n=T_{n,0} \subset T_{n,1} \subset \ldots \subset T_{n,p-1} \subset H$ is a good hierarchy of 
$H$ under $\mu_5$, with the same $\Gamma$-sets as $T_{n,p}$ under $\varphi_n$. 

By (5.1) and Theorem~\ref{thm:tech10}(vi), $H$ satisfies MP under $\mu_5$.  Set $T'=H$. Let us grow $T'$
by using the following algorithm: 

(11) While there exists $f\in\partial(T')$ with $\mu_5(f)\in\overline{\mu}_5(T')$, do: set $T'=T'+f$  if 
the resulting $T'$ satisfies $\Gamma^{p-1}_h \cap \mu_5 \langle T'(v_{\eta_h})-T_{n,p-1}  \rangle =\emptyset$ for 
all $\eta_h \in D_{n,p-1}$. 

Note that this algorithm is exactly the same as Step 2 in Algorithm 5.6. From (11) we see that

(12)  $T'$ is $(\cup_{\eta_h\in D'_{n,p}}\Gamma^{p-1}_h)^-$-closed with respect to $\mu_5$, where $D_{n,p}'=\cup_{h\leq n}S_h -\overline{\mu}_5(T')$ (so $D_{n,p}' \subseteq D_{n,p-1}$). 

In view of (10) and (11), we conclude that

(13) $T_n=T_{n,0} \subset T_{n,1} \subset \ldots \subset T_{n,p-1} \subset T'$ is a good hierarchy of $T'$ under 
$\mu_5$, with the same $\Gamma$-sets as $T_{n,p}$ under $\varphi_n$. 

Clearly, $T'$ satisfies MP under $\mu_5$. By (13), (6.1), and Theorem \ref{hierarchy}, $V(T')$ is elementary with 
respect to $\mu_5$. Observe that none of $a_1,a_2,a_3$ is contained in $T'$, for otherwise, let $a_i \in V(T_2)$ for 
some $i$ with $1\le i \le 3$. Since $\{\beta,\gamma\}\cap \overline{\mu}_5(a_i) \ne \emptyset$ and 
$\beta \in \overline{\mu}_5(b_3)$, we obtain $\gamma \in \overline{\varphi}_2(a_i)$. Recall that $\beta,\gamma\notin\Gamma^{p-1}$. Hence from (11) we see that 
$P_1,P_2,P_3$ are all entirely contained in $G[T']$, which in turn implies  $\gamma \in \overline{\varphi}_2(a_j)$ 
for $j=1,2,3$. So $V(T')$ is not elementary with respect to $\mu_5$, a contradiction. Therefore, each $P_i$ contains a subpath 
$L_i$, which is a $T'$-exit path with respect to $\mu_5$. Since $f_1$ is not contained in $L_1$, we obtain $|L_1|+|L_2|+|L_3|<|P_1|+|P_2|+|P_3|$. 
Thus, in view of (12), the existence of the counterexample $(\mu_5, T', \gamma, 
\beta, L_1, L_2,L_3)$ violates the minimality assumption on $(\varphi_n, T_{n,p}, \alpha, \beta, P_1,P_2,P_3)$. 

{\bf Subcase 2.2}. $b_3\in V(T_{n,p-1}^*)$.	
	
The proof in this subcase is essentially the same as that in Subcase 1.2.  Let $\theta$ be a color as described in
(7). Consider $\mu_3=\varphi_n/(G-T_{n,p}, \alpha,\theta)$. Then we can verify that $\theta$ and $\beta$ are not $T_{n,0}^*$-interchangeable under $\mu_3$ if $p=1$ and not $T_{n,p-1}$-interchangeable under $\mu_3$ if $p\ge 2$, which contradicts 
Lemma \ref{interchange}(iii) or the minimality assumption on $p$; for the omitted details, see the proof in Subcase 1.2. \qed \\

Let us make some further preparations before proving Theorem \ref{hierarchy}. Let $T_n=T_{n,0} \subset T_{n,1} \subset 
\ldots \subset T_{n,q+1}=T$ be a good hierarchy of $T$ (see (5.5) and Definition \ref{R2}). Recall that $T_{n,0}^*= 
T_{n}\vee R_n$ if $\Theta_n=PE$ and $T_{n,0}^*= T_{n}$ otherwise, $T_{n,0}^* \subset T_{n,1}$ by (5.5), and $T_{n,i}^*=T_{n,i}$ 
if $i \ge 1$. Let $T$ be constructed from $T_{n,q}^*$ using TAA by recursively adding edges 
$e_1,e_2, \ldots, e_p$ and vertices $y_1,y_2, \ldots, y_p$, where $y_i$ is the end of $e_i$ outside 
$T(y_{i-1})$ for $i\ge 1$, with $T(y_0)=T_{n,q}^*$. Write $T=T_{n,q}^* \cup\{e_1,y_1,e_2,...,e_p,y_p\}$.  The {\em path number}\label{pathnumber} 
of $T$, denoted by $p(T)$, is defined to be the smallest subscript $i\in \{1,2,...,p\}$ such that the sequence $(y_i,e_{i+1},...,e_p,y_p)$ corresponds 
to a path in $G$. Note that $p(T)=p$ if this path contains the vertex $y_p$ only.

A coloring $\sigma_n \in {\cal C}^k(G-e)$ is called a $(T_{n,0}^*, D_n,\varphi_n)$-{\em weakly stable coloring}\label{t0strong}
if it is a $(T_n\oplus R_n,D_n,\varphi_n)$-stable coloring when $\Theta_n=PE$ and is a $(T_n,D_n,\varphi_n)$-stable 
coloring when $\Theta_n \ne PE$.  By Lemma \ref{hku}(iv) and (5.12), every $(T_{n,0}^*, D_n,\varphi_n)$-weakly stable 
coloring is $(T_{n,0}^*, \varphi_n)$-invariant.    

A coloring $\sigma_n \in {\cal C}^k(G-e)$ is called a $(T_{n,i}^*, D_n, \varphi_n)$-{\em weakly stable coloring}\label{tistrong}, with
$1 \le i \le q$, if it is both a $(T_{n,0}^*, D_n, \varphi_n)$-weakly stable and a $(T_{n,i}^*,\varphi_n)$-invariant coloring. 
By Lemma \ref{hku}(iv), every $(T_{n,i}^*, D_n, \varphi_n)$-stable coloring is $(T_{n,i}^*, D_n, \varphi_n)$-weakly stable. 
From Theorem~\ref{thm:tech10}(vi) it is also clear that, under a $(T_{n,i}^*, D_n, \varphi_n)$-weakly stable coloring $\sigma_n$, 
$T^*_{n,i}$ is an ETT satisfying MP (this statement will frequently be used directly in subsequent proofs without even 
citing Theorem~\ref{thm:tech10}(vi)). 

As stated before, our proof of Theorem \ref{hierarchy} proceeds by induction on $q$ (see (6.1)). The induction step will be 
carried out by contradiction. Throughout the remainder of this section and Subsection 7.1, $(T, \varphi_n)$ stands for 
a minimum counterexample to Theorem \ref{hierarchy}; that is,

{\bf (6.2)} $T$ is an ETT that admits a good hierarchy $T_n=T_{n,0} \subset T_{n,1} \subset \ldots \subset
T_{n,q} \subset T_{n,q+1}$

\hskip 10mm  $=T$  and satisfies MP with respect to the generating coloring $\varphi_n$;

{\bf (6.3)} subject to (6.2), $V(T)$ is not elementary with respect to $\varphi_n$;

{\bf (6.4)} subject to (6.2) and (6.3), $p(T)$ is minimum; and

{\bf (6.5)} subject to (6.2)-(6.4), $|T|-|T_{n,q}|$ is minimum.

\noindent Our objective is to find another counterexample $(T', \sigma_n)$ to Theorem \ref{hierarchy}, which violates 
the minimality assumption (6.4) or (6.5) on $(T, \varphi_n)$.

The following fact will be used frequently in subsequent proof. 

{\bf (6.6)} $V(T(y_{p-1}))$ is elementary with respect to $\varphi_n$.

\vskip 2mm

Let us exhibit some basic properties satisfied by the minimum counterexample $(T, \varphi_n)$ as specified above. 

\vskip 3mm

\begin{lemma}\label{9n}

For $0 \le i \le p-1$, the inequality
$$|\overline{\varphi}_n(T(y_i))- \phibar_n(T_{n,0}^*-V(T_n)) - \varphi_n \langle T(y_i) - T_{n,q}^* \rangle | \geq 2n+11$$
holds, where $T(y_0)=T_{n,q}^*$. Furthermore, if 
$$|\overline{\varphi}_n(T(y_i))- \phibar_n(T_{n,0}^*-V(T_n)) - \varphi_n \langle T(y_i) - T_{n,q}^* \rangle
-(\Gamma^{q}\cup D_{n,q}) |\le 4,$$
then there exist $7$ distinct colors $\eta_{h}\in D_{n,q}\cap \phibar_n(T(y_i))$ such that  $(\Gamma^{q}_h \cup \{\eta_{h}\})
\cap \varphi_n \langle T(y_i) - T_{n,q}^* \rangle=\emptyset$, where $\Gamma^{q}$ and
$\Gamma^{q}_h$ are introduced in Definition \ref{R2}.
\end{lemma}

{\bf Proof.} By (6.6), $T(y_{p-1})$ is elementary with respect to $\varphi_n$. Since the number of vertices in $T(y_i)-V(T_{n,q}^*)$ is $i$, and the number of edges in 
$T(y_i) - T_{n,q}^*$ is also $i$, we obtain $|\overline{\varphi}_n (T(y_i) - V(T_{n,q}^*))| \ge | \varphi_n \langle 
T(y_i) - T_{n,q}^* \rangle|$. Hence 

\begin{equation*}
\begin{aligned}
& \hskip 2mm |\overline{\varphi}_n(T(y_i))- \phibar_n(T_{n,0}^*-V(T_n)) - \varphi_n \langle T(y_i) - T_{n,q}^* \rangle |\\
   \ge  & \hskip 2mm |\overline{\varphi}_n(T(y_i))|- |\phibar_n(T_{n,0}^*-V(T_n))| - |\varphi_n \langle T(y_i) - T_{n,q}^* \rangle |\\
 \ge & \hskip 2mm |\overline{\varphi}_n(T(y_i))| - |\phibar_n(T_{n,0}^*-V(T_n))| - |\overline{\varphi}_n (T(y_i) - V(T_{n,q}^*))| \\
 = & \hskip 2mm |\overline{\varphi}_n(T_{n,q}^*)| - |\phibar_n(T_{n,0}^*-V(T_n))| \\
 \ge & \hskip 2mm |\overline{\varphi}_n(T_{n,0}^*)| - |\phibar_n(T_{n,0}^*)- \phibar_n(T_n)| \\
  = & \hskip 2mm |\overline{\varphi}_n(T_n)|\\
  \ge & \hskip 2mm 2n+11, 
\end{aligned}
\end{equation*}
where the last inequality can be found in the proof of Theorem \ref{good} (see (3) therein).  So the first inequality is
established.

Suppose the second inequality also holds. Then these two inequalities guarantee the existence of at least $2n+7$ colors in
the intersection $C$ of $\overline{\varphi}_n(T(y_i))- \phibar_n(T_{n,0}^*-V(T_n)) - \varphi_n \langle T(y_i) - T_{n,q}^* \rangle 
$ and $\Gamma^{q}\cup D_{n,q}$. By (5.6), we have $|D_{n,q}|\le  |D_{n}| \le n$
and $|\Gamma^{q}| \le 2 |D_{n,q}|\le 2n$. So $|\Gamma^{q}\cup D_{n,q}|\le 3n$. Since $|C|\le |\Gamma^{q}\cup D_{n,q}|$,
it follows that $2n+7 \le 3n$, which implies $n \ge 7$. Note that $C=\cup_{\eta_{h}\in D_{n,q}} (\Gamma^{q}_h \cup \{\eta_{h}\})
\cap C$ and $|(\Gamma^{q}_h \cup \{\eta_{h}\}) \cap C|\le 3$ for any $\eta_{h}$ in $D_{n,q}$. Since $|C|\ge 2n+7$ and $n\ge 7$, 
by the Pigeonhole Principle, there exist at least $7$ distinct colors $\eta_{h}$ in $D_{n,q}$, such that $| (\Gamma^{q}_h \cup 
\{\eta_{h}\}) \cap C|=3$, or equivalently, $\Gamma^{q}_h \cup \{\eta_{h}\} \subseteq C$.  For each of these $\eta_{h}$, clearly 
$\eta_{h}\in D_{n,q}\cap \phibar_n(T(y_i))$ and $(\Gamma^{q}_h \cup \{\eta_{h}\}) \cap \varphi_n \langle T(y_i) - T_{n,q}^* \rangle
=\emptyset$ by the definition of $C$. \qed

\vskip 2mm
Let $v$ be a vertex of $T$ and $T'\subset T$. By $T' \prec v$ we mean that $u \prec v$ for any $u \in V(T')$. Given a 
color $\alpha \in [k]$, we use $v_{\alpha}$ to denote the first vertex $u$ of $T$ in the order $\prec$ for which $\alpha \in 
\overline{\varphi}_n(u)$, if any, and defined to be the last vertex of $T$ in the order $\prec$ otherwise. 

\begin{lemma}\label{rutcor}
Suppose $q\geq 1$ and $\alpha\in \phibar_n(T_{n,q})$. If there exists a subscript $i$ with $0\le i \le q$, 
such that $\alpha$ is closed in $T^*_{n,i}$ with respect to $\varphi_n$,  then $\alpha \notin \varphi_n \langle 
T_{n,q}- T_{n,r}^* \rangle$, where $r$ is the largest such $i$. If there is no such subscript $i$, then 
$\alpha\in \cup_{\eta_h\in D_{n,j}} \Gamma^{j-1}_h \subseteq \Gamma^{j-1}$ for $1 \le j \le q$, 
$\Theta_n=PE$, $v_{\alpha}\in V(T_n)-V(R_n)$, and $\alpha \notin \varphi_n \langle T_{n,q}- T_n \rangle$.  
\end{lemma}

{\bf Proof.}
Recall that $T$ has a good hierarchy by (6.2).   Let us first assume the existence of a subscript $i$ with $0\le i \le q$, 
such that $\alpha$ is closed in $T^*_{n,i}$ with respect to $\varphi_n$. By definition, $r$ is the largest such $i$. 
As the statement holds trivially when $r=q$, we may assume that $r<q$. Let $s$ be an arbitrary index with $r+1\leq s\leq q$. 
From the definition of $r$, we see that $\alpha$ is not closed in $T_{n,s}$ with respect to $\varphi_n$. It follows from Definition \ref{R2}(v)
that $\alpha\in\Gamma^{s-1}_h$ for some $\eta_h\in D_{n,s}\subseteq D_{n,s-1}$. By definition, $D_{n,s}=\cup_{h\leq n}S_h-\phibar_n(T_{n,s})$, 
so $\eta_h\notin\phibar_n(T_{n,s})$ and hence $T_{n,s}(v_{\eta_h})= T_{n,s}$ (see paragraphs above Definition \ref{R2} for the notation 
$T_{n,s}(v_{\eta_h})$). Since $\alpha\in\Gamma^{s-1}_h$, Definition \ref{R2}(i) (with $j=s-1$)
implies $\alpha \notin \phiv_n \langle T_{n,s}(v_{\eta_h})-T_{n,s-1}^* \rangle=\phiv_n \langle T_{n,s}-T_{n,s-1}^* \rangle$. As this property holds 
for all $s$ with $r+1\leq s\leq q$, we get $\alpha \notin \varphi_n \langle T_{n,q}- T^*_{n,r} \rangle$.

Next we assume that there exists no subscript $i$ with $0\le i \le q$, such that $\alpha$ is closed in $T^*_{n,i}$ with 
respect to $\varphi_n$. Since $\alpha\in \phibar_n(T_{n,q})$, it follows from (5.10) that $\alpha  \in \phibar_n(T_{n,0}^*)$.
By Definition \ref{R2}(v), we obtain 

(1) $\alpha\in \cup_{\eta_h\in D_{n,j}}\Gamma^{j-1}_h \subseteq \Gamma^{j-1}$ for $1 \le j \le q$. 

\noindent Hence $\alpha\in\Gamma^j$ for all $0 \le j \le q-1$. From the definition of $\Gamma^0$, we see that $v_{\alpha} 
\in V(T_n)$. If $\Theta_n\neq PE$, then $\alpha$ would be closed in $T_{n}=T_{n,0}^*$ under $\varphi_n$, a contradiction. 
So $\Theta_n=PE$. Moreover, since $\alpha$ is not closed in $T^*_{n,0}$, by (5.4), we have
$v_{\alpha} \in V(T_n)-V(R_n)$. Since $R_n$ is a closure of $T_n(v_n)$ under $\varphi_n$, using (6.6) and TAA we obtain 

(2) $\alpha \notin \overline{\varphi}_n (R_n-V(T_n))$ and $\alpha \notin \varphi_n \langle R_n-T_n \rangle$.

(3) $\alpha \notin \varphi_n \langle T_{n,q}-T_{n,0}^* \rangle$.

Let $s$ be an arbitrary index with $1\le s \le q$.  By (1), we have $\alpha\in\Gamma^{s-1}_h$ for some $\eta_h \in D_{n, s}\subseteq D_{n,s-1}$. 
As $D_{n,s}=\cup_{h\leq n}S_h-\phibar_n(T_{n,s})$, there holds $\eta_h\notin\phibar_n(T_{n,s})$. So $T_{n,s}(v_{\eta_h})= T_{n,s}$.
From Definition \ref{R2}(i) (with $j=s-1$),  we deduce that $\alpha \notin \phiv_n \langle T_{n,s}(v_{\eta_h})-T_{n,s-1}^* \rangle=
\phiv_n \langle T_{n,s}-T_{n,s-1}^* \rangle$. Since this property is valid for all $s$ with $1\le s \le q$, we establish (3).

Combining (2) and (3), we conclude that $\alpha \notin \varphi_n \langle T_{n,q}- T_n \rangle$.  \qed

\vskip 3mm
Our proof of Theorem \ref{hierarchy} relies heavily on the following two technical lemmas.

\begin{lemma}\label{change}
Let $\alpha$ and $\beta$ be two colors in $\phibar_n(T(y_{p-1}))$. Suppose both $v_{\alpha} \prec v_{\beta}$ and $\alpha
\notin \varphi_n \langle T(v_{\beta}) -  T_{n,q}^* \rangle$ hold if $\{\alpha, \beta\}-\overline{\varphi}_n(T_{n,q}^*)
\ne \emptyset$. Then $P_{v_{\alpha}}(\alpha,\beta,\varphi_n)=P_{v_{\beta}}(\alpha,\beta,\varphi_n)$ if one of the 
following cases occurs:
\begin{itemize}
\vspace{-1.5mm}
\item[(i)] $q\geq 1$, and $\alpha\in \phibar_n(T_{n,q})$ or $\{\alpha, \beta\}\cap D_{n,q}=\emptyset$;
\vspace{-2mm}		
\item[(ii)] $q=0$, and $\alpha\in \phibar_n(T_{n})$ or $\{\alpha, \beta\} \cap D_{n}=\emptyset$; and
\vspace{-1.5mm}
\item[(iii)] $\alpha \in \phibar_n(T_{n,q}^*)$ and is closed in $T_{n,q}^*$ with respect to $\varphi_n$.
\vspace{-1.5mm}
\end{itemize}
Furthermore, in Case (iii), $P_{v_{\alpha}}(\alpha,\beta,\varphi_n)=P_{v_{\beta}}(\alpha,\beta,\varphi_n)$
is the only $(\alpha,\beta)$-path with respect to $\varphi_n$ intersecting $T_{n,q}^*$.
\end{lemma}

{\bf Proof.} Let $a=v_{\alpha}$ and $b=v_{\beta}$. We distinguish among three cases according to 
the locations of $a$ and $b$.

{\bf Case 1}. 	$\{a,b\} \subseteq V(T_{n,q}^*)$. 
	
By (6.6), $V(T_{n,q}^*)$ is elementary with respect to $\varphi_n$. So $a$ (resp. $b$) is the only vertex in $T_{n,q}^*$ 
missing $\alpha$ (resp. $\beta$). If both $\alpha$ and $\beta$ are closed in $T_{n,q}^*$ with respect to
$\varphi_n$, then no boundary edge of $T_{n,q}^*$ is colored by $\alpha$ or $\beta$. Hence  $P_{a}(\alpha,\beta,\varphi_n)=P_{b}(\alpha,\beta,\varphi_n)$ is the only path intersecting $T_{n,q}^*$. 
So we may assume that $\alpha$ or $\beta$ is not closed in $T_{n,q}^*$ with respect to $\varphi_n$. 
It follows that if $q=0$, then $\Theta_n=PE$, for otherwise, Algorithm 3.1 would imply that
both $\alpha$ and $\beta$ are closed in $T_{n}=T_{n,0}^*$, a contradiction. Therefore 

(1) $T_{n,0}^*= T_{n}\vee R_n$ if $q=0$.  

Let us first assume that precisely one of $\alpha$ and $\beta$ is closed in $T_{n,q}^*$ with respect to $\varphi_n$. In
this subcase, by Lemma \ref{step1} if $q \ge 1$ and by (1) and Lemma~\ref{interchange}(iii) if $q=0$, colors  
$\alpha$ and $\beta$ are $T_{n,q}^*$-interchangeable under $\varphi_n$, so $P_{a}(\alpha,\beta,\varphi_n)=P_{b}(\alpha,
\beta,\varphi_n)$ is the only path intersecting $T_{n,q}^*$.  

Next we assume that neither $\alpha$ nor $\beta$ is closed in $T_{n,q}^*$ with respect to $\varphi_n$. In this subcase,
we only need to show that $P_{a}(\alpha,\beta,\varphi_n)=P_{b}(\alpha,\beta,\varphi_n)$. Symmetry allows us to assume
that $a \prec b$. Let $r$ be the subscript with $\beta\in\phibar_n(T_{n,r}^*-V(T_{n,r-1}^*))$, where $0\leq r \le q$ 
and $T_{n,-1}^*=\emptyset$. Then $a, b \in V(T_{n,r}^*)$. By (6.2),  $T_n=T_{n,0} \subset T_{n,1} \subset \ldots \subset 
T_{n,q} \subset T_{n,q+1}=T$ is a good hierarchy of $T$. If $r \ge 1$, then $\beta$ is closed in $T_{n,r}$ with respect 
to $\varphi_n$ by Definition \ref{R2} (see (5.10)). From the above discussion about $T_{n,q}^*$ (with $r$ in place of $q$), 
we similarly deduce that $P_{a}(\alpha,\beta,\varphi_n)=P_{b}(\alpha,\beta,\varphi_n)$.  So we may assume that $r=0$.  If 
$\Theta_n\ne PE$, then both $\alpha$ and $\beta$ are closed in $T_n$ with respect to $\varphi_n$ (see Algorithm 3.1), so $P_{a}(\alpha,\beta,\varphi_n)=P_{b}(\alpha,\beta,\varphi_n)$ by (6.6). If $\Theta_n =PE$, then it follows from Lemma~\ref{interchange}(i), (ii) and (iv) that $P_{a}(\alpha,\beta,\varphi_n)=P_{b}(\alpha,\beta,\varphi_n)$. 
	
{\bf Case 2}. $\{a,b\} \cap V(T_{n,q}^*)=\emptyset$.
	
By the hypotheses of the present case and the present lemma, we have $\{\alpha, \beta\}\cap D_{n,q}=\emptyset$ if 
$q \ge 1$ and $\{\alpha, \beta\} \cap D_{n}=\emptyset$ if $q=0$. So 

(2) $\alpha, \beta \notin D_{n,q}\cup \phibar_n(T_{n,q}^*)$ if $q \ge 1$ and 
$\alpha, \beta \notin D_{n}\cup \phibar_n(T_{n,0}^*)$ if $q=0$. 

\noindent By the definitions of $D_n$ and $D_{n,q}$, we have $D_n \cup \phibar_n(T_n) \subseteq D_{n,q}\cup \phibar_n(T_{n,q}^*)$. By Lemma~\ref{hku}(iv) 
and Algorithm 3.1, we obtain $\varphi_n\langle T(a) \rangle\subseteq D_n\cup \phibar_n(T(a)-a)$ and $\varphi_n\langle T(b) \rangle\subseteq D_n\cup \phibar_n(T(b)-b)$. 
Since $T_{n,q}^*\prec a \prec b$ and $\alpha \notin \varphi_n \langle T(b) -  T_{n,q}^* \rangle$ (by the hypotheses of the present case and the present lemma), 
from (2) we see that 

(3) $\alpha,\beta\notin\varphi_n \langle T(b) \rangle$. 

\noindent Suppose on the contrary that $P_a(\alpha,\beta,\varphi_n)\neq P_{b}(\alpha,\beta,\varphi_n)$. 
Consider $\sigma_n=\varphi_n/P_{b}(\alpha,\beta,\varphi_n)$. By (2), $\sigma_n$ is a $(T_{n,q}^*, D_n,\varphi_n)$-stable coloring, so it is 
also a $(T_{n,q}^*, D_n,\varphi_n)$-weakly stable coloring. Thus, by (3) and (6.6),  every edge of $T(b)$ is colored under $\sigma_n$ the same as under 
$\varphi_n$. So $T(b)$ is still an ETT satisfying MP with respect to $\sigma_n$. Moreover, from (2), (3), and (6.6), we deduce that $T_n=T_{n,0} 
\subset T_{n,1}  \subset \ldots \subset T_{n,q} \subset T(b)$ is still a good hierarchy of $T(b)$ under $\sigma_n$, with the same $\Gamma$-sets 
as $T$ under $\varphi_n$ (see Definition \ref{R2}). As $\alpha\in\overline{\sigma}_n(a)\cap \overline{\sigma}_n(b)$, the pair $(T(b), \sigma_n)$ 
is a counterexample to Theorem \ref{hierarchy}, which contradicts the minimality assumption (6.5) on $(T, \varphi_n)$.

{\bf Case 3}. $a \in V(T_{n,q}^*)$ and $b \notin V(T_{n,q}^*).$ 

By the hypotheses of the present case and the present lemma, (6.6) and TAA, we obtain

(4) $\alpha \notin \varphi_n \langle T(b) -  T_{n,q}^* \rangle$ and $\beta \notin \overline{\varphi}_n (T(b)-b)$.
So $\beta$ is not used by any edge in $T(b) -  T_{n,q}^*$, except possibly $e_1$ when $q=0$ and $T_{n,0}^*=T_n$ 
(now $e_1=f_n$ in Algorithm 3.1 and $\beta \in D_n$).

Let us first assume that $\alpha$ is closed in $T_{n,q}^*$ with respect to $\varphi_n$. By Lemma \ref{step1} if $q \ge 1$ 
and by Lemma~\ref{interchange}(iii) or Theorem \ref{thm:tech10}(ii) (see (5.1)) if $q=0$, colors $\alpha$ and $\beta$ are $T_{n,q}^*$-interchangeable under $\varphi_n$. So $P_{a}(\alpha,\beta,\varphi_n)$ is the only $(\alpha,\beta)$-path 
intersecting $T_{n,q}^*$.  Suppose on the contrary that $P_{a}(\alpha,\beta,\varphi_n) \ne P_{b}(\alpha,\beta,\varphi_n)$. 
Then  $P_{b}(\alpha,\beta,\varphi_n)$ is vertex-disjoint from $T_{n,q}^*$ and hence contains no edge incident to $T_{n,q}^*$.

Consider $\sigma_n=\varphi_n/P_{b}(\alpha,\beta,\varphi_n)$.  It is routine to check 
that $\sigma_n$ is a $(T_{n,q}^*, D_n,\varphi_n)$-weakly stable coloring, and $T(b)$ is an ETT 
satisfying MP with respect to $\sigma_n$. Moreover, $T_n=T_{n,0} \subset T_{n,1}  \subset \ldots \subset T_{n,q} 
\subset T(b)$ is a good hierarchy of $T(b)$ under $\sigma_n$, with the same $\Gamma$-sets as $T$ under $\varphi_n$,
by (4). As $\alpha\in\overline{\sigma}_n(a)\cap \overline{\sigma}_n(b)$, the pair $(T(b), \sigma_n)$ 
is a counterexample to Theorem \ref{hierarchy}, which contradicts the minimality assumption (6.5) on $(T, \varphi_n)$. 

So we assume hereafter that 

(5) $\alpha$ is not closed in $T_{n,q}^*$ with respect to $\varphi_n$. 

Hence our objective is to show that $P_{a}(\alpha,\beta,\varphi_n)=P_{b}(\alpha,\beta,\varphi_n)$. Assume the contrary: $P_{a}(\alpha,\beta,\varphi_n)\ne P_{b}(\alpha,\beta,\varphi_n)$. We distinguish between two subcases according to 
the value of $q$.

{\bf Subcase 3.1.} $q=0$. 

By the hypothesis of the present lemma, $\alpha\in \phibar_n(T_{n})$ or $\{\alpha, \beta\} \cap D_{n}=\emptyset$.
So $\alpha\notin D_n$. From (5) and Algorithm 3.1 we deduce that $T_{n,0}^* \ne T_n$. Hence

(6) $\Theta_n=PE$, which together with (5) and (5.4) yields $a \notin V(T_n)\cap V(R_n)$.

Consider $\sigma_n=\varphi_n/P_{b}(\alpha,\beta,\varphi_n)$. We claim that

(7) $\sigma_n$ is a $(T_{n,0}^*, D_n,\varphi_n)$-weakly stable coloring. 

To justify this, note that if $a\in V(T_n)-V(R_n)$, then $\alpha,\beta\notin \phibar_n (R_n)$ by (6.6) and the 
hypothesis of the present case. By definition, $\sigma_n$ is $(R_n,\emptyset,\varphi_n)$-stable. In view of Lemma~\ref{interchange}(ii), $P_{b}(\alpha,\beta,\varphi_n)$ is disjoint from $T_n$ and hence contains no edge 
incident to $T_n$. So $\sigma_n$ is $(T_n,D_n,\varphi_n)$-stable. Hence (7) holds.   Suppose $a\in V(R_n)-V(T_n)$. 
By the hypothesis of the present lemma, $\{\alpha, \beta\} \cap D_{n}=\emptyset$. By (6.6), we also have 
$\alpha,\beta\notin \phibar_n(T_n)$. Thus $\alpha,\beta\notin \phibar_n(T_n) \cup D_n$.  By definition,
$\sigma_n$ is $(T_n,D_n,\varphi_n)$-stable. Using Lemma~\ref{interchange}(i), $P_{b}(\alpha,\beta,\varphi_n)$ 
is disjoint from $R_n$ and hence contains no edge incident to $R_n$. By definition, $\sigma_n$ is $(R_n,\emptyset,
\varphi_n)$-stable. Therefore (7) is true.

From (4), (7) and (6.6) we see that $\sigma_n(f)=\varphi_n(f)$ for each $f\in E(T(b))$ and $\overline{\sigma}_n(u)=\phibar_n(u)$ for each $u\in V(T(b)-b)$ 
(recall that every $(T_{n,0}^*, D_n,\varphi_n)$-weakly stable coloring is $(T_{n,0}^*, \varphi_n)$-invariant). Furthermore, $T(b)$ is an ETT 
satisfying MP with respect to $\sigma_n$, and $T_n=T_{n,0} \subset T(b)$ is a good hierarchy of $T(b)$ under $\sigma_n$, 
with the same $\Gamma$-sets as $T$ under $\varphi_n$. As $\alpha\in\overline{\sigma}_n(a)\cap \overline{\sigma}_n(b)$, 
the pair $(T(b), \sigma_n)$ is a counterexample to Theorem \ref{hierarchy}, which contradicts the minimality assumption 
(6.5) on $(T, \varphi_n)$. 

{\bf Subcase 3.2.} $q\ge 1$. 	

Let us first assume that $\alpha$ is closed in $T^*_{n,i}$ with respect to $\varphi_n$ for some $i$ with
$0\le i \le q$. Let $r$ be the largest subscript $i$ with this property. Then $r\le q-1$ by (5). Note that $\alpha\in\phibar_n(T^*_{n,r})$ since $\alpha$ is closed in $T^*_{n,r}$.
By Lemma \ref{rutcor}, we have $\alpha \notin \varphi_n \langle T_{n,q}- T_{n,r}^* \rangle$,  which together 
with (4) yields 

(8) $\alpha\notin \phiv_n \langle T(b)-T_{n,r}^* \rangle$. 

By Lemma \ref{step1} if $r \ge 1$ and by Theorem \ref{thm:tech10}(ii) or Lemma~\ref{interchange}(iii) if $r=0$, colors  
$\alpha$ and $\beta$ are $T_{n,r}^*$-interchangeable under $\varphi_n$. So $P_{a}(\alpha,\beta,\varphi_n)$ is the only $(\alpha,\beta)$-path with respect to $\varphi_n$ intersecting $T_{n,r}^*$. Hence $P_{b}(\alpha,\beta,\varphi_n)$ is 
vertex-disjoint from $T_{n,r}^*$ and therefore contains no edge incident to  $T_{n,r}^*$. Let $\sigma_n=\varphi_n/P_{b}(\alpha,\beta,\varphi_n)$. By Lemma \ref{LEM:Stable}, $\sigma_n$ is a $(T_{n,r}^*, D_n,\varphi_n)$-weakly stable coloring, and $T_{n,r}^*$ is an ETT having a good hierarchy and satisfying MP with respect to $\sigma_n$. By (4)
and TAA, $\beta$ is not used by any edge in $T(b) -  T_{n,r}^*$, except possibly $e_1$ when $r=0$ and $T_{n,0}^*=T_n$ 
(now $e_1=f_n$ in Algorithm 3.1 and $\beta \in D_n$). Since $\sigma_n$ is $(T_n, D_n,\varphi_n)$-stable,
it follows from (8) and (6.6) that $\sigma_n(f)=\varphi_n(f)$ for each $f\in E(T(b))$ and $\overline{\sigma}_n(u)=
\phibar_n(u)$ for each $u\in V(T(b)-b)$. So $T(b)$ is an ETT satisfying MP with respect to $\sigma_n$. Moreover, 

(9) $T_n=T_{n,0} \subset T_{n,1}  \subset \ldots \subset T_{n,q} \subset T(b)$ is a good hierarchy of $T(b)$ under $\sigma_n$, 
with the same $\Gamma$-sets as $T$ under $\varphi_n$.

Since $\sigma_n(f)=\varphi_n(f)$ for each $f\in E(T(b))$ and $\overline{\sigma}_n(u)=\phibar_n(u)$ for each $u\in V(T(b)-b)$, 
to justify (9), it suffices to verify that Definition \ref{R2}(v) is satisfied with respect to $\sigma_n$; that is,
$T_{n,j}$ is $(\cup_{\eta_h\in D_{n,j}}\Gamma^{j-1}_h)^-$-closed with respect to $\sigma_n$ for $1\le j \le q$. 
As the statement holds trivially if $P_{b}(\alpha,\beta,\varphi_n)$ is vertex-disjoint from $T_{n,j}$, we may
assume that $P_b(\alpha, \beta, \varphi_n)$ intersects $T_{n,j}$. Thus $r+1\le j\le q$. Observe that $\alpha \in \cup_{\eta_h\in D_{n,j}} \Gamma^{j-1}_h$, for otherwise, $\alpha$ is closed in $T_{n,j}$ with respect to $\varphi_n$ by Definition \ref{R2}(v),
contradicting the definition of $r$. By (6.6), we also obtain $\beta \notin \overline{\varphi}_n(T_{n,j})$. Consequently,
$T_{n,j}$ is $(\cup_{\eta_h\in D_{n,j}}\Gamma^{j-1}_h)^-$-closed with respect to $\sigma_n$. (Note that $\alpha$ may 
become closed in $T_{n,j}$ with respect to $\sigma_n$. Yet, even in this situation the desired statement is true.)
This proves (9).  

As $\alpha\in\overline{\sigma}_n(a)\cap \overline{\sigma}_n(b)$, the existence of $(T(b), \sigma_n)$ contradicts the minimality assumption (6.5) on $(T, \varphi_n)$.  

Next we assume that $\alpha$ is not closed in $T^*_{n,i}$ with respect to $\varphi_n$ for any $i$ with $0\le i \le q$. By the hypothesis of the present
subcase, $q\geq 1$. In view of Lemma \ref{rutcor}, we obtain

(10) $\alpha\in \cup_{\eta_h\in D_{n,j}}\Gamma^{j-1}_h \subseteq \Gamma^{j-1}$ for $1 \le j \le q$,  $\Theta_n=PE$,  
$a \in V(T_n)-V(R_n)$, and $\alpha \notin \varphi_n \langle T_{n,q}- T_n \rangle$.  

It follows from (4), (10) and TAA that
 
(11) $\alpha\notin \phiv_n \langle T(b)-T_n \rangle$ and $\beta \notin \phiv_n \langle T(b)- T_{n,0}^* \rangle$. 

\noindent Since $R_n$ is a closure of $T_n(v_n)$ under $\varphi_n$, using (10), (6.6) and TAA we obtain 

(12) $\alpha, \beta \notin \overline{\varphi}_n (R_n)$ and $\beta \notin \varphi_n \langle R_n-T_n \rangle$.

By Lemma~\ref{interchange}(ii), colors $\alpha$ and $\beta$ are $T_n$-interchangeable under $\varphi_n$. So $P_{a}(\alpha,\beta,\varphi_n)$ is the only $(\alpha,\beta)$-path with respect to $\varphi_n$ intersecting 
$T_n$. Hence $P_{b}(\alpha,\beta,\varphi_n)$ is vertex-disjoint from $T_n$ and therefore contains no edge incident 
to  $T_n$.  Consider $\sigma_n=\varphi_n/P_{b}(\alpha,\beta,\varphi_n)$. By Lemma \ref{LEM:Stable}, 
$\sigma_n$ is a $(T_n, D_n,\varphi_n)$-stable coloring, and $T_n$ is an ETT satisfying MP with respect to $\sigma_n$.  
From (11) and (12) we further deduce that $\sigma_n$ is a $(T_{n,0}^*, D_n, \varphi_n)$-weakly stable coloring,
$\sigma_n(f)=\varphi_n(f)$ for each $f\in E(T(b))$, and $\overline{\sigma}_n(u)=\phibar_n(u)$ for each $u\in V(T(b)-b)$. 
So $T(b)$ is an ETT satisfying MP with respect to $\sigma_n$. Moreover, $T_n=T_{n,0} \subset T_{n,1}  \subset \ldots 
\subset T_{n,q} \subset T(b)$ is a good hierarchy of $T(b)$ under $\sigma_n$, with the same $\Gamma$-sets as $T$ 
under $\varphi_n$ (see (10) and the proof of (9) for omitted details). As $\alpha\in\overline{\sigma}_n(a)\cap \overline{\sigma}_n(b)$, the existence of $(T(b), \sigma_n)$ contradicts the minimality assumption (6.5) on 
$(T, \varphi_n)$.  \qed

\vskip 3mm

\begin{lemma}\label{stablechange}
Let $\alpha$ and $\beta$ be two colors in $\phibar_n(T(y_{p-1}))$, let $Q$ be an $(\alpha,\beta)$-chain with
respect to $\varphi_n$, and let $\sigma_n=\varphi_n/Q$. Suppose one of the following cases occurs:
\begin{itemize}
\vspace{-2mm}
\item[1)] $q\geq 1$, $\alpha\in\phibar_n(T_{n,q})$, and $Q$ is an $(\alpha,\beta)$-path disjoint from $P_{v_{\alpha}}(\alpha,\beta,\varphi_n)$;
\vspace{-2mm}
\item[2)] $q=0$, $\alpha\in\phibar_n(T_{n})$, or $\alpha\in\phibar_n(T_{n,0}^*)$ with $\alpha,\beta\notin D_n$, and 
$Q$ is an $(\alpha,\beta)$-path disjoint from $P_{v_{\alpha}}(\alpha,\beta,\varphi_n)$; and
\vspace{-2mm}
\item[3)] $T_{n,q}^*\prec v_{\alpha} \prec v_{\beta}$, $\alpha,\beta\notin D_{n,q}$, $\alpha\notin \phiv_n \langle T(v_{\beta})-T(v_{\alpha}) \rangle$, and $Q$ is an arbitrary $(\alpha,\beta)$-chain.

\vspace{-2mm}
\end{itemize} 
Then the following statements hold:
\begin{itemize}
\vspace{-1.5mm}
\item[(i)] $\sigma_n$ is a $(T_{n,q}^*, D_n,\varphi_n)$-weakly stable coloring;
\vspace{-2mm} 
\item[(ii)] $T_{n,q}^*$ is an ETT satisfying MP with respect to $\sigma_n$; and
\vspace{-2mm}
\item[(iii)] if $q \ge 1$, then  $T_n=T_{n,0} \subset T_{n,1} \subset \ldots \subset T_{n,q}$ is a good hierarchy 
of $T_{n,q}$ under $\sigma_n$, with the same $\Gamma$-sets (see Definition \ref{R2}) as $T$ under $\varphi_n$,
and $T_{n,q}$ is $(\cup_{\eta_h\in D_{n,q}}\Gamma^{q-1}_h)^-$-closed with respect to $\sigma_n$.
\vspace{-2mm}
\end{itemize}
Furthermore, in Case 3, $T$ is also an ETT satisfying MP with respect to $\sigma_n$, and $T_n=T_{n,0} \subset T_{n,1}
\subset \ldots \subset T_{n,q}\subset T_{n,q+1}=T$ remains to be a good hierarchy of $T$ under $\sigma_n$, with the 
same $\Gamma$-sets (see Definition \ref{R2}) as $T$ under $\varphi_n$. 
\end{lemma}

\noindent {\bf Remark.} In the proof of Theorem \ref{hierarchy}, frequently we need to check whether a ``smaller" 
counterexample $T'$ with $T_{n,q} \subset T'$ has a good hierarchy with the same $\Gamma$-sets under $\sigma_n$ as 
$T$ under $\varphi_n$. Lemma \ref{stablechange} is established to fulfill such needs: We shall use the above Statement 
(iii) to ensure that Definition \ref{R2}(i)-(v) are satisfied by $T_{n,q}$ and that Definition \ref{R2}(v) is satisfied by 
$T'$. Since the $\Gamma$-sets used under $\sigma_n$ are the same as those under $\varphi_n$, Definition \ref{R2}(ii)-(iv) 
are automatically satisfied by $T'$. One technical question remains unanswered: How can we verify that Definition 
\ref{R2}(i) is satisfied by $T'$?  It is only a straightforward matter, as we shall see.

\vskip 3mm  

{\bf Proof of Lemma \ref{stablechange}}. 	Write $a=v_{\alpha}$ and $b=v_{\beta}$. Let us consider the three cases 
described in the lemma separately.

{\bf Case 1.} $q\geq 1$, $\alpha\in\phibar_n(T_{n,q})$, and $Q$ is an $(\alpha,\beta)$-path disjoint from $P_{v_{\alpha}}(\alpha,\beta,\varphi_n)$. 

We distinguish between two subcases according to the location of $b$.

{\bf Subcase 1.1.} $b \in V(T_{n,q})$.  

Let us first assume that there exists a subscript $i$ with $0\le i \le q$, such that $\alpha$ or $\beta$ is closed 
in $T^*_{n,i}$ with respect to $\varphi_n$. Let $r$ be the largest such $i$. By (5.10) and Lemma \ref{rutcor}, we have

(1) $\{a, b\} \subseteq V(T_{n,r}^*)$ and $\alpha, \beta \notin \varphi_n \langle T_{n,q}- T_{n,r}^* \rangle$. 

(2) $\alpha$ and $\beta$ are $T_{n,r}^*$-interchangeable under $\varphi_n$. So $P_a(\alpha,\beta,
\varphi_n)=P_b(\alpha,\beta,\varphi_n)$. 

To justify this, note that if $r \geq 1$, then (2) holds by Lemma~\ref{step1}. So we assume that $r=0$. 
Then $\alpha$ or $\beta$ is closed in $T_{n,0}^*$ with respect to $\varphi_n$.  Hence, by Lemma~\ref{interchange}(iii) 
if $\Theta_n= PE$ and by (5.1) and Theorem \ref{thm:tech10}(ii) otherwise, $\alpha$ and $\beta$ are $T_{n,0}^*$-interchangeable 
under $\varphi_n$. This proves (2). 

It follows from (2) that $Q$ is vertex-disjoint from $T_{n,r}^*$ and hence contains no edge incident to $T_{n,r}^*$. 
By Lemma \ref{LEM:Stable}, $\sigma_n=\varphi_n/Q$ is a $(T_{n,r}^*,D_n,\varphi_n)$-weakly stable coloring, and $T_{n,r}^*$ 
is an ETT satisfying MP with respect to $\sigma_n$. By (1) and (6.6), we obtain $\sigma_n(f)=\varphi_n(f)$ for each edge $f$ of $T_{n,q}$ and $\overline{\sigma}_n(u)=\phibar_n(u)$ for each vertex $u$ of $T_{n,q}$.  Therefore $\sigma_n$ is a $(T_{n,q}, D_n,\varphi_n)$-weakly stable coloring. By the definition of $r$, for any $r+1\le j \le q$ and $\theta \in \{\alpha,\beta\}$, we have $\partial_{\varphi_n, \theta}(T_{n,j})\ne \emptyset$, so $\theta \in \cup_{\eta_h\in D_{n,j}}\Gamma^{j-1}_h$ by
Definition \ref{R2}(v). It is then routine to check that $T_n=T_{n,0} \subset T_{n,1} \subset \ldots \subset T_{n,q}$ is a good 
hierarchy of $T_{n,q}$  under $\sigma_n$, with the same $\Gamma$-sets as $T$ under $\varphi_n$\footnote{See the justification
of (9) in the proof of Lemma \ref{change} for omitted details. Note that $\alpha$ or $\beta$ may become closed in $T_{n,j}$ 
with respect to $\sigma_n$ for some $j$ with $r+1\le j \le q$. Yet, even in this situation Definition \ref{R2}(v) remains valid with respect to $\sigma_n$.}, and $T_{n,q}$ is $(\cup_{\eta_h\in D_{n,q}}\Gamma^{q-1}_h)^-$-closed with respect to $\sigma_n$. 

Next we assume that there exists no subscript $i$ with $0\le i \le q$, such that $\alpha$ or $\beta$ is closed 
in $T^*_{n,i}$ with respect to $\varphi_n$. By Lemma \ref{rutcor}, we have

(3) $\alpha, \beta\in \cup_{\eta_h\in D_{n,j}}\Gamma^{j-1}_h \subseteq \Gamma^{j-1}$ for $1 \le j \le q$, 
$\Theta_n=PE$, $a, b\in V(T_n)-V(R_n)$, and $\alpha, \beta \notin \varphi_n \langle T_{n,q}- T_n \rangle$.  

Since $R_n$ is a closure of $T_n(v_n)$ under $\varphi_n$, using (6.6) and TAA we obtain 

(4) $\alpha, \beta \notin \overline{\varphi}_n (R_n)$.

By Lemma~\ref{interchange}(ii), colors $\alpha$ and $\beta$ are $T_n$-interchangeable under $\varphi_n$. So $P_{a}(\alpha,\beta,\varphi_n)$ is the only $(\alpha,\beta)$-path with respect to $\varphi_n$ intersecting 
$T_n$. Hence $Q$ is vertex-disjoint from $T_n$ and therefore contains no edge incident to  $T_n$.  
By Lemma \ref{LEM:Stable}, $\sigma_n=\varphi_n/Q$ is a $(T_n, D_n,\varphi_n)$-stable coloring, and 
$T_n$ is an ETT satisfying MP with respect to $\sigma_n$. By (3), (4) and (6.6), we further deduce that
$\sigma_n$ is a $(T_{n,0}^*, D_n,\varphi_n)$-stable coloring, $\sigma_n(f)=\varphi_n(f)$ for each 
edge $f$ of $T_{n,q}$, and $\overline{\sigma}_n(u)=\phibar_n(u)$ for each vertex $u$ of $T_{n,q}$. It is then routine 
to check that the desired statements hold.  

{\bf Subcase 1.2.} $b \notin V(T_{n,q})$.  

Let us first assume that there exists a subscript $i$ with $0\le i \le q$, such that $\alpha$ is closed 
in $T^*_{n,i}$ with respect to $\varphi_n$. Let $r$ be the largest such $i$. By (5.10), Lemma \ref{rutcor} and TAA, we have

(5) $a \subseteq V(T_{n,r}^*)$ and $\alpha \notin \varphi_n \langle T_{n,q}- T_{n,r}^* \rangle$. Furthermore,
no edge in $T_{n,q}- T_{n,r}^*$ is colored by $\beta$, except possibly $e_1$ when $r=0$ and $T_{n,0}^*=T_n$ 
(now $e_1=f_n$ in Algorithm 3.1 and $\beta \in D_n$). 

Using the same argument as that of (2), we obtain 

(6) $\alpha$ and $\beta$ are $T_{n,r}^*$-interchangeable under $\varphi_n$. 

It follows from (6) that $Q$ is vertex-disjoint from $T_{n,r}^*$ and hence contains no edge incident to $T_{n,r}^*$. 
By Lemma \ref{LEM:Stable}, $\sigma_n=\varphi_n/Q$ is a $(T_{n,r}^*,D_n,\varphi_n)$-weakly stable coloring, and $T_{n,r}^*$ 
is an ETT satisfying MP with respect to $\sigma_n$. Using (5), we obtain $\sigma_n(f)=\varphi_n(f)$ for each edge $f$ 
of $T_{n,q}$ and $\overline{\sigma}_n(u)=\phibar_n(u)$ for each vertex $u$ of $T_{n,q}$.  Therefore $\sigma_n$ is a $(T_{n,q}, D_n,\varphi_n)$-weakly stable coloring, $T_n=T_{n,0} \subset T_{n,1} \subset \ldots \subset T_{n,q}$ is a good 
hierarchy of $T_{n,q}$  under $\sigma_n$, with the same $\Gamma$-sets as $T$ under $\varphi_n$, and $T_{n,q}$ is 
$(\cup_{\eta_h\in D_{n,q}}\Gamma^{q-1}_h)^-$-closed with respect to $\sigma_n$ (see the justification
of (9) in the proof of Lemma \ref{change} for omitted details). 

Next we assume that there exists no subscript $i$ with $0\le i \le q$, such that $\alpha$ is closed 
in $T^*_{n,i}$ with respect to $\varphi_n$. By Lemma \ref{rutcor}, we have

(7) $\alpha \in \cup_{\eta_h\in D_{n,j}}\Gamma^{j-1}_h \subseteq \Gamma^{j-1}$ for $1 \le j \le q$, 
$\Theta_n=PE$, $a \in V(T_n)-V(R_n)$, and $\alpha \notin \varphi_n \langle T_{n,q}- T_n \rangle$.  

It follows that (4) also holds. By Lemma~\ref{interchange}(ii), colors $\alpha$ and $\beta$ are $T_n$-interchangeable under $\varphi_n$. So $P_{a}(\alpha,\beta,\varphi_n)$ is the only $(\alpha,\beta)$-path with respect to $\varphi_n$ intersecting 
$T_n$. Hence $Q$ is vertex-disjoint from $T_n$ and therefore contains no edge incident to  $T_n$. By Lemma \ref{LEM:Stable}, $\sigma_n=\varphi_n/Q$ is a $(T_n, D_n,\varphi_n)$-stable coloring, and $T_n$ is an ETT satisfying MP with respect to 
$\sigma_n$.  Since $b \notin V(T_{n,q})$, no edge in $T_{n,q}- T_{n,0}^*$ is colored by $\beta$ by TAA, because
$T_{n,0}^*=T_n \vee R_n$ by (7). Using (4) and (7), it is routine to check that 
the desired statements hold.

{\bf Case 2.}	$q=0$, $\alpha\in\phibar_n(T_{n})$, or $\alpha\in\phibar_n(T_{n,0}^*)$ with $\alpha,\beta\notin D_n$, and 
$Q$ is an $(\alpha,\beta)$-path disjoint from $P_{v_{\alpha}}(\alpha,\beta,\varphi_n)$.

Let us first assume that $\alpha$ or $\beta$ is closed in $T_{n,0}^*$ with respect to $\varphi_n$. By Lemma~\ref{interchange}(iii) 
or Theorem \ref{thm:tech10}(ii) (see (5.1)), colors $\alpha$ and $\beta$ are $T_{n,0}^*$-interchangeable under $\varphi_n$. 
So $P_{a}(\alpha,\beta,\varphi_n)$ is the only $(\alpha,\beta)$-path intersecting $T_{n,0}^*$, and hence $Q$ is vertex-disjoint from
$T_{n,0}^*$. It is then routine to check that $\sigma_n=\varphi_n/Q$ is a $(T_{n,0}^*, D_n,\varphi_n)$-weakly stable 
coloring, and $T_{n,0}^*$ is an ETT satisfying MP with respect to $\sigma_n$ by Theorem \ref{thm:tech10}(vi). So we assume 
hereafter that 

(8) neither $\alpha$ nor $\beta$ is closed in $T_{n,0}^*$ with respect to $\varphi_n$. 

By the hypothesis of the present case, $\alpha\in \phibar_n(T_{n})$ or $\{\alpha, \beta\} \cap D_{n}=\emptyset$.
So $\alpha\notin D_n$. From (8) and Algorithm 3.1 we deduce that $T_{n,0}^* \ne T_n$. Hence

(9) $\Theta_n=PE$, which together with (5.4) yields $a, b \notin V(T_n)\cap V(R_n)$.

Let us show that

(10) $\sigma_n=\varphi_n/Q$ is a $(T_{n,0}^*, D_n,\varphi_n)$-weakly stable coloring. 

To justify this, note that if one of $a$ and $b$ is contained in $V(T_n)-V(R_n)$ and the other is contained in
$V(R_n)-V(T_n)$, then $\alpha$ and $\beta$ are $T_{n,0}^*$-interchangeable under $\varphi_n$ by 
Lemma~\ref{interchange}(iv). So $Q$ is vertex-disjoint from $T_{n,0}^*$ and hence (10) holds. 
In view of (9), we may assume that 

(11) if $a, b \in V(T_{n,0}^*)$, then either $a, b \in V(T_n)-V(R_n)$ or $a, b \in V(R_n)-V(T_n)$.    

Let us first assume that $a\in V(T_n)-V(R_n)$. Then $\alpha \notin \phibar_n (R_n)$ by (6.6) and 
$b \in V(T_n)-V(R_n)$ if $b \in V(T_{n,0}^*)$ by (11). So $\alpha$ and $\beta$ are $T_n$-interchangeable under 
$\varphi_n$ by Lemma~\ref{interchange}(ii) and $\beta \notin \phibar_n (R_n)$ by (6.6). It follows that $Q$ is 
vertex-disjoint from $T_{n}$ and that $\sigma_n(f)=\varphi_n(f)$ for any edge $f$ incident to $R_n$ with $\varphi_n(f)\in\phibar_n(R_n)$. Hence (10) holds.   

Next we assume that $a\in V(R_n)-V(T_n)$. Then $\alpha \notin \phibar_n (T_n)$ by (6.6) and 
$b \in V(R_n)-V(T_n)$ if $b \in V(T_{n,0}^*)$ by (11). So $\alpha$ and $\beta$ are $R_n$-interchangeable under 
$\varphi_n$ by Lemma~\ref{interchange}(i) and $\beta \notin \phibar_n (T_n)$ by (6.6). It follows that $Q$ is 
vertex-disjoint from $R_{n}$. By the hypothesis of the present case, $\{\alpha, \beta\} \cap D_{n} =\emptyset$. 
So $\alpha, \beta \notin \overline{\varphi}_n(T_n) \cup D_{n}$ and hence (10) holds.   

From (10) we deduce that $T_{n,0}^*$ is an ETT satisfying MP with respect to $\sigma_n$.
	  
{\bf Case 3}. $T_{n,q}^*\prec v_{\alpha} \prec v_{\beta}$, $\alpha,\beta\notin D_{n,q}$, $\alpha\notin \phiv_n \langle T(v_{\beta})-T(v_{\alpha}) \rangle$, and $Q$ is an arbitrary $(\alpha,\beta)$-chain.

By (6.6), $V(T(y_{p-1}))$ is elementary with respect to $\varphi_n$. So $\alpha,\beta\notin\phibar_n(T_{n,q}^*)$.
By hypothesis, $\alpha,\beta\notin D_{n,q}$. Hence 

(12) $\alpha,\beta\notin\phibar_n(T_{n,q}^*) \cup D_{n,q}$.  

\noindent By the definitions of $D_n$ and $D_{n,q}$, we have $ D_n \cup \phibar_n(T_n) \subseteq D_{n,q}\cup 
\phibar_n(T_{n,q}^*)$. So $\alpha,\beta \notin\phibar_n(T_n)\cup D_n$. From Lemma \ref{hku}(iv), TAA and the 
hypothesis of the present case, we further deduce that

(13) $\alpha, \beta \notin \varphi_n\langle T(b)\rangle$. 

\noindent In view of Lemma~\ref{change}, we obtain

(14) $P_a(\alpha,\beta,\varphi)=P_b(\alpha,\beta,\varphi)$. (Possibly $Q$ is this path.)  

Since $T_{n,q}^*\prec a \prec b$, using (12)-(14), it is straightforward to verify that $\sigma_n=\varphi_n/Q$ is a 
$(T_{n,q}^*, D_n,\varphi_n)$-stable coloring, so $\sigma_n$ is also $(T_{n,q}^*, D_n,\varphi_n)$-weakly stable. 

From (12) and (13) we also see that $T(b)$ can be obtained from $T_{n,q}^*$ by using TAA, no matter whether $Q=
P_a(\alpha,\beta,\varphi)$. Thus $T$ is an ETT corresponding to $(\sigma_n,T_n)$. As neither $\alpha$ nor $\beta$ is contained in any $\Gamma$-set,
it is clear that $T$ also satisfies MP under $\sigma_n$, and $T_n=T_{n,0} \subset T_{n,1} \subset \ldots \subset T_{n,q}\subset 
T_{n,q+1}=T$ remains to be a good hierarchy of $T$ under $\sigma_n$, with the same $\Gamma$-sets as $T$ under 
$\varphi_n$. \qed

\section{Elementariness and Interchangeability} 

In Section 5 we have developed a control mechanism over Kempe changes; that is, a good hierarchies of an ETT.  In Section 6 
we have derived some properties satisfied by such hierarchies. Now we are ready to present a proof of Theorem \ref{hierarchy} 
by using Kempe changes based on these hierarchies, whose origin can be traced back to Tashkinov's proof of  
Theorem \ref{TashTree} \cite{T} (see Stiebitz et al. \cite{SSTF} for an English version). 

\subsection{Proof of Theorem \ref{hierarchy}}

By hypothesis, $T$ is an ETT constructed from a $k$-triple $(G,e, \varphi)$ by using the Tashkinov series 
$\TT=\{(T_i, \phiv_{i-1}, S_{i-1}, F_{i-1}, \Theta_{i-1}): 1\le i \le n+1\}$. Furthermore, $T$ admits a good hierarchy 
$T_n=T_{n,0} \subset T_{n,1} \subset \ldots \subset T_{n,q+1}=T$ and satisfies MP with respect to $\varphi_n$. Our 
objective is to show that $V(T)$ is elementary with respect to $\varphi_n$. 

As introduced in the preceding section,  $T=T_{n,q}^* \cup\{e_1,y_1,e_2,...,e_p,y_p\}$, where $y_i$ is the end 
of $e_i$ outside $T(y_{i-1})$ for $i\ge 1$, with $T(y_0)=T_{n,q}^*$. Suppose on the contrary that $V(T)$ is not 
elementary with respect to $\varphi_n$.  Then

{\bf (7.1)} $\overline{\varphi}_n(T(y_{p-1}))\cap\overline{\varphi}_n(y_p) \ne \emptyset$ by (6.6). 

For ease of reference, recall that (see (3) in the proof of Theorem \ref{good})

{\bf (7.2)} $|\phibar_n(T_{n})|\ge 2n+11$ and $|D_{n,j}|\le |D_n|\le n$ for $0 \le j \le q$. 

In our proof, by $A\cap B=\emptyset$ we mean $A$ and $B$ are vertex-disjoint, provided that $A$ is a path and $B$ is a tree. 
We shall frequently make use of a coloring $\sigma_n \in {\cal C}^k(G-e)$ with properties (i)-(iii) as described 
in Lemma \ref{stablechange}; that is,

{\bf (7.3)} $\sigma_n$ is a $(T_{n,q}^*, D_n,\varphi_n)$-weakly stable coloring, and $T_{n,q}^*$ is an ETT 
satisfying MP with respect to $\sigma_n$. Furthermore, if $q \ge 1$, then $T_{n,q}$ admits a good hierarchy 
$T_n=T_{n,0} \subset T_{n,1} \subset \ldots \subset T_{n,q}$ under $\sigma_n$, with the same $\Gamma$-sets 
(see Definition \ref{R2}) as $T$ under $\varphi_n$, and $T_{n,q}$ is $(\cup_{\eta_h\in D_{n,q}}\Gamma^{q-1}_h)^-$-closed 
with respect to $\sigma_n$ (see the remark succeeding Lemma \ref{stablechange}).

\begin{claim}\label{cla-p>0}
$p \ge  2.$
\end{claim}

\vspace{-1mm}
Assume the contrary: $p=1$; that is, $T=T_{n,q}^* \cup\{e_1,y_1\}$. Then 

(1) there exists a color $\alpha$ in $\overline{\varphi}_n(T_{n,q}^*)\cap\overline{\varphi}_n(y_1)$ by (7.1). 

\noindent We consider two cases according to the value of $q$. 

{\bf Case 1.} $q=0$. In this case, from (1) and Algorithm 3.1 we see that $\Theta_n\ne SE$. Let us first assume
that $\Theta_n=RE$. Let $\delta_n, \gamma_n$ be as specified in RE of Algorithm 3.1.
Since $\alpha, \delta_n \in \overline{\varphi}_n(T_n)$, both of them are closed in $T_n$ with
respect to $\varphi_n$. Hence $P_{y_1}(\alpha, \delta_n, \varphi_n)$ is vertex-disjoint from $T_n$.  Let
$\sigma_n=\varphi_n/P_{y_1}(\alpha, \delta_n, \varphi_n)$. Then $\delta_n \in \overline{\sigma}_n(T_n) \cap 
\overline{\sigma}_n(y_1)$. By Lemma \ref{LEM:Stable}, $\sigma_n$ is a $(T_n, D_n, \varphi_n)$-stable coloring
and hence, by Theorem \ref{thm:tech10}(vi), it is a $\phiv_n\bmod T_n$ coloring. In view of Definition \ref{wz2}, 
$f_n=e_1$ is still an RE connecting edge under $\sigma_n$. From Algorithm 3.1 we see that
$q\geq 1$ and $e_1$ is contained in a $(\delta_n,\gamma_n)$-cycle under $\sigma_n$, which is 
impossible because $\delta_n \in \overline{\sigma}_n(y_1)$. 

So we may assume that $\Theta_n=PE$. Let $\beta=\varphi_n(e_1)$. From TAA we see that $\beta\in\phibar_n(T_{n,0}^*)$. 
Let $\theta\in \phibar_n(T_n) \cap \phibar_n(R_n)$. Then $\theta$ is closed in $T_{n,0}^*$ under $\varphi_n$ by (5.4).
By Lemma~\ref{interchange}(iii), $P_{v_{\theta}}(\alpha,\theta,\varphi_n)$ 
is the only $(\alpha,\theta)$-path intersecting $T_{n,0}^*$. Thus $P_{y_1}(\alpha,\theta,\varphi_n)\cap T_{n,0}^*=
\emptyset$. Let $\sigma_n=\varphi_n/P_{y_1}(\alpha,\theta,\varphi_n)$. Then
$\theta$ is also closed in $T_{n,0}^*$ with respect to $\sigma_n$, and $\sigma_n$ is a $(T_{n,0}^*, D_n, \varphi_n)$-weakly stable coloring by Lemma~\ref{LEM:Stable}. In view of
Lemma~\ref{interchange}(iii), $\beta$ and $\theta$ are $T_{n,0}^*$-interchangeable under $\sigma_n$.  
As $P_{y_1}(\theta,\beta,\sigma_n)\cap T_{n,0}^*\neq\emptyset$ and $\theta,\beta\in\overline{\sigma}_n(T_{n,0}^*)$, there are at least two $(\theta,\beta)$-paths with
respect to $\sigma_n$ intersecting $T_{n,0}^*$, a contradiction.   

{\bf Case 2.} $q\ge 1$. In this case, by Definition \ref{R2}(v), we have	

(2) $T_{n,q}$ is $(\cup_{\eta_h\in D_{n,q}}\Gamma^{q-1}_h)^-$-closed with respect to $\varphi_n$ 

\noindent So $e_1$ is colored by some color $\gamma$ in $\cup_{\eta_h\in D_{n,q}}\Gamma^{q-1}_h$. By Definition \ref{R2}(i) 
and (5.9), we have $\gamma \notin\Gamma^{q}$. Let $\theta\in\phibar_n(T_{n,q})-\phibar_n(T_{n,q-1}^*)$.
Then $\theta \notin \Gamma^{q-1}$ (so $\theta \neq\gamma$) by Definition \ref{R2}(i). Furthermore, 
$\theta$ is closed in $T_{n,q}$ under $\varphi_n$ by (2). In view of Lemma \ref{step1}, $\alpha$ and $\theta$ are $T_{n,q}$-interchangeable under $\varphi_n$. So  $P_{v_{\theta}}(\alpha,\theta,\varphi_n)=P_{v_{\alpha}}(\alpha,\theta,\varphi_n)$ 
is the unique $(\alpha,\theta)$-path intersecting $T_{n,q}$. Hence $P_{y_1}(\alpha,\theta,\varphi_n)\cap T_{n,q}=\emptyset$. 
Let $\sigma_n=\varphi_n/P_{y_1}(\alpha,\theta,\varphi_n)$. Then $\sigma_n$ satisfies all the properties described in (7.3) by Lemma~\ref{stablechange}. Since $e_1$ is still colored by $\gamma\in\Gamma^{q-1}$ under $\sigma_n$ and 
$\gamma\notin \Gamma^{q}$, we can obtain $T$ from $T_{n,q}$ by TAA under $\sigma_n$, so $T$ is an ETT satisfying MP under $\sigma_n$.  Moreover, $T_n=T_{n,0} \subset T_{n,1} \subset \ldots \subset T_{n,q+1}=T$ remains to be a good hierarchy of $T$ 
under $\sigma_n$, with the same $\Gamma$-sets as those under $\varphi_n$.  Hence $(T, \sigma_n)$ is also a minimum counterexample 
to Theorem \ref{hierarchy} (see (6.2)-(6.5)). As $P_{y_1}(\theta,\gamma,\sigma_n)\cap T_{n,q}\neq \emptyset$ and $\theta,\gamma\in\overline{\sigma}_n(T_{n,q})$, there are 
at least two $(\theta,\gamma)$-paths with respect to $\sigma_n$ intersecting $T_{n,q}$, contradicting Lemma~\ref{change}(iii) 
(with $\sigma_n$ in place of $\varphi_n$), because $\theta, \gamma\in \overline{\sigma}_n(T_{n,q})$ and $\theta$ is also 
closed in $T_{n,q}$ under $\sigma_n$ by (2). Hence Claim \ref{cla-p>0} is justified. 

\vskip 3mm
Recall that the path number $p(T)$ of $T$ is the smallest subscript $i\in \{1,2,...,p\}$, such that the sequence $(y_i,e_{i+1},...,e_p,y_p)$ corresponds to a path in $G$, where $p\ge 2$ by Claim \ref{cla-p>0}. Depending on 
the value of $p(T)$, we distinguish among three situations, labeled as Situation 7.1, Situation 7.2, and Situation 7.3. \\

\noindent {\bf Situation 7.1.} $p(T)=1$. Now $T - V(T_{n,q}^*)$ is a path obtained by using TAA under $\varphi_n$.

\begin{claim}\label{claim9}
We may assume that $\overline{\varphi}_n(y_i)\cap\overline{\varphi}_n(y_p) \ne \emptyset$ for some $i$ with $1\le i \le p-1$.
\end{claim} 	
\vspace{-1mm}

To justify this, let $\alpha \in \overline{\varphi}_n(T(y_{p-1})) \cap\overline{\varphi}_n(y_p)$ (see 
(7.1)). If $\alpha\in \overline{\varphi}_n(y_i)\cap\overline{\varphi}_n(y_p)$ for some $i$ with $1\le i \le p-1$,
we are done. So we assume that  

(1) $\alpha\in\overline{\varphi}_n(T_{n,q}^*) \cap\overline{\varphi}_n(y_p)$ and $\alpha \notin \overline{\varphi}_n(y_i)$
for all $1 \le i \le p-1$.

(2) If $\Theta_n=PE$ and $q=0$, then we may further assume that $\alpha\in \phibar_n(T_n)$. 

Let us justify (2).
By (1), we have $\alpha\in \phibar_n(T_{n,0}^*)$. Suppose $\alpha\in \phibar_n(R_n-V(T_n))$. Then $\alpha\notin\Gamma^{0}$ by Definition \ref{R2}(i). In view of (7.2), we have $|\phibar_n(T_n)|\geq 11+2n$ and $|\Gamma^{0}|\leq 2 |D_{n,0}| \le 2n$. 
So there exists $\beta\in\phibar_n(T_n)-\Gamma^{0}$. By Lemma \ref{interchange}(iv), $\alpha$ and $\beta$ are $T_{n,0}^*$-interchangeable under $\varphi_n$. Thus $P_{v_{\alpha}}(\alpha,\beta,\varphi_n)=P_{v_{\beta}}(\alpha,\beta,\varphi_n)$ 
and $P_{y_p}(\alpha,\beta,\varphi_n)$ is disjoint from $T_{n,0}^*$. Let $\sigma_n=\varphi_n/P_{y_p}(\alpha,\beta,\varphi_n)$.
By Lemma~\ref{stablechange} (the second case), $\sigma_n$ is a $(T_{n,0}^*, D_n,\varphi_n)$-weakly stable coloring, 
and $T_{n,0}^*$ is an ETT satisfying MP with respect to $\sigma_n$. Note that $T$ can also be obtained from $T_{n,0}^*$ by 
TAA under $\sigma_n$, because $\alpha,\beta \in \overline{\sigma}_n(T_{n,0}^*)$. Hence $T$ is an ETT satisfying MP under $\sigma_n$ as 
well. Since $\alpha,\beta\notin \Gamma^{0}$ and $\alpha,\beta\notin \phibar_n(T(y_{p-1})-V(T_{n,0}^*))$, the hierarchy
$T_n=T_{n,0} \subset T_{n,1}=T$ remains to be good under $\sigma_n$, with the same $\Gamma$-sets
as those under $\varphi_n$. Therefore $(T, \sigma_n)$ is also a minimum counterexample to Theorem \ref{hierarchy} (see 
(6.2)-(6.5)). As $\beta\in \overline{\sigma}_n(T_n) \cap \overline{\sigma}_n(y_p)$, replacing $\varphi_n$ by $\sigma_n$ 
and $\alpha$ by $\beta$ if necessary, we see that (2) holds.

Depending on whether $\alpha$ is used by edges in $T -  T_{n,q}^*$, we consider two cases.

{\bf Case 1.} $\alpha\notin \varphi_n \langle T -  T_{n,q}^* \rangle$. In this case, let $\beta\in\overline{\varphi}_n(y_{p-1})$.
Then $\beta$ is not used by any edge in $T -  T_{n,q}^*$, except possibly $e_1$ when $q=0$ and $T_{n,0}^*=T_n$ (now $e_1=f_n$
in Algorithm 3.1 and $\varphi_n(e_1)=\beta \in D_n$). By (1) and (2), we have $\alpha\in\phibar_n(T_{n,q})$ if $q \ge 1$ and $\alpha\in\phibar_n(T_n)$ if $q=0$. It follows from Lemma~\ref{change} that $P_{v_{\alpha}}(\alpha,\beta,\varphi_n)=P_{y_{p-1}}(\alpha,\beta,\varphi_n)$. So $P_{y_{p}}(\alpha,\beta,\varphi_n)$ is disjoint from $P_{v_{\alpha}}(\alpha,\beta,\varphi_n)$.  Let $\sigma_n = \varphi_n/P_{y_{p}}(\alpha,\beta,\varphi_n)$. By Lemma~\ref{stablechange},  $\sigma_n$ satisfies all the properties described in (7.3). In particular, if $e_1=f_n$ and $\varphi_n(e_1)=\beta \in D_n$,
then $\sigma_n(e_1)=\varphi_n(e_1)$, which implies that $e_1$ is outside $P_{y_{p}}(\alpha,\beta,\varphi_n)$.
So $\sigma_n(f)=\varphi_n(f)$ for each $f\in E(T)$ and $\overline{\sigma}_n(u)=\phibar_n(u)$ for each $u\in V(T(y_{p-1}))$. 
Thus $T$ can be obtained from $T_{n,q}^*+e_1$ by TAA and is an ETT satisfying MP under $\sigma_n$. Furthermore, $T_n=T_{n,0} \subset T_{n,1} 
\subset \ldots \subset T_{n,q} \subset T_{n,q+1}=T$ remains to be a good hierarchy of $T$ under $\sigma_n$, with the same $\Gamma$-sets as those under $\varphi_n$. Therefore, $(T, \sigma_n)$ is also a minimum counterexample to Theorem \ref{hierarchy} 
(see (6.2)-(6.5)). As $\beta\in \overline{\sigma}_n(y_{p-1}) \cap \overline{\sigma}_n(y_p)$, replacing $\varphi_n$ by 
$\sigma_n$ if necessary, we see that Claim~\ref{claim9} is true.

{\bf Case 2.} $\alpha\in \varphi_n \langle T -  T_{n,q}^* \rangle$. In this case, let $e_j$ be the edge with the smallest
subscript in $T -T_{n,q}^*$ such that $\varphi(e_j)=\alpha$. We distinguish between two subcases according to the value of 
$j$.

{\bf Subcase 2.1.} $j\geq 2$. In this subcase, let $\beta \in \overline{\varphi}_n(y_{j-1})$. Then $\beta$ is not used 
by any edge in $T(y_j) -  T_{n,q}^*$, except possibly $e_1$ when $q=0$ and $T_{n,0}^*=T_n$ (now $e_1=f_n$ in Algorithm 3.1 
and $\varphi_n(e_1)=\beta \in D_n$). By (1) and (2), we have $\alpha\in\phibar_n(T_{n,q})$ if $q \ge 1$ and $\alpha\in\phibar_n(T_n)$ 
if $q=0$. It follows from Lemma \ref{change} that  $P_{v_{\alpha}}(\alpha,\beta,\varphi_n)=P_{y_{j-1}}(\alpha,\beta,\varphi_n)$. 
So $P_{y_{p}}(\alpha,\beta,\varphi)$ is disjoint from $P_{v_{\alpha}}(\alpha,\beta,\varphi_n)$. Let $\sigma_n = \varphi_n/P_{y_{p}}(\alpha,\beta,\varphi_n)$. By Lemma~\ref{stablechange}, $\sigma_n$ satisfies all the properties 
described in (7.3). In particular, if $e_1=f_n$ and $\varphi_n(e_1)=\beta \in D_n$, then $\sigma_n(e_1)=\varphi_n(e_1)$, 
which implies that $e_1$ is outside $P_{y_{p}}(\alpha,\beta,\varphi_n)$. So $T$ can be obtained from $T_{n,q}^*+e_1$ by TAA 
under $\sigma_n$ and hence is an ETT satisfying MP under $\sigma_n$. 

Note that $\beta\notin\Gamma^{q}$ by Definition \ref{R2}(i) and that $\overline{\sigma}_n(u)=\phibar_n(u)$ for each $u\in V(T(y_{p-1}))$ by (6.6).  If $\alpha\notin \Gamma^{q}$, then clearly $T_n=T_{n,0} \subset T_{n,1} \subset \ldots \subset T_{n,q} 
\subset T_{n,q+1}=T$ is a good hierarchy of $T$ under $\sigma_n$, with the same $\Gamma$-sets as those under $\varphi_n$. 
If $\alpha\in \Gamma^{q}$, say $\alpha\in \Gamma^{q}_{h}$ for some $\eta_h\in D_{n,q}$, then Definition \ref{R2}(i)
implies that $\eta_h \in\overline{\varphi}_n(w)$ for some $w\preceq  y_{j-1}$. Since only edges outside $T(w)$ may 
change colors between $\alpha$ and $\beta$ as we transform $\varphi_n$ into $\sigma_n$, it follows that $T_n=T_{n,0} 
\subset T_{n,1} \subset \ldots \subset T_{n,q} \subset T_{n,q+1}=T$ remains to be a good hierarchy of $T$ under $\sigma_n$, 
with the same $\Gamma$-sets as those under $\varphi_n$. Hence $(T, \sigma_n)$ is also a minimum counterexample to Theorem 
\ref{hierarchy} (see (6.2)-(6.5)). Since $\beta\in \overline{\sigma}_n(y_{j-1}) \cap \overline{\sigma}_n(y_p)$, replacing 
$\varphi_n$ by $\sigma_n$ if necessary, we see that Claim~\ref{claim9} holds.

{\bf Subcase 2.2.} $j=1$. In this subcase, $\alpha=\varphi(e_1)$. Note that $\alpha\notin\Gamma^{q}$ by Definition \ref{R2}(i)
and (5.9). We propose to show that 

(3) there exists a color $\gamma$ in $\phibar_n(T_{n,q})-\Gamma^{q}$ if $q\ge 1$ and in $\phibar_n(T_{n})-\Gamma^{0}$ 
if $q=0$, such that $\gamma$ is closed in $T_{n,q}^*$ with respect to $\varphi_n$. 

Let us first assume that $q\geq 1$. By (7.2), we obtain $|\phibar_n(T_{n,q})|\ge |\phibar_n(T_{n})|\ge 2n+11$ and 
$|\Gamma^{q-1}|\le 2 |D_{n,q-1}|\le 2n$. So $|\phibar_n(T_{n,q})-\Gamma^{q-1}|\ge 11$. By Definition \ref{R2}(iii), we 
have $|\Gamma^{q}-\Gamma^{q-1}|= 2$. So $|\phibar_n(T_{n,q})-(\Gamma^{q-1} \cup \Gamma^{q}) |\ge 9$. Let $\gamma$ be a 
color in $\phibar_n(T_{n,q})-(\Gamma^{q-1} \cup \Gamma^{q})$. By Definition \ref{R2}(v), $\gamma$ is closed in $T_{n,q}$ 
with respect to $\varphi_n$.  

Next we assume that $q=0$. Again, by (7.2), we have $|\phibar_n(T_{n})|\ge 2n+11$ and $|\Gamma^{0}|\le 2 |D_{n,0}|\le 
2 |D_n| \le 2n$. Let $\gamma$ be a color in $\phibar_n(T_{n})-\Gamma^{0}$ if $\Theta_n\ne PE$ and a color in 
$\phibar_n(T_n) \cap \phibar_n(R_n)-\Gamma^0$ if $\Theta_n=PE$ (see Definition \ref{R2}(iv)). By Algorithm 3.1 
and (5.4),  $\gamma$ is closed in $T_{n,0}^*$ with respect to $\varphi_n$. So (3) holds.

By (3) and Lemma \ref{change}, $P_{v_{\alpha}}(\alpha,\gamma,\varphi_n)=P_{v_{\gamma}}(\alpha,\gamma,
\varphi_n)$ is the only  $(\alpha,\gamma)$-path intersecting $T_{n,q}^*$. So $P_{y_{p}}(\alpha, \gamma,\varphi_n)$ is 
disjoint from $T_{n,q}^*$ and hence it does not contain $e_1$. Let $\sigma_n = \varphi_n/P_{y_{p}}(\alpha, \gamma,\varphi_n)$. 
Then $\sigma_n$ satisfies all the properties described in (7.3) by Lemma~\ref{stablechange}. Moreover, $\overline{\sigma}_n(u)=\phibar_n(u)$ 
for all $u\in V(T(y_{p-1}))$ by (6.6). Since $\alpha, \gamma\in \phibar_n(T_{n,q}^*)$, we have $\alpha, \gamma\in \overline{\sigma}_n(T_{n,q}^*)$. 
Hence we can obtain $T$ from $T_{n,q}^*+e_1$ by using TAA under $\sigma_n$, so $T$ is an ETT satisfying MP under $\sigma_n$. Since 
$\alpha, \gamma \notin\Gamma^{q}$, the hierarchy $T_n=T_{n,0} \subset T_{n,1} \subset \ldots \subset T_{n,q} 
\subset T_{n,q+1}=T$ remains to be good under $\sigma_n$,  with the same $\Gamma$-sets
as those under $\varphi_n$. Therefore, $(T, \sigma_n)$ is also a minimum counterexample to Theorem \ref{hierarchy} 
(see (6.2)-(6.5)). Since $e_1$ is outside $P_{y_p}(\alpha,\gamma,\varphi_n)$, we have $\sigma_n(e_1)=\alpha$. As 
$\gamma\in\overline{\sigma}_n(y_p)\cap \overline{\sigma}_n(v)$ for some $v\in V(T_{n,q})$ and $\alpha\neq\gamma$, 
the present subcase reduces to Case 1 if $\gamma\notin \sigma_n \langle T-T_{n,q}^* \rangle$ or to Subcase 2.1 if 
$\gamma\in \sigma_n \langle T-T_{n,q}^* \rangle$. This proves Claim \ref{claim9}.

\vskip 3mm

\begin{claim}\label{p-1}
We may assume that $\overline{\varphi}_n(y_{p-1})\cap\overline{\varphi}_n(y_p) \ne \emptyset$.
\end{claim}

\vspace{-1mm}
To justify this, let ${\cal K}$ be the set of all minimum counterexamples $(T, \varphi_n)$ to Theorem \ref{hierarchy} 
(see (6.2)-(6.5)), and let $i$ be the largest subscript with $1 \le i \le p-1$, such that there exists a 
member $(T, \mu_n)$ of ${\cal K}$ with $\overline{\mu}_n(y_{i})\cap\overline{\mu}_n(y_p) \ne \emptyset$;
this $i$ exists by Claim \ref{claim9}. We aim to show that $i=p-1$. Thus Claim \ref{p-1} follows by replacing 
$\varphi_n$ with $\mu_n$, if necessary.

With a slight abuse of notation, we assume that $\overline{\varphi}_n(y_{i})\cap\overline{\varphi}_n(y_p) \ne 
\emptyset$ and assume, on the contrary, that $i\le p-2$. Let $\alpha \in \overline{\varphi}_n(y_{i})\cap\overline{\varphi}_n(y_p)$.
Using (6.6) and TAA, we obtain 

(1) $\alpha \notin \overline{\varphi}_n(T(y_{i-1}))$, where $T(y_0)=T_{n,q}^*$. So $\alpha$ is not used by any edge in 
$T(y_{i+1}) -  T_{n,q}^*$, except possibly $e_1$ when $q=0$ and $T_{n,0}^*=T_n$ (now $e_1=f_n$ in Algorithm 3.1 and 
$\varphi_n(e_1)=\alpha \in D_n$). 

Recall that Definition \ref{R2} involves $\Gamma^{q}_h=\{\gamma^{q}_{h_1}, \gamma^{q}_{h_2}\}$ for each $\eta_h\in 
D_{n,q}$. Nevertheless, the proof of this claim only involves one $\eta_h\in D_{n,q}$. For simplicity, we abbreviate 
its corresponding $\gamma^{q}_{h_j}$ to $\gamma_{j}$ for $j=1,2$. By Definition \ref{R2}(i) and (5.9), we have 

(2)  $\gamma_{j}\in \phibar_n (T_{n,q})$ if $q\geq 1$ and $\gamma_{j}\in \phibar_n(T_{n})$ if $q=0$.  
Moreover, if $\eta_h \in \overline{\varphi}_n(y_t)$ for some $t\ge 1$, then $\gamma_{j}\notin \varphi_n 
\langle T(y_t)-T_{n,q}^* \rangle$ for $j=1,2$. 

Depending on whether $\alpha \in D_{n,q}$, we consider two cases.

{\bf Case 1.} $\alpha\notin D_{n,q}$. In this case, let $\theta\in\overline{\varphi}_n(y_{i+1})$. From TAA and (6.6)
it follows that 

(3) $\theta \notin \overline{\varphi}_n(T(y_i))$, so $\theta$ is not used by any edge in $T(y_{i+1}) -  T_{n,q}^*$, 
except possibly $e_1$ when $q=0$ and $T_{n,0}^*=T_n$ (now $e_1=f_n$ in Algorithm 3.1 and $\varphi_n(e_1)=\theta \in D_n$).  

If $\theta\notin D_{n,q}$,  then $\{\alpha, \theta\}\cap D_{n,q} = \emptyset$. By the definitions of $D_n$ and
$D_{n,q}$, we have $\overline{\varphi}_n(T_n) \cup D_n \subseteq  \overline{\varphi}_n(T_{n,q}^*) \cup D_{n,q}$, 
which together with (1) and (3) implies $\{\alpha, \theta\}\cap D_{n} = \emptyset$.  Hence
$P_{y_{i}}(\alpha,\theta,\varphi_n)=P_{y_{i+1}}(\alpha,\theta,\varphi_n)$ by Lemma~\ref{change}. 
Let $\sigma_n=\varphi_n/P_{y_{p}}(\alpha,\theta,\varphi_n)$. Since both $y_i$ and $y_{i+1}$ are contained in $T-V(T_{n,q}^*)$
and (1) holds, by Lemma~\ref{stablechange} (the third case), $\sigma_n$ satisfies all the properties described in
(7.3).  Furthermore, $T$ is also an ETT satisfying MP with respect to $\sigma_n$, and $T_n=T_{n,0} \subset T_{n,1}
\subset \ldots \subset T_{n,q}\subset T_{n,q+1}=T$ remains to be a good hierarchy of $T$ under $\sigma_n$, with the 
same $\Gamma$-sets as those under $\varphi_n$. Hence $(T, \sigma_n)$ is also a minimum counterexample to Theorem 
\ref{hierarchy} (see (6.2)-(6.5)). Since $\theta\in\overline{\sigma}_n(y_p)\cap\overline{\sigma}_n(y_{i+1})$, 
we reach a contradiction to the maximality assumption on $i$.

So we may assume that $\theta\in D_{n,q}$. Let $\theta=\eta_h\in D_{n,q}$. In view of (2) and Lemma \ref{change}, we obtain $P_{v_{\gamma_1}}(\alpha,\gamma_{1},\varphi_n)=P_{y_{i}}(\alpha,\gamma_{1},\varphi_n)$, which is disjoint from $P_{y_{p}}(\alpha,\gamma_{1},\varphi_n)$.  Let $\sigma_n=\varphi_n/P_{y_{p}}(\alpha,\gamma_{1},\varphi_n)$. 
By Lemma~\ref{stablechange}, $\sigma_n$ satisfies all the properties described in (7.3). In particular,
if $e_1=f_n$ and $\varphi_n(e_1)=\alpha \in D_n$, then $\sigma_n(e_1)=\varphi_n(e_1)$, which implies that $e_1$ is outside 
$P_{y_{p}}(\alpha,\gamma_{1},\varphi_n)$. By (6.6), (1) and (2), we have $\overline{\sigma}_n(u)=\phibar_n(u)$ for each 
$u \in V(T(y_{p-1}))$ and $\sigma_n(f)=\varphi_n(f)$ for each edge $f$ in $T(y_{i+1})$. So $T$ can be obtained from 
$T_{n,q}^*+e_1$ by TAA under $\sigma_n$, and hence is an ETT satisfying MP under $\sigma_n$. Furthermore, $T_n=T_{n,0} \subset T_{n,1} 
\subset \ldots \subset T_{n,q}\subset T_{n,q+1}=T$ remains to be a good hierarchy of $T$ under $\sigma_n$, with the same 
$\Gamma$-sets as those under $\varphi_n$. Hence $(T, \sigma_n)$ is also a minimum counterexample to Theorem \ref{hierarchy} 
(see (6.2)-(6.5)), with $\gamma_{1}\in \overline{\sigma}_n(y_p) \cap \overline{\sigma}_n(T_{n,q})$.

Using (2) and Lemma \ref{change}, we obtain $P_{v_{\gamma_{1}}}(\eta_h, \gamma_{1},\sigma_n) =P_{y_{i+1}}(\eta_h,\gamma_{1}, \sigma_n)$, which is disjoint from $P_{y_{p}}(\eta_h, \gamma_{1},\sigma_n)$. Let $\sigma_n'=\sigma_n/P_{y_{p}}(\eta_h,\gamma_{1},\sigma_n)$.  By Lemma~\ref{stablechange}, $\sigma_n'$ satisfies all the 
properties described in (7.3) (with $\sigma_n'$ in place of $\sigma_n$). In particular, if $e_1=f_n$ and $\sigma_n(e_1)
=\eta_h \in D_n$, then $\sigma_n'(e_1)=\sigma_n(e_1)$, which implies that $e_1$ is outside $P_{y_{p}}(\eta_h,\gamma_{1},\sigma_n)$. 
By (6.6), (2) and (3), we have $\overline{\sigma}'_n(u)=\overline{\sigma}_n(u)$ for each $u \in V(T(y_{p-1}))$ and $\sigma_n'(f)=\sigma_n(f)$ for each edge $f$ in $T(y_{i+1})$.  So $T$ can be obtained from $T_{n,q}^*+e_1$ by TAA under 
$\sigma_n'$, and hence is an ETT satisfying MP under $\sigma_n'$. Furthermore, since $\eta_h \in \overline{\sigma}'_n(y_{i+1})$,
the hierarchy $T_n=T_{n,0} \subset T_{n,1} \subset \ldots \subset T_{n,q}\subset T_{n,q+1}=T$ remains 
to be good under $\sigma_n'$, with the same $\Gamma$-sets as those under $\varphi_n$. Therefore 
$(T, \sigma_n')$ is also a minimum counterexample to Theorem \ref{hierarchy} (see (6.2)-(6.5)). Since $\eta_h \in\overline{\sigma}_n'(y_p)\cap\overline{\sigma}_n'(y_{i+1})$, we reach a contradiction to the maximality assumption on $i$.

{\bf Case 2.} $\alpha\in D_{n,q}$. In this case, let $\alpha=\eta_h\in D_{n,q}$. Then $\Gamma^{q}_h=\{\gamma_{1},
\gamma_{2}\}$ (see the paragraph above (2)). Renaming subscript if necessary, we may assume that $\varphi_n(e_{i+1}) 
\neq \gamma_{1}$. By (1) and (2), we have

(4) $\gamma_{1}\notin \varphi_n \langle T(y_{i+1})-T_{n,q}^* \rangle$ and $\eta_h$ is not used by any edge in 
$T(y_{i+1}) -  T_{n,q}^*$, except possibly $e_1$ when $q=0$ and $T_{n,0}^*=T_n$ (now $e_1=f_n$ in Algorithm 3.1 
and $\varphi_n(e_1)=\eta_h \in D_{n,q} \subseteq D_n$).

By (4) and Lemma~\ref{change}, we obtain $P_{v_{\gamma_1}}(\eta_h,\gamma_{1},\varphi_n)=P_{y_i}(\eta_h,\gamma_{1}, \varphi_n)$,
which is disjoint from the path $P_{y_p}(\eta_h,\gamma_{1},\varphi_n)$. Let $\sigma_n=\varphi_n/P_{y_p}(\eta_h,\gamma_{1},
\varphi_n)$. By Lemma~\ref{stablechange}, $\sigma_n$ satisfies all the properties described in (7.3). 
In particular, if $e_1=f_n$ and $\varphi_n(e_1)=\eta_h \in D_n$, then $\sigma_n(e_1)=\varphi_n(e_1)$, which implies that 
$e_1$ is outside $P_{y_p}(\eta_h,\gamma_{1}, \varphi_n)$. By  (6.6) and 
(4), we have $\overline{\sigma}_n(u)=\phibar_n(u)$ for each $u \in V(T(y_{p-1}))$ and $\sigma_n(f)=\varphi_n(f)$ for
each edge $f$ in $T(y_{i+1})$. So $T$ can be obtained from $T_{n,q}^*+e_1$ by TAA under $\sigma_n$, and hence is an ETT 
satisfying MP under $\sigma_n$. Furthermore, $T_n=T_{n,0} \subset T_{n,1} \subset \ldots \subset T_{n,q}\subset T_{n,q+1}=T$ 
remains to be a good hierarchy of $T$ under $\sigma_n$, with the same $\Gamma$-sets as those under $\varphi_n$. Therefore,
$(T, \sigma_n)$ is also a minimum counterexample to Theorem \ref{hierarchy} (see (6.2)-(6.5)), with $\gamma_{1}\in \overline{\sigma}_n(y_p) \cap \overline{\sigma}_n(T_{n,q})$.  Let $\theta\in\overline{\sigma}_n(y_{i+1})$. From TAA
we see that 

(5) $\theta$ is not used by any edge in $T(y_{i+1}) -  T_{n,q}^*$ under $\sigma_n$, except possibly $e_1$ when $q=0$ and $T_{n,0}^*=T_n$ (now $e_1=f_n$ in Algorithm 3.1 and $\sigma_n(e_1)=\theta \in D_n$).

By (6.6), we have $\theta\neq\gamma_1$. Using (4) and Lemma~\ref{change}, we get $P_{v_{\gamma_1}}(\theta, \gamma_{1},\sigma_n)=P_{y_{i+1}}(\theta, \gamma_{1}, \sigma_n)$.  Let $\sigma_n'=\sigma_n/P_{y_p} (\theta,\gamma_{1},\sigma_n)$.
By Lemma~\ref{stablechange}, $\sigma_n'$ satisfies all the properties described in (7.3) (with $\sigma_n'$ in place of $\sigma_n$). 
In particular, if $e_1=f_n$ and $\sigma_n(e_1)=\theta \in D_n$, then $\sigma_n'(e_1)=\sigma_n(e_1)$, which implies that 
$e_1$ is outside $P_{y_p} (\theta,\gamma_{1},\sigma_n)$. 
From (6.6) and (4) we deduce that $\overline{\sigma}'_n(u)=\overline{\sigma}_n(u)$ for each $u \in V(T(y_{p-1}))$,  
and $\sigma_n'(f)=\sigma_n(f)$ for each edge $f$ in $T(y_{i+1})$. So $T$ can also be obtained from $T_{n,q}^*+e_1$ 
by TAA under $\sigma_n'$, and hence is an ETT satisfying MP under $\sigma_n'$. Furthermore, $T_n=T_{n,0} \subset T_{n,1} \subset 
\ldots \subset T_{n,q}\subset T_{n,q+1}=T$ remains to be a good hierarchy of those under $\sigma_n'$, with the same 
$\Gamma$-sets as those under $\varphi_n$. Therefore, $(T, \sigma_n')$ is also a minimum counterexample to Theorem 
\ref{hierarchy} (see (6.2)-(6.5)). Since $\theta \in\overline{\sigma}_n'(y_p)\cap\overline{\sigma}_n'(y_{i+1})$, 
we reach a contradiction to the maximality assumption on $i$. Hence Claim \ref{p-1} is established.

\vskip 2mm

By Claim \ref{cla-p>0},  $p \ge  2$.  By Claim \ref{p-1}, $\overline{\varphi}_n(y_{p-1})\cap\overline{\varphi}_n(y_p) 
\ne \emptyset$. Let $\alpha \in \overline{\varphi}_n(y_{p-1})\cap\overline{\varphi}_n(y_p)$ and $\beta=\varphi_n(e_{p})$. Let
$\sigma_n$ be obtained from $\varphi_n$ by recoloring $e_p$ with $\alpha$ and let $T'=T(y_{p-1})$. Then
$\beta \in \overline{\sigma}_n(y_{p-1})\cap\overline{\sigma}_n(T'(y_{p-2}))$ and  $T_n=T_{n,0} \subset T_{n,1} \subset 
\ldots \subset T_{n,q}\subset T'$ is a good hierarchy of $T'$ under $\sigma_n$. 
So $(T', \sigma_n)$ is a counterexample to Theorem \ref{hierarchy} (see (6.2)-(6.4)), which violates the 
minimality assumption (6.5) on $(T, \varphi_n)$.  This completes our discussion about Situation 7.1.\\

\noindent {\bf Situation 7.2.} $p(T)=p$. Now $e_p$ is not incident to $y_{p-1}$.  

\vskip 2mm

By (7.1), there exists a color $\alpha\in \overline{\varphi}_n(T(y_{p-1}))\cap\overline{\varphi}_n(y_p)$. 
We divide this situation into $3$ cases and further into $6$ subcases (see Figure 4), depending on whether 
$v_{\alpha}=y_{p-1}$ and $\alpha \in D_{n,q}$. Our proof of Subcase 1.1 is self-contained. Yet, in our discussion 
Subcase 1.2 may be redirected to Subcase 1.1 and Subcase 2.1, and Subcase 2.1 may be redirected to Subcase 1.1, etc. 
Figure 4 illustrates such redirections (note that no cycling occurs). 
\begin{figure}[htpb]
\vspace{-1mm}
\centerline{\includegraphics[width=11cm]{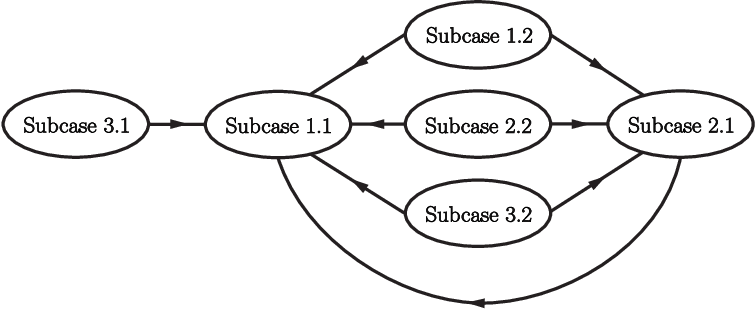}}
\vspace{-1mm}\caption{Redirections}
\end{figure}

Throughout this situation we reserve the symbol $\theta$ for $\varphi_n(e_p)$. Clearly, $\theta\neq \alpha$.

\vskip 2mm
{\bf Case 1.} $\alpha\in\overline{\varphi}_n(y_p)\cap\overline{\varphi}_n(y_{p-1})$ and  $\alpha \in D_{n,q}$.

\vskip 1mm
Let $\alpha=\eta_{m}\in D_{n,q}$. For simplicity, we abbreviate the two colors $\gamma^{q}_{m_1}$ and $\gamma^{q}_{m_2}$ in $\Gamma^{q}_m$ (see Definition \ref{R2}) to $\gamma_1$ and $\gamma_2$, respectively. 
Since $\eta_m \in \overline{\varphi}_n(y_p)\cap\overline{\varphi}_n(y_{p-1})$, from TAA and Definition \ref{R2}(i) we 
see that

(1) $\gamma_{1},\gamma_{2}\notin \varphi_n \langle T(y_{p-1}) -  T_{n,q}^* \rangle$ and $\eta_m$ is not used by 
any edge in $T -  T_{n,q}^*$, except possibly $e_1$ when $q=0$ and $T_{n,0}^*=T_n$ (now $e_1=f_n$ in Algorithm 3.1 
and $\varphi_n(e_1)=\eta_m \in  D_{n,q} \subseteq D_n$). 

By (1) and Lemma \ref{change} (with respect to $(T, \varphi_n)$), we have 

(2) $P_{v_{\gamma_j}}(\eta_m, \gamma_{j},\varphi_n)=P_{y_{p-1}}(\eta_m,\gamma_{j},\varphi_n)$ for $j=1,2$. 

\vskip 2mm

Let us consider two subcases according to whether $\theta\in\overline{\varphi}_n(y_{p-1})$.

{\bf Subcase 1.1.} $\theta\notin\overline{\varphi}_n(y_{p-1})$.

In our discussion about this subcase, we shall appeal to the following two tree-sequences:

\vskip 1mm
$\bullet$ $T^-=(T_{n,q}^*, e_1, y_1, e_2, \ldots, e_{p-2}, y_{p-2}, e_p, y_p)$ and

\vskip 1mm
$\bullet$ $T^* \hskip 0.6mm =(T_{n,q}^*, e_1, y_1, e_2, \ldots ,y_{p-2}, e_p, y_p, e_{p-1}, y_{p-1})$.

\vskip 1mm
\noindent Note that $T^-$ is obtained from $T$ by deleting $y_{p-1}$ and $T^*$ arises from $T$ by interchanging
the order of $(e_{p-1}, y_{p-1})$ and $(e_p, y_p)$. We propose to show that both $T^-$ and $T^*$ are ETTs corresponding
to $\varphi_n$. Indeed, if $T(y_{p-2}) \ne T_n$, then both $T^-$ and $T^*$ can be obtained from  $T(y_{p-2})$ 
by using TAA under $\varphi_n$.  So we assume that $T(y_{p-2})= T_n$.  By the hypothesis of the present subcase, $\varphi_n(e_p)=\theta \notin\overline{\varphi}_n(y_{p-1})$. From Algorithm 3.1 we deduce that now $\Theta_n=PE$. 
Hence both $T^-$ and $T^*$ can be obtained from  $T(y_{p-2})$ by using TAA under $\varphi_n$ as well.  Therefore
both $T^-$ and $T^*$ are ETTs corresponding to $\varphi_n$. In view of the maximum property enjoyed by $T$, we
further conclude that both $T^-$ and $T^*$ are ETTs satisfying MP with respect to $\varphi_n$. 

Let us first assume that $\theta\notin\Gamma^{q}$. Now it is easy to see that $T_n=T_{n,0} \subset T_{n,1} \subset 
\ldots \subset T_{n,q} \subset T^-$ is a good hierarchy of $T^-$ under $\varphi_n$, with the same $\Gamma$-sets 
(see Definition \ref{R2}) as $T$. (If $\theta\in\Gamma^{q}$, say $\theta \in \Gamma^{q}_h$, and 
$\eta_h \in \overline{\varphi}_n(y_{p-1})$, then $T^-$ no longer satisfies Definition \ref{R2}(i).)  Observe that $\gamma_{1} \notin 
\overline{\varphi}_n(y_{p})$, for otherwise, $\gamma_{1}$ is missing at two vertices in $T^-$. Thus $(T^-, \varphi_n)$ is 
a counterexample to Theorem \ref{hierarchy} (see (6.2) and (6.3)), which violates the minimality assumption
$(6.4)$ or $(6.5)$ on $(T, \varphi_n)$. Let us turn to considering $T^*$. Since $\theta \notin\overline{\varphi}_n(y_{p-1})$ and $\theta\notin\Gamma^{q}$, it is clear that $T_n=T_{n,0} \subset T_{n,1} \subset \ldots \subset T_{n,q} \subset 
T^*$ is a good hierarchy of $T^*$ under $\varphi_n$, with the same $\Gamma$-sets as $T$. 
Moreover, by (1), we have $\gamma_{1} \notin \varphi_n \langle T^*(y_{p}) -  T_{n,q}^* \rangle$.
It follows from Lemma \ref{change} (with respect to $(T^*, \varphi_n)$) that
$P_{v_{\gamma_1}}(\eta_m, \gamma_{1},\varphi_n)=P_{y_{p}}(\eta_m,\gamma_{1},\varphi_n)$, contradicting (2).

Next we assume that $\theta\in\Gamma^{q}$. Then $\theta\in \Gamma^{q}_h$ for some $\eta_h \in D_{n,q}$. 
If $\eta_h\notin \overline{\varphi}_n(y_{p-1})$, then $\eta_h \in\overline{\varphi}_n(T(y_{p-2}))$ by
Definition \ref{R2}(i). So we can still ensure that both $T^-$ and $T^*$ have good 
hierarchies under $\varphi_n$. Thus, using the same argument as employed in the preceding paragraph, we can  
reach a contradiction. Hence we may assume that $\eta_h\in \overline{\varphi}_n(y_{p-1})$.

Clearly, $\theta \ne \gamma_{1}$ or $\gamma_{2}$. Renaming subscripts if necessary, we may assume that 

(3) $\theta \ne \gamma_2$.

\noindent Since $P_{v_{\gamma_2}}(\eta_m, \gamma_{2},\varphi_n)=P_{y_{p-1}}(\eta_m,\gamma_{2},\varphi_n)$ by (2), this 
path is disjoint from $P_{y_{p}}(\eta_m,\gamma_{2},\varphi_n)$. Let $\mu_1 = \varphi_n/P_{y_{p}}(\eta_m,\gamma_2,\varphi_n)$.
By Lemma~\ref{stablechange},  $\mu_1$ satisfies all the properties described in (7.3) (with $\mu_1$ in place of
$\sigma_n$).  In particular, if $e_1=f_n$ and $\varphi_n(e_1)=\eta_m \in D_n$, then $\mu_1(e_1)=\varphi_n(e_1)$, which 
implies that $e_1$ is outside $P_{y_{p}}(\eta_m,\gamma_2,\varphi_n)$. 
By (1) and (3), we have $\mu_1(f)=\varphi_n(f)$ for each $f\in E(T)$ and $\overline{\mu}_1(u)= 
\phibar_n(u)$ for each $u\in V(T(y_{p-1}))$. So we can obtain $T$ from $T_{n,q}^*+e_1$ by using TAA under $\mu_1$; thereby 
$T$ is an ETT satisfying MP under $\mu_1$. Furthermore, $T_n=T_{n,0} \subset T_{n,1} \subset \ldots \subset T_{n,q} 
\subset T_{n,q+1}=T$ remains to be a good hierarchy of $T$ under $\mu_1$, with the same $\Gamma$-sets as those under $\varphi_n$. Therefore, $(T, \mu_1)$ is also a minimum counterexample to Theorem \ref{hierarchy} (see (6.2)-(6.5)), in which $\gamma_2$ is missing at two vertices.

By Lemma \ref{9n}, we have 
$|\overline{\mu}_1(T(y_{p-2}))- \overline{\mu}_1(T_{n,0}^*-V(T_n)) - \mu_1 \langle T(y_{p-2}) - T_{n,q}^* \rangle | 
\geq 2n+11$, where $T(y_0)=T_{n,q}^*$. It follows that $|\overline{\mu}_1(T(y_{p-2}))- \overline{\mu}_1(T_{n,0}^*-V(T_n)) 
- \mu_1 \langle T - T_{n,q}^* \rangle | \geq 2n+9$. As $|\Gamma^{q}|\le 2 |D_{n,q}| \le 2|D_n| \le 2n$ by Lemma \ref{Dnzang},
using (6.6) we obtain 

(4) there exists a color $\beta$ in $\overline{\mu}_1(T(y_{p-2}))- \overline{\mu}_1(T_{n,0}^*-V(T_n))- 
\mu_1 \langle T - T_{n,q}^* \rangle- \Gamma^{q}$.  

\noindent By Lemma \ref{change} (with $\gamma_2$ in place of $\alpha$), $P_{v_{\gamma_{2}}}(\beta, \gamma_{2},\mu_1)=P_{v_{\beta}}(\beta, \gamma_{2}, \mu_1)$, so it is disjoint from $P_{y_{p}}(\beta, \gamma_{2},\mu_1)$. 
Let $\mu_2=\mu_1/P_{y_{p}}(\beta,\gamma_{2}, \mu_1)$.  By Lemma~\ref{stablechange},  $\mu_2$ satisfies all the properties 
described in (7.3) (with $\mu_2$ in place of $\sigma_n$). By (1), (3) and (4), we have $\beta, \gamma_2 \notin \mu_1 \langle T(y_{p})-T_{n,q}^* \rangle$. So $\mu_2(f)=\mu_1(f)$ for each $f\in E(T)$ and $\overline{\mu}_2(u)=\overline{\mu}_1(u)$ 
for each $u\in V(T(y_{p-1}))$. Hence we can obtain $T$ from $T_{n,q}^*+e_1$ by using TAA under $\mu_2$; thereby $T$ satisfies 
MP under $\mu_2$. Furthermore, $T_n=T_{n,0} \subset T_{n,1} \subset \ldots \subset T_{n,q} \subset T_{n,q+1}=T$ remains to 
be a good hierarchy of $T$ under $\mu_2$, with the same $\Gamma$-sets as those under $\mu_1$.
Therefore, $(T, \mu_2)$ is also a minimum counterexample to Theorem \ref{hierarchy} 
(see (6.2)-(6.5)), in which $\beta$ is missing at two vertices. Since $\theta\in \Gamma^{q}_h$ and $\eta_h\in \overline{\varphi}_n(y_{p-1})=\overline{\mu}_1(y_{p-1})=\overline{\mu}_2(y_{p-1})$, we obtain 

(5) $\theta \notin \mu_2 \langle T(y_{p-1}) - T_{n,q}^* \rangle$. 

\noindent By (4), we also have 

(6) $\beta \notin \mu_2 \langle T - T_{n,q}^* \rangle$. 

\noindent It follows from (5) and Lemma~\ref{change} (with $\theta$ in place of $\alpha$) that $P_{v_{\theta}}(\beta,\theta,\mu_2)=P_{v_{\beta}}(\beta, \theta,\mu_2)$, so it is disjoint from 
$P_{y_{p}}(\beta,\theta,\mu_2)$. Finally, set $\mu_3=\mu_2/P_{y_{p}}(\beta,\theta,
\mu_2)$. By Lemma~\ref{stablechange},  $\mu_3$ satisfies all the properties described in (7.3) (with $\mu_3$ in 
place of $\sigma_n$). From (5) and (6) we see that $T$ can be obtained from $T_{n,q}^*+e_1$ by using TAA under $\mu_3$. Hence 
$T$ is an ETT satisfying MP under $\mu_3$. Note that $\mu_3(f)=\mu_2(f)$ for each $f\in E(T(y_{p-1}))$,
$\mu_3(e_p)=\beta$, and $\overline{\mu}_3(u)=\overline{\mu}_2(u)$ for each $u\in V(T(y_{p-1}))$. Moreover, $\beta \notin 
\Gamma^q$ by (4). It is a routine matter to check that $T_n=T_{n,0} \subset T_{n,1} \subset \ldots \subset T_{n,q} \subset T_{n,q+1}=T$ remains to be a good hierarchy of $T$ under $\mu_3$, with the same $\Gamma$-sets as those under $\mu_2$. 
Since $\mu_3(e_p)=\beta \notin \Gamma^q$ and $v_{\beta} \prec y_{p-1}$,  we see that 
$T^-$ has a good hierarchy and is an ETT satisfying MP with respect to $\mu_3$. As $\theta$ is missing at two vertices
in $T^-$, we conclude that $(T^-, \mu_3)$ is a counterexample to Theorem \ref{hierarchy} (see (6.2) and (6.3)), 
which contradicts the minimality assumption $(6.4)$ or $(6.5)$ on  $(T, \varphi_n)$.

{\bf Subcase 1.2.} $\theta \in\overline{\varphi}_n(y_{p-1})$. 

In this subcase, from (6.6) and TAA we see that

(7) $\theta \notin \overline{\varphi}_n(T(y_{p-2}))$, so $\theta \notin \Gamma^q$ and hence $\theta\ne \gamma_1, \gamma_2$.
Furthermore, $\theta$ is not used by any edge in $T(y_{p-1}) -  T_{n,q}^*$, except possibly $e_1$ when $q=0$ and $T_{n,0}^*=T_n$ 
(now $e_1=f_n$ in Algorithm 3.1 and $\varphi_n(e_1)=\theta \in D_n$). 

Since $P_{v_{\gamma_1}}(\eta_m, \gamma_{1},\varphi_n)=P_{y_{p-1}}(\eta_m,\gamma_{1},\varphi_n)$ by (2), this
path is disjoint from $P_{y_{p}}(\eta_m,\gamma_{1},\varphi_n)$. Let $\mu_1 = \varphi_n/P_{y_{p}}(\eta_m,\gamma_1,\varphi_n)$.
By Lemma~\ref{stablechange},  $\mu_1$ satisfies all the properties described in (7.3) (with $\mu_1$ in place of
$\sigma_n$). In particular, if $e_1=f_n$ and $\varphi_n(e_1)=\eta_m \in D_n$, then $\mu_1(e_1)=\varphi_n(e_1)$, which 
implies that $e_1$ is outside $P_{y_{p}}(\eta_m,\gamma_1,\varphi_n)$.
By (1) and (6.6), we have $\mu_1(f)=\varphi_n(f)$ for each $f\in E(T)$ and $\overline{\mu}_1(u)=\phibar_n(u)$ 
for each $u\in V(T(y_{p-1}))$. So we can obtain $T$ from $T_{n,q}^*+e_1$ by using TAA under $\mu_1$, and hence $T$ satisfies 
MP under $\mu_1$. Furthermore, $T_n=T_{n,0} \subset T_{n,1} \subset \ldots \subset T_{n,q} \subset T_{n,q+1}=T$ 
remains to be a good hierarchy of $T$ under $\mu_1$, with the same $\Gamma$-sets as those under $\varphi_n$.
Therefore, $(T, \mu_1)$ is also a minimum counterexample to Theorem \ref{hierarchy} (see (6.2)-(6.5)), in which $\gamma_1$ is missing at two vertices.

From (1) and the definition of $\mu_1$, we see that

(8) $\gamma_1 \notin \mu_1 \langle T - T_{n,q}^* \rangle$. 

\noindent From (8) and Lemma~\ref{change} (with $\gamma_1$ in place of $\alpha$), we deduce that  $P_{v_{\gamma_{1}}}(\theta,\gamma_{1},\mu_1)=P_{y_{p-1}}(\theta,\gamma_{1},\mu_1)$, which is disjoint from 
$P_{y_{p}}(\theta,\gamma_{1},\mu_1)$. Let $\mu_2=\mu_1/P_{y_{p}}(\theta,\gamma_{1},\mu_1)$. By Lemma~\ref{stablechange},  
$\mu_2$ satisfies all the properties described in (7.3) (with $\mu_2$ in place of $\sigma_n$). 
In particular, if $e_1=f_n$ and $\mu_1(e_1)=\theta \in D_n$, then $\mu_2(e_1)=\mu_1(e_1)$, which 
implies that $e_1$ is outside $P_{y_{p}}(\theta,\gamma_{1},\mu_1)$. In view of (7), (8) and (6.6),
we have $\mu_2(f)=\mu_1(f)$ for each $f\in E(T(y_{p-1}))$, $\mu_2(e_p)=\gamma_1$, and $\overline{\mu}_2(u)=
\overline{\mu}_1(u)$ for each $u\in V(T(y_{p-1}))$. Moreover, $\theta \notin \Gamma^q$. So $T$ can be obtained 
from $T_{n,q}^*+e_1$ by using TAA under $\mu_2$, and hence is an ETT satisfying MP under $\mu_2$. It is a routine 
matter to check that $T_n=T_{n,0} \subset T_{n,1} \subset \ldots \subset T_{n,q} \subset T_{n,q+1}=T$ remains to 
be a good hierarchy of $T$ under $\mu_2$, with the same $\Gamma$-sets as those under $\mu_1$. 
Therefore, $(T, \mu_2)$ is also a minimum counterexample to Theorem 
\ref{hierarchy} (see (6.2)-(6.5)). Since $\theta \in \overline{\mu}_2(y_p)\cap \overline{\mu}_2(y_{p-1})$
and $\mu_2(e_p)= \gamma_1 \notin \overline{\mu}_2(y_{p-1})$, the present subcase reduces to Subcase 1.1 if 
$\theta\in D_{n,q}$ and reduces to Subcase 2.1 (to be discussed below) if $\theta\notin D_{n,q}$. 

\vskip 2mm

{\bf Case 2.} $\alpha\in \phibar_n(y_p)\cap \phibar_n(y_{p-1})$ and $\alpha\notin D_{n,q}$. 

\vskip 1mm
By the definitions of $D_n$ and $D_{n,q}$, we have $\overline{\varphi}_n(T_n) \cup D_n \subseteq  
\overline{\varphi}_n(T_{n,q}^*) \cup D_{n,q}$. Using (6.6) and this set inclusion, we obtain 

(9) $\alpha \notin  \overline{\varphi}_n(T(y_{p-2}))$ and $\alpha \notin D_n$. So $\alpha \notin\varphi_n \langle 
T-T_{n,q}^* \rangle$ by TAA (see, for instance, (1)).

Recall that $T(y_0)=T_{n,q}^*$ and $\theta=\varphi_n(e_p)$. We consider two subcases according to whether $\theta  \in \overline{\varphi}_n(y_{p-1})$.

{\bf Subcase 2.1.} $\theta  \notin\overline{\varphi}_n(y_{p-1})$.

In our discussion about this subcase, we shall also appeal to the following two tree-sequences:

\vskip 1mm
$\bullet$ $T^-=(T_{n,q}^*, e_1, y_1, e_2, \ldots, e_{p-2}, y_{p-2}, e_p, y_p)$ and

\vskip 1mm
$\bullet$ $T^* \hskip 0.6mm =(T_{n,q}^*, e_1, y_1, e_2, \ldots ,y_{p-2}, e_p, y_p, e_{p-1}, y_{p-1})$.

\vskip 1mm
\noindent As stated in Subcase 1.1, $T^-$ is obtained from $T$ by deleting $y_{p-1}$ and $T^*$ arises from $T$ by 
interchanging the order of $(e_{p-1}, y_{p-1})$ and $(e_p, y_p)$.  Furthermore, both $T^-$ and $T^*$ 
are ETTs satisfying MP with respect to $\varphi_n$. Observe that 

(10) $T_n=T_{n,0} \subset T_{n,1} \subset \ldots \subset T_{n,q} \subset T_{n,q+1}=T^*$ is a good hierarchy 
of $T^*$ under $\varphi_n$, unless $\theta\in \Gamma^{q}_h$ for some $\eta_h\in D_{n,q}$ such that $\eta_h 
\in\phibar_n(y_{p-1})$. 

Let us first assume that the exceptional case in (10) does not occur; that is, there exists no $\eta_h\in D_{n,q}$ such 
that $\eta_h \in\phibar_n(y_{p-1})$ and $\theta\in \Gamma^{q}_h$. It is easy to see that now
$T_n=T_{n,0} \subset T_{n,1} \subset \ldots \subset T_{n,q} \subset T^-$ is a good hierarchy of $T^-$ under $\varphi_n$. 

By Lemma \ref{9n}, we have 
$|\overline{\varphi}_1(T(y_{p-2}))- \overline{\varphi}_n(T_{n,0}^*-V(T_n)) - \varphi_1 \langle T(y_{p-2}) - T_{n,q}^* \rangle | 
\geq 2n+11$ holds, where $T(y_0)=T_{n,q}^*$. Since $|\Gamma^{q}|\le 2 |D_{n,q}| \le 2|D_n| \le 2n$ by 
Lemma \ref{Dnzang}, using (6.6) we obtain

(11) there exists a color $\beta$ in $\overline{\varphi}_n(T(y_{p-2}))- \overline{\varphi}_n(T_{n,0}^*-V(T_n))- 
\varphi_n \langle T - T_{n,q}^* \rangle-\Gamma^{q}$.  

\noindent Note that $\beta \notin \overline{\varphi}_n(y_p)$, for otherwise, $(T^-, \varphi_n)$  would be 
a counterexample to Theorem \ref{hierarchy} (see (6.2) and (6.3)), which violates the minimality assumption
$(6.4)$ or $(6.5)$ on $(T, \sigma_n)$. Since $\alpha, \beta \notin \varphi_n \langle T - T_{n,q}^* \rangle$ by
(9) and (11), applying Lemma \ref{change} to $(T, \varphi_n)$ and $(T^*, \varphi_n)$, respectively, we obtain 
$P_{v_{\beta}}(\alpha,\beta,\varphi_n)=P_{y_{p-1}}(\alpha,\beta,\varphi_n)$
and $P_{v_{\beta}}(\alpha,\beta,\varphi_n)=P_{y_{p}}(\alpha,\beta,\varphi_n)$, a contradiction.

So we assume that the exceptional case in (10) occurs; that is, there exists $\eta_h\in D_{n,q}$ such 
that $\eta_h \in\phibar_n(y_{p-1})$ and $\theta\in \Gamma^{q}_h$. For simplicity, we abbreviate the two colors 
$\gamma^{q}_{h_1}$ and $\gamma^{q}_{h_2}$ in $\Gamma^{q}_h$ (see Definition \ref{R2}) to $\gamma_1$ and $\gamma_2$, 
respectively. Renaming subscripts if necessary, we may assume that $\theta=\gamma_1$. By Definition \ref{R2}(i) and
TAA, we have

(12)  $\gamma_2 \notin\varphi_n \langle T-T_{n,q}^* \rangle$ and  $\eta_h$ is not used by 
any edge in $T -  T_{n,q}^*$, except possibly $e_1$ when $q=0$ and $T_{n,0}^*=T_n$ (now $e_1=f_n$ in Algorithm 3.1 
and $\varphi_n(e_1)=\eta_h \in D_{n,q} \subseteq D_n$). 

By (12) and Lemma \ref{change} (with $\alpha$ in place of $\beta$), we obtain $P_{v_{\gamma_2}}(\alpha,\gamma_{2},\varphi_n)=P_{y_{p-1}}(\alpha,\gamma_{2},
\varphi_n)$, which is disjoint from $P_{y_{p}}(\alpha,\gamma_{2},\varphi_n)$.  Let $\mu_1=\varphi_n/P_{y_{p}}(\alpha,\gamma_{2},
\varphi_n)$. By Lemma~\ref{stablechange},  $\mu_1$ satisfies all the properties described in (7.3) (with $\mu_1$ in place of
$\sigma_n$).  Since $\alpha, \gamma_2 \notin \varphi_n \langle T(y_{p})-T_{n,q}^* \rangle$ by (9) and (12), we have 
$\mu_1(f)=\varphi_n(f)$ for each $f\in E(T)$ and $\overline{\mu}_1(u)= \phibar_n(u)$ for each $u\in V(T(y_{p-1}))$. 
So we can obtain $T$ from $T_{n,q}^*+e_1$ by using TAA under $\mu_1$, and hence $T$  
is an ETT satisfying MP under $\mu_1$. Furthermore, $T_n=T_{n,0} \subset T_{n,1} \subset \ldots \subset T_{n,q} \subset T_{n,q+1}=T$ 
remains to be a good hierarchy of $T$ under $\mu_1$, with the same $\Gamma$-sets as those under $\varphi_n$. 
Therefore, $(T, \mu_1)$ is also a minimum counterexample to Theorem \ref{hierarchy} (see (6.2)-(6.5)), in which 
$\gamma_2$ is missing at two vertices.

If $\eta_h\in\overline{\mu}_1(y_{p})$, then $\eta_h \in\overline{\mu}_1(y_p)\cap\overline{\mu}_1(y_{p-1})$,
$\eta_h\in D_{n,q}$, and $\mu_1(e_p)=\gamma_{1} \notin \phibar_n(y_{p-1})$. Thus the present subcase reduces 
to Subcase 1.1. So we may assume that $\eta_h\notin\overline{\mu}_1(y_{p})$. By (12) and the definition of $\mu_1$, we have

(13) $\gamma_2 \notin\mu_1 \langle T-T_{n,q}^* \rangle$ and $\eta_h$ is not used by any edge in $T -  T_{n,q}^*$ 
under $\mu_1$, except possibly $e_1$ when $q=0$ and $T_{n,0}^*=T_n$ (now $e_1=f_n$ in Algorithm 3.1 and 
$\mu_1(e_1)=\eta_h \in D_n$). 

By (13) and Lemma \ref{change} (with $\gamma_2$ in place of $\alpha$), we obtain $P_{v_{\gamma_2}}(\eta_h,\gamma_{2},\mu_1)=P_{y_{p-1}}(\eta_h,\gamma_{2},
\mu_1)$, which is disjoint from $P_{y_{p}}(\eta_h,\gamma_{2},\mu_1)$.  Let $\mu_2=\mu_1/P_{y_{p}}(\eta_h,\gamma_{2},
\mu_1)$. By Lemma~\ref{stablechange},  $\mu_2$ satisfies all the properties described in (7.3) (with $\mu_2$ in place of
$\sigma_n$).  In particular, if $e_1=f_n$ and $\mu_1(e_1)=\eta_h \in D_n$, then $\mu_2(e_1)=\mu_1(e_1)$, which 
implies that $e_1$ is outside $P_{y_{p}}(\eta_h,\gamma_{2},\mu_1)$. By (13), we have $\mu_2(f)=\mu_1(f)$ for each 
$f\in E(T)$ and $\overline{\mu}_2(u)= \overline{\mu}_1(u)$ for each $u\in V(T(y_{p-1}))$. So we can obtain $T$ from
$T_{n,q}^*+e_1$ by using TAA under $\mu_2$, and hence $T$ is an ETT satisfying MP under $\mu_2$. Furthermore,
$T_n=T_{n,0} \subset T_{n,1} \subset \ldots \subset T_{n,q} \subset T_{n,q+1}=T$ 
remains to be a good hierarchy of $T$ under $\mu_2$, with the same $\Gamma$-sets as those under $\mu_1$. 
Therefore, $(T, \mu_2)$ is also a minimum counterexample to Theorem \ref{hierarchy} (see (6.2)-(6.5)), in which $\eta_h \in\overline{\mu}_2(y_p)\cap\overline{\mu}_2(y_{p-1})$,  $\eta_h\in D_{n,q}$, and $\mu_2(e_p)=\gamma_1 \notin \overline{\mu}_2(y_{p-1})$. Thus the present subcase reduces to Subcase 1.1.

{\bf Subcase 2.2.} $\theta\in\overline{\varphi}_n(y_{p-1})$. 

Let us first assume that $\theta \in D_{n,q}$; that is, $\theta=\eta_{m}$ for some $\eta_m\in D_{n,q}$. For 
simplicity, we use $\varepsilon_1$ and $\varepsilon_2$ to denote the two colors $\gamma^{q}_{m_1}$ and 
$\gamma^{q}_{m_2}$ in $\Gamma^{q}_m$ (see Definition \ref{R2}), respectively. By Definition \ref{R2}(i) and TAA, we have

(14) $\varepsilon_1, \varepsilon_2 \notin \varphi_n \langle T-T_{n,q}^* \rangle$ and $\eta_m$ is not used by 
any edge in $T(y_{p-1}) -  T_{n,q}^*$, except possibly $e_1$ when $q=0$ and $T_{n,0}^*=T_n$ (now $e_1=f_n$ in 
Algorithm 3.1 and $\varphi_n(e_1)=\eta_m \in D_n$). 

By (14) and Lemma \ref{change}, we obtain $P_{v_{\varepsilon_1}}(\alpha, \varepsilon_1, \varphi_n)
=P_{y_{p-1}}(\alpha, \varepsilon_1, \varphi_n)$, which is disjoint from $P_{y_{p}}(\alpha, \varepsilon_1, \varphi_n)$.
Let $\mu_1=\varphi_n/P_{y_{p}}(\alpha, \varepsilon_1, \varphi_n)$. By Lemma~\ref{stablechange}, $\mu_1$ satisfies all 
the properties described in (7.3) (with $\mu_1$ in place of $\sigma_n$).  By (9) and (14), we have 

(15) $\alpha, \varepsilon_1 \notin \mu_1 \langle T-T_{n,q}^* \rangle$ and $\eta_m$ is not used by any edge in 
$T(y_{p-1}) -  T_{n,q}^*$ under $\mu_1$, except possibly $e_1$ when $q=0$ and $T_{n,0}^*=T_n$ (now $e_1=f_n$ in 
Algorithm 3.1 and $\mu_1(e_1)=\eta_m \in D_n$). 

So $\mu_1(f)=\varphi_n(f)$ for each $f\in E(T)$ and $\overline{\mu}_1(u)= \overline{\varphi}_n(u)$ for each 
$u\in V(T(y_{p-1}))$. Thus $T$ can be obtained from $T_{n,q}^*+e_1$ by using TAA under $\mu_1$, and hence is an ETT 
satisfying MP under $\mu_1$. Furthermore, $T_n=T_{n,0} \subset T_{n,1} \subset 
\ldots \subset T_{n,q} \subset T_{n,q+1}=T$ remains to be a good hierarchy of $T$ under $\mu_1$, with the same 
$\Gamma$-sets as those under $\varphi_n$. Therefore, $(T, \mu_1)$ is also a minimum counterexample to Theorem 
\ref{hierarchy} (see (6.2)-(6.5)), in which $\varepsilon_1$ is missing at two vertices. 

By (15) and Lemma \ref{change} (with $\varepsilon_1$ in place of $\alpha$), we obtain $P_{v_{\varepsilon_1}}(\eta_m, 
\varepsilon_1, \mu_1) =P_{y_{p-1}}(\eta_m, \varepsilon_1, \mu_1)$, which is disjoint from $P_{y_{p}}(\eta_m, \varepsilon_1, 
\mu_1)$. Let $\mu_2=\mu_1/P_{y_{p}}(\eta_m, \varepsilon_1, \mu_1)$. By Lemma~\ref{stablechange}, $\mu_2$ satisfies all 
the properties described in (7.3) (with $\mu_2$ in place of $\sigma_n$).  
In particular, if $e_1=f_n$ and $\mu_1(e_1)=\eta_m \in D_n$, then $\mu_2(e_1)=\mu_1(e_1)$, which 
implies that $e_1$ is outside $P_{y_{p}}(\eta_m, \varepsilon_1, \mu_1)$.  In view of (15), we have 
$\mu_2(f)=\mu_1(f)$ for each $f\in E(T(y_{p-1}))$, $\mu_2(e_p)=\varepsilon_1$, and $\overline{\mu}_2(u)= \overline{\mu}_1(u)$ 
for each $u\in V(T(y_{p-1}))$. So $T$ can be obtained from $T_{n,q}^*+e_1$ by using TAA under $\mu_2$, and hence satisfies 
MP under $\mu_2$. Furthermore, $T_n=T_{n,0} \subset T_{n,1} \subset \ldots \subset T_{n,q} \subset T_{n,q+1}=T$ 
remains to be a good hierarchy of $T$ under $\mu_2$, with the same $\Gamma$-sets as those under $\mu_1$. 
Therefore, $(T, \mu_2)$ is also a minimum counterexample to Theorem \ref{hierarchy} (see (6.2)-(6.5)), in which $\eta_m \in\overline{\mu}_2(y_p)\cap\overline{\mu}_2(y_{p-1})$,  $\eta_m\in D_{n,q}$, and $\mu_2(e_p)=\varepsilon_1 \notin \overline{\mu}_2(y_{p-1})$. Thus the present subcase reduces to Subcase 1.1. 

Next we assume that $\theta\notin D_{n,q}$. Set $T(y_0)=T_{n,q}^*$. We propose to show that 

(16) there exists a color $\beta \in \overline{\varphi}_n(T(y_{p-2}))- \phibar_n(T_{n,0}^*-V(T_n)) - 
\varphi_n \langle T - T_{n,q}^* \rangle- D_{n,q}$, such that either $\beta \notin \Gamma^{q}$
or $\beta \in \Gamma^{q}_h$ for some $\eta_h \in D_{n,q}\cap \phibar_n(T(y_{p-2}))$.

To justify this, note that if $|\overline{\varphi}_n(T(y_{p-2}))- \phibar_n(T_{n,0}^*-V(T_n)) - \varphi_n \langle 
T(y_{p-2}) - T_{n,q}^* \rangle -(\Gamma^{q}\cup D_{n,q}) |\ge  5$, then 
$|\overline{\varphi}_n(T(y_{p-2}))- \phibar_n(T_{n,0}^*-V(T_n)) - \varphi_n \langle 
T - T_{n,q}^* \rangle - (\Gamma^{q}\cup D_{n,q}) |\ge 3$, because $T-T(y_{p-2})$ contains precisely two edges.
Thus there exists a color $\beta \in \overline{\varphi}_n(T(y_{p-2}))- \phibar_n(T_{n,0}^*-V(T_n)) - 
\varphi_n \langle T - T_{n,q}^* \rangle- D_{n,q}$, such that $\beta \notin \Gamma^{q}$.

So we assume that $|\overline{\varphi}_n(T(y_{p-2}))- \phibar_n(T_{n,0}^*-V(T_n)) - \varphi_n \langle 
T(y_{p-2}) - T_{n,q}^* \rangle - (\Gamma^{q}\cup D_{n,q}) |\le 4 $. By Lemma \ref{9n}, there exist $7$ distinct 
colors $\eta_{h}\in D_{n,q}\cap \phibar_n(T(y_{p-2}))$ such that  $(\Gamma^{q}_h \cup \{\eta_{h}\})
\cap \varphi_n \langle T(y_{p-2}) - T_{n,q}^* \rangle=\emptyset$. Let $\beta$ be an arbitrary color in such a
$\Gamma^{q}_h$. From Definition \ref{R2}, we see that $\Gamma^{q}_h \subseteq \overline{\varphi}_n(T_{n,q}^*)
\subseteq \overline{\varphi}_n(T(y_{p-2}))$,  $\Gamma^{q}_h \cap \phibar_n(T_{n,0}^*-V(T_n)) =\emptyset$,
and $\Gamma^{q}_h \cap D_{n,q} =\emptyset$ (see (5.7)). So $\beta \in \overline{\varphi}_n(T(y_{p-2}))- \phibar_n(T_{n,0}^*-V(T_n)) - 
\varphi_n \langle T(y_{p-2}) - T_{n,q}^* \rangle- D_{n,q}$. Since $T-T(y_{p-2})$ contains precisely two edges, there
exists $\beta \in \overline{\varphi}_n(T(y_{p-2}))- \phibar_n(T_{n,0}^*-V(T_n)) - \varphi_n \langle T - T_{n,q}^* 
\rangle- D_{n,q}$, such that $\beta \in \Gamma^{q}_h$ for some $\eta_h \in D_{n,q}\cap \phibar_n(T(y_{p-2}))$. 
Hence (16) is established.

By the definitions of $D_n$ and $D_{n,q}$, we have $\overline{\varphi}_n(T_n) \cup D_n \subseteq  
\overline{\varphi}_n(T_{n,q}^*) \cup D_{n,q}$. By (16), $\beta \notin \phibar_n(T_{n,0}^*-V(T_n)) \cup D_{n,q}$.
It follows from these two observations that 

(17) if $q\ge 1$, then $\beta \in \overline{\varphi}_n(T_{n,q}^*)$ or $\beta \notin D_n$; if $q=0$, then
$\beta \in \overline{\varphi}_n(T_n)$ or $\beta \notin D_n$.  

By (9), (16), (17) and Lemma \ref{change}, we obtain $P_{v_{\beta}}(\alpha,\beta,\varphi_n)=P_{y_{p-1}}(\alpha,\beta,\varphi_n)$,
which is disjoint from $P_{y_{p}}(\alpha,\beta,\varphi_n)$. Let $\mu_3=\varphi_n/P_{y_{p}}(\alpha,\beta,\varphi_n)$.
By Lemma~\ref{stablechange},  $\mu_3$ satisfies all the properties described in (7.3) (with $\mu_3$ in place of
$\sigma_n$). By (9) and (16), we have $\alpha, \beta \notin \varphi_n \langle T-T_{n,q}^* \rangle$. So

(18) $\alpha, \beta \notin \mu_3 \langle T-T_{n,q}^* \rangle$,

\noindent $\mu_3(f)=\varphi_n(f)$ for each $f\in E(T)$, and $\overline{\mu}_3(u)= \phibar_n(u)$ for each $u\in V(T(y_{p-1}))$. 
Thus we can obtain $T$ from $T_{n,q}^*+e_1$ by using TAA under $\mu_3$, and hence $T$ is an ETT satisfying MP under $\mu_3$. 
Furthermore, $T_n=T_{n,0} \subset T_{n,1} \subset \ldots \subset T_{n,q} \subset T_{n,q+1}=T$ remains to be a good hierarchy 
of $T$ under $\mu_3$, with the same $\Gamma$-sets as those under $\varphi_n$. Therefore, $(T, \mu_3)$ is also a minimum 
counterexample to Theorem \ref{hierarchy} (see (6.2)-(6.5)), in which $\beta$ is missing at two vertices.

Since $\theta \in \overline{\varphi}_n(y_{p-1})$, it follows from (6.6) that $\theta\notin \overline{\varphi}_n(T_{n,q}^*)$.
By assumption, $\theta\notin D_{n,q}$. As $\overline{\varphi}_n(T_n) \cup D_n \subseteq  \overline{\varphi}_n(T_{n,q}^*) \cup D_{n,q}$, we obtain

(19) $\theta\notin D_n$ and hence $\theta \notin \mu_3 \langle T(y_{p-1})-T_{n,q}^* \rangle$ by TAA.

By (17)-(19) and Lemma \ref{change}, we obtain $P_{v_{\beta}}(\theta,\beta,\mu_3)=P_{y_{p-1}}(\theta,\beta,\mu_3)$,
which is disjoint from $P_{y_{p}}(\theta,\beta,\mu_3)$. Let $\mu_4=\mu_3/P_{y_{p}}(\theta,\beta,\mu_3)$.
By Lemma~\ref{stablechange},  $\mu_4$ satisfies all the properties described in (7.3) (with $\mu_4$ in place of
$\sigma_n$).  By (18) and (19), we have $\mu_4(f)=\mu_3(f)$ for each $f\in E(T(y_{p-1}))$ and $\overline{\mu}_4(u)= \overline{\mu}_3(u)$ for each $u\in V(T(y_{p-1}))$. So we can obtain $T$ from $T_{n,q}^*+e_1$ by using TAA under $\mu_4$, and hence 
$T$ is an ETT satisfying MP under $\mu_4$. Since either $\beta \notin \Gamma^{q}$ or $\beta \in \Gamma^{q}_h$ for some $\eta_h 
\in D_{n,q}\cap \overline{\mu}_3(T(y_{p-2}))$ by (16), it follows that $T_n=T_{n,0} \subset T_{n,1} \subset \ldots \subset 
T_{n,q} \subset T_{n,q+1}=T$ remains to be a good hierarchy of $T$ under $\mu_4$, with the same $\Gamma$-sets as those under 
$\mu_3$. Therefore, $(T, \mu_4)$ is also a minimum counterexample to Theorem \ref{hierarchy} (see (6.2)-(6.5)), in which 
$\theta \in\overline{\mu}_4(y_p)\cap\overline{\mu}_4(y_{p-1})$, $\theta \notin D_{n,q}$, and $\mu_4(e_p)=\beta \notin \overline{\mu}_4(y_{p-1})$. Thus the present subcase reduces to Subcase 2.1.

\vskip 2mm

{\bf Case 3.} $\alpha\in \overline{\varphi}_n(y_p)\cap\overline{\varphi}_n(v)$ for some vertex $v\prec y_{p-1}$.
\vskip 1mm
Set $T(y_0)=T_{n,q}^*$. Let us first impose some restrictions on $\alpha$.

(20) We may assume that $\alpha\in\phibar_n(T(y_{p-2}))-\varphi_n \langle T - T_{n,q}^* \rangle$, such that 
either $\alpha\notin D_{n,q}\cup\Gamma^{q}$ if $q \ge 1$ and $\alpha\notin D_{n}\cup\Gamma^{0}$ if $q=0$, 
or $\alpha$ is some $\eta_{h}\in D_{n,q}$ satisfying $\Gamma^{q}_h \cap \varphi_n \langle T - T_{n,q}^* \rangle=\emptyset$.

To justify this, note that if $|\overline{\varphi}_n(T(y_{p-2}))- \phibar_n(T_{n,0}^*-V(T_n)) - \varphi_n \langle 
T(y_{p-2}) - T_{n,q}^* \rangle - (\Gamma^{q}\cup D_{n,q}) |\ge  5$, then 
$|\overline{\varphi}_n(T(y_{p-2}))- \phibar_n(T_{n,0}^*-V(T_n)) - \varphi_n \langle 
T - T_{n,q}^* \rangle -(\Gamma^{q}\cup D_{n,q}) |\ge 3$, because $T-T(y_{p-2})$ contains precisely two edges.
Thus there exists a color $\beta \in \overline{\varphi}_n(T(y_{p-2}))- \phibar_n(T_{n,0}^*-V(T_n)) - 
\varphi_n \langle T - T_{n,q}^* \rangle- (\Gamma^{q}\cup D_{n,q})$. Clearly,
$\beta\in\phibar_n(T(y_{p-2}))-\varphi_n \langle T - T_{n,q}^* \rangle$ and
$\beta\notin D_{n,q}\cup\Gamma^{q}$ if $q \ge 1$ and $\beta \notin D_{n}\cup\Gamma^{0}$ if $q=0$ 
(note that $\beta \notin D_{n}$ because $\beta \notin \phibar_n(T_{n,0}^*-V(T_n)) \cup D_{n,0}$).

If $|\overline{\varphi}_n(T(y_{p-2}))- \phibar_n(T_{n,0}^*-V(T_n)) - \varphi_n \langle 
T(y_{p-2}) - T_{n,q}^* \rangle -(\Gamma^{q}\cup D_{n,q}) |\le 4 $, then, by Lemma \ref{9n}, there exist $7$ distinct 
colors $\eta_{h}\in D_{n,q}\cap \phibar_n(T(y_{p-2}))$ such that  $(\Gamma^{q}_h \cup \{\eta_{h}\})
\cap \varphi_n \langle T(y_{p-2}) - T_{n,q}^* \rangle=\emptyset$. Since $T-T(y_{p-2})$ contains precisely two edges,
there exists one of these $\eta_{h}$, denoted by $\beta$, such that $(\Gamma^{q}_h \cup \{\eta_{h}\})
\cap \varphi_n \langle T - T_{n,q}^* \rangle=\emptyset$.

Combining the above observations, we conclude that

(21) there exists $\beta \in\phibar_n(T(y_{p-2}))-\varphi_n \langle T - T_{n,q}^* \rangle$, such that 
either $\beta \notin D_{n,q}\cup\Gamma^{q}$ if $q \ge 1$ and $\beta \notin D_{n}\cup\Gamma^{0}$ if $q=0$, 
or $\beta$ is some $\eta_{h}\in D_{n,q}$ satisfying $\Gamma^{q}_h \cap \varphi_n \langle T - T_{n,q}^* \rangle=\emptyset$.

If $\beta \in \overline{\varphi}_n(y_p)$, then (20) holds by replacing $\alpha$ with $\beta$ (recall the hypothesis of
the present case). So we assume hereafter that $\beta \notin \overline{\varphi}_n(y_p)$. Let $Q=P_{y_p}(\alpha, \beta, \varphi_n)$ 
and let $\sigma_n=\varphi_n/Q$. We propose to show that one of the following statements (a) and (b) holds:
\begin{itemize}
\vspace{-2mm}
\item[(a)] $\sigma_n$ is a $(T_{n,q}^*, D_n,\varphi_n)$-weakly stable coloring, 
 $T$ is also an ETT satisfying MP with respect to $\sigma_n$, and $T_n=T_{n,0} \subset T_{n,1}
\subset \ldots \subset T_{n,q}\subset T_{n,q+1}=T$ remains to be a hierarchy of $T$ under $\sigma_n$, with the 
same $\Gamma$-sets (see Definition \ref{R2}) as those under $\varphi_n$. Moreover, (20) holds with respect to 
$(T, \sigma_n)$. 
\vspace{-2mm}
\item[(b)] There exists an ETT $T'$ satisfying MP with respect to $\varphi_n$, such that $T_n=T_{n,0} \subset T_{n,1}
\subset \ldots \subset T_{n,q}\subset T'$ is a good hierarchy of $T'$ under $\varphi_n$, with the 
same $\Gamma$-sets as $T$ under $\varphi_n$. Moreover, $V(T')$ is not elementary with respect to $\varphi_n$
and $p(T') < p(T)$.
\vspace{-2mm}
\end{itemize}

\noindent Note that if (b) holds, then $(T', \varphi_n)$ would be a counterexample to Theorem \ref{hierarchy} (see (6.2) 
and (6.3)), which violates the minimality assumption $(6.4)$ on $(T, \varphi_n)$.

Let us first assume that $Q$ is vertex-disjoint from $T(y_{p-1})$. By Lemma \ref{LEM:Stable}, $\sigma_n$ is both $(T(y_{p-1}),D_n,\varphi_n)$-stable and $(T(y_{p-1}),\varphi_n)$-invariant. If $\Theta_n=PE$, then $\sigma_n$ is also
$(T_n\oplus R_n,D_n,\varphi_n)$-stable. Furthermore, $T(y_{p-1})$ is an ETT satisfying MP with respect to $\sigma_n$,
and $T_n=T_{n,0} \subset T_{n,1} \subset \ldots \subset T_{n,q} \subset T(y_{p-1})$ is a good hierarchy of 
$T(y_{p-1})$, with the same $\Gamma$-sets as $T$ under $\sigma_n$. By definition, $\sigma_n$ is a $(T_{n,q}^*, D_n,\varphi_n)$-weakly stable coloring. By the hypothesis of Case 3 and assumption on $\beta$, we have $\varphi_n(e_p) \ne \alpha, 
\beta$. Thus it is clear that (a) is true, and (20) follows if we replace $\varphi_n$ by $\sigma_n$ and $\alpha$ by $\beta$.

Next we assume that $Q$ and $T(y_{p-1})$ have vertices in common. Let $u$ be the first vertex of $Q$ contained in
$T(y_{p-1})$ as we traverse $Q$ from $y_p$. Define $T'=T(y_{p-1}) \cup Q[u,y_p]$ if $u= y_{p-1}$
and $T'=T(y_{p-2})\cup Q[u,y_p]$ otherwise.  By the hypothesis of Case 3 and (21), we have $\alpha,\beta \in
\phibar_n(T(y_{p-2}))$. So $T'$ can be obtained from $T(y_{p-2})$ by using TAA under $\varphi_n$, with $p(T')<p(T)$. It 
follows that $T'$ is an ETT satisfying MP with respect to $\varphi_n$.

By Definition \ref{R2}, we have $D_{n,q} \cap \Gamma^q=\emptyset$ (see (5.7)). Thus

(22)  $\beta\notin \Gamma^{q}$ by (21).

Let us proceed by considering three possibilities for $\alpha$.
  
$\bullet$ $\alpha\notin \Gamma^{q}$. Since both $\alpha$ and $\beta$ are outside $\Gamma^{q}$ (see (22)),
it is easy to see that $T_n=T_{n,0} \subset T_{n,1} \subset \ldots \subset T_{n,q}\subset T'$ is a good hierarchy of 
$T'$ under $\varphi_n$, with the same $\Gamma$-sets as $T$ under $\varphi_n$. Hence (b) holds.

$\bullet$ $\alpha\in \Gamma^{q}\cap \varphi_n \langle T-T_{n,q}^* \rangle$.  Let $\alpha \in \Gamma^{q}_h$ for some
$\eta_h \in D_{n,q}$. Since $\varphi(e_p)\neq\alpha$, we have $\alpha\in \varphi_n \langle T(y_{p-1})-T_{n,q}^* \rangle$. 
Hence  $\eta_h\in \phibar_n(T(y_{p-2}))$ by Definition \ref{R2}(i).  Furthermore, $\beta\in\phibar_n(T(y_{p-2}))$ 
and $\beta\notin\Gamma^q$ by (21) and (22). Therefore, $T_n=T_{n,0} \subset T_{n,1} \subset \ldots \subset T_{n,q}\subset T'$ 
is a good hierarchy of $T'$ under $\varphi_n$, with the same $\Gamma$-sets as $T$ under $\varphi_n$. Hence (b) holds.

$\bullet$ $\alpha \in \Gamma^{q} - \varphi_n \langle T-T_{n,q}^* \rangle$. By the definition of $\Gamma^q$, we have
$\alpha \in \overline{\varphi}_n(T_{n,q})$ if $q \ge 1$ and $\alpha \in \overline{\varphi}_n(T_n)$ if $q=0$. 
It follows from Lemma~\ref{change} that $P_{v_{\alpha}}(\alpha,\beta,\varphi_n)=P_{v_{\beta}}(\alpha,\beta,\varphi_n)$, 
which is disjoint from $Q$. By Lemma~\ref{stablechange}, $\sigma_n=\varphi_n/Q$ satisfies all the properties described 
in (7.3).  Since $\alpha, \beta \notin \varphi_n \langle T-T_{n,q}^* \rangle$ by the assumption on $\alpha$, (21) and (6.6), 
we have $\sigma_n(f)=\varphi_n(f)$ for each $f\in E(T)$ and $\overline{\sigma}_n(u)= \phibar_n(u)$ for each $u\in V(T(y_{p-1}))$.
So we can obtain $T$ from $T_{n,q}^*+e_1$ by using TAA under $\sigma_n$, and hence $T$ is an ETT satisfying MP under $\sigma_n$. 
Furthermore, $T_n=T_{n,0} \subset T_{n,1} \subset \ldots \subset T_{n,q} \subset T_{n,q+1}=T$ remains to be a good hierarchy of $T$ 
under $\sigma_n$, with the same $\Gamma$-sets as those under $\varphi_n$. Therefore, $(T, \sigma_n)$ is also a minimum counterexample to Theorem \ref{hierarchy} (see (6.2)-(6.5)), in which $\beta$ is missing at two vertices. So (a) holds and therefore (20) is established by replacing $\varphi_n$ with $\sigma_n$ and $\beta$ with $\alpha$.

\vskip 3mm

Let $\alpha$ be a color as specified in (20). Recall that $\theta=\varphi_n(e_p)$. We consider 
two subcases according to whether $\theta\in \phibar_n(y_{p-1})$.

{\bf Subcase 3.1.} $\theta \notin\phibar_n(y_{p-1})$.

Consider the tree-sequence $T^-=(T_{n,q}^*, e_1, y_1, e_2, \ldots, e_{p-2}, y_{p-2}, e_p, y_p)$. As stated in 
Subcase 1.1, $T^-$ arises from $T$ by deleting $y_{p-1}$, and $T^-$ is an ETT satisfying MP with respect to $\varphi_n$. 
Observe that 

(23) $T_n=T_{n,0} \subset T_{n,1} \subset \ldots \subset T_{n,q} \subset T^-$ is a good hierarchy 
of $T^-$ under $\varphi_n$, unless $\theta\in \Gamma^{q}_m$ for some $\eta_m\in D_{n,q}$ such that $\eta_m 
\in\phibar_n(y_{p-1})$. 

It follows that the exceptional case stated in (23) must occur, for otherwise, $(T^-, \varphi_n)$ would be a counterexample
to Theorem \ref{hierarchy} (see (6.2) and (6.3)), which violates the minimality assumption $(6.4)$ or $(6.5)$ on 
$(T, \varphi_n)$. So $\theta\in \Gamma^{q}_m$ for some $\eta_m\in D_{n,q}$ such that $\eta_m \in\phibar_n(y_{p-1})$. 

Since $\alpha\in\phibar_n(T(y_{p-2}))$, we have $\alpha\neq \eta_m$ by (6.6). From Definition \ref{R2}(i),
we see that 

(24) $\theta \notin \varphi_n \langle T(y_{p-1}) - T_{n,q}^* \rangle$. 

\noindent By the definition of $\Gamma^q$, we have $\theta \in \overline{\varphi}_n(T_{n,q})$ if $q \ge 1$ and $\theta \in \overline{\varphi}_n(T_n)$ if $q=0$. Thus, by (20), (24) and Lemma~\ref{change}, we obtain $P_{v_{\alpha}} (\alpha, \theta, 
\varphi_n) = P_{v_{\theta}}(\alpha, \theta, \varphi_n)$, which is disjoint from $P_{y_p}(\alpha, \theta, \varphi_n)$. Let $\mu_1=\varphi_n/P_{y_{p}}(\alpha,\theta, \varphi_n)$. By Lemma~\ref{stablechange},  $\mu_1$ satisfies all the 
properties described in (7.3) (with $\mu_1$ in place of $\sigma_n$).  Using (20) and (24), we get 

(25) $\alpha, \theta \notin \mu_1 \langle T(y_{p-1})-T_{n,q}^* \rangle$,  

\noindent $\mu_1(f)=\varphi_n(f)$ for each $f\in E(T(y_{p-1}))$, $\mu_1(e_p)=\alpha \notin \Gamma^q$ (see (20)), 
and $\overline{\mu}_1(u)= \phibar_n(u)$ for each $u\in V(T(y_{p-1}))$. So we can obtain $T$ from $T_{n,q}^*+e_1$ by 
using TAA under $\mu_1$ and hence $T$ is an ETT satisfying MP under $\mu_1$. Furthermore, $T_n=T_{n,0} \subset T_{n,1} \subset 
\ldots \subset T_{n,q} \subset T_{n,q+1}=T$ remains to be a good hierarchy of $T$ under $\mu_1$,
with the same $\Gamma$-sets as those under $\varphi_n$. Therefore, $(T, \mu_1)$ is also a minimum counterexample to Theorem \ref{hierarchy} (see (6.2)-(6.5)), in which $\theta$ is missing at two vertices.

By (25) and  Lemma~\ref{change}, we obtain $P_{v_{\theta}} (\eta_m, \theta, \mu_1) = P_{y_{p-1}}(\eta_m, \theta, 
\mu_1)$, which is disjoint from $P_{y_p}(\eta_m, \theta, \mu_1)$. Let $\mu_2=\mu_1/P_{y_{p}}(\eta_m,\theta, 
\mu_1)$. By Lemma~\ref{stablechange},  $\mu_2$ satisfies all the properties described in (7.3) (with $\mu_2$ in 
place of $\sigma_n$). Note that $\eta_m$ is not used by any edge in $T -  T_{n,q}^*$ under $\mu_1$, except possibly 
$e_1$ when $q=0$ and $T_{n,0}^*=T_n$ (now $e_1=f_n$ in Algorithm 3.1 and $\mu_1(e_1)=\eta_m \in D_n$). So
$e_1$ is outside $P_{y_p}(\eta_m, \theta, \mu_1)$. Hence 
$\mu_2(f)=\mu_1(f)$ for each $f\in E(T)$, and $\overline{\mu}_2(u)= \overline{\mu}_1(u)$ for each $u\in V(T(y_{p-1}))$.
It follows that $T$ can be obtained from $T_{n,q}^*+e_1$ by using TAA and hence is an ETT satisfying MP under $\mu_2$. Furthermore, 
$T_n=T_{n,0} \subset T_{n,1} \subset \ldots \subset T_{n,q} \subset T_{n,q+1}=T$ remains to be a good hierarchy of $T$ under $\mu_2$, with the same $\Gamma$-sets as those under $\mu_1$. Therefore, $(T, \mu_2)$ is also a minimum counterexample to Theorem \ref{hierarchy} (see (6.2)-(6.5)). Since $\eta_m \in \overline{\mu}_2(y_p) \cap \overline{\mu}_2(y_{p-1})$, $\eta_m\in D_{n,q}$, and $\mu_2(e_p)=\alpha \notin \overline{\mu}_2(y_{p-1})$, the present subcase reduces to Subcase 1.1.

{\bf Subcase 3.2.} $\theta \in\overline{\varphi}_n(y_{p-1})$.

We first assume that $\theta \in D_{n,q}$. Let $\theta=\eta_{m}\in D_{n,q}$. For simplicity, we abbreviate the 
two colors $\gamma^{q}_{m_1}$ and $\gamma^{q}_{m_2}$ in $\Gamma^{q}_m$ (see Definition \ref{R2}) to $\gamma_1$ and 
$\gamma_2$, respectively. By (20) and Definition \ref{R2}(i), we have

(26) $\{\alpha, \gamma_1, \gamma_2\}\cap \varphi_n \langle T - T_{n,q}^* \rangle=\emptyset$.

By (26) and  Lemma~\ref{change}, we obtain $P_{v_{\alpha}} (\alpha, \gamma_1, \varphi_n) = P_{v_{\gamma_1}}(\alpha, \gamma_1, \varphi_n)$, which is disjoint from $P_{y_p}(\alpha, \gamma_1, \varphi_n)$. Let $\mu_1=\varphi_n/P_{y_{p}}(\alpha, \gamma_1, 
\varphi_n)$. By Lemma~\ref{stablechange},  $\mu_1$ satisfies all the properties described in (7.3) (with $\mu_1$ in 
place of $\sigma_n$). Since $\mu_1(f)=\varphi_n(f)$ for each $f\in E(T)$, and $\overline{\mu}_1(u)= 
\overline{\varphi}_n(u)$ for each $u\in V(T(y_{p-1}))$, we can obtain $T$ from $T_{n,q}^*+e_1$ by using TAA under 
$\mu_1$ and hence $T$ is an ETT satisfying MP under $\mu_1$. Furthermore, $T_n=T_{n,0} \subset T_{n,1} \subset \ldots \subset 
T_{n,q} \subset T_{n,q+1}=T$ remains to be a good hierarchy of $T$ under $\mu_1$, with the same $\Gamma$-sets as those under $\mu_1$. Therefore, $(T, \mu_1)$ is also a minimum counterexample to Theorem \ref{hierarchy} (see (6.2)-(6.5)),
in which $\gamma_1$ is missing at two vertices. In view of (26) and Definition \ref{R2}(i), we get

(27) $\{\alpha, \gamma_1, \gamma_2\}\cap \mu_1 \langle T - T_{n,q}^* \rangle=\emptyset$, and $\eta_m$ is not used by 
any edge in $T -  T_{n,q}^*$ under $\mu_1$, except possibly $e_1$ when $q=0$ and $T_{n,0}^*=T_n$ (now $e_1=f_n$ in 
Algorithm 3.1 and $\mu_1(e_1)=\eta_m \in D_{n,q} \subseteq D_n$). 

By (27) and  Lemma~\ref{change}, we obtain $P_{v_{\gamma_1}} (\gamma_1, \eta_m, \mu_1) = P_{y_{p-1}}(\gamma_1,\eta_m,  
\mu_1)$, which is disjoint from $P_{y_p}(\gamma_1,\eta_m, \mu_1)$. Let $\mu_2=\mu_1/P_{y_p}(\gamma_1,\eta_m, \mu_1)$.
By Lemma~\ref{stablechange},  $\mu_2$ satisfies all the properties described in (7.3) (with $\mu_2$ in 
place of $\sigma_n$). In particular, if $e_1=f_n$ and $\mu_1(e_1)=\eta_m \in D_n$, then $\mu_2(e_1)=\mu_1(e_1)$, which 
implies that $e_1$ is outside $P_{y_p}(\gamma_1,\eta_m, \mu_1)$.
Since $\mu_2(f)=\mu_1(f)$ for each $f\in E(T(y_{p-1}))$ by (27), and $\overline{\mu}_2(u)= 
\overline{\mu}_1(u)$ for each $u\in V(T(y_{p-1}))$, we can obtain $T$ from $T_{n,q}^*+e_1$ by using TAA under 
$\mu_2$ and hence $T$ is an ETT satisfying MP under $\mu_2$. Furthermore, $T_n=T_{n,0} \subset T_{n,1} \subset \ldots \subset 
T_{n,q} \subset T_{n,q+1}=T$ remains to be a good hierarchy of $T$ under $\mu_2$, with the same $\Gamma$-sets as those under $\mu_1$. Therefore, $(T, \mu_2)$ is also a minimum counterexample to Theorem \ref{hierarchy} (see (6.2)-(6.5)). Since $\eta_m \in \overline{\mu}_2(y_p) \cap \overline{\mu}_2(y_{p-1})$, $\eta_m\in D_{n,q}$,  and $\mu_2(e_p)=\gamma_1 \notin \overline{\mu}_2(y_{p-1})$, the present subcase reduces to Subcase 1.1.

Next we assume that $\theta \notin D_{n,q}$. By (6.6) and the hypothesis of the present subcase, we have $\theta \notin \overline{\varphi}_n(T_{n,q}^*)$. So $\theta \notin \overline{\varphi}_n(T_{n,q}^*) \cup D_{n,q}$, which implies
$\theta \notin \overline{\varphi}_n(T_n) \cup D_n$. In particular,

(28) $\theta \notin D_{n,q}\cup \Gamma^q$ if $q\ge 1$ and $\theta \notin D_n \cup \Gamma^0$ if $q=0$. Furthermore,
$\theta$ is not used by any edge in $T(y_{p-1}) -  T_{n,q}^*$ by TAA (see, for instance, (1)).

We proceed by considering two possibilities for $\alpha$.

$\bullet$ $\alpha\notin D_{n,q}$. Now it follows from (20) that  

(29) $\alpha \notin D_{n,q}\cup \Gamma^q$ if $q\ge 1$ and $\alpha \notin D_n \cup \Gamma^0$ if $q=0$.

By (20) and Lemma~\ref{change}, we obtain $P_{v_{\alpha}} (\alpha, \theta, \varphi_n) = P_{y_{p-1}}(\alpha, \theta, \varphi_n)$, 
which is disjoint from $P_{y_p}(\alpha, \theta, \varphi_n)$. Let $\sigma_n=\varphi_n/P_{y_p}(\alpha, \theta, \varphi_n)$. By Lemma~\ref{stablechange},  $\sigma_n$ satisfies all the properties described in (7.3). Since $\sigma_n(f)=\varphi_n(f)$ for 
each $f\in E(T(y_{p-1}))$ by (20) and (28), and $\overline{\sigma}_n(u)=\overline{\varphi}_n(u)$ for each $u\in V(T(y_{p-1}))$, 
we can obtain $T$ from $T_{n,q}^*+e_1$ by using TAA under $\sigma_n$ and hence $T$ is an ETT satisfying MP under $\sigma_n$. In view of (28) and (29), $T_n=T_{n,0} \subset T_{n,1} \subset \ldots \subset T_{n,q} \subset T_{n,q+1}=T$ remains to be a good hierarchy of $T$ under $\sigma_n$, with the same $\Gamma$-sets as those under $\varphi_n$. Therefore, $(T, \sigma_n)$ is also a minimum counterexample 
to Theorem \ref{hierarchy} (see (6.2)-(6.5)). Since $\theta \in \overline{\sigma}_n(y_p) \cap \overline{\sigma}_n(y_{p-1})$, 
$\theta \notin D_{n,q}$, and $\sigma_n(e_p)=\alpha \notin \overline{\sigma}_n(y_{p-1})$, the present subcase reduces to Subcase 2.1.

$\bullet$ $\alpha\in D_{n,q}$. Let $\alpha=\eta_{h}\in D_{n,q}$. For simplicity, we use $\varepsilon_1$ and $\varepsilon_2$
to denote the two colors $\gamma^{q}_{h_1}$ and $\gamma^{q}_{h_2}$ in $\Gamma^{q}_h$ (see Definition \ref{R2}), respectively. 
By (20), we have

(30) $\{\alpha, \varepsilon_1, \varepsilon_2\}\cap \varphi_n \langle T - T_{n,q}^* \rangle=\emptyset$.

By (30) and Lemma~\ref{change}, we obtain $P_{v_{\alpha}} (\alpha, \varepsilon_1, \varphi_n) = P_{v_{\varepsilon_1}}
(\alpha, \varepsilon_1, \varphi_n)$, which is disjoint from $P_{y_p}(\alpha, \varepsilon_1, \varphi_n)$. 
Let $\mu_1=\varphi_n/P_{y_p}(\alpha, \varepsilon_1, \varphi_n)$. By Lemma~\ref{stablechange},  $\mu_1$ satisfies 
all the properties described in (7.3) (with $\mu_1$ in place of $\sigma_n$). Since $\mu_1(f)=\varphi_n(f)$ for each $f\in E(T)$
by (30), and $\overline{\mu}_1(u)=\overline{\varphi}_n(u)$ for each $u\in V(T(y_{p-1}))$, we can obtain $T$ from $T_{n,q}^*+e_1$ by 
using TAA under $\mu_1$ and hence $T$ is an ETT satisfying MP under $\mu_1$. Furthermore, $T_n=T_{n,0} \subset T_{n,1} \subset \ldots \subset T_{n,q} \subset T_{n,q+1}=T$ remains to be a good hierarchy of $T$ under $\mu_1$, with the same $\Gamma$-sets as those under $\varphi_n$. Therefore, $(T, \mu_1)$ is also a minimum counterexample to Theorem \ref{hierarchy} (see (6.2)-(6.5)),
in which $\varepsilon_1$ is missing at two vertices. From (30) and Definition \ref{R2}(i) we see that 

(31) $\varepsilon_1 \notin \mu_1 \langle T - T_{n,q}^* \rangle$.

By (31) and Lemma~\ref{change}, we obtain $P_{v_{\varepsilon_1}} (\theta, \varepsilon_1, \mu_1) = P_{y_{p-1}}
(\theta, \varepsilon_1, \mu_1)$, which is disjoint from $P_{y_p}(\theta, \varepsilon_1, \mu_1)$. 
Let $\mu_2=\mu_1/P_{y_p}(\theta, \varepsilon_1, \mu_1)$. By Lemma~\ref{stablechange},  $\mu_2$ satisfies 
all the properties described in (7.3) (with $\mu_2$ in place of $\sigma_n$). In view of (28) and (31),
we have $\mu_2(f)=\mu_1(f)$ for each $f\in E(T(y_{p-1}))$ and $\overline{\mu}_1(u)=\overline{\varphi}_n(u)$ 
for each $u\in V(T(y_{p-1}))$. So $T$ can be obtained from $T_{n,q}^*+e_1$ by using TAA and hence is an ETT satisfying MP under 
$\mu_2$. Furthermore, $T_n=T_{n,0} \subset T_{n,1} \subset \ldots \subset T_{n,q} \subset T_{n,q+1}=T$ remains to be a good hierarchy of $T$ under $\mu_2$, with the same $\Gamma$-sets as those under $\mu_1$. Therefore, $(T, \mu_2)$ is also a minimum counterexample to Theorem \ref{hierarchy} (see (6.2)-(6.5)). Since $\theta \in \overline{\mu}_2(y_p) \cap \overline{\mu}_2(y_{p-1})$,
$\theta \notin D_{n,q}$, and $\mu_2(e_p)= \varepsilon_1 \notin \overline{\mu}_2(y_{p-1})$, the present subcase reduces to 
Subcase 2.1. This completes our discussion about Situation 7.2.\\

\noindent {\bf Situation 7.3.} $2\le p(T)\le p-1$. 

\vskip 2mm
Recall that $T=T_{n,q}^* \cup\{e_1,y_1,e_2,...,e_p,y_p\}$, and the path number $p(T)$ of $T$ is the smallest subscript 
$t \in \{1,2,...,p\}$ such that the sequence $(y_t,e_{t+1},...,e_p,y_p)$ corresponds to a path in $G$. Set $I_{\varphi_n}
=\{1\le t \le p-1: \, \overline{\varphi}_n(y_p)\cap\overline{\varphi}_n(y_t)\ne \emptyset \}$. We use $\max(I_{\varphi_n})$
to denote the maximum element of $I_{\varphi_n}$ if $I_{\varphi_n} \ne \emptyset$. For convenience, set 
$\max(I_{\varphi_n})=-1$ if $I_{\varphi_n}=\emptyset$. 

If $\max(I_{\varphi_n}) \ge p(T)$, then we may assume that $\max(I_{\varphi_n})=p-1$ (the proof is exactly the same as
that of Claim \ref{p-1}). Let $\alpha \in \overline{\varphi}_n(y_{p-1})\cap\overline{\varphi}_n(y_p)$ and 
$\beta=\varphi_n(e_{p})$. Let $\sigma_n$ be obtained from $\varphi_n$ by recoloring $e_p$ with $\alpha$ and let 
$T'=T(y_{p-1})$. Then $\beta \in \overline{\sigma}_n(y_{p-1})\cap\overline{\sigma}_n(T')$ and 
$T_n=T_{n,0} \subset T_{n,1} \subset \ldots \subset T_{n,q}\subset T'$ is a good hierarchy of $T'$ under $\sigma_n$. 
So $(T', \sigma_n)$ is a counterexample to Theorem \ref{hierarchy} (see (6.2) and (6.3)), which violates the minimality 
assumption (6.4) or (6.5) on $(T, \varphi_n)$. 

So we may assume hereafter that $\max(I_{\varphi_n})< p(T)$. Let $i=\max(I_{\varphi_n})$ if $I_{\varphi_n} \ne \emptyset$,
and let $j=p(T)$. Then $e_j$ is not incident to $y_{j-1}$. In our proof we reserve $y_0$ for the maximum vertex 
(in the order $\prec$) in $T_{n,q}^*$.  

\begin{claim}\label{j-1}
We may assume that there exists $\alpha\in\overline{\varphi}_n(y_p)\cap\overline{\varphi}_n(T(y_{j-2}))$, such that 
either $\alpha\notin\Gamma^q\cup \phibar_n(T_{n,0}^*-V(T_n))$ or $\alpha \in\Gamma^{q}_m$ for some $\eta_m\in D_{n,q}$ 
with $v_{\eta_m} \preceq y_{j-2}$.
\end{claim}

To establish this statement, we consider two cases, depending on whether $I_{\varphi}$ is nonempty. 

{\bf Case 1.} $I_{\varphi}\neq\emptyset$. 

By assumption, $\max(I_{\varphi_n})< p(T)$. So $i \le j-1$. Let $\alpha\in\overline{\varphi}_n(y_p)\cap
\overline{\varphi}_n(y_i)$. By (6.6), we obtain 

(1) $\alpha \notin \phibar_n(T_{n,q}^*)$. So $\alpha \notin \Gamma^{q}\cup \phibar_n(T_{n,0}^*-V(T_n))$. 

\noindent If $i\le j-2$, then $\alpha\in\phibar_n(T(y_{j-2}))$, as desired. Thus we may assume that $i=j-1$. 

(2) There exists a color $\beta \in \overline{\varphi}_n(T(y_{j-2}))- \phibar_n(T_{n,0}^*-V(T_n)) - \varphi_n \langle 
T(y_{j-1}) - T_{n,q}^* \rangle -(\Gamma^{q}\cup D_{n,q})$ or a color $\beta \in\Gamma^{q}_m$ for some $\eta_m\in D_{n,q}$ 
with $v_{\eta_m} \preceq y_{j-2}$ and $(\Gamma^{q}_m \cup \{\eta_m\}) \cap \varphi_n \langle T(y_{j-1}) - T_{n,q}^* 
\rangle=\emptyset$.

To justify this, note that if $|\overline{\varphi}_n(T(y_{j-2}))- \phibar_n(T_{n,0}^*-V(T_n)) - \varphi_n \langle 
T(y_{j-2}) - T_{n,q}^* \rangle -(\Gamma^{q}\cup D_{n,q}) |\ge  5$, then there exists a color $\beta$ in 
$\overline{\varphi}_n(T(y_{j-2}))- \phibar_n(T_{n,0}^*-V(T_n)) - \varphi_n \langle T(y_{j-1}) - T_{n,q}^* \rangle 
-(\Gamma^{q}\cup D_{n,q})$, because  $T(y_{j-1})-T(y_{j-2})$ contains only one edge.

If $|\overline{\varphi}_n(T(y_{j-2}))- \phibar_n(T_{n,0}^*-V(T_n)) - \varphi_n \langle 
T(y_{j-2}) - T_{n,q}^* \rangle -(\Gamma^{q}\cup D_{n,q}) |\le 4 $, then, by Lemma \ref{9n}, there exist $7$ distinct 
colors $\eta_{h}\in D_{n,q}\cap \phibar_n(T(y_{j-2}))$ such that  $(\Gamma^{q}_h \cup \{\eta_{h}\})
\cap \varphi_n \langle T(y_{j-2}) - T_{n,q}^* \rangle=\emptyset$. Since $T(y_{j-1})-T(y_{j-2})$ contains only one edge,
there exists at least one of these $\eta_{h}$, say $\eta_m$, such that $(\Gamma^{q}_m \cup \{\eta_{m}\})
\cap \varphi_n \langle T(y_{j-1}) - T_{n,q}^* \rangle=\emptyset$. So (2) is true. 

Depending on whether $\alpha$ is contained in $D_{n,q}$, we distinguish between two subcases.

{\bf Subcase 1.1.} $\alpha \in D_{n,q}$. In this subcase, let $\alpha=\eta_{h}\in D_{n,q}$. For simplicity, we abbreviate the two colors $\gamma^{q}_{h_1}$ and $\gamma^{q}_{h_2}$ in $\Gamma^{q}_h$ (see Definition \ref{R2}) to $\gamma_1$ and $\gamma_2$, 
respectively. Since $\eta_{h} \in\overline{\varphi}_n(y_{j-1})$, by Definition \ref{R2}(i) and TAA, we have

(3) $\gamma_1, \gamma_2 \notin \varphi_n \langle T(y_{j-1}) - T_{n,q}^* \rangle$, and $\eta_h$ is not used by 
any edge in $T(y_{j-1}) -  T_{n,q}^*$, except possibly $e_1$ when $q=0$ and $T_{n,0}^*=T_n$ (now $e_1=f_n$ in 
Algorithm 3.1 and $\varphi_n(e_1)=\eta_h \in  D_{n,q} \subseteq  D_n$). 

By (3) and Lemma \ref{change}, we obtain $P_{v_{\gamma_1}}(\gamma_1,\eta_h,\varphi_n)=P_{y_{j-1}}(\gamma_1,\eta_h,
\varphi_n)$, which is disjoint from $P_{y_{p}}(\gamma_1,\eta_h,\varphi_n)$. Let $\mu_1=\varphi_n/P_{y_{p}}(\gamma_1,
\eta_h,\varphi_n)$. By Lemma~\ref{stablechange},  $\mu_1$ satisfies all the properties described in (7.3) (with $\mu_1$ 
in place of $\sigma_n$). In particular, if $e_1=f_n$ and $\varphi_n(e_1)=\eta_h \in D_n$, then $\mu_1(e_1)=\varphi_n(e_1)$, 
which implies that $e_1$ is outside $P_{y_{p}}(\gamma_1,\eta_h,\varphi_n)$. 
Using (3) and (6.6), we get $\mu_1(f)=\varphi_n(f)$ for each $f\in E(T(y_{j-1}))$ and 
$\overline{\mu}_1(u)= \overline{\varphi}_n(u)$ for each $u\in V(T(y_{p-1}))$. So we can obtain $T$ from $T_{n,q}^*+e_1$ by 
using TAA under $\mu_1$, and hence $T$ is an ETT satisfying MP under $\mu_1$. Furthermore, since $\eta_h \in \overline{\mu}_1(y_{j-1})$,
the hierarchy $T_n=T_{n,0} \subset T_{n,1} \subset \ldots \subset T_{n,q} \subset T_{n,q+1}=T$ remains to be good 
under $\mu_1$, with the same $\Gamma$-sets as those under $\mu_1$. Therefore, $(T, \mu_1)$ is also a minimum counterexample 
to Theorem \ref{hierarchy} (see (6.2)-(6.5)), in which $\gamma_1$ is missing at two vertices.

From (3) we see that

(4) $\gamma_1, \gamma_2 \notin \mu_1 \langle T(y_{j-1}) - T_{n,q}^* \rangle$, and $\eta_h$ is not used by 
any edge in $T(y_{j-1}) -  T_{n,q}^*$ under $\mu_1$, except possibly $e_1$ when $q=0$ and $T_{n,0}^*=T_n$ 
(now $e_1=f_n$ in Algorithm 3.1 and $\mu_1(e_1)=\eta_h \in  D_{n,q} \subseteq  D_n$). 

Let $\beta$ be a color as specified in (2). Note that 

(5) $\beta\notin \mu_1 \langle T(y_{j-1}) - T_{n,q}^* \rangle$,  $\beta\notin D_{n,q}$, and $\beta\neq\eta_h=\alpha$. 

Since $\gamma_1 \in \overline{\mu}_1(T_{n,q})$ if $q \ge 1$ and $\gamma_1 \in \overline{\mu}_1(T_n)$ if $q=0$,
from (4) and Lemma~\ref{change} we deduce that $P_{v_{\gamma_1}} (\gamma_1, \beta, \mu_1) = P_{v_{\beta}}
(\gamma_1, \beta, \mu_1)$, which is disjoint from $P_{y_p}(\gamma_1, \beta, \mu_1)$. Let $\mu_2=\mu_1/P_{y_p}(\gamma_1, 
\beta, \mu_1)$. By Lemma~\ref{stablechange},  $\mu_2$ satisfies all the properties described in (7.3) (with $\mu_2$ 
in place of $\sigma_n$). By (4), (5) and (6.6), we have $\mu_2(f)=\mu_1(f)$ for each $f\in E(T(y_{j-1}))$, and $\overline{\mu}_2(u)=\overline{\mu}_1(u)$ for each $u\in V(T(y_{p-1}))$. So we can obtain $T$ from $T_{n,q}^*+e_1$ by 
using TAA under $\mu_2$ and hence $T$ is an ETT satisfying MP under $\mu_2$. If $\beta \notin \Gamma^q$, then clearly 
$T_n=T_{n,0} \subset T_{n,1} \subset \ldots \subset T_{n,q} \subset T_{n,q+1}=T$ remains to be a good hierarchy 
of $T$ under $\mu_2$, with the same $\Gamma$-sets as those under $\mu_1$. So we assume that $\beta \in \Gamma^q$.
By (2), we have $\beta \in\Gamma^{q}_m$ for some $\eta_m\in D_{n,q}$ with $v_{\eta_m} \preceq y_{j-2}$ and 
$(\Gamma^{q}_m \cup \{\eta_m\}) \cap \varphi_n \langle T(y_{j-1}) - T_{n,q}^* \rangle=\emptyset$. It follows
that $T_n=T_{n,0} \subset T_{n,1} \subset \ldots \subset T_{n,q} \subset T_{n,q+1}=T$ is still a good hierarchy 
of $T$ under $\mu_2$, with the same $\Gamma$-sets as those under $\mu_1$.  Therefore, $(T, \mu_2)$ is also a minimum 
counterexample to Theorem \ref{hierarchy} (see (6.2)-(6.5)), in which $\beta \in \overline{\mu}_2(y_p) \cap 
\overline{\mu}_2(T(y_{j-2}))$. From (2) and the definitions of $\mu_1$ and $\mu_2$, we see that 
either $\beta\notin\Gamma^q\cup \phibar_n(T_{n,0}^*-V(T_n))$ or $\beta \in\Gamma^{q}_m$ for some $\eta_m\in D_{n,q}$ 
with $v_{\eta_m} \preceq y_{j-2}$. Thus Claim~\ref{j-1} holds by replacing $\varphi_n$ with $\mu_2$ and $\alpha$ 
with $\beta$. 

{\bf Subcase 1.2.} $\alpha \notin D_{n,q}$. In this subcase, using (1) and the set inclusion $\overline{\varphi}_n(T_n) \cup D_n \subseteq  \overline{\varphi}_n(T_{n,q}^*) \cup D_{n,q}$, we get

(6) $\alpha\notin D_n$. So $\alpha$ is not used by any edge in $T(y_{j-1}) -  T_{n,q}^*$ by TAA.

Let $\beta$ be a color as specified in (2). Then there are two possibilities for $\beta$.

$\bullet$  $\beta \in \overline{\varphi}_n(T(y_{j-2}))- \phibar_n(T_{n,0}^*-V(T_n)) - \varphi_n \langle 
T(y_{j-1}) - T_{n,q}^* \rangle -(\Gamma^{q}\cup D_{n,q})$.  Now it follows from Lemma~\ref{change} that $P_{v_{\beta}} 
(\alpha, \beta, \varphi_n) = P_{y_{j-1}}(\alpha, \beta, \varphi_n)$, so this path is disjoint from $P_{y_p}(\alpha, \beta, \varphi_n)$.  Let $\sigma_n=\varphi_n/P_{y_p}(\alpha, \beta, \varphi_n)$. By Lemma~\ref{stablechange},  $\sigma_n$ 
satisfies all the properties described in (7.3). By (6), the assumption on $\beta$ and (6.6), we have $\sigma_n(f)=
\varphi_n(f)$ for each $f\in E(T(y_{j-1}))$, and $\overline{\sigma}_n(u)=\overline{\varphi}_n(u)$ for each 
$u\in V(T(y_{p-1}))$.  So we can obtain $T$ from $T_{n,q}^*+e_1$ by using TAA under $\sigma_n$ and hence $T$ satisfies 
MP under $\sigma_n$. Since $\alpha, \beta \notin \Gamma^q$ (see (1)), the hierarchy $T_n=T_{n,0} \subset 
T_{n,1} \subset \ldots \subset T_{n,q} \subset T_{n,q+1}=T$ remains to be good under $\sigma_n$, with the 
same $\Gamma$-sets as those under $\varphi_n$. Therefore, $(T, \sigma_n)$ is also a minimum counterexample to Theorem 
\ref{hierarchy} (see (6.2)-(6.5)), in which $\beta \in \overline{\sigma}_n(y_p)\cap\overline{\sigma}_n(T(y_{j-2}))$. 
Thus Claim~\ref{j-1} holds by replacing $\varphi_n$ with $\sigma_n$ and $\alpha$ with $\beta$. 

$\bullet$ $\beta \in\Gamma^{q}_m$ for some $\eta_m\in D_{n,q}$ with $v_{\eta_m} \preceq y_{j-2}$ and $(\Gamma^{q}_m 
\cup \{\eta_m\}) \cap \varphi_n \langle T(y_{j-1}) - T_{n,q}^* \rangle=\emptyset$. Note that $\eta_m \in \overline{\varphi}_n(T(y_{j-2}))$ and hence $\alpha\neq \eta_m$ by (6.6).  In view of Lemma~\ref{change}, we obtain
$P_{v_{\beta}} (\alpha, \beta, \varphi_n) = P_{y_{j-1}} (\alpha, \beta, \varphi_n)$, which is disjoint from 
$P_{y_p} (\alpha, \beta, \varphi_n)$. Let $\sigma_n=\varphi_n/P_{y_p} (\alpha, \beta, \varphi_n)$.
By Lemma~\ref{stablechange},  $\sigma_n$ satisfies all the properties described in (7.3). By (6), the assumption 
on $\beta$ and (6.6), we have $\sigma_n(f)=\varphi_n(f)$ for each $f\in E(T(y_{j-1}))$, and $\overline{\sigma}_n(u)=\overline{\varphi}_n(u)$ for each $u\in V(T(y_{p-1}))$.  So we can obtain $T$ from $T_{n,q}^*+e_1$ 
by using TAA under $\sigma_n$ and hence $T$ is an ETT satisfying MP under $\sigma_n$. Since $\alpha \notin \Gamma^q$ 
(see (1)) and $\eta_m \in \overline{\varphi}_n(T(y_{j-2}))$, the hierarchy $T_n=T_{n,0} \subset T_{n,1} \subset 
\ldots \subset T_{n,q} \subset T_{n,q+1}=T$ remains to be good under $\sigma_n$, with the same $\Gamma$-sets as 
those under $\varphi_n$. Therefore, $(T, \sigma_n)$ is also a minimum counterexample to Theorem \ref{hierarchy} 
(see (6.2)-(6.5)), in which $\beta \in \overline{\sigma}_n(y_p)\cap\overline{\sigma}_n(T(y_{j-2}))$. Thus Claim~\ref{j-1} holds by replacing $\varphi_n$ with $\sigma_n$ and $\alpha$ with $\beta$.

{\bf Case 2.} $I_{\varphi}=\emptyset$. 

Let $\alpha\in \overline{\varphi}_n(y_p)\cap\overline{\varphi}_n(T(y_{p-1}))$. By the hypothesis of the present
case, we have $\alpha\in\phibar_n(T_{n,q}^*)$. If $\alpha\notin\Gamma^q\cup \phibar_n(T_{n,0}^*-V(T_n))$, 
we are done. So we assume that $\alpha\in \Gamma^q\cup \phibar_n(T_{n,0}^*-V(T_n))$.  
 
{\bf Subcase 2.1.} $\alpha\in \phibar_n(T_{n,0}^*-V(T_n))-\Gamma^q$. Let us first show that

(7) there exists a color $\beta\in\phibar_n(T_{n,q}^*)-\phibar_n(T_{n,0}^*-V(T_n))-\Gamma^q$. 

Indeed, since $V(T_{n,q}^*)$ is elementary with respect to $\varphi_n$, we have $|\phibar_n(T_{n,q}^*)-\phibar_n(T_{n,0}^*-V(T_n))-
\Gamma^q| \ge |\phibar_n(T_{n,0}^*)- \phibar_n(T_{n,0}^*-V(T_n))-\Gamma^q| = |\phibar_n(T_n)- \Gamma^q|$.
In view of (7.2), we obtain $|\phibar_n(T_{n})|\ge 2n+11$ and $|\Gamma^q| \le 2 |D_{n,q}|\le 2n$. So 
$|\phibar_n(T_{n,q}^*)-\phibar_n(T_{n,0}^*-V(T_n))-\Gamma^q| \ge 11$, which implies (7).

By (7) and Lemma~\ref{change}, we obtain $P_{v_{\alpha}} (\alpha, \beta, \varphi_n) = P_{v_{\beta}}(\alpha, \beta, \varphi_n)$, 
which is disjoint from $P_{y_p}(\alpha, \beta, \varphi_n)$.  Let $\sigma_n=\varphi_n/P_{y_p}(\alpha, \beta, \varphi_n)$. By Lemma~\ref{stablechange},  $\sigma_n$ satisfies all the properties described in (7.3). Since $\alpha, \beta\in
\phibar_n(T_{n,q}^*)$, we have $\sigma_n(f)=\varphi_n(f)$ for each $f\in E(T_{n,q}^*)$, and $\overline{\sigma}_n(u)=\overline{\varphi}_n(u)$ for each $u\in V(T(y_{p-1}))$.  So we can obtain $T$ from 
$T_{n,q}^*+e_1$ by using TAA under $\sigma_n$ and hence $T$ is an ETT satisfying MP under $\sigma_n$. As $\alpha, \beta \notin \Gamma^q$, 
the hierarchy $T_n=T_{n,0} \subset T_{n,1} \subset \ldots \subset T_{n,q} \subset T_{n,q+1}=T$ remains to 
be good under $\sigma_n$, with the same $\Gamma$-sets as those under $\varphi_n$. Therefore, $(T, \sigma_n)$ is also a 
minimum counterexample to Theorem \ref{hierarchy} (see (6.2)-(6.5)), in which $\beta \in \overline{\sigma}_n(y_p)\cap\overline{\sigma}_n(T(y_{j-2}))$. 
Thus Claim~\ref{j-1} holds by replacing $\varphi_n$ with $\sigma_n$ and $\alpha$ with $\beta$.

{\bf Subcase 2.2.}  $\alpha\in \Gamma^q$. Let $\alpha \in\Gamma^{q}_m$ for some $\eta_{m}\in D_{n,q}$. Depending on
whether $\eta_m$ is contained in $\phibar_n(T(y_{p-1}))$, we consider two possibilities.

$\bullet$ $\eta_m\notin\phibar_n(T(y_{p-1}))$. By Definition \ref{R2}(i), we have $\alpha \notin\varphi_n
\langle T-T_{n,q}^* \rangle$. Since $T-T(y_{p-2})$ contains precisely two edges, Lemma~\ref{9n} guarantees
the existence of a color $\beta$ in $\overline{\varphi}_n(T(y_{p-2}))- \phibar_n(T_{n,0}^*-V(T_n)) - \varphi_n 
\langle T - T_{n,q}^* \rangle -(\Gamma^{q}\cup D_{n,q} )$ or a color $\beta=\eta_{h}\in D_{n,q}\cap 
\phibar_n(T(y_{p-2}))$ such that $(\Gamma^{q}_h \cup \{\eta_{h}\}) \cap \varphi_n \langle T - T_{n,q}^* 
\rangle=\emptyset$. Note that $\beta \in \overline{\varphi}_n(T(y_{p-2}))-\varphi_n \langle T - 
T_{n,q}^* \rangle$. By  Lemma~\ref{change}, we obtain $P_{v_{\alpha}} (\alpha, \beta, \varphi_n) = 
P_{v_{\beta}}(\alpha, \beta, \varphi_n)$, which is disjoint from $P_{y_p}(\alpha, \beta, \varphi_n)$.  Let $\sigma_n=\varphi_n/P_{y_p}(\alpha, \beta, \varphi_n)$. By Lemma~\ref{stablechange},  $\sigma_n$ satisfies all 
the properties described in (7.3). Since $\alpha, \beta \notin\varphi_n \langle T-T_{n,q}^* \rangle$
and $\alpha, \beta\in \phibar_n(T(y_{p-2}))$, we have $\sigma_n(f)=\varphi_n(f)$ 
for each $f\in E(T)$, and $\overline{\sigma}_n(u)=\overline{\varphi}_n(u)$ for each $u\in V(T(y_{p-1}))$.  
So we can obtain $T$ from $T_{n,q}^*+e_1$ by using TAA under $\sigma_n$ and hence $T$ is an ETT satisfying MP under $\sigma_n$. 
Furthermore, the hierarchy $T_n=T_{n,0} \subset T_{n,1} \subset \ldots \subset T_{n,q} \subset T_{n,q+1}=T$ remains 
to be good under $\sigma_n$, with the same $\Gamma$-sets as those under $\varphi_n$. Therefore, $(T, \sigma_n)$ is 
also a minimum counterexample to Theorem \ref{hierarchy} (see (6.2)-(6.5)), in which $\beta \in \overline{\sigma}_n(y_p)\cap 
\overline{\sigma}_n(v_{\beta})$. Thus Claim~\ref{j-1} holds if $v_{\beta}\preceq y_{j-2}$, the present subcase
reduces to the case when $\max(I_{\sigma_n}) \ge p(T)$ if $y_j \preceq v_{\beta}$ (see the paragraphs above
Claim \ref{j-1}), and the present subcase reduces to Case 1 (where $I_{\sigma_n}\neq \emptyset$) if $y_{j-1}=v_{\beta}$.
 
$\bullet$ $\eta_m\in\phibar_n(T(y_{p-1}))$. Note that $\eta_m\notin\phibar_n(T_{n,q}^*)$ because $\eta_m\in D_{n,q}$. 
So $\eta_{m}\in\phibar_n(y_t)$ for some $1\leq t \leq p-1$. If $t\le j-2$, then Claim~\ref{j-1} holds. Thus we may assume
that $t \ge j-1$.  Since $\eta_m \in \phibar_n(y_t)$, it is not used by any edge in $T(y_t) -  T_{n,q}^*$, except possibly 
$e_1$ when $q=0$ and $T_{n,0}^*=T_n$ (now $e_1=f_n$ in Algorithm 3.1 and $\varphi_n(e_1)=\eta_m \in  D_{n,q} \subseteq  
D_n$). Since $\alpha \in\Gamma^{q}_m$, by Definition \ref{R2}(i), $\alpha$ is not used by any edge in $T(y_t) -  T_{n,q}^*$. 
It follows from Lemma~\ref{change} that $P_{v_{\alpha}} (\alpha, \eta_m, \varphi_n) = P_{y_t} (\alpha, \eta_m, \varphi_n)$,
which is disjoint from $P_{y_p} (\alpha,\eta_m, \varphi_n)$. Let $\sigma_n=\varphi_n/P_{y_p} (\alpha,\eta_m, \varphi_n)$. By Lemma~\ref{stablechange},  $\sigma_n$ satisfies all the properties described in (7.3). 
In particular, if $e_1=f_n$ and $\varphi_n(e_1)=\eta_m \in D_n$, then $\sigma_n(e_1)=\varphi_n(e_1)$, which 
implies that $e_1$ is outside $P_{y_p} (\alpha,\eta_m, \varphi_n)$. Since $\sigma_n(f)=\varphi_n(f)$ 
for each $f\in E(T(y_t))$ and $\overline{\sigma}_n(u)=\overline{\varphi}_n(u)$ for each $u\in V(T(y_{p-1}))$, we can 
obtain $T$ from $T_{n,q}^*+e_1$ by using TAA under $\sigma_n$, so $T$ is an ETT satisfying MP under $\sigma_n$. Furthermore,
As $\alpha, \eta_m \in \overline{\sigma}_n(T(y_t))$, the hierarchy $T_n=T_{n,0} \subset T_{n,1} \subset \ldots 
\subset T_{n,q} \subset T_{n,q+1}=T$ remains to be good under $\sigma_n$, with the same $\Gamma$-sets as those under 
$\varphi_n$. Therefore, $(T, \sigma_n)$ is also a minimum counterexample to Theorem \ref{hierarchy} (see (6.2)-(6.5)), 
in which $\eta_m \in \overline{\sigma}_n(y_p)\cap\overline{\sigma}_n(y_t)$. Thus the present subcase
reduces to the case when $\max(I_{\sigma_n}) \ge p(T)$ if $j \preceq t$ (see the paragraphs above
Claim \ref{j-1}), and reduces to Case 1 (where $I_{\sigma_n}\neq \emptyset$) if $t=j-1$. This proves 
Claim~\ref{j-1}. 

\vskip 3mm

Let $\alpha$ be a color as specified in Claim \ref{j-1}; that is, $\alpha\in\overline{\varphi}_n(y_p)\cap\overline{\varphi}_n(T(y_{j-2}))$, such that either $\alpha\notin\Gamma^q\cup \phibar_n(T_{n,0}^*-V(T_n))$ or $\alpha \in\Gamma^{q}_m$ for some $\eta_m\in D_{n,q}$ with $v_{\eta_m} \preceq y_{j-2}$.
Since $T(y_j)-T(y_{j-2})$ contains precisely two edges, Lemma~\ref{9n} guarantees the existence of a color $\beta$ in $\overline{\varphi}_n(T(y_{j-2}))- \phibar_n(T_{n,0}^*-V(T_n)) - \varphi_n \langle T(y_j) - T_{n,q}^* \rangle -
(\Gamma^{q}\cup D_{n,q} )$ or a color $\beta=\eta_{h}\in D_{n,q}\cap \phibar_n(T(y_{j-2}))$ such that $(\Gamma^{q}_h \cup \{\eta_{h}\}) \cap \varphi_n \langle T(y_j) - T_{n,q}^* \rangle=\emptyset$. Note that 

(8) $\beta \notin \varphi_n \langle T(y_j) - T_{n,q}^* \rangle \cup \Gamma^{q}$. 

\noindent Let $Q=P_{y_{p}}(\alpha,\beta,\varphi_n)$. We consider two cases, depending on whether $Q$ intersects  $T(y_{j-1})$.

{\bf Case 1.} $Q$ and $T(y_{j-1})$ have vertices in common. Let $u$ be the first vertex of $Q$ contained in
$T(y_{j-1})$ as we traverse $Q$ from $y_p$. Define $T'=T(y_{j-1}) \cup Q[u,y_p]$ if $u= y_{j-1}$
and $T'=T(y_{j-2})\cup Q[u,y_p]$ otherwise.  By the choices of $\alpha$ and $\beta$, we have $\alpha,\beta \in
\phibar_n(T(y_{j-2}))$. So $T'$ can be obtained from $T(y_{j-2})$ by using TAA under $\varphi_n$. It follows
that $T'$ is an ETT satisfying MP with respect to $\varphi_n$, with $p(T')<p(T)$. If $\alpha\notin \Gamma^{q}$, then 
both $\alpha$ and $\beta$ are outside $\Gamma^{q}$ (see (8)), so $T_n=T_{n,0} \subset T_{n,1} \subset 
\ldots \subset T_{n,q}\subset T'$ is a good hierarchy of $T'$ under $\varphi_n$, with the same $\Gamma$-sets 
as $T$ under $\varphi_n$. If $\alpha\in \Gamma^{q}$, then $\alpha \in\Gamma^{q}_m$ for some $\eta_m\in D_{n,q}$ 
with $v_{\eta_m} \preceq y_{j-2}$ by Claim \ref{j-1}. Since $\alpha,\eta_m \in \phibar_n(T(y_{j-2}))$ and 
$\beta \notin \Gamma^{q}$, it is clear that $T_n=T_{n,0} \subset T_{n,1} \subset \ldots \subset T_{n,q}\subset T'$ 
is also a good hierarchy of $T'$ under $\varphi_n$, with the same $\Gamma$-sets as $T$ under $\varphi_n$. 
So $(T', \varphi_n)$ is a counterexample to Theorem \ref{hierarchy} (see (6.2) and (6.3)), which violates the 
minimality assumption $(6.4)$ on $(T, \varphi_n)$.

{\bf Case 2.} $Q$ is vertex-disjoint from $T(y_{j-1})$. Let $\sigma_n=\varphi_n/Q$. By Lemma \ref{LEM:Stable}, $\sigma_n$ is  $(T(y_{j-1}),D_n,\varphi_n)$-stable. In particular, $\sigma_n$ is  $(T(y_{j-1}),\varphi_n)$-invariant. If $\Theta_n=PE$, 
then $\sigma_n$ is also $(T_n\oplus R_n,D_n,\varphi_n)$-stable. Furthermore, $T(y_{j-1})$ is an ETT satisfying MP with respect to $\sigma_n$, and $T_n=T_{n,0} \subset T_{n,1} \subset \ldots \subset T_{n,q} \subset T(y_{j-1})$ is a good hierarchy of 
$T(y_{j-1})$ under $\sigma_n$, with the same $\Gamma$-sets as $T$ under $\varphi_n$. By definition, $\sigma_n$ is a $(T_{n,q}^*, D_n,\varphi_n)$-weakly stable coloring. If $\alpha\notin \Gamma^{q}$, then both $\alpha$ 
and $\beta$ are outside $\Gamma^{q}$ (see (8)), so $T_n=T_{n,0} \subset T_{n,1} \subset 
\ldots \subset T_{n,q}\subset T$ is a good hierarchy of $T$ under $\sigma_n$, with the same $\Gamma$-sets 
as $T$ under $\varphi_n$. If $\alpha\in \Gamma^{q}$, then $\alpha \in\Gamma^{q}_m$ for some $\eta_m\in D_{n,q}$ 
with $v_{\eta_m} \preceq y_{j-2}$ by Claim \ref{j-1}. Since $\alpha,\eta_m \in \phibar_n(T(y_{j-2}))$ and 
$\beta \notin \Gamma^{q}$, it is clear that $T_n=T_{n,0} \subset T_{n,1} \subset \ldots \subset T_{n,q}\subset T$ 
is also a good hierarchy of $T$ under $\sigma_n$, with the same $\Gamma$-sets as $T$ under $\varphi_n$. 
So $(T, \sigma_n)$ is a counterexample to Theorem \ref{hierarchy}, in which $\beta$ is missing at two vertices. 

From the choice of $\beta$ above $(8)$ and the definition of $\sigma_n$, we see that

(9) either $\beta \notin \overline{\sigma}_n(T_{n,0}^*-V(T_n)) \cup \sigma_n \langle T(y_j) - T_{n,q}^* \rangle \cup
(\Gamma^{q}\cup D_{n,q} )$ or $\beta=\eta_{h}\in D_{n,q}\cap \overline{\sigma}_n(T(y_{j-2}))$, such that $(\Gamma^{q}_h 
\cup \{\eta_{h}\}) \cap \sigma_n \langle T(y_j) - T_{n,q}^* \rangle=\emptyset$.  

Let $\theta \in\overline{\sigma}_n(y_j)$. Then $\theta \notin\Gamma^{q}$. We proceed by considering two subcases.

{\bf Subcase 2.1.} $\theta \notin D_{n,q}$. In this subcase, using (6.6) and the set inclusion $\overline{\varphi}_n(T_n) 
\cup D_n \subseteq  \overline{\varphi}_n(T_{n,q}^*) \cup D_{n,q}$, we obtain 

(10) $\theta \notin \overline{\sigma}_n(T(y_{j-1}))$ and $\theta \notin D_n$. So $\theta$ is not assigned to any
edge in $T(y_{j})-T_{n,q}^*$ by TAA. 

As described in (9), there are two possibilities for $\beta$.

$\bullet$ $\beta \notin \overline{\sigma}_n(T_{n,0}^*-V(T_n)) \cup \sigma_n \langle T(y_j) - T_{n,q}^* \rangle \cup
(\Gamma^{q}\cup D_{n,q} )$. Observe that $\beta\notin D_{n}$ if $q=0$. By Lemma~\ref{change}, we obtain 
$P_{v_{\beta}} (\beta, \theta, \sigma_n) = P_{y_j}(\beta, \theta, \sigma_n)$, which is disjoint from $P_{y_p}(\beta,
\theta, \sigma_n)$.  Let $\mu_1=\sigma_n/P_{y_p}(\beta, \theta, \sigma_n)$. By Lemma~\ref{stablechange},  $\mu_1$ 
satisfies all the properties described in (7.3). By (10), the assumption on $\beta$ and (6.6), we have
$\mu_1(f)=\sigma_n(f)$ for each $f\in E(T(y_j))$ and $\overline{\mu}_1(u)=\overline{\sigma}_n(u)$ for each 
$u\in V(T(y_{p-1}))$.  So we can obtain $T$ from $T_{n,q}^*+e_1$ by using TAA under $\mu_1$ and hence $T$ is an ETT satisfying MP under 
$\mu_1$. As $\beta, \theta \notin \Gamma^q$, the hierarchy $T_n=T_{n,0} \subset T_{n,1} \subset \ldots \subset 
T_{n,q} \subset T_{n,q+1}=T$ remains to be good under $\mu_1$, with the same $\Gamma$-sets as those under $\sigma_n$. 
Therefore, $(T, \mu_1)$ is also a minimum counterexample to Theorem \ref{hierarchy} (see (6.2)-(6.5)), in which 
$\theta \in \overline{\mu}_1(y_p)\cap \overline{\mu}_1(y_j)$. Thus the present subcase reduces to the case when 
$\max(I_{\mu_1}) \ge p(T)$ (see the paragraphs above Claim \ref{j-1}).

$\bullet$ $\beta=\eta_{h}\in D_{n,q}\cap \overline{\sigma}_n(T(y_{j-2}))$, such that $(\Gamma^{q}_h \cup \{\eta_{h}\}) 
\cap \sigma_n \langle T(y_j) - T_{n,q}^* \rangle=\emptyset$. For simplicity, we abbreviate the two colors $\gamma^{q}_{h_1}$ and $\gamma^{q}_{h_2}$ in $\Gamma^{q}_h$ (see Definition \ref{R2}) to $\gamma_1$ and $\gamma_2$, respectively.  
By Lemma~\ref{change}, we obtain $P_{v_{\beta}} (\beta, \gamma_1, \sigma_n) = P_{v_{\gamma_1}}(\beta, \gamma_1, \sigma_n)$, 
which is disjoint from $P_{y_p}(\beta, \gamma_1, \sigma_n)$.  Let $\mu_2=\sigma_n/P_{y_p}(\beta, \gamma_1, \sigma_n)$.
By Lemma~\ref{stablechange},  $\mu_2$ satisfies all the properties described in (7.3). By the assumption on $\beta$,
neither $\beta$ nor $\gamma_1$ is used by any edge in $T(y_j) - T_{n,q}^*$. So $\mu_2(f)=\sigma_n(f)$ for each $f\in E(T(y_j))$. 
By (6.6), we get $\overline{\mu}_2(u)=\overline{\sigma}_n(u)$ for each $u\in V(T(y_{p-1}))$. It follows that $T$ can 
be obtained from $T_{n,q}^*+e_1$ by using TAA under $\mu_2$ and hence $T$ is an ETT satisfying MP under $\mu_2$. Furthermore, the 
hierarchy $T_n=T_{n,0} \subset T_{n,1} \subset \ldots \subset T_{n,q} \subset T_{n,q+1}=T$ remains to be good under 
$\mu_2$, with the same $\Gamma$-sets as those under $\sigma_n$. Therefore, $(T, \mu_2)$ is also a minimum counterexample 
to Theorem \ref{hierarchy} (see (6.2)-(6.5)), in which $\gamma_1$ is missing at both $y_p$ and $v_{\gamma_1}$. 

From the assumption on $\beta$ and the definition of $\mu_2$, we deduce that 

(11) $\beta=\eta_{h}\in D_{n,q}\cap \overline{\mu}_2(T(y_{j-2}))$, such that $(\Gamma^{q}_h \cup \{\eta_{h}\}) 
\cap \mu_2 \langle T(y_j) - T_{n,q}^* \rangle=\emptyset$. 

By (11) and Lemma~\ref{change}, we obtain $P_{v_{\gamma_1}} (\theta, \gamma_1, \mu_2) = P_{y_j}(\theta, \gamma_1, 
\mu_2)$, which is disjoint from $P_{y_p}(\theta, \gamma_1, \mu_2)$. Let $\mu_3=\mu_2/P_{y_p}(\theta, \gamma_1, \mu_2)$.
By Lemma~\ref{stablechange},  $\mu_3$ satisfies all the properties described in (7.3). By (10), (11) and (6.6),
we have $\mu_3(f)=\mu_2(f)$ for each $f\in E(T(y_j))$ and $\overline{\mu}_3(u)=\overline{\mu}_2(u)$ for each 
$u\in V(T(y_{p-1}))$.  So we can obtain $T$ from $T_{n,q}^*+e_1$ by using TAA under $\mu_3$ and hence $T$ is ETT satisfying MP 
under $\mu_3$. Furthermore, the hierarchy $T_n=T_{n,0} \subset T_{n,1} \subset \ldots \subset T_{n,q} \subset 
T_{n,q+1}=T$ remains to be good under $\mu_3$, with the same $\Gamma$-sets as those under $\mu_2$. Therefore, 
$(T, \mu_3)$ is also a minimum counterexample to Theorem \ref{hierarchy} (see (6.2)-(6.5)), in which $\theta$ is 
missing at both $y_p$ and $y_j$. Thus the present subcase reduces to the case when $\max(I_{\mu_3}) \ge p(T)$ (see 
the paragraphs above Claim \ref{j-1}).

{\bf Subcase 2.2.} $\theta \in D_{n,q}$. Let $\theta=\eta_t \in D_{n,q}$. For simplicity, we use $\varepsilon_1$ and 
$\varepsilon_2$ to denote the two colors $\gamma^{q}_{t_1}$ and $\gamma^{q}_{t_2}$ in $\Gamma^{q}_t$ (see Definition 
\ref{R2}), respectively. Then

(12) $\varepsilon_1, \varepsilon_2 \notin \sigma_n \langle T(y_j) - T_{n,q}^* \rangle$ and  $\eta_t$ is not used by 
any edge in $T(y_j) -  T_{n,q}^*$ under $\sigma_n$, except possibly $e_1$ when $q=0$ and $T_{n,0}^*=T_n$ (now $e_1=f_n$ 
in Algorithm 3.1 and $\sigma_n(e_1)=\eta_t  \in  D_{n,q} \subseteq D_n$). 

By (12) and Lemma~\ref{change} (with $\varepsilon_1$ in place of $\alpha$), we obtain $P_{v_{\varepsilon_1}} (\varepsilon_1, 
\beta, \sigma_n) = P_{v_{\beta}}(\varepsilon_1, \beta, \sigma_n)$, which is disjoint from $P_{y_p}(\varepsilon_1, \beta, 
\sigma_n)$. Let $\mu_4=\sigma_n/P_{y_p}(\varepsilon_1, \beta, \sigma_n)$. By Lemma~\ref{stablechange},  $\mu_4$ satisfies all 
the properties described in (7.3). By (9), we have $\beta\notin \sigma_n \langle T(y_j) - T_{n,q}^* \rangle$,
which together with (12) and (6.6) implies $\mu_4(f)=\sigma_n(f)$ for each $f\in E(T(y_j))$ and $\overline{\mu}_4(u)=\overline{\sigma}_n(u)$ for each $u\in V(T(y_{p-1}))$.  So we can obtain $T$ from $T_{n,q}^*+e_1$ 
by using TAA under $\mu_4$ and hence $T$ is an ETT satisfying MP under $\mu_4$. Since $\beta \notin  \Gamma^q$ by (9)
and $\eta_t \in\overline{\mu}_4(y_j)$, the hierarchy $T_n=T_{n,0} \subset T_{n,1} \subset \ldots \subset 
T_{n,q} \subset T_{n,q+1}=T$ remains to be good under $\mu_4$, with the same $\Gamma$-sets as those under $\sigma_n$. 
Therefore, $(T, \mu_4)$ is also a minimum counterexample to Theorem \ref{hierarchy} (see (6.2)-(6.5)), in which 
$\varepsilon_1$ is missing at both $y_p$ and $v_{\varepsilon_1}$. 

From (12) and (6.6) it can be seen that 

(13) $\varepsilon_1, \varepsilon_2 \notin \mu_4 \langle T(y_j) - T_{n,q}^* \rangle$ and $\eta_t \notin 
\overline{\mu}_4(T(y_{j-1}))$. So $\eta_t$ is not used by any edge in $T(y_j)-T_{n,q}^*$ under $\mu_4$, except 
possibly $e_1$ when $q=0$ and $T_{n,0}^*=T_n$ (now $e_1=f_n$ in Algorithm 3.1 and $\mu_4(e_1)=\eta_t \in D_{n,q} \subseteq D_n$).  

By (13) and Lemma~\ref{change}, we obtain $P_{v_{\varepsilon_1}} (\varepsilon_1, \eta_t, \mu_4) = 
P_{y_j}(\varepsilon_1, \eta_t, \mu_4)$, which is disjoint from $P_{y_p}(\varepsilon_1, \eta_t, \mu_4)$.
Let $\mu_5=\mu_4/P_{y_p}(\varepsilon_1, \eta_t, \mu_4)$. By Lemma~\ref{stablechange},  $\mu_5$ satisfies all 
the properties described in (7.3). In particular, if $e_1=f_n$ and $\mu_4(e_1)=\eta_t  \in D_n$, then $\mu_5(e_1)=\mu_4(e_1)$, 
which implies that $e_1$ is outside $P_{y_p}(\varepsilon_1, \eta_t, \mu_4)$. By (13) and (6.6), we have $\mu_5(f)=\mu_4(f)$ 
for each $f\in E(T(y_j))$ and $\overline{\mu}_5(u)=\overline{\mu}_4(u)$ for each $u\in V(T(y_{p-1}))$.  So we can obtain 
$T$ from $T_{n,q}^*+e_1$ by using TAA under $\mu_5$ and hence $T$ is an ETT satisfying MP under $\mu_5$. Since $\eta_t, \varepsilon_1 \in
\overline{\mu}_5(T(y_j))$, the hierarchy $T_n=T_{n,0} \subset T_{n,1} \subset \ldots \subset 
T_{n,q} \subset T_{n,q+1}=T$ remains to be good under $\mu_5$, with the same $\Gamma$-sets as those under $\mu_4$. 
Therefore, $(T, \mu_5)$ is also a minimum counterexample to Theorem \ref{hierarchy} (see (6.2)-(6.5)), in which 
$\theta=\eta_t$ is missing at both $y_p$ and $y_j$. Thus the present subcase reduces to the case when $\max(I_{\mu_5}) 
\ge p(T)$ (see the paragraphs above Claim \ref{j-1}).

This completes our discussion about Situation 7.3 and hence our proof of Theorem \ref{hierarchy}}. \qed

\subsection{Proof of Theorem \ref{thm:tech10}(ii)}

In the preceding subsection we have proved Theorem \ref{hierarchy}} and hence Theorem \ref{thm:tech10}(i).
To complete the proof of Theorem \ref{thm:tech10}, we still need to establish the interchangeability property
as described in Theorem \ref{thm:tech10}(ii).

\begin{lemma}\label{rutgers}
Suppose Theorem~\ref{thm:tech10}(i), (iv), and (vi) hold for all ETTs with $n$ rungs and satisfying MP, and
suppose Theorem~\ref{thm:tech10}(ii) holds for all ETTs with $n-1$ rungs and satisfying MP. Then  
Theorem~\ref{thm:tech10}(ii) holds for all ETTs with $n$ rungs and satisfying MP; that is, 
$T_{n+1}$ has the interchangeability property with respect to $\varphi_n$.
\end{lemma}

{\bf Proof.} Let $T=T_{n+1}$, let $\sigma_n$ be a $(T, D_n, \phiv_{n})$-stable coloring, and let   
$\alpha$ and $\beta$ be two colors in $[k]$ with $\alpha\in \overline{\sigma}_n(T)$ (equivalently 
$\alpha\in \phibar_{n}(T)$). We aim to prove that $\alpha$ and $\beta$ are $T$-interchangeable under
$\sigma_n$. Assume the contrary: there are at least two $(\alpha,\beta)$-paths $Q_1$ and $Q_2$ with respect to $\sigma_n$ intersecting 
$T$. By Theorem~\ref{thm:tech10}(i), $V(T)$ is elementary with respect to $\varphi_n$, so it is also elementary with 
respect to $\sigma_n$. Since $T=T_{n+1}$ is closed with respect to $\varphi_n$, it is also closed with respect to $\sigma_n$. 
As $\alpha\in \overline{\sigma}_n(T)$, it follows that $|V(T)|$ is odd and $\beta \notin \overline{\sigma}_n(T)$. 
From the existence of $Q_1$ and $Q_2$, we see that $G$ contains at least three $(T, \sigma_n, \{\alpha, \beta\})$-exit paths $P_1,P_2,P_3$.

We call the tuple $(\sigma_n, T, \alpha, \beta, P_1,P_2,P_3)$ a {\em counterexample} and use 
${\cal K}$ to denote the set of all such counterexamples.
With a slight abuse of notation, let $(\sigma_n, T, \alpha, \beta, P_1,P_2,P_3)$ be a counterexample in ${\cal K}$ with the minimum $|P_1|+|P_2|+|P_3|$. For $i=1,2,3$, let $a_i$ and $b_i$ be the 
ends of $P_i$ with $b_i \in V(T)$, and $f_i$ be the edge of $P_i$ incident to $b_i$.  Renaming subscripts if necessary, 
we may assume that $b_1\prec b_2 \prec b_3$.  We propose to show that

(1) $b_2 \notin V(T_n)$ 

Otherwise, $b_2 \in V(T_n)$. Let $\gamma$ be a color in $\overline{\sigma}_n(T_n)-\{\delta_n\}$ if $\Theta_n=PE$ and a color in $\overline{\sigma}_n(T_n)$ 
otherwise.  Since $T=T_{n+1}$ is closed with respect to $\sigma_n$, both $\alpha$ and $\gamma$ are closed in 
$T$ with respect to $\sigma_n$. Let $\mu_1=\sigma_n/ (G-T,\alpha,\gamma)$.  Then $P_1$ and $P_2$ are two $(T_n, \mu_1, \{\gamma, \beta\})$-exit paths. 
By Lemma \ref{LEM:Stable}, $\mu_1$ is a $(T,D_n,\sigma_n)$-stable coloring, so it is also $(T,D_n,\varphi_n)$-stable. As $T_n\subset T$, 
$\mu_1$ is a $(T_n,D_n,\varphi_n)$-stable coloring. 

If $\Theta_n=SE$ or $RE$ then, by Algorithm 3.1 and Lemma \ref{hku}(i), $\mu_1$ is $(T_n,D_{n-1},\varphi_{n-1})$-stable and hence is 
$(T_{n-1},D_{n-1},\varphi_{n-1})$-stable. By Theorem~\ref{thm:tech10}(vi) and TAA, $T_n$ is an ETT corresponding to $\mu_1$ (see Definition 
\ref{wz2}) and satisfies MP under $\mu_1$, with $n-1$ rungs. Since $P_1$ and $P_2$ are two $(T_n, \mu_1, \{\gamma, \beta\})$-exit paths and 
$\gamma\in\overline{\mu}_1(T_n)=\overline{\sigma}_n(T_n)$, there are at least two $(\gamma, \beta)$-paths with respect to $\mu_1$ intersecting 
$T_n$. Hence $\gamma$ and $\beta$ are not $T_n$-interchangeable under $\mu_1$, contradicting Theorem~\ref{thm:tech10}(ii) because $T_n$ has $n-1$ rungs.  

So we assume that $\Theta_n=PE$. Since $T_n$ is an ETT under $\varphi_{n-1}$, $V(T_n)$ is elementary under $\varphi_{n-1}$ by 
Theorem~\ref{thm:tech10}(i), and hence $\delta_n\in\phibar_n(T_n)$ and $\gamma_n\notin\phibar_n(T_n)$ by Algorithm 3.1. Since 
$\sigma_n$ is $(T, D_n, \phiv_{n})$-stable, $\phibar_n(T_n)=\overline{\sigma}_n(T_n)$ and  $\partial_{\varphi_n, \gamma_n}(T_n)=
\partial_{\sigma_n, \gamma_n}(T_n)$. As $\gamma\in\overline{\sigma}_n(T_n)-\delta_n$, $\delta_n\in\overline{\sigma}_n(T_n)$, 
and $\gamma_n, \beta \notin \overline{\sigma}_n(T_n)$,  we have $\gamma\notin \{\gamma_n,\delta_n\}$ and $\beta\neq\delta_n$. 
In view of Lemma \ref{hku}(v), we obtain $|\partial_{\varphi_n, \gamma_n}(T_n)|=1$. So $|\partial_{\sigma_n, \gamma_n}(T_n)|=1$, 
which implies $\beta\neq\gamma_n$, because $\{f_1,f_2\} \subseteq  \partial_{\sigma_n, \beta}(T_n)$ (as $b_2\in V(T_n)$). Therefore 
$\{\beta, \gamma\}\cap \{\gamma_n,\delta_n\}=\emptyset$. 
Since $\mu_1$ is $(T_n,D_n,\varphi_n)$-stable, $P_{v_n}(\gamma_n,\delta_n,\mu_1) \cap T_n=\{v_n\}$ by Theorem~\ref{thm:tech10}(iv). 
Let $\mu_2=\mu_1/P_{v_n}(\gamma_n,\delta_n,\mu_1)$. Then $\mu_2$ is $(T_n,D_{n-1},\varphi_{n-1})$-stable by Lemma~\ref{stablezang}.
As $\{\beta, \gamma\}\cap \{\gamma_n,\delta_n\}=\emptyset$,  we see that $P_1$ and $P_2$ are two $(T_n, \mu_2, \{\gamma, \beta\})$-exit 
paths and $\gamma\in\overline{\mu}_2(T_n)$. So there are at least two $(\gamma, \beta)$-paths with respect to $\mu_2$ intersecting $T_n$. 
Thus $\gamma$ and $\beta$ are not $T_n$-interchangeable under $\mu_2$, contradicting Theorem~\ref{thm:tech10}(ii) because $T_n$ has $n-1$ rungs. 
Therefore (1) is established.

Let $\gamma\in\overline{\sigma}_n(b_3)$ and let $\mu_3=\sigma_n/(G-T,\alpha,\gamma)$. By Lemma \ref{LEM:Stable}, 
$\mu_3$ is $(T,D_n,\varphi_n)$-stable and hence is $(T_n,D_n,\varphi_n)$-stable. By Theorem~\ref{thm:tech10}(vi), $T$
is an ETT corresponding to $\mu_3$ and satisfies MP under $\mu_3$. Furthermore, $f_i$ is colored by $\beta$ under both $\mu_3$ 
and $\sigma_n$, for $i=1,2,3$, and $P_3=P_{b_3}(\beta, \gamma, \mu_3)$.

Consider $\mu_4=\mu_3/P_{b_3}(\beta, \gamma, \mu_3)$. Clearly, $\beta \in \overline{\mu}_4(b_3)$. Since $P_{b_3}(\beta, \gamma, \mu_3)\cap T=\{b_3\}$,
by Lemma \ref{LEM:Stable}, $\mu_4$ is $(T(b_3)-b_3, D_n,\varphi_n)$-stable and $T(b_3)-b_3$ is an ETT corresponding to $\mu_4$ and satisfies MP 
under $\mu_4$. Since $b_2\prec b_3$, it is contained in $T(b_3)-b_3$. So (1) implies that $b_3$ is not the first vertex added to $T_n$ 
in the construction of $T$. According to Algorithm 3.1, $b_3$ is added to $T(b_3)-b_3$ by TAA under 
$\varphi_{n}$. Since colors on the edges of $T(b_3)$ are not affected under this Kempe change and $\mu_3$ is $(T,D_n,\varphi_n)$-stable, 
$b_3$ can still be added to $T(b_3)-b_3$ by TAA under $\mu_4$. Hence $T(b_3)$ is still an ETT satisfying MP under $\mu_4$ by Theorem~\ref{thm:tech10}(vi). 
Let $T'$ be a closure of $T(b_3)$ under $\mu_4$. Then $T'$ is an ETT satisfying MP under $\mu_4$. Since both $f_1$ and $f_2$ are colored by $\beta$ 
under $\mu_4$ and $\beta \in \overline{\mu}_4(b_3)$, the ends of $f_1$ and $f_2$ are all contained in $T'$.  By Theorem~\ref{thm:tech10}(i), $V(T')$ 
is elementary with respect to $\mu_4$, because $T'$ has $n$ rungs. 

Observe that none of $a_1,a_2,a_3$ is contained in $T'$, for otherwise, let $a_i \in V(T')$ for some $i$
with $1\le i \le 3$. Since $\{\beta,\gamma\}\cap \overline{\mu}_4(a_i) \ne \emptyset$ and $\beta \in \overline{\mu}_4(b_3)$, 
we obtain $\gamma \in \overline{\mu}_4(a_i)$. Hence from TAA we see that $P_1,P_2,P_3$ are all entirely 
contained in $G[T']$, which in turn implies $\gamma \in \overline{\mu}_4(a_j)$ for $j=1,2,3$. So $V(T')$ is not 
elementary with respect to $\mu_4$, a contradiction. Thus each $P_i$ contains a subpath $L_i$, which is a $T'$-exit path 
with respect to $\mu_4$. Since both ends of $f_1$ are contained in $T'$, $f_1$ is outside $L_1$. It follows that $|L_1|+|L_2|+|L_3|<|P_1|+|P_2|+|P_3|$.
Therefore the existence of the counterexample $(\mu_4, T', \gamma, \beta, L_1,L_2,L_3)$ violates the minimality assumption on 
$(\sigma_n, T, \alpha, \beta, P_1,P_2,P_3)$. This completes our proof of Lemma \ref{rutgers} and hence the whole proof of Theorem \ref{thm:tech10}. \qed

\newpage

\noindent {\Large \bf Subject Index}

\vskip 5mm

\noindent $(\alpha, \beta)$-chain,~\pageref{alphabetapath}

\noindent $(\alpha, \beta)$-path,~\pageref{alphabetapath}

\noindent $C$-closed subgraph with respect to $\varphi_n$,~\pageref{cclosed}

\noindent $C^-$-closed subgraph with respect to $\varphi_n$,~\pageref{cclosed}

\noindent chromatic index,~\pageref{cindex}

\noindent closed set,~\pageref{closedsets}

\noindent color class,~\pageref{colorclass}

\noindent coloring sequence~\pageref{coloringsequence}

\noindent connecting color,~\pageref{connectingcolor}

\noindent connecting edge,~\pageref{connectingedge}

\noindent critical multigraph,~\pageref{critical}

\noindent defective color,~\pageref{defectivecolor}

\noindent defective edge,~\pageref{defectiveedge}

\noindent defective vertex,~\pageref{defectivevertex}

\noindent density,~\pageref{density}

\noindent edge-coloring problem (ECP),~\pageref{ECP}

\noindent elementary multigraph,~\pageref{elementarymulti}

\noindent elementary set,~\pageref{elementaryset}

\noindent ETT corresponding to $(\sigma_n,T_n)$ or to $\sigma_n$,~\pageref{correspondingett}

\noindent extended Tashkinov tree (ETT),~\pageref{ETT}

\noindent extension vertex,~\pageref{extensionvertex}

\noindent fractional chromatic index,~\pageref{findex}

\noindent fractional edge-coloring problem (FECP),~\pageref{FECP}

\noindent $\Gamma$-set,~\pageref{gammasets}

\noindent generating coloring,~\pageref{generatingcoloring}

\noindent Goldberg-Seymour conjecture,~\pageref{GS}

\noindent good hierarchy,~\pageref{goodhierarchy}

\noindent hierarchy,~\pageref{hierar}

\noindent interchangeability property,~\pageref{interchangeproperty}

\noindent Kempe change,~\pageref{kempe}

\noindent $k$-critical multigraph,~\pageref{kcritical}

\noindent $k$-edge coloring,~\pageref{kedgecoloring}

\noindent $k$-triple,~\pageref{ktriple}

\noindent ladder,~\pageref{ladder}

\noindent maximum defective vertex,~\pageref{maximumdefective}

\noindent maximum property (MP),~\pageref{maximumproperty}

\noindent missing color,~\pageref{missingcolor}

\noindent parallel extension (PE),~\pageref{PEextension}

\noindent path number,~\pageref{pathnumber}

\noindent revisiting extension (RE),~\pageref{REextension}

\noindent rung number,~\pageref{rung}

\noindent segment of tree-sequence,~\pageref{segment}

\noindent series extension (SE),~\pageref{SEextension}

\noindent strongly closed set,~\pageref{sclosed} 

\noindent supporting vertex,~\pageref{supportingvertex}

\noindent Tashkinov's augmentation algorithm (TAA),~\pageref{TAAalgorithm}

\noindent Tashkinov tree,~\pageref{tashkinovtree}

\noindent tree-sequence,~\pageref{tsequence}

\noindent $T$-exit path,~\pageref{exitpath}

\noindent $T$-interchangeability,~\pageref{interchangable}



\noindent $(T,C,\varphi)$-stable coloring,~\pageref{tcstable}

\noindent $(T, \varphi_n)$-invariant coloring,~\pageref{tinvariant}

\noindent $(T,\varphi,\{\alpha,\beta\})$-exit,~\pageref{exit}

\noindent $(T,\varphi,\{\alpha,\beta\})$-exit path,~\pageref{tphiexitpath}  

\noindent $(T_n\oplus R_n,D_n,\varphi_n)$-stable coloring,~\pageref{trstable}

\noindent $(T_{n,0}^*, D_n,\varphi_n)$-weakly stable coloring,~\pageref{t0strong}

\noindent $(T_{n,i}^*, D_n,\varphi_n)$-weakly stable coloring,~\pageref{tistrong}

\noindent $\varphi_n \bmod T_n$-coloring,~\pageref{modcoloring}

\end{document}